\documentclass[12pt]{amsart}         
\usepackage{amscd}                   
\usepackage{amssymb}                 

\newtheorem{thm}{Theorem}[section]
\newtheorem{lem}[thm]{Lemma}

\newtheorem{cor}[thm]{Corollary}
\newtheorem{claim}[thm]{Claim}
\renewcommand{\thestep}{}
\newtheorem{prop}[thm]{Proposition}

\theoremstyle{definition}

\newtheorem{case}{Case}\renewcommand{\thecase}{}

\newtheorem{defn}[thm]{Definition}

\newtheorem{rmk}[thm]{Remark}

\theoremstyle{remark}

\makeatletter
\def\alphenumi{
  \def\theenumi{\alph{enumi}}
  \def\p@enumi{\theenumi}
  \def\labelenumi{(\@alph\c@enumi)}}
\makeatother

\makeatletter
\def\thecase{\@arabic\c@case}
\makeatother

\numberwithin{equation}{section}

\makeatletter
\def\thestep{\@arabic\c@step}
\makeatother

\newenvironment{pf}{\begin{proof}[\proofname]}{\end{proof}}
\newenvironment{pf*}[1]{\begin{proof}[#1]}{\end{proof}}

\newcommand\barc{{{\bar c}}}
\newcommand\barB{{\bar{B}}}

\newcommand\barM{{\bar{M}}}

\newcommand\barOm{{\bar\Omega}}

\newcommand\barV{{\overline{V}}}

\newcommand\barY{{\overline{Y}}}

\newcommand\barW{{\overline{W}}}

\newcommand\fs{{\mathfrak{S}}}
\newcommand\fS{{\underline{\mathfrak{S}}}}
\newcommand\ubarfV{{\underline{\mathfrak{V}}}}

\newcommand\AAA{\mathbb{A}}

\newcommand\CC{\mathbb{C}}

\newcommand\NN{\mathbb{N}}

\newcommand\RR{\mathbb{R}}

\newcommand\ZZ{\mathbb{Z}}

\newcommand\bPsi{{\boldsymbol{\Psi}}}

\newcommand\bB{{\mathbf{B}}}

\newcommand\bE{{\mathbf{E}}}

\newcommand\bF{{\mathbf{F}}}

\newcommand\bK{{\mathbf{K}}}

\newcommand\bM{{\mathbf{M}}}

\newcommand\bp{{\mathbf{p}}}

\newcommand\bS{{\mathbf{S}}}

\newcommand\bU{{\mathbf{U}}}

\newcommand\bx{{\mathbf{x}}}
\newcommand\by{{\mathbf{y}}}

\newcommand\ssG{{{}^\circ\mathcal{G}}}

\newcommand{\cov}{\nabla}

\newcommand{\rd}{\partial}

\newcommand\half{{\textstyle{\frac{1}{2}}}}

\newcommand\quarter{{\textstyle{\frac{1}{4}}}}
\newcommand\threehalf{{\textstyle{\frac{3}{2}}}}
\newcommand\threequarter{{\textstyle{\frac{3}{4}}}}
\newcommand\fivequarter{{\textstyle{\frac{5}{4}}}}

\newcommand\eighth{{\textstyle{\frac{1}{8}}}}

\newcommand\sixtyfourth{{\textstyle{\frac{1}{{64}}}}}

\newcommand\fB{{\mathfrak{B}}}

\newcommand\ff{{\mathfrak{f}}}
\newcommand\fg{{\mathfrak{g}}}

\newcommand\fh{{\mathfrak{h}}}
\newcommand\fH{{\mathfrak{H}}}

\newcommand\fm{{\mathfrak{m}}}
\newcommand\fM{{\mathfrak{M}}}

\newcommand\fV{{\mathfrak{V}}}

\newcommand\al{\alpha}
\newcommand\be{\beta}
\newcommand\De{\Delta}
\newcommand\de{\delta}
\newcommand\eps{\varepsilon}
\newcommand\ga{\gamma}
\newcommand\Ga{\Gamma}
\newcommand\la{\lambda}
\newcommand\La{\Lambda}
\newcommand\ka{\kappa}
\newcommand\om{\omega}
\newcommand\Om{\Omega}
\newcommand\si{\sigma}

\newcommand\hatA{{\hat A}}
\newcommand\hatB{{\hat B}}

\newcommand\hatPhi{{\hat\Phi}}
\newcommand\hatPsi{{\hat\Psi}}
\newcommand\hatu{{\hat u}}

\newcommand\gl{{\mathfrak{g}\mathfrak{l}}}
\newcommand\fsl{{\mathfrak{s}\mathfrak{l}}}

\newcommand\so{{\mathfrak{s}\mathfrak{o}}}

\newcommand\su{{\mathfrak{s}\mathfrak{u}}}
\newcommand\fu{{\mathfrak{u}}}

\newcommand\PU{\operatorname{PU}}

\newcommand\SO{\operatorname{SO}}

\newcommand\SU{\operatorname{SU}}
\newcommand\U{\operatorname{U}}

\newcommand\less{\setminus}

\newcommand{\8}{\infty}

\newcommand\ad{{\operatorname{ad}}}
\newcommand\Ad{{\operatorname{Ad}}}
\newcommand\asd{{\operatorname{asd}}}

\newcommand\Aut{\operatorname{Aut}}

\newcommand\ch{\operatorname{ch}}

\newcommand\Coker{\operatorname{Coker}}

\newcommand\dist{\operatorname{dist}}
\newcommand\divg{\operatorname{div}}
\newcommand\Dom{\operatorname{Dom}}

\newcommand\End{\operatorname{End}}

\newcommand\Hol{\operatorname{Hol}}
\newcommand\Hom{\operatorname{Hom}}

\newcommand\Ind{\operatorname{Ind}}
\newcommand\Imag{\operatorname{Im}}
\newcommand\Isom{\operatorname{Isom}}
\newcommand\Ker{\operatorname{Ker}}

\newcommand\loc{\operatorname{loc}}
\newcommand\Map{\operatorname{Map}}
\newcommand\Met{\operatorname{Met}}

\newcommand\Real{\operatorname{Re}}

\newcommand\Rank{\operatorname{Rank}}
\newcommand\red{{\operatorname{red}}}

\newcommand\Spec{\operatorname{Spec}}

\newcommand\Stab{\operatorname{Stab}}

\newcommand\supp{\operatorname{supp}}

\newcommand\Sym{\operatorname{Sym}}

\newcommand\tr{\operatorname{tr}}

\newcommand\vol{\operatorname{vol}}

\newcommand\id{{\mathrm{id}}}

\newcommand\reg{{\mathrm{reg}}}

\newcommand\spinc{\text{$\text{spin}^c$ }\allowbreak}

\newcommand\sA{{\mathcal{A}}}
\newcommand\sB{{\mathcal{B}}}
\newcommand\sC{{\mathcal{C}}}
\newcommand\sD{{\mathcal{D}}}

\newcommand\sG{{\mathcal{G}}}

\newcommand\sP{{\mathcal{P}}}
\newcommand\sQ{{\mathcal{Q}}}

\newcommand\sU{{\mathcal{U}}}

\newcommand\sX{{\mathcal{X}}}

\newcommand\sZ{{\mathcal{Z}}}

\newcommand\tA{{\tilde A}}

\newcommand\tsC{{\tilde\sC}}

\newcommand\tE{{\widetilde E}}

\newcommand\tOm{{\widetilde \Om}}

\newcommand\tPhi{{\tilde\Phi}}

\newcommand\tu{{\tilde u}}
\newcommand\tU{{\widetilde U}}

\newcommand\vecfm{{\vec{\fm}}}
\newcommand\vectau{{\vec{\tau}}}
\newcommand\vecvartheta{{\vec{\vartheta}}}

\begin{document}
\title[PU(2) Monopoles, I]
{PU(2) Monopoles, I: Regularity, Uhlenbeck Compactness, and
Transversality$\text{}^1$}  
\author[Paul M. N. Feehan]{Paul M. N. Feehan}
\address{Department of Mathematics\\
Harvard University\\
Cambridge, MA 02138}
\email{feehan@math.harvard.edu}
\author[Thomas G. Leness]{Thomas G. Leness}
\address{Department of Mathematics\\
Michigan State University\\
East Lansing, MI 48824}
\email{leness@math.msu.edu}
\thanks{The first author was supported in part by an NSF Mathematical 
Sciences Postdoctoral Fellowship under grant DMS 9306061.}
\footnotetext[1]{Journal of Differential Geometry, to appear.
Submitted September 11, 1996, and in revised form, May 5, 1997.
This version: October 29, 1997. \texttt{dg-ga/9710032}}
\maketitle


\section{Introduction} 
At seminars at Harvard and MIT, during October 1994, Edward Witten
introduced the $\U(1)$ monopole equations and the Seiberg-Witten invariants
to smooth four-manifold topology and conjectured their relationship with
Donaldson invariants on the basis of new developments in quantum field
theory \cite{DonSW,Witten}. The conjecture, recently extended in
\cite{MooreWitten},  
has been verified for all four-manifolds whose Donaldson and Seiberg-Witten
invariants have been independently computed.  Within two months of Witten's
announcement, a program was outlined by V. Pidstrigach and A. Tyurin and
others which should lead to a mathematical proof of the relationship
between these two invariants \cite{OTQuaternion,PTCambridge,PTLocal}.  This
approach is unrelated to the quantum field-theoretic arguments of
\cite{MooreWitten,Witten} and uses a moduli space of $\PU(2)$ monopoles to
construct a cobordism between links of Seiberg-Witten moduli spaces of
$\U(1)$ monopoles and the Donaldson moduli space of anti-self-dual
connections, which appear as singularities in this larger stratified moduli
space. 

It was soon recognized, however, that despite the appeal and elegance of
the $\PU(2)$ monopole program, its implementation 
involves substantial technical difficulties due to the contributions of
moduli spaces of $\U(1)$ monopoles in the lower levels of the Uhlenbeck
compactification of the moduli space of $\PU(2)$ monopoles.  Many of
these difficulties had never been resolved even in the case of Donaldson
theory where similar problems
arise, albeit in a rather simpler form, in attempts to
prove the Kotschick-Morgan conjecture for Donaldson invariants of
four-manifolds $X$ with $b^+(X)=1$. That conjecture asserts that the
Donaldson invariants computed using metrics lying in different chambers of
the positive cone of $H^2(X;\RR)/\RR^*$ differ by terms depending only the
homotopy type of $X$
\cite{KoM}. In the case of the Kotschick-Morgan conjecture, the heart of
the problem lies in describing the links of the lower-level reducibles via
gluing and then in calculating the pairings of the Donaldson cohomology
classes with those links. Thus far, such links have been described and
their pairings with cohomology classes computed only in certain special
cases \cite{DonApplic,DonConn,DonHCobord,DK,EG,FrQ,L,Y}. The methods used
to obtain these special cases fall very far short of the kind of general
analysis needed to prove the Kotschick-Morgan conjecture. On the other
hand, by assuming the Kotschick-Morgan conjecture, L. G\"ottsche computed
the coefficients of the wall-crossing formula in \cite{KoM} in terms of
modular forms by exploiting the presumed
homotopy invariance of the coefficients
\cite{Goettsche}. A related approach to the Witten conjecture was
proposed by Pidstrigach and Tyurin \cite{PTOberwolfach}. Certain aspects of
the $\PU(2)$ monopole program have been considered from a quantum-field
theoretic viewpoint in \cite{RConst,HPP1,HPP2,LabLoz,LabMar1,LabMar2,LabMar3}.

In the present article and its sequels \cite{FL2,FL3,FL4} we address
the analytical problems associated with constructing the
links of lower-level Seiberg-Witten moduli spaces and in establishing the
analogues of the Kotschick-Morgan conjecture for $\PU(2)$ monopoles needed
to compute the pairings of cohomology classes with these links. We hope to
return to the actual computations and a verification of Witten's conjecture
\cite{MooreWitten, Witten} in a subsequent paper.  In this article we
describe the basic regularity, Uhlenbeck compactness, and transversality
results we need for the moduli space of $\PU(2)$ monopoles and in the
sequels \cite{FL3,FL4} we develop the gluing theory required to construct
the links of the lower-level Seiberg-Witten moduli spaces. An announcement
of the main results of the this article appeared in \cite{FLGeorgia}.
 
\subsection{Statement of results}
\subsubsection{$\PU(2)$ monopoles and holonomy perturbations}
We consider Hermitian two-plane bundles $E$ over $X$ whose determinant line
bundles $\det E$ are isomorphic to a fixed Hermitian line bundle over $X$
endowed with a fixed $C^\8$, unitary connection.  Let
$(\rho,W^+,W^-)$ be a \spinc structure on $X$, where
$\rho:T^*X\to\End W$ is the Clifford map, and the Hermitian
four-plane bundle $W=W^+\oplus W^-$ is endowed with a $C^\8$ \spinc
connection. 

Let $k\ge 3$ be an integer and let $\sA_E$ be the space of $L^2_k$
connections $A$ on the $\U(2)$ bundle $E$ all inducing the fixed
determinant connection on $\det E$.  Equivalently, following \cite[\S
2(i)]{KMStructure}, we may view $\sA_E$ as the space of $L^2_k$ connections
$A$ on the $\PU(2)=\SO(3)$ bundle $\su(E)$.  We pass back and forth between
these viewpoints, via the fixed connection on $\det E$, and rely on the
context to make the distinction clear. Given a connection $A$ on $\su(E)$
with curvature $F_A\in L^2_{k-1}(\La^2\otimes\so(\su(E)))$, then
$\ad^{-1}(F_A^+) \in L^2_{k-1}(\La^+\otimes\su(E))$ is its self-dual
component, viewed as a section of $\La^+\otimes\su(E)$ via the isomorphism
$\ad:\su(E)\to\so(\su(E))$. When no confusion can arise, the isomorphism
$\ad:\su(E)\to\so(\su(E))$ will be implicit and so we regard $F_A$ as a
section of $\La^+\otimes\su(E)$ when $A$ is a connection on $\su(E)$.  Let
$D_A:L^2_k(W^+\otimes E)\to L^2_{k-1}(W^-\otimes E)$ be the corresponding
Dirac operator.

For an $L^2_k$ section $\Phi$ of $W^+\otimes E$, let $\Phi^*$ be its
pointwise Hermitian dual and let $(\Phi\otimes\Phi^*)_{00}$ be the
component of the Hermitian endomorphism $\Phi\otimes\Phi^*$ of $W^+\otimes
E$ which lies in $\su(W^+)\otimes\su(E)$. The Clifford multiplication
$\rho$ defines an isomorphism $\rho:\La^+\to\su(W^+)$ and thus an
isomorphism $\rho=\rho\otimes\id_{\su(E)}$ of $\La^+\otimes\su(E)$ with
$\su(W^+)\otimes\su(E)$. Then
\begin{align}
F_A^+ - \rho^{-1}(\Phi\otimes\Phi^*)_{00} &= 0, 
\label{eq:IntroUnpertPT}\\
D_A\Phi &= 0, \notag
\end{align}
are the unperturbed $\PU(2)$ monopole equations considered in
\cite{OTVortex,OTQuaternion,PTCambridge,PTLocal}, with a slightly differing
trace conditions (see below), for a pair $(A,\Phi)$ consisting of a
connection on $\su(E)$ and a section $\Phi$ of $W^+\otimes E$.

Donaldson's proof of the connected-sum theorem for his polynomial
invariants \cite[Theorem B]{DonPoly} makes use of certain `extended
anti-self-dual equations' \cite[Equation (4.24)]{DonPoly}
to which the Freed-Uhlenbeck generic metrics
theorem does not apply \cite[\S 4(v)]{DonPoly}. These extended equations
model a neighborhood of the product connection appearing in the Uhlenbeck
compactification of the moduli space of anti-self-dual $\SU(2)$ connections.
To obtain transversality for the zero locus of
the extended equations, he employs holonomy perturbations which give
gauge-equivariant $C^\8$ maps $\sA_E^*(X)\to
\Om^+(\su(E))$ which are then used to perturb the extended anti-self-dual
equations \cite[\S 2]{DonOrient}, \cite[pp. 282--287]{DonPoly}. 
These perturbations are continuous with respect to Uhlenbeck limits and
yield transversality not only for the top-level moduli space, but also for
all lower-level moduli spaces and for all intersections of the geometric
representatives defining the Donaldson invariants.

In \S \ref{subsubsec:HolonomyPerturbations} and in the Appendix we describe
a generalization of Donaldson's idea which we use to prove transversality
for the moduli space of solutions to a perturbed version of the $\PU(2)$
monopole equations \eqref{eq:IntroUnpertPT}.  Unfortunately, in the case of
the moduli space of $\PU(2)$ monopoles, the analysis is considerably more
intricate and the method we employ here is rather different to the one
developed in \cite{DonPoly}.  We use an infinite sequence of holonomy
sections defined on the infinite-dimensional configuration space of pairs;
when restricted to small enough open balls in the configuration space, away
from reducibles, only finitely many of these perturbing sections are
non-zero and they vanish along the reducibles.

Let $\sG_E$ be the Hilbert Lie group of $L^2_{k+1}$ unitary gauge
transformations of $E$ with determinant one. Let $S_Z^1$ denote the center
of $\U(2)$ and set
$$
\ssG_E := S_Z^1\times_{\{\pm\id_E\}}\sG_E,
$$ 
which we may view as the group of $L^2_{k+1}$ unitary gauge transformations
of $E$ with constant determinant.  The stabilizer of a unitary connection
$A$ on $E$ in $\ssG_E$ (which coincides with its stabilizer in the full
group $\Aut E$ of unitary automorphisms of $E$) always contains the center
$S^1_Z\subset\U(2)$, corresponding to the constant, central, unitary 
automorphisms of $E$. We call $A$ {\em irreducible\/} if its stabilizer
is exactly $S_Z^1$ and {\em reducible\/} otherwise. 

It is also possible, as in \cite{PTCambridge,PTLocal}, to fix a smooth
representative $\omega\in\Om^2(X,\RR)$ for $c_1(E)$ and instead consider
the space of unitary connections $A$ on $E$ satisfying the trace condition
$\tr F_A = -2\pi i\omega$, modulo the action of the full group $\Aut E$ of
unitary autommorphisms of $E$. The resulting moduli space of nonabelian
monopoles is then a torus bundle over the moduli space we define below, with
fibers $H^1(X;\RR)/H^1(X;\ZZ)$. These tori complicate the analysis
of the links of singularities when $b^1(X)>0$ and do not contain any
additional information, so we choose to eliminate them by instead imposing the
stronger fixed-determinant-connection trace condition and working 
with a compatible group of gauge transformations. A similar framework is
used in \cite{OTQuaternion,TelemanMonopole}. 
 
We refer to \S \ref{subsec:MonopoleEqns} for a detailed account
of the construction of our holonomy perturbations.  The
large number of technical points involving regularity and uniform estimates
for these perturbations (which still allow us to obtain an Uhlenbeck
compactification) are discussed in the Appendix. 
We fix $r\ge k+1$ and define gauge-equivariant $C^\8$ maps (see
\S \ref{subsec:MonopoleEqns}),
\begin{align}
&\sA_E(X) \ni A\mapsto \vectau\cdot\vecfm(A)
\in L^2_{k+1}(X,\gl(\La^+)\otimes_\RR\so(\su(E))), 
\label{eq:GaugeEquivariantMap}\\
&\sA_E(X) \ni A\mapsto \vecvartheta\cdot\vecfm(A)
\in L^2_{k+1}(X,\Hom(W^+,W^-)\otimes_\CC\fsl(E)), \notag
\end{align}
where $\vectau := (\tau_{j,l,\alpha})$ is a suitably convergent
sequence in $C^r(X,\gl(\La^+))$
and $\vecvartheta := (\vartheta_{j,l,\alpha})$ is a 
suitably convergent sequence in
$C^r(X,\La^1\otimes\CC)$, while $\vecfm(A) := (\fm_{j,l,\alpha}(A))$ is a
sequence in $L^2_{k+1}(X,\su(E))$ of holonomy sections constructed by
extending the method of
\cite{DonOrient,DonPoly}, and 
\begin{align*}
\vectau\cdot\vecfm(A)
&:= 
\sum_{j,l,\alpha}\tau_{j,l,\alpha}\otimes_\RR
\ad(\fm_{j,l,\alpha}(A)), \\
\vecvartheta\cdot\vecfm(A)
&:=
\sum_{j,l,\alpha}\rho(\vartheta_{j,l,\alpha})
\otimes_\CC\fm_{j,l,\alpha}(A).
\end{align*}
We call a point $(A,\Phi)$ in the pre-configuration space of $L^2_k$
pairs $\tsC_{W,E}:=\sA_E\times L^2_k(W^+\otimes E)$ a {\em $\PU(2)$
monopole\/} if it solves
\begin{align}
F_A^+ - (\id+\tau_0\otimes\id_{\su(E)}+\vectau\cdot\vecfm(A))
\rho^{-1}(\Phi\otimes\Phi^*)_{00} &= 0, 
\label{eq:IntroPT} \\
D_A\Phi + \rho(\vartheta_0)\Phi
+ \vecvartheta\cdot\vecfm(A)\Phi &= 0. \notag
\end{align}
We let $M_{W,E}$ be the moduli space 
of solutions cut out of the configuration space of pairs
$\sC_{W,E} := \tsC_{W,E}/\ssG_E$ by the equations
\eqref{eq:IntroPT}. We let $M_{W,E}^{*,0}\subset M_{W,E}$ be the subspace
of pairs $[A,\Phi]$ such that $A$ is irreducible and the section $\Phi$ is
not identically zero.
The sections $\vectau\cdot\vecfm(A)$ and $\vecvartheta\cdot\vecfm(A)$
vanish at reducible connections $A$ by construction; plainly, the terms
in \eqref{eq:IntroPT} involving the perturbations $\vectau\cdot\vecfm(A)$ and
$\vecvartheta\cdot\vecfm(A)$ are zero when $\Phi$ is zero. 

\subsubsection{Uhlenbeck compactness}
The holonomy-perturbation maps in \eqref{eq:GaugeEquivariantMap}
are continuous with respect to the Uhlenbeck
topology (see \S \ref{subsec:DefnUhlCompact}), 
just as are those of \cite{DonPoly}. Suppose $\{A_\beta\}$ is a
sequence in $\sA_E(X)$ which converges 
in the Uhlenbeck topology to a limit $(A,\bx)$ in
$\sA_{E_{-\ell}}(X)\times\Sym^\ell(X)$.
The sections
$\vectau\cdot\vecfm(A_\beta)$ and $\vecvartheta\cdot\vecfm(A_\beta)$ then
converge in $L^2_{k+1}(X)$ to a section $\vectau\cdot\vecfm(A,\bx)$ of
$\gl(\La^+)\otimes_\RR\so(\su(E_{-\ell}))$ and a section
$\vecvartheta\cdot\vecfm(A,\bx)$ of $\Hom(W^+,W^-)\otimes_\CC\fsl(E_{-\ell})$,
respectively.  For each $\ell\ge 0$, the maps of
\eqref{eq:GaugeEquivariantMap} extend continuously to gauge-equivariant
maps
\begin{align}
&\sA_{E_{-\ell}}(X)\times\Sym^\ell(X) 
\to L^2_{k+1}(X,\gl(\La^+)\otimes_\RR\so(\su(E_{-\ell}))), 
\label{eq:GaugeEquivariantExtendedMap}\\
&\sA_{E_{-\ell}}(X)\times\Sym^\ell(X) 
\to L^2_{k+1}(X,\Hom(W^+,W^-)\otimes_\CC\fsl(E_{-\ell})), \notag
\end{align}
given by $(A,\bx)\mapsto \vectau\cdot\vecfm(A,\bx)$ and 
$(A,\bx)\mapsto \vecvartheta\cdot\vecfm(A,\bx)$, 
which are $C^\8$ on each $C^\8$ stratum determined by $\Sym^\ell(X)$.

Our construction of the Uhlenbeck compactification for
$M_{W,E}$ requires us to consider moduli spaces
$$
\bM_{W,E_{-\ell}} \subset \sC_{W,E_{-\ell}}\times\Sym^\ell(X)
$$
of triples $[A,\Phi,\bx]$ given by the zero locus of the 
$\ssG_{E_{-\ell}}$-equivariant map
$$
\fs:\tsC_{W,E_{-\ell}}\times\Sym^\ell(X) \to 
L^2_k(\La^+\otimes\su(E_{-\ell}))\oplus L^2_k(W^-\otimes E_{-\ell})
$$
defined as in \eqref{eq:IntroPT} except using the perturbing sections
$\vectau\cdot\vecfm$ and $\vecvartheta\cdot\vecfm$ in
\eqref{eq:GaugeEquivariantExtendedMap} instead of those in
\eqref{eq:GaugeEquivariantMap}. We call
$\bM_{W,E_{-\ell}}$ a {\em lower-level\/} moduli space if $\ell>0$ and
call $\bM_{W,E_{-0}} = M_{W,E}$ the {\em top\/} or {\em highest level\/}.

In the more familiar case of the unperturbed $\PU(2)$ monopole equations
\eqref{eq:IntroUnpertPT}, the spaces $\bM_{W,E_{-\ell}}$
would simply be products $M_{W,E_{-\ell}}\times \Sym^\ell(X)$.  In general,
though, the spaces $\bM_{W,E_{-\ell}}$ are not products when $\ell > 0$ due
to the slight dependence of the section $\fs(A,\Phi,\bx)$ on the points
$\bx\in\Sym^\ell(X)$ through the perturbations $\vectau\cdot\vecfm$
and $\vecvartheta\cdot\vecfm$. A similar phenomenon is encountered in
\cite[\S 4(iv)--(vi)]{DonPoly} for the case of the extended anti-self-dual
equations, where holonomy perturbations are also employed in order to
achieve transversality.

We define $\barM_{W,E}$ to be the Uhlenbeck closure of $M_{W,E}$ in the
space of ideal $\PU(2)$ monopoles,
$$
\bigcup_{\ell=0}^N \bM_{W,E_{-\ell}} 
\subset 
\bigcup_{\ell=0}^N \left(\sC_{W,E_{-\ell}}\times\Sym^\ell(X)\right),
$$
for any integer $N\ge N_p$, where $N_p$ is a sufficiently large integer to
be specified below.  

\begin{thm}
\label{thm:Compactness}
Let $X$ be a closed, oriented, smooth four-manifold with $C^\8$ Riemannian
metric, \spinc structure $(\rho,W^+,W^-)$ with 
\spinc connection on $W=W^+\oplus W^-$, and a Hermitian two-plane
bundle $E$ with unitary connection on $\det E$.  Then there is a positive
integer $N_p$, depending at most on the curvatures of the fixed connections
on $W$ and $\det E$ together with $c_2(E)$, such that for all $N\ge N_p$
the topological space $\barM_{W,E}$ is compact, second countable,
Hausdorff, and is given by the closure of $M_{W,E}$ in
$\cup_{\ell=0}^{N}\bM_{W,E_{-\ell}}$.
\end{thm}

\begin{rmk}
The existence of an Uhlenbeck compactification for the moduli space of
solutions to the unperturbed $\PU(2)$ monopole equations
\eqref{eq:IntroUnpertPT} was announced by Pidstrigach
\cite{PTCambridge} and an argument was outlined in \cite{PTLocal}.
A similar argument for the equations 
\eqref{eq:IntroUnpertPT} was outlined by Okonek and Teleman in
\cite{OTQuaternion}. Theorem \ref{thm:Compactness} yields the standard
Uhlenbeck compactification for the system \eqref{eq:IntroUnpertPT} and for
the perturbations of \eqref{eq:IntroUnpertPT} described in
\cite{FeehanGenericMetric,TelemanGenericMetric} --- see Remark
\ref{rmk:TransversalityThm}. 
An independent proof of Uhlenbeck compactness for
\eqref{eq:IntroUnpertPT} and certain perturbations of these equations
is given in \cite{TelemanMonopole,TelemanGenericMetric}.
\end{rmk}

\subsubsection{Transversality}
The space $\Sym^\ell(X)$ is smoothly stratified, the strata being
enumerated by partitions of $\ell$. If $\Sigma\subset\Sym^\ell(X)$ is a
smooth stratum, we define
$$
\bM_{W,E_{-\ell}}|_\Sigma
:=\{[A,\Phi,\bx]\in\bM_{W,E_{-\ell}}: \bx\in\Sigma \},
$$
with $\bM_{W,E_{-0}} := M_{W,E}$ when $\ell=0$. We then have the following
transversality result for the equations \eqref{eq:IntroPT}, at least away
from the solutions where the connection is reducible or the spinor vanishes.

\begin{thm}
\label{thm:Transversality}
Let $X$ be a closed, oriented, smooth four-manifold with $C^\8$ Riemannian
metric, \spinc structure $(\rho,W^+,W^-)$ with 
\spinc connection on $W=W^+\oplus W^-$, and a Hermitian
line bundle $\det E$ with unitary connection.  Then there is a
first-category subset of the space of $C^\8$ perturbation parameters such
that the following holds. For each $4$-tuple
$(\tau_0,\vartheta_0,\vectau,\vecvartheta)$ in the complement of this
first-category subset, integer $\ell\ge 0$, and smooth stratum
$\Sigma\subset\Sym^\ell(X)$, the moduli space
$\bM^{*,0}_{W,E_{-\ell}}|_\Sigma(\tau_0,\vartheta_0,\vectau,\vecvartheta)$
is a smooth manifold of the expected dimension,
\begin{align*} 
\dim\bM^{*,0}_{W,E_{-\ell}}|_\Sigma
&= \dim M^{*,0}_{W,E_{-\ell}} + \dim\Sigma 
\\
&= -2p_1(\su(E_{-\ell}))-\threehalf(e(X)+\sigma(X))  + \dim\Sigma\\
&\quad + \half p_1(\su(E_{-\ell}))+\half((c_1(W^+)+c_1(E))^2-\sigma(X))-1,
\end{align*}
where $\det(E_{-\ell})=\det E$ and $c_2(E_{-\ell})=c_2(E)-\ell$.
\end{thm}

\begin{rmk}
\label{rmk:TransversalityThm}
Different approaches to the question of transversality for the $\PU(2)$
monopole equations \eqref{eq:IntroUnpertPT} with generic perturbation
parameters have been considered by Pidstrigach and Tyurin in \cite{PTLocal}
and by Teleman in \cite{TelemanMonopole}.
More recently, a new approach to transversality for
\eqref{eq:IntroUnpertPT} has been discovered independently by the first
author \cite{FeehanGenericMetric} and by Teleman
\cite{TelemanGenericMetric}: the method uses only the perturbations
$(\tau_0,\vartheta_0)$ together with perturbations of the Riemannian metric
$g$ on $X$ and compatible Clifford map $\rho$.
\end{rmk}

A choice of generic Riemannian metric on $X$ ensures that the moduli space
$M_E^{\asd,*}$ of irreducible anti-self-dual connections on $\su(E)$ is
smooth and of the expected dimension \cite{DK,FU}, although the points of
$M_E^{\asd,*}$ need not be regular points of $M^*_{W,E}$ as the
linearization of \eqref{eq:IntroPT} may not be surjective there.
A choice of generic parameter $\tau_0$ ensures that the moduli spaces
$M^{\red,0}_{W,E,L_1}$ of non-zero-section solutions to
\eqref{eq:IntroPT} which are reducible with respect to the splitting
$E=L_1\oplus (\det E)\otimes L_1^*$ are smooth and of the expected
dimension \cite{FL2}. Again, the points  
of $M^{\red,0}_{W,E,L_1}$ need not be regular points of $M^0_{W,E}$ since the
linearization of \eqref{eq:IntroPT} may not be surjective there.

We note that related transversality
and compactness issues have been recently considered in approaches to
defining Gromov-Witten invariants for general symplectic manifolds
\cite{LiTian,RuanGW,Siebert}.

\subsection{Outline}
We indicate how the remainder of our article is organized.
In \S \ref{subsec:ConfigSpace} we prove a slice result for the
configuration space $\sC_{W,E}$ (Proposition
\ref{prop:Slice}) while in \S \ref{subsec:DeformationComplex} we
describe the elliptic deformation complex for \eqref{eq:IntroPT} and
compute the expected dimension of $M_{W,E}$.  We develop the regularity
theory for a generalization of the $\PU(2)$ monopole equations
\eqref{eq:IntroPT} (obtained by allowing additional inhomogeneous terms) in
\S \ref{sec:regularity}: the main technical result there is that an $L^2_1$
solution to an inhomogeneous version of
\eqref{eq:IntroPT} and the Coulomb gauge equation is $C^\8$ (Corollary
\ref{cor:L2_1InhomoReg}). By combining this with the slice result of
Proposition \ref{prop:Slice}, we then show that any $\PU(2)$ monopole
$(A,\Phi)$ in $L^2_k$ is $L^2_{k+1}$ gauge-equivalent to a $\PU(2)$
monopole in $C^\8$ (Proposition \ref{prop:GlobalReg}).    

We establish local estimates for $L^2_1$ solutions to the inhomogeneous
version of \eqref{eq:IntroPT} in \S \ref{subsec:L2_1InhomoLocal} and
\S \ref{subsec:Uhlenbeck}. We use the sharp
$L^2_1$ regularity result of Corollary
\ref{cor:L2_1InhomoReg} in \S \ref{sec:compact} to prove the
removability of point singularities for $\PU(2)$ monopoles (Theorem
\ref{thm:PTRemovSing}). In the sequels \cite{FL3,FL4}, these regularity
results and estimates are needed to prove that $L^2_1$ gluing solutions to
\eqref{eq:IntroPT} are $C^\8$ and to analyse the asymptotic behavior of
Taubes' gluing maps and their differentials near the lower strata of the
Uhlenbeck compactification.

The proof of Theorem \ref{thm:Compactness} relies heavily on both the
regularity theory of \S \ref{sec:regularity} and the fact that solutions
$(A,\Phi)$ to \eqref{eq:IntroPT} satisfy a `universal energy bound', with
constant depending only on the data in the hypotheses of
Theorem \ref{thm:Compactness}. The section $\Phi$ also satisfies a
universal $C^0$ bound: these bounds are the analogues for $\PU(2)$
monopoles of the now well-known {\em a priori} estimates for Seiberg-Witten
monopoles \cite{KrM,MorganSWNotes,Salamon,Witten} and follow, in much the
same way, from the maximum principle and the Bochner-Weitzenb\"ock formula
for $D_A$ provided $k\ge 3$ (see \S \ref{subsec:BW}).  In \S
\ref{sec:compact} we prove the removability of point singularities for
solutions to \eqref{eq:IntroPT} (Theorem
\ref{thm:PTRemovSing}) using our {\em a priori} bounds and regularity results.
While the $\PU(2)$ monopole equations are not
conformally invariant, they are invariant under constant rescalings of the
metric (in the sense of \S
\ref{subsec:ScaleInv}) and, as in the case of anti-self-dual connections,
this scale invariance is exploited in the proof of Theorem
\ref{thm:Compactness}, whose proof is completed in \S
\ref{subsec:SeqCompact}.

Theorem \ref{thm:Transversality} is initially established in \S
\ref{sec:transv} for $C^r$ perturbations for any fixed $3\le r<\8$, in
order to avail of the Sard-Smale theorem for Fredholm maps of Banach manifolds
\cite{Smale}, while in \S
\ref{subsec:SmoothPert} we show that generic $C^\8$ perturbation parameters 
are sufficient for transversality. (See Corollary
\ref{cor:GenericSmoothAnalyticProjection} for the special case $\ell=0$ and
\S \ref{subsubsec:LowerStratumTransv} for its extension to the general case
$\ell\ge 0$.)

As we shall explain in \S \ref{sec:transv}, our proof of Theorem
\ref{thm:Transversality} ultimately hinges on the fact that if $(A,\Phi)$
is a $\PU(2)$ monopole and $A$ is reducible on a non-empty open subset
containing the support of the 
holonomy perturbations depending on $A$, then $A$ is
reducible over $X$. The proof of this result (Theorem
\ref{thm:LocalToGlobalReducible}) occupies \S
\ref{subsec:LocalToGlobal}; the result follows from the 
Agmon-Nirenberg unique continuation theorem after the system
\eqref{eq:IntroPT} has been transformed into an ordinary differential
equation for a one-parameter family of pairs over $S^3$. The unique
continuation property holds for both the perturbed $\PU(2)$ monopole equations 
\eqref{eq:IntroPT}, when the initial open subset of $X$ contains
the balls in $X$ supporting holonomy perturbations, and the 
unperturbed equations
\eqref{eq:IntroUnpertPT} for any initial open subset. 

\subsection{Other approaches to transversality}
As we noted in Remark \ref{rmk:TransversalityThm}, transversality for the
$\PU(2)$ monopole equations \eqref{eq:IntroUnpertPT}
has also been proved very recently using only the
perturbations $(\tau_0,\vartheta_0)$, together with perturbations of the
Riemannian metric $g$ on $X$ and compatible Clifford map $\rho$
\cite{FeehanGenericMetric,TelemanGenericMetric}. The transversality proof
given in \cite{FeehanGenericMetric} is considerably more delicate
than the method we employ in the present article.  When using holonomy
perturbations, the main technical difficulties are due to the
noncompactness of the moduli space of $\PU(2)$ monopoles and the core
transversality argument itself is more straightforward, whereas in
\cite{FeehanGenericMetric} the situation is entirely reversed. To place
these various transversality results in a suitable context, we briefly
discuss some other approaches to transversality, both for the $\PU(2)$
monopole equations and the closely related `\spinc-ASD' equations of
\cite{PTDirac}, which Pidstrigach and Tyurin used to define spin polynomial
invariants. 

Suppose $(A,\Phi)$ is a solution to either the
equations \eqref{eq:IntroUnpertPT}, the equations
\eqref{eq:IntroUnpertPT} with perturbations $(\tau_0,\vartheta_0)$, or the
holonomy-perturbed equations \eqref{eq:IntroPT}.  If $A$ is 
reducible, then $\Phi$ has rank less than or equal to one
(see Lemma \ref{lem:ReducibleStabilizers}).  However, as
observed by Teleman \cite{OTSurvey,TelemanMonopole}, if $\Phi$ is rank one
then $A$ is not necessarily reducible and he describes a simple
counterexample for the equations \eqref{eq:IntroUnpertPT} when $X$ is a
K\"ahler manifold with its canonical \spinc structure.

It is not too difficult to prove that $M_{W,E}^{*,0}$ is a smooth manifold
of the expected dimension away from the locus of irreducible, rank-one
solutions using the perturbation parameters $(\tau_0,\vartheta_0)$ alone.
However, as irreducible, rank-one solutions to
\eqref{eq:IntroPT} could be present in $M_{W,E}^{*,0}$, it appears impossible
to prove that the entire space $M_{W,E}^{*,0}$ is a smooth manifold of the
expected dimension using only the parameters $(\tau_0,\vartheta_0)$. A
similar problem arises in the proof of transversality for the \spinc-ASD
equations given in
\cite[Proposition I.3.5]{PTDirac}; a version of these equations can be
obtained from the equations 
\eqref{eq:IntroUnpertPT} by omitting the quadratic term
$\rho^{-1}(\Phi\otimes\Phi^*)_{00}$. In the proof of \cite[Proposition
I.3.5]{PTDirac} it is claimed that if $D_A\Phi=0$ and $\Phi$ is rank one,
then $A$ is reducible [p. 277]: Teleman's counterexample shows that this
claim is incorrect and he points out an error in their argument
\cite{OTSurvey,TelemanMonopole}. On the other hand,
the possible presence or absence of irreducible, rank-one solutions to the
$\PU(2)$ monopole equations makes no difference to the transversality
argument we describe in \S \ref{sec:transv} using holonomy perturbations,
as these perturbations are strong enough to yield transversality without a
separate analysis of the locus of irreducible, rank-one solutions. 

\subsection{Applications}
In \cite{FL2} we discuss the singularities of the moduli space $M_{W,E}$.
We introduce cohomology classes and geometric representatives on
$M_{W,E}^*$, construct the links of the top-level singular strata of
anti-self-dual and Seiberg-Witten moduli spaces, compute the Chern
characters of their normal bundles in suitably defined ambient manifolds,
and define the orientations of these links. We thus obtain a relation
between Donaldson and Seiberg-Witten invariants when the reducible
solutions appear only in the top level of the Uhlenbeck compactification.
In \cite{FL3,FL4}, we develop the gluing theory for $\PU(2)$ monopoles:
this is used to construct links of the lower-level Seiberg-Witten moduli
spaces, to show that the pairing of the cohomology classes with this link
is well-defined, and to eliminate the requirement of
\cite{FL2} that the reducible solutions appear only in the top level. A
survey of some of the results contained in the present article and its
sequels \cite{FL2,FL3,FL4} is provided in \cite{FLGeorgia}.

\subsubsection*{Acknowledgements}
The first author is especially grateful to Peter Kronheimer, Tom Mrowka,
and Cliff Taubes for helpful conversations and for their enthusiasm and
encouragement. He is particularly grateful to Tom Mrowka for his
corrections to the manuscript.  He would also like to thank the Harvard
Mathematics Department and the National Science Foundation for their
generous support. The second author is very grateful to Selman Akbulut, Ron
Fintushel, and Tom Parker for helpful conversations and the Mathematics
Departments at Harvard and Michigan State Universities for their
hospitality and support. We warmly thank the referee for his corrections,
careful reading of the manuscript, and thoughtful comments.


\section{The $\PU(2)$ monopole equations}
\label{sec:prelim}
We lay out a framework for gauge theory for pairs in
\S \ref{subsec:FuncSpace} and \S \ref{subsec:ConfigSpace}:
In \S \ref{subsec:FuncSpace} we define the Hilbert
spaces of pairs and gauge transformations, while in \S
\ref{subsec:ConfigSpace} we establish the main slice result we require for the
configuration space of pairs modulo gauge transformations, paying
particular attention to the structure near `reducibles'. The $\PU(2)$
monopole equations and their holonomy perturbations are introduced in \S
\ref{subsec:MonopoleEqns}. The role of the perturbations in obtaining
transversality results will be explained in \S
\ref{sec:transv}. To preserve continuity a detailed discussion of the 
technical points which need to be addressed when using holonomy
perturbations in the present context is deferred to the Appendix.
The moduli space
and the elliptic deformation complex for solutions to the $\PU(2)$
monopole equations are described in \S
\ref{subsec:DeformationComplex}.  

\subsection{Configuration spaces of connections and pairs}
\label{subsec:FuncSpace}
In this section we define the spaces of connections,
pairs, gauge transformations, and
configuration spaces we will use throughout our work.

\subsubsection{Connections on Hermitian two-plane bundles}
We consider Hermitian two-plane bundles $E$ over $X$ whose determinant line
bundles $\det E$ are isomorphic to a fixed Hermitian line bundle over $X$
endowed with a fixed $C^\8$, unitary connection $A_e$. {\em The Hermitian
line bundle over $X$ and its unitary connection $A_e$ are fixed for the
remainder of this article.\/}

Let $k\ge 2$ be an integer and let $\sA_E$ be the space of $L^2_k$
connections $A$ on the $\U(2)$ bundle $E$ all inducing the fixed
determinant connection $A_e$ on $\det E$.  Equivalently, following \cite[\S
2(i)]{KMStructure}, we may view $\sA_E$ as the space of $L^2_k$ connections
$A$ on the $\PU(2)=\SO(3)$ bundle $\su(E)$.  We shall pass back and
forth between these viewpoints, via the fixed connection on $\det E$,
relying on the context to make the distinction clear.  
Explictly, if $A$ is a unitary connection on $E$ and $s\in\Om^0(X,E)$, then
$\cov_A s\in\Om^1(E)$. For $\zeta\in\Om^0(X,\fu(E))$, then
$\cov_A\zeta\in\Om^1(X,\fu(E))$ is determined by
$$
(\cov_A\zeta)s = \cov_A (\zeta s) - \zeta(\cov_A s),
$$
so, if $A^\tau\in\Om^1(U,\fu(2))$ 
denotes the local connection matrix defined by a
choice of local frame for $E$ over an open subset $U\subset X$, then
\begin{align*}
\cov_As &= ds + A^\tau s, \\
\cov_A\zeta &= d\zeta + [A^\tau,\zeta] = d\zeta + (\ad A^\tau)\zeta.
\end{align*}
Similarly, the connection $A$ on $E$ induces connections on $\su(E)$ and
$\det E$. For example,
if $\xi\in\Om^0(X,\su(E))$, then $\cov_A\xi\in\Om^1(X,\su(E))$ is also
given locally by 
$$
\cov_A\xi = d\xi + [A^\tau,\xi] = d\xi + [(A^\tau)_0,\xi]
 = d\zeta_0 + (\ad A^\tau)\zeta_0,
$$
while if $\la\in\Om^0(X,\det E)$ then $\cov_A\la\in\Om^1(X,\det E)$
is given by 
$$
\cov_A \la = d\la + (\tr A^\tau)\la. 
$$
The above local connection matrices are related by
$$
A^\tau = (A^\tau)_0 + \half (\tr A^\tau)\,\id_{\CC^2}
\in \Om^1(U,\fu(2)),
$$
where $(A^\tau)_0 = A^\tau - \half (\tr A^\tau)\,\id_{\CC^2}
\subset\Om^1(U,\su(2))$ is the {\em traceless\/} component of $A^\tau$ while
$\tr A^\tau\subset\Om^1(U,\fu(1))=\Om^1(U,i\RR)$ is the induced local
connection form for $\det E$ given by the {\em trace\/} component of
$A^\tau$.  Note that $(A^\tau)_0\in \Om^1(U,\su(2))$ and that $\ad
(A^\tau)_0 = \ad A^\tau \in \Om^1(U,\so(3))$ is the induced local
connection matrix for the $\SO(3)$ bundle $\su(E)$, where we use the
standard identification $\so(3) = \ad(\su(2))$.

Conversely, given a unitary connection $A_e$ for $\det E$ and a Riemannian
connection $\hatA$ for $\su(E)$, we obtain a unitary connection for $E$
given in terms of local connection matrices by
$$
A^\tau = \ad^{-1}(\hatA^\tau) + \half A_e^\tau\,\id_{\CC^2}
\in \Om^1(U,\fu(2)).
$$
The curvatures of these connection matrices are related by
$$
F_{A^\tau} = \ad^{-1}(F_{\hatA^\tau}) 
+ \half F_{A_e^\tau}\,\id_{\CC^2} \in \Om^2(U,\fu(2)),
$$
with $F_{\hatA^\tau}\in \Om^2(U,\so(3))$ 
and $F_{A_e^\tau}\in \Om^2(U,\fu(1))$. Thus,
$(F_{A^\tau})_0 = \ad^{-1}(F_{\hatA^\tau})$ and $\tr F_{A^\tau} 
= F_{A_e^\tau} = F(\tr A^\tau)$. Therefore, globally we have
$(F_A)_0 = \ad^{-1}F_{\hatA}\in\Om^2(X,\su(E))$ and 
$\tr F_A = F_{A_e}\in\Om^2(X,\fu(1))=\Om^2(X,i\RR)$, so that
$$
F_A = \ad^{-1}(F_{\hatA}) + \half F_{A_e}\,\id_E \in \Om^2(X,\fu(E)),
$$
as $\fu(E) = \su(E)\oplus (i\RR)\,{\id_E}$. 

When we are not explicitly discussing connections which
are reducible with respect to some splitting of $E$, it is generally more
convenient to view our connections as being defined on $\su(E)$ rather than
$E$ as we can then avoid explicit mention of the otherwise unimportant
choice of fixed connection on $\det E$. Of course, these viewpoints are
equivalent via the choice of this fixed determinant connection. The
isomorphism $\ad:\su(E)\to\so(\su(E))$ will remain implicit when no
confusion can arise so that, given a connection $A$ on $\su(E)$, we write
its curvature as $F_A\in\Om^2(\su(E))$ and associated deformation complexes
on $\Om^\bullet(\su(E))$ rather than $F_A\in\Om^2(\so(\su(E)))$, 
with deformation complexes on $\Om^\bullet(\so(\su(E)))$
(see \cite[Chapter 10]{FU}, for comparison). 

\subsubsection{Spin${}^c$ structures}
The minimal, axiomatic approach to the definition of \spinc structures and
the Dirac operator employed by Kronheimer-Mrowka \cite{KrMContact} and
Mrowka-Ozsv\'ath-Yu \cite{MrowkaOzsvathYu} is extremely useful for our
purposes, so this is the method we shall follow here.

Recall that a real-linear map $\rho_+:T^*X\to\Hom(W^+,W^-)$ defines a
Clifford map $\rho:\Lambda^\bullet(T^*X)\otimes\CC\to\End(W)$, with $W :=
W^+\oplus W^-$, if and only if
\cite{LM,Salamon}
\begin{equation}
\rho_+(\alpha)^\dagger\rho_+(\alpha) = g(\alpha,\alpha)\id_{W^+},
\qquad \alpha\in C^\8(T^*X),
\label{eq:MinimalCliffordMap}
\end{equation}
where $g$ denotes the Riemannian metric on $T^*X$.
The real-linear map $\rho:T^*X\to\End_\CC(W^+\oplus W^-)$ is obtained by
defining a real-linear map 
$$
\rho_-:T^*X\to\Hom(W^-,W^+),
\qquad
\alpha \mapsto \rho_-(\alpha) := -\rho_+(\alpha)^\dagger,
$$
and setting
\begin{equation}
\rho(\alpha) 
:=
\begin{pmatrix}
0 & \rho_-(\alpha) \\ \rho_+(\alpha) & 0
\end{pmatrix}.
\label{eq:DefnFullCliffordMap}
\end{equation}
This extends to a linear map
$\rho:\Lambda^\bullet(T^*X)\otimes\CC\to\End(W^+\oplus W^-)$ and
satisfies
\begin{equation}
\rho(\alpha)^\dagger = - \rho(\alpha)
\quad\text{and}\quad
\rho(\alpha)^\dagger\rho(\alpha) = g(\alpha,\alpha)\id_W,
\qquad \alpha\in C^\8(T^*X).
\label{eq:CliffordMap}
\end{equation}

A unitary connection $\cov$ on $W$ and $g$-compatible Clifford map $\rho$
induces a unique $\SO(4)$ connection $\cov^g$ on $T^*X$ by requiring that
\begin{equation}
[\cov_\eta,\rho(\alpha)] = -\rho(\cov_\eta^g\alpha),
\label{eq:RhoCompatibility}
\end{equation}
for all $\eta \in C^\8(TX)$ and $\alpha \in \Omega^1(X,\RR)$. 
(As in \cite{MrowkaOzsvathYu}, we use a minus sign on the right-hand side
of \eqref{eq:RhoCompatibility} since our basic Clifford map and Dirac
operator are defined using Clifford 
multiplication by one-forms rather than the more
traditional tangent vectors of \cite{LM} or \cite{Salamon}.) The unitary
connection $\cov$ on $W$ uniquely determines a unitary connection on $\det W^+
\simeq \det W^-$ in the standard way \cite{Kobayashi}. Conversely, the
preceding data uniquely determines a unitary connection $\cov$ on $W$. The
connection $\cov$ on $W$ is called a {\em \spinc\/} connection if the
connection $\cov^g$ on $T^*X$ is also torsion free, that is, if $\cov^g$ is the
Levi-Civita connection for the metric $g$.

Given a unitary connection $A$ on an auxiliary Hermitian two-plane bundle
$E$, we let $\cov_A$ denote the induced unitary connection on $W\otimes E$.
The corresponding Dirac operator $D_A$ is defined by 
$$
D_A := \sum_{\mu=1}^4\rho(v^\mu)\cov_{A,v_\mu},
$$
where $\{v_\mu\}$ is a local frame for $TX$ and $\{v^\mu\}$ is the dual
frame for $T^*X$ defined by $v^\mu(v_\mu) = \delta_{\mu\nu}$.

\subsubsection{Pre-configuration spaces and automorphism groups}
We denote $\Lambda^\bullet := \Lambda^\bullet(T^*X)$.
For any fixed $L^2_k$ connection $A_0\in \sA_E$ we therefore write
$$
\sA_E = A_0 + L^2_k(\La^1\otimes\su(E)),
$$
for either fixed-determinant, unitary connections on $E$ or
determinant-one, orthogonal connections on $\su(E)$. Our {\em
pre-configuration space\/} of pairs is given by
$$
\tsC_{W,E}
:= 
\sA_E\times L^2_k(W^+\otimes E),
$$ 
and for any fixed $L^2_k$ pair $(A_0,\Phi_0)\in \tsC_{W,E}$, we have
\begin{equation}
\tsC_{W,E}
= 
(A_0,\Phi_0) + L^2_k(\La^1\otimes\su(E))\oplus L^2_k(W^+\otimes E),
\label{eq:TangentPairSpace}
\end{equation}
for the cases of either fixed-determinant, unitary connections on $E$ or
determinant-one, orthogonal connections on $\su(E)$. 

Given any $C^\8$ connection $A_0\in\sA_E$, our Sobolev norms are defined in
the usual way: for example, if $a\in\Om^1(\su(E))$, we write
$$
\|a\|_{L^p_{\ell,A_0}(X)} 
:= \left(\sum_{j=0}^\ell\|\cov_{A_0}^j a\|_{L^p(X)}^p\right)^{1/p},
$$
and if $(a,\phi)\in \Om^1(\su(E))\oplus\Om^0(W^+\otimes E)$,
we write 
$$
\|(a,\phi)\|_{L^p_{\ell,A_0}(X)} 
:= \left(\|a\|_{L^p_{\ell,A_0}(X)}^2 
+ \|\phi\|_{L^p_{\ell,A_0}(X)}^p\right)^{1/p}, 
$$
for any $1\le p\le \8$ and integer $\ell\ge 0$.

For convenience, we let $i\RR_Z = (i\RR)\,\id_{\CC^2}\subset\fu(2)$ denote
the center of the Lie algebra $\fu(2)$ and let
$S^1_Z=\exp(i\RR_Z)\subset\U(2)$ denote the center $Z(\U(2))$ of the Lie
group $\U(2)$ given by
\begin{equation}
i\RR_Z=\left\{\begin{pmatrix} i\theta & 0 \\ 0 & i\theta\end{pmatrix}
: \theta\in\RR\right\}
\quad\text{and}\quad
S^1_Z = \left\{\begin{pmatrix}e^{i\theta} &0\\ 0 &e^{i\theta}\end{pmatrix}
: \theta\in\RR\right\}. 
\label{eq:DefnS1ZandRZ}
\end{equation}
Following \cite[\S 2]{KMStructure} we consider the induced action for
connections on $\su(E)$ by the group $\sG_E$ of determinant-one, unitary
$L^2_{k+1}$ automorphisms of $E$ rather than the group $\sG_{\su(E)}$ of
determinant-one, orthogonal $L^2_{k+1}$ automorphisms of $\su(E)$ even to
define quotient spaces of connections on $\su(E)$. The reasons for this
choice are explained further in \cite{FL2}; see also \S
\ref{subsec:GaugeGroup}. We have
$$
\sG_E := \{u\in L^2_{k+1}(\gl(E)): u^\dagger u = \id_E\text{ and }
\det u = 1 \text{ a.e.}\},
$$
and 
$$
\ssG_E := S_Z^1\times_{\{\pm\id_E\}}\sG_E.
$$
One can define the action of $\ssG_E$ on $\sA_E$ either by push-forward or
pullback and one obtains equivalent quotient spaces in either case, but the
choice does affect the orientation of moduli spaces so we specify the
action here: for $u\in\ssG_E$ and $(A,\Phi)\in
\tsC_{W,E}$, the action of $\ssG_E$ is defined
by
$$
u(A,\Phi) := (u(A),u\Phi) = (u_*A,u\Phi) = ((u^{-1})^*A,u\Phi).
$$
The push-forward action, $u(A)=u_*A$, agrees with the conventions of 
\cite{BradlowDask,DK}. For the associated covariant derivative
$\cov_A$ on $E$ and $u\in\Aut E$, we have
$$
\cov_{u(A)} = u\circ\cov_A\circ u^{-1},
$$
so that $u(A) = A - (\cov_Au)u^{-1}$.  In terms of a local connection
matrix $A^\tau\in\Om^1(U,\fu(2))$, this gives $u(A^\tau) = uA^\tau u^{-1} -
(du)u^{-1}$, where we use $u$ to denote both the gauge transformation and
the element of $\Om^0(U,\U(2))$.

In order to define quotients by the action of $\ssG_E$, we need to
choose $k\ge 2$, so gauge transformations are at least continuous. 
The proof of the following proposition is a standard application of the
Sobolev multiplication theorem: see \cite[Props. A.2, A.3, A.9]{FU}.

\begin{prop}\label{prop:GaugeGroup}
Let $X$ be a closed Riemannian four-manifold, let $E$ be a Hermitian vector
bundle over $X$, and let $k\ge 2$ be an integer. Then the following hold.
\begin{enumerate}
\item The space $\ssG_E$ is a Hilbert Lie group with
Lie algebra $T_{\id}\ssG_E=L^2_{k+1}(\su(E))\oplus i\RR_Z$;
\item The action of $\ssG_E$ on $\tsC_{W,E}$ is smooth;
\item For $(A,\Phi)\in \tsC_{W,E}$,
the differential, at the identity $\id_E\in \ssG_E$, 
of the map $\ssG_E\to \tsC_{W,E}$
given by $u\mapsto u(A,\Phi)=(A-(d_Au)u^{-1},u\Phi)$ is 
$\zeta\mapsto -d_{A,\Phi}^0\zeta := (-d_A\zeta,\zeta\Phi)$ as a map
$$
L^2_{k+1}(\su(E))\oplus i\RR_Z
\to 
L^2_k(\La^1\otimes\su(E)) \oplus L^2_k(W^+\otimes E).
$$
\end{enumerate}
\end{prop}

We denote $\sB_E = \sA_E/\sG_E = \sA_E/\ssG_E$.  The {\em configuration
space of pairs\/} is given by $\sC_{W,E} := \tsC_{W,E}/\ssG_E$, that is,
$$
\sC_{W,E} 
:= \left(\sA_E\times 
L^2_k(W^+\otimes E)\right)/\ssG_E
=\sA_E\times_{\ssG_E}L^2_k(W^+\otimes E),
$$
and is endowed with the quotient $L^2_k$ topology. 

If $\pi$ denotes the projection $\tsC_{W,E}\to
\sC_{W,E}$, 
a base for the quotient $L^2_k$ topology of $\sC_{W,E}$ is
given by open subsets $\pi^{-1}\pi(B_{(A,\Phi)}(\eps))\supset
B_{(A,\Phi)}(\eps)$, where $B_{(A,\Phi)}(\eps)$ is the $L^2_k$-ball of
radius $\eps$ and center $(A,\Phi)$ in
$\tsC_{W,E}$ given by
$\{(A_1,\Phi_1)\in \tsC_{W,E}:
\|(A_1,\Phi_1)-(A,\Phi)\|_{L^2_{k,A}} < \eps\}$.

\subsubsection{Stabilizers}
The space $\sC_{W,E}$ is not a manifold --- it has singularities at points
$[A,\Phi]$ with non-trivial stabilizer $\Stab_{A,\Phi}\subset \ssG_E$.
Recall that the stabilizer subgroup (of the group of bundle automorphisms)
for a connection on a $G$ bundle always contains the center $Z(G)\subset G$
\cite[\S 4.2.2]{DK}. We let $S_Z^1\subset \ssG_E$ denote the constant,
central automorphisms of $E$. 

\begin{defn}
Suppose $(A,\Phi)$ is a point in $\tsC_{W,E}$. The stabilizer
$\Stab_{A,\Phi}\subset\ssG_E$ of $(A,\Phi)$ in $\ssG_E$ is 
given by $\{\ga\in \ssG_E:\ga(A,\Phi) = (A,\Phi)\}$, while the
stabilizers of $A$ and $\Phi$ are denoted by $\Stab_A$ and
$\Stab_\Phi$, respectively. 
\begin{enumerate}
\item The point $(A,\Phi)$ is a {\em zero-section pair\/} 
if $\Phi\equiv 0$;
\item The point $(A,\Phi)$ is an {\em irreducible pair\/} if the connection
$A$ has minimal stabilizer $\Stab_A = Z(\U(2))=S_Z^1\subset\ssG_E$
and is {\em reducible} otherwise, that is, if $\Stab_A \supsetneq S_Z^1$.
\end{enumerate}
The point $(A,\Phi)$ is a {\em reducible, zero-section pair\/} if
$(A,\Phi)$ is both a reducible and zero-section pair.
\end{defn}

As usual, the stabilizer subgroup $\Stab_A\subset\ssG_E$ may be identified
with a closed subgroup of $\U(E_{x_0})\simeq\U(2)$ for any point $x_0\in X$ by
parallel translation with respect to the connection $A$ \cite[\S 4.2.2]{DK},
\cite{MorganGTNotes}. The following lemma implies that the stabilizer
in $\Aut E$ of a unitary connection on $E$ coincides with its stabilizer on
$\ssG_E$. 

\begin{lem}\label{lem:ConstDet}
Let $E$ be a Hermitian two-plane bundle over a connected four-manifold
$X$, let $A$ be a unitary connection on $E$, and let $u$ be a unitary
automorphism of $E$ such that $u\in\Stab_A\subset\Aut E$.  Then $\det
u:X\to S^1$ is a constant map.
\end{lem}

\begin{pf}
The gauge transformation $u$ may be viewed as an $\Ad$-equivariant map
$u:P\to\U(2)$, where $E=P\times_{\U(2)}\CC^2$,
so $\det u$ may be viewed as a map from $P$ to $S^1$.
Since $\det u$ is constant on the fibers of $P$, it descends to a map on $X$.
Over a small enough open set $U\subset X$, we may
write $u=\exp\zeta$, 
where $\zeta\in\Om^0(\ad P)$ and $\ad P :=
P\times_{\ad}\fg \simeq \fu(E)$, with $G=\U(2)$ and
$\fg=\fu(2)$. Differentiating the action of $u$ on $A$ we see that
$\cov_A\zeta = 0$, so $\zeta$ is covariantly constant over $U$. Then
$\det u = \exp(\tr\zeta)$ and $\cov(\det u) = \tr(\cov_A\zeta) = 0$,
so $\det u$ is constant on $U$ and thus on $X$, since $X$ is connected. 
\end{pf}

For a pair $(A,\Phi)$ on $(E,W^+\otimes
E)$ or $(\su(E),W^+\otimes E)$, the Lie algebras of $\Stab_A$ and
$\Stab_{A,\Phi}$ are given by
\begin{align*}
H_A^0 &:= \Ker\{d_A:L^2_{k+1}(\su(E))\oplus i\RR_Z
\to L^2_k(\La^1\otimes\su(E))\}, \\
H_{A,\Phi}^0 &:= \Ker\{d_{A,\Phi}^0:L^2_{k+1}(\su(E))\oplus i\RR_Z\to 
L^2_k(\La^1\otimes\su(E))\oplus L^2_k(\La^1\otimes W^+\otimes E)\}.
\end{align*}
The first identification is standard (see \cite[\S 4.2.2]{DK},
\cite[Chapter 3]{FU}, or \cite{MorganGTNotes}) and the second follows by the
same argument.   

{}From \cite[\S 4.2.2]{DK} and \cite[Theorem 10.8]{FU} we have the following
characterization of the stabilizer $\Stab_B\subset
\sG_V = \Aut V$ of a connection $B$ on an $\SO(3)$ bundle $V$. Note that the
four-manifold $X$ below need not be 
simply-connected or closed and that $H_B^0 =
\Ker\{d_B:L^2_{k+1}(\so(V))\to L^2_k(\La^1\otimes\so(V))\}$. 

\begin{lem}\label{lem:AdConnStab}
Let $B$ be an orthogonal connection on an $\SO(3)$ bundle $V$ over a
connected smooth four-manifold $X$.  The following are
equivalent provided $B$ is not flat:
\begin{enumerate}
\item $\Stab_B\simeq \SO(2)\simeq S^1$;
\item $H_B^0\ne 0$;
\item $B$ is reducible with respect to 
an orthogonal splitting $V =
N\oplus \underline{\RR}$, where $N$ is a complex line bundle over $X$ and
$\underline{\RR} = X\times\RR$ (that is, $B$ reduces to an $\SO(2)$
connection); 
\item $\Stab_B\ne \{\id_V\}$.
\end{enumerate}
Finally, $B$ is flat if $\Stab_B=\SO(3)$.
\end{lem}

\begin{rmk}
Note that the `twisted reducible' (that is, locally reducible) connections
on $\su(E)$ discussed by Kronheimer and Mrowka in \cite[\S
2(i)]{KMStructure} are globally irreducible in the sense that they have
minimal stabilizer $\{\id_{\su(E)}\}$. The condition that $B$ is not
flat is used in the proof that (4) implies (1) \cite[p. 48]{FU}: The
stabilizer $\Stab_B$ is isomorphic to the centralizer of the
holonomy group $\Hol_B$ so, if $\Stab_B\supsetneq S^1$,
then $\Hol_B\subsetneq S^1$ and hence $\Hol_B$ is
discrete. Since the Lie algebra of $\Hol_B$ vanishes, the
Ambrose-Singer holonomy theorem implies that the connection $B$ is
flat. In particular, we only need $X$ to be a connected four-manifold for
these arguments to hold.
\end{rmk} 

As customary, we call an $\SO(3)$ connection 
$B$ on $V$ {\em irreducible} if $\Stab_B =
\{\id_V\}$. We say that a unitary connection $A$ on $E$ is {\em
projectively flat\/} if the induced connection $A^{\ad}$ on $\su(E)$ is
flat. By modifying the proofs of Theorem 3.1 in \cite{FU} and Proposition
II.8.10 \cite{Lawson}, we obtain the following relation between (i) the
stabilizers $\Stab_A$ in $\ssG_E$ of $\U(2)$ connections $A$ on $E$ or
their induced $\SO(3)$ connections $A^{\ad}$ on $\su(E)$ and (ii) the
stabilizers $\Stab_{A^{\ad}}$ in $\sG_{\su(E)}$ of $\SO(3)$ connections
$A^{\ad}$ on $\su(E)$.

\begin{lem}\label{lem:ConnStab}
Let $A$ be a $\U(2)$ connection on a Hermitian two-plane bundle $E$ over a
smooth, connected four-manifold $X$ and let $A^{\ad}$ be the induced
$\SO(3)$ connection on the bundle $\su(E)$.  Then the following
are equivalent:
\begin{enumerate}
\item $\Stab_{A^{\ad}}\simeq S^1$;
\item $\Stab_A\simeq S^1\times S_Z^1$;
\item $H_A^0 \simeq i\RR\oplus i\RR_Z$;
\item $H_{A,\Phi}^0\supsetneq i\RR_Z$;
\item $A$ is reducible with respect to an orthogonal splitting $E=L_1\oplus
L_2$, where $L_1$ and $L_2$ are complex line bundles over $X$ (that is, $A$
reduces to a $T^2 = S^1\times S^1$ connection); 
\item $\Stab_A\supsetneq S_Z^1$.
\end{enumerate}
Finally, the connection $A$ on $E$ is projectively flat (or the connection
$A^{\ad}$ on $\su(E)$ is flat) if $\Stab_A=\U(2)$.
\end{lem}

\begin{pf}
$(1)\implies (2)$: Let $\pi:\U(2)\to\SO(3)=\U(2)/S_Z^1$ be the
projection $u\mapsto \pm(\det u)^{-1/2}u$. Then $\pi^{-1}(\Stab_{A^{\ad}})
\simeq \Stab_{A^{\ad}}\times S_Z^1$ and so $\Stab_{A^{\ad}}\simeq S^1$
implies that $\Stab_A\simeq S^1\times S_Z^1$.

$(2)\implies(3)$: Immediate from the identification of $H_A^0$ as the Lie
algebra of $\Stab_A$.

$(3)\implies(4)$: Trivial.

$(4)\implies (5)$: We argue as in the proof of Theorem 3.1 in
\cite{FU}. Note that
$\Ker\{d_A:L^2_{k+1}(X,\fu(E))\to \Om^1(X,\fu(E))\}$ is isomorphic to
\begin{align*}
&\Ker\{d_A:L^2_{k+1}(X,\su(E))\} 
\oplus \Ker\{d:L^2_{k+1}(X,i\RR_Z)\} \\
&\qquad =
\Ker\{d_{A^{\ad}}:L^2_{k+1}(X,\su(E))\}\oplus i\RR_Z. 
\end{align*}
Given (4), we may choose $0\not\equiv\zeta\in \Ker\{d_A:L^2_{k+1}
\Om^0(X,\su(E))\}$. The pointwise traceless, skew-Hermitian endomorphism
$\zeta$ of $E$ has eigenvalues $i\la_1,i\la_2 = \pm i\la$ for some
$0\not\equiv\la\in L^2_{k+1}(X,\RR)$. On the open subset of $X$ where
$\la\ne 0$, let $\xi_j$, $j=1,2$, be $L^2_{k+1}$ eigenvectors of $\zeta$ such
that $\zeta(\xi_j) = i\la_j\xi_j$ and $\langle\xi_j,\xi_k\rangle =
\de_{jk}$. Just as in
\cite[p. 47]{FU}, we find that $\la$ is constant, the eigenvectors $\xi_j$ are
globally defined and $d_A\xi_j = 0$ for
$j=1,2$. In particular, we have an orthogonal splitting $E=L_1\oplus L_2$,
where $\xi_j\in L^2_{k+1}(X,E)$ is a section of $L_j$  
and $d_A = d_{A_1}\oplus d_{A_2}$ with respect to
this splitting, so this gives (5).

$(5)\implies(6)$: Given (5), the connection $A$ reduces to $A_1\oplus A_2$ with
respect to the splitting $E=L_1\oplus L_2$ and so has stabilizer 
$S_{L_1}^1\times S_{L_2}^1$, where we identify the constant maps in
$\Map^{k+1}(X,S^1)$ with $S_{L_i}^1$ for $i=1,2$. Since $S_{L_1}^1\times
S_{L_2}^1 \simeq S^1\times S_Z^1$, this gives (2).

$(6)\implies (1)$: Given $u\in\Stab_A$ and $u\notin S_Z^1$, we obtain $u_1
= \pm(\det u)^{-1/2}u\in\Stab_{A^{\ad}}$. (By Lemma
\ref{lem:ConstDet}, the determinant $\det u$ is a constant map.) 
If $u_1 = \id_{\su(E)}$ then we would have $u_1
= \pm\id_E$ and $u=\pm(\det u)^{1/2}u_1\in S_Z^1$,  so $u_1\ne \id_{\su(E)}$.
Therefore, $\Stab_{A^{\ad}}\supsetneq \{\id_{\su(E)}\}$ and
Lemma \ref{lem:AdConnStab} implies that $\Stab_{A^{\ad}}\simeq S^1$.  

Lastly, $A$ is projectively flat if
$\Stab_{A^{\ad}}=\SO(3)$ (by Lemma
\ref{lem:AdConnStab}) and $\Stab_{A^{\ad}}=\SO(3)\simeq\U(2)/S_Z^1$ if and
only if $\Stab_A=\U(2)$.  
\end{pf}

The equivalence of (1) and (2) above can alternatively be seen by noting
that $H_{A^{\ad}}^0\simeq i\RR$ if and only if $H_A^0 \simeq i\RR\oplus
i\RR_Z$, so $\Stab_{A^{\ad}} \simeq S^1$ if and only if $\Stab_A \simeq
S^1\times S_Z^1$.
  
{}From Lemmas \ref{lem:AdConnStab} and \ref{lem:ConnStab} we see that $A$
is an irreducible $\U(2)$ connection if and only if $A^{\ad}$ is an
irreducible $\SO(3)$ connection.  Note that $\Stab_{A,\Phi} =
\Stab_A\cap\Stab_\Phi$. The stabilizer $\Stab_{A,\Phi}$ of a zero-section
pair $(A,\Phi)$ contains $S^1_Z$.  The stabilizer $\Stab_{A,\Phi}$ of a
reducible pair $(A,\Phi)$ need not contain the stabilizer $\Stab_A$
since $\Stab_A$ may not fix $\Phi$. 

We write $\tsC_{W,E}^*$ (respectively, $\tsC^0_{W,E}$) for the complement
of the reducible pairs (respectively, zero section pairs) in $\tsC_{W,E}$
and let $\tsC^{*,0}_{W,E}=\tsC_{W,E}^*\cap \tsC^0_{W,E}$.  The quotients
$\sC_{W,E}^*$, $\sC_{W,E}^0$, and $\sC_{W,E}^{*,0}$ are similarly defined.

If $(A,\Phi)\in\tsC^{*,0}_{W,E}$, then
$\Stab_{A,\Phi}=\{\id_E\}$, as the stabilizer of an irreducible connection
$A$ will be $S^1_Z$ and if $S^1_Z$ stabilizes the section $\Phi$ then
$e^{i\theta}\Phi = 0$ for all $\theta\in\RR$ 
and so $\Phi \equiv 0\in L^2_k(X,W^+\otimes E)$ (note that $\Phi$ need not
be continuous). Therefore, $\ssG_E$ acts freely on
$\tsC^{*,0}_{W,E}$ and, as we shall see in the next section, the
quotient $\sC^{*,0}_{W,E}$ is a Banach manifold.
Conversely, if $(A,\Phi)$ is a $\PU(2)$ monopole, $A$ is
reducible, and $\Phi\not\equiv 0$, then Lemma
\ref{lem:ReducibleStabilizers} implies that $\Stab_{A,\Phi}\simeq S^1$.

\subsection{A slice for the action of the group of gauge transformations}
\label{subsec:ConfigSpace}
Let $k\ge 2$ be an integer. The {\em slice\/} $\bS_{A,\Phi}\subset
\tsC_{W,E}$ through a pair $(A,\Phi)$ is given by $\bS_{A,\Phi} :=
(A,\Phi) + \bK_{A,\Phi}$, where
\begin{equation}
\bK_{A,\Phi} := \Ker d^{0,*}_{A,\Phi} 
\subset 
L^2_k(\La^1\otimes\su(E))\oplus L^2_k(W^+\otimes E)
\label{eq:Slice}
\end{equation}
If $\pi$ is the projection from $\tsC_{W,E}$ onto $\sC_{W,E}
= \tsC_{W,E}/\ssG_E$, denoted by $(A,\Phi)\mapsto
[A,\Phi]$, we let 
$$
\bB_{A,\Phi}(\eps)=\pi^{-1}B_{A,\Phi}(\eps)\cap \bS_{A,\Phi},
$$
be the open $L^2_k$ ball in $\bS_{A,\Phi}$ with center $(A,\Phi)$ and radius
$\eps$, so that
\begin{align*}
\bB_{A,\Phi}(\eps)
&:= \{(A_1,\Phi_1)\in \bS_{A,\Phi}:\|(A_1,\Phi_1)-(A,\Phi)\|_{L^2_{k,A}}
< \eps\} \\
&= (A,\Phi) + \{(a,\phi)\in \bK_{A,\Phi}:\|(a,\phi)\|_{L^2_{k,A}}
< \eps\}.
\end{align*}
The Hilbert Lie group $\ssG_E$ has Lie algebra 
$L^2_{k+1}(\su(E))\oplus i\RR_Z\subset L^2_{k+1}(\su(E))$ 
and exponential map $\exp:L^2_{k+1}(\su(E))\oplus i\RR_Z\to
\ssG_E$ given by $\zeta\mapsto u=\exp\zeta$.
Let $\Stab^\perp_{A,\Phi}=\exp((\Ker d^0_{A,\Phi})^\perp) 
\subset \ssG_E$, where
$$
(\Ker(d^0_{A,\Phi}|_{L^2_{k+1}}))^\perp
= \Imag(d^{0,*}_{A,\Phi}|_{L^2_{k+2}}) 
\subset L^2_{k+1}(\su(E)),
$$
noting that 
$$
d_{A,\Phi}^{0,*}:L^2_{k+2}\Om^1(\su(E))\to L^2_{k+1}(\su(E))
\oplus L^2_k(W^+\otimes E)
$$ 
has closed range.
Recall that $\Stab_{A,\Phi}\subset\ssG_E$ is given by $\{\ga\in
\ssG_E:\ga(A,\Phi) = (A,\Phi)\}$ and has Lie algebra
$H_{A,\Phi}^0 = \Ker(d^0_{A,\Phi}|_{L^2_{k+1}})
\subset L^2_{k+1}(\su(E))\oplus i\RR_Z$, so that
\begin{equation}
L^2_{k+1}(\su(E))\oplus i\RR_Z = 
(\Ker(d^0_{A,\Phi}|_{L^2_{k+1}}))^\perp \oplus H_{A,\Phi}^0. 
\label{eq:GaugeLieAlgOrthogDecomp}
\end{equation}
The subspace $\Stab^\perp_{A,\Phi}\subset \ssG_E$ is closed and is a
Banach submanifold of $\ssG_E$ with codimension $\dim H_{A,\Phi}^0$. 

The map $d_{A,\Phi}^0:L^2_{k+1}(\su(E))\oplus i\RR_Z
\to L^2_k(\La^1\otimes\su(E))
\oplus L^2_k(W^+\otimes E)$ has closed range and so
we have an $L^2$-orthogonal decomposition 
\begin{equation}
T_{A,\Phi}\tsC_{W,E}
= \Imag(d^0_{A,\Phi}|_{L^2_{k+1}}) \oplus \bK_{A,\Phi}
\label{eq:TangentPairOrthogDecomp}
\end{equation}
of the tangent space to the space of $L^2_k$ pairs at the point
$(A,\Phi)$. 

The proof that the quotient space $\sC_{W,E}$ is Hausdorff and
our later proof of removable singularities for $\PU(2)$ monopoles make
use of the following well-known technical result \cite[Proposition 2.3.15]{DK},
\cite[Proposition A.5]{FU}, \cite[Theorem II.7.11]{Lawson}. Note that 
the space of $L^2_2$ unitary automorphisms of $E$ is
neither a Hilbert Lie group nor does it act smoothly on 
the space of $L^2_2$ unitary connections on $E$.

\begin{lem}\label{lem:GaugeSequence}
Let $E$ be a Hermitian bundle over a Riemannian manifold $X$ and let $k\ge
2$ be an integer.  Suppose $\{A_\al\}$ and $\{B_\al\}$ are sequences of
$L^2_k$ unitary connections on $E$ and that $\{u_\al\}$ is a sequence of
unitary automorphisms of $E$ such that $u_\al(A_\al) = B_\al$. Then the
following hold.
\begin{enumerate}
\item The sequence $\{u_\al\}$ is in $L^2_{k+1}$;
\item If $\{A_\al\}$ and $\{B_\al\}$ converge in 
$L^2_k$ to limits $A_\8$, $B_\8$, then  
there is a subsequence $\{\al'\}\subset\{\al\}$ such that
$\{u_{\al'}\}$ converges in $L^2_{k+1}$ to $u_\8$ and $B_\8 = u_\8(A_\8)$.
\end{enumerate}
\end{lem}

The following slice result was established by Parker \cite{Parker}, but only
for pairs which are neither zero-section nor reducible; see also
\cite{BradlowDask}. It is, of course, the
analogue of the usual result for the topology and manifold structure of the
configuration space $\sB_E=\sA_E/\sG_E$ of connections.  The
proof we give here is modelled on the corresponding arguments for
connections given in \cite[Proposition 2.3.4]{DK}, \cite[Theorem 3.2 \&
4.4]{FU}, and \cite[Theorem II.10.4]{Lawson}. We will ultimately need a
rather more involved version of this method in order to show that our
gluing maps are diffeomorphisms, so we give the argument in the simpler
model case below in some detail and establish some of the notation and
conventions we later require.

\begin{prop}\label{prop:Slice}
Let $X$ be a closed, oriented, Riemannian \spinc four-manifold, let $E$
be a Hermitian two-plane bundle over $X$, and let
$k\ge 2$ be an integer.  Then the following hold.
\begin{enumerate}
\item The space $\sC_{W,E}$ is Hausdorff;
\item The subspace $\sC_{W,E}^{*,0}\subset \sC_{W,E}$ is
open and is a $C^\8$ Hilbert manifold with local parametrizations given by
$\pi:\bB_{A_0,\Phi_0}(\eps)\to \sC_{W,E}^{*,0}$ for sufficiently
small $\eps=\eps(A_0,\Phi_0,k)$;
\item The projection $\pi:\tsC_{W,E}^{*,0} \to
\sC_{W,E}^{*,0}$ is a $C^\8$ principal $\ssG_E$ bundle;
\item For $(A_0,\Phi_0)\in \tsC_{W,E}$, the projection  
$\pi:\bB_{A_0,\Phi_0}(\eps)/\Stab_{A_0,\Phi_0}\to \sC_{W,E}$ is a
homeomorphism onto an open neighborhood of $[A_0,\Phi_0]\in
\sC_{W,E}$ and a diffeomorphism on the complement of the
set of points in $\bB_{A_0,\Phi_0}(\eps)$ with non-trivial stabilizer.
\end{enumerate}
\end{prop}

\begin{pf}
The stabilizer $\Stab_{A_0,\Phi_0}$ (which we can identify with a Lie
subgroup $\U(2)$) acts freely on $\ssG_E$ and thus
on the Hilbert manifold $\ssG_E\times 
\bS_{A_0,\Phi_0}$ by $(u,A,\Phi)\mapsto \ga\cdot(u,A,\Phi) = 
(u\ga^{-1},\ga(A,\Phi))$ and so the
quotient $\ssG_E\times_{\Stab_{A_0,\Phi_0}}\bS_{A_0,\Phi_0}$ is
again a Hilbert manifold. We define a smooth map   
$$
\bPsi: \ssG_E\times_{\Stab_{A_0,\Phi_0}} \bS_{A_0,\Phi_0}
\rightarrow \tsC_{W,E}, \qquad
[u,A,\Phi]  \mapsto u(A,\Phi).
$$
Our first task is to show that the map $\bPsi$ is a 
diffeomorphism onto its image upon restriction to a
sufficiently small neighborhood $\ssG_E\times_{\Stab_{A_0,\Phi_0}}
\bB_{A_0,\Phi_0}(\eps)$. 
Given $\de>0$, let $B_{\id}(\de)$ be the ball
$\{u\in \ssG_E: \|u-\id_E\|_{L^2_{k+1,A_0}}<\de\}$ 
and let $B_{\id}^\perp(\de)= B_{\id}(\de)\cap\Stab^\perp_{A_0,\Phi_0}$.

\begin{claim}\label{claim:GaugeGroupSlice}
For small enough $\de=\de(A_0,\Phi_0,k)$, the ball $B_{\id}(\de)$ is
diffeomorphic to an open neighborhood in $B_{\id}^\perp(\de)\times
\Stab_{A_0,\Phi_0}$, with inverse map given by $(u_0,\ga)\mapsto u =
u_0\gamma$.
\end{claim}  

\begin{pf}
The differential of the multiplication map
$$
\Stab^\perp_{A_0,\Phi_0}\times \Stab_{A_0,\Phi_0}\rightarrow \ssG_E, 
\qquad (u_0,\ga)\mapsto u_0\ga,
$$ 
at $(\id_E,\id_E)$ is given by 
$$
\Ker(d^0_{A_0,\Phi_0}|_{L^2_{k+1}})^\perp
\oplus H_{A_0,\Phi_0}^0\rightarrow L^2_{k+1}(\su(E))
\qquad (\zeta,\chi)\mapsto \zeta + \chi,
$$
and is just the identity map with respect to the 
$L^2$-orthogonal decomposition \eqref{eq:GaugeLieAlgOrthogDecomp}
of the range. Hence, the
Hilbert space implicit function theorem implies that there is
a diffeomorphism from an open 
neighborhood of $(\id_E,\id_E)$ onto an open neighborhood of $\id_E\in
\ssG_E$. For small enough $\de$, we may suppose that if
$u\in B_{\id}(\de)$, then $u$ can
be written uniquely as $u=u_0\gamma$ with 
$u_0\in  B_{\id}^\perp(\de)$ and
$\gamma\in\Stab_{A_0,\Phi_0}$.  
\end{pf}

\begin{claim}
For small enough $\eps=\eps(A_0,\Phi_0,k)$, 
the map $\bPsi$ is a 
diffeomorphism from $\ssG_E\times_{\Stab_{A_0,\Phi_0}}
\bB_{A_0,\Phi_0}(\eps)$ onto its image in $\tsC_{W,E}$.
\end{claim}

\begin{pf}
We first restrict the map $\bPsi$ to a neighborhood 
$B_{\id}(\de)\times_{\Stab_{A_0,\Phi_0}}\bS_{A_0,\Phi_0}$, which is
diffeomorphic to the neighborhood
$B_{\id}^\perp(\de)\times\bS_{A_0,\Phi_0}$ in  
$\Stab^\perp_{A_0,\Phi_0}\times \bS_{A_0,\Phi_0}$ by Claim
\ref{claim:GaugeGroupSlice}. The differential of the induced map
$$
\bPsi:\Stab^\perp_{A_0,\Phi_0} \times \bS_{A_0,\Phi_0}
\to \tsC_{W,E}, \qquad
(u,A,\Phi) \mapsto u(A,\Phi),
$$
at $(\id_E,A_0,\Phi_0)$ is given by
\begin{align*}
(D\bPsi)_{(\id,A_0,\Phi_0)}: &T_{\id}\Stab^\perp_{A_0,\Phi_0}
\oplus T_{A_0,\Phi_0}\bS_{A_0,\Phi_0}
\to T_{A_0,\Phi_0}\tsC_{W,E}, \\
&(\zeta,a,\phi) \mapsto -d^0_{A_0,\Phi_0}\zeta + (a,\phi),
\end{align*}
where we recall that $T_{A_0,\Phi_0}\bS_{A_0,\Phi_0} = \bK_{A_0,\Phi_0}$ and
$$
T_{\id}\Stab^\perp_{A_0,\Phi_0} = (\Ker(d^0_{A_0,\Phi_0}|_{L^2_{k+1}}))^\perp
= \Imag(d^{0,*}_{A_0,\Phi_0}|_{L^2_{k+2}}).
$$
Using the $L^2$-orthogonal
decomposition \eqref{eq:TangentPairOrthogDecomp} of the range we see that
the map 
$$
-d^0_{A_0,\Phi_0}\oplus \id_E: (\Ker(d^0_{A_0,\Phi_0}|_{L^2_{k+1}}))^\perp
\oplus \bK_{A_0,\Phi_0} \to \Imag(d^0_{A_0,\Phi_0}|_{L^2_{k+1}})
\oplus \bK_{A_0,\Phi_0}
$$
given by $(\zeta,b,\psi) \mapsto -d^0_{A_0,\Phi_0}\zeta +
(b,\psi)$ is a Hilbert space isomorphism.  So, by the
Hilbert space implicit function theorem, there are
positive constants $\eps=\eps(A_0,\Phi_0,k)$ and $\de=\de(A_0,\Phi_0,k)$ 
and an open neighborhood $\sU_{A_0,\Phi_0}\subset 
\tsC_{W,E}$ such that the map
$$
\bPsi: B_{\id}^\perp(\de) 
\times \bB_{A_0,\Phi_0}(\eps)\to \sU_{A_0,\Phi_0}, \qquad
(u,A,\Phi)\mapsto u(A,\Phi),
$$ 
gives a diffeomorphism from an open neighborhood of
$(\id_E,A_0,\Phi_0)$ onto an open neighborhood of $(A_0,\Phi_0)$. 
In particular,
we obtain a map $\sU_{A_0,\Phi_0}\to \Stab^\perp_{A_0,\Phi_0}$, given by
$(A,\Phi) \mapsto u = u_{A,\Phi}$, such that
$$
\bPsi^{-1}(A,\Phi) = \left(u,u^{-1}(A,\Phi)\right)
\in B_{\id}^\perp(\de) \times \bB_{A_0,\Phi_0}(\eps)
\subset \Stab^\perp_{A_0,\Phi_0}\times \bB_{A_0,\Phi_0}(\eps).
$$
Hence, for any $(A,\Phi)\in \sU_{A_0,\Phi_0}$ there is a unique
$u\in  B_{\id}^\perp(\de)$ such that
$u^{-1}(A,\Phi)-(A_0,\Phi_0)\in\bK_{A_0,\Phi_0}$: 
\begin{equation}
d^{0,*}_{A_0,\Phi_0}(u^{-1}(A_1,\Phi_1)-(A_0,\Phi_0)) = 0.
\label{eq:Coulomb}
\end{equation}
The neighborhood $\bB_{A_0,\Phi_0}(\eps)$ is
$\Stab_{A_0,\Phi_0}$-invariant: if $(A,\Phi)\in
\bB_{A_0,\Phi_0}(\eps)$ and 
$\ga\in\Stab_{A_0,\Phi_0}$, then 
\begin{align*}
\|\ga(A_1,\Phi_1)-(A,\Phi)\|_{L^2_{k,A_0}}
&= \|(A_1,\Phi_1)-\ga^{-1}(A,\Phi)\|_{L^2_{k,\ga^{-1}(A_0)}} \\
&= \|(A_1,\Phi_1)-(A,\Phi)\|_{L^2_{k,A_0}} <\eps,
\end{align*} 
and
\begin{align*}
d^{0,*}_{A_0,\Phi_0}\left(\ga(A,\Phi)-(A_0,\Phi_0)\right) 
&= \ga\left(d^{0,*}_{\ga^{-1}(A_0,\Phi_0)}
\left((A,\Phi)-(A_0,\Phi_0)\right)\right) \\
&= \ga\left(d^{0,*}_{A_0,\Phi_0}\left((A,\Phi)-(A_0,\Phi_0)\right)\right)
= 0,
\end{align*} 
so $\ga(A,\Phi)\in \bB_{A_0,\Phi_0}(\eps)$. 

The group $\ssG_E$ acts on $\ssG_E\times
\bS_{A_0,\Phi_0}$ by $(u,A,\Phi)\mapsto (vu,A,\Phi)$, and so gives a
diffeomorphism  
$$
B_{\id}(\de) \times \bB_{A_0,\Phi_0}(\eps)
\to B_v(\de) \times \bB_{A_0,\Phi_0}(\eps), \quad
(u,A,\Phi) \to (vu,A,\Phi),
$$
and as this action commutes with the given action of
$\Stab_{A_0,\Phi_0}$, it descends to a diffeomorphism 
$$
B_{\id}(\de) \times_{\Stab_{A_0,\Phi_0}} \bB_{A_0,\Phi_0}(\eps)
\to B_v(\de) \times_{\Stab_{A_0,\Phi_0}} \bB_{A_0,\Phi_0}(\eps), \quad
[u,A,\Phi] \to [vu,A,\Phi],
$$
for each $v\in \ssG_E$. Consequently, the 
$\ssG_E$-equivariant map
$$
\ssG_E\times_{\Stab_{A_0,\Phi_0}}
\bB_{A_0,\Phi_0}(\eps)\to \tsC_{W,E}
$$
is a diffeomorphism onto its image, as desired.
\end{pf}

Plainly, $[\ga(A,\Phi)] =
[A,\Phi]$ for each $\ga\in\Stab_{A_0,\Phi_0}$ and $(A,\Phi) \in
\bB_{A_0,\Phi_0}(\eps)$  
and hence, the projection $\pi:\bB_{A_0,\Phi_0}(\eps)\to
\sC_{W,E}$ factors through
$\bB_{A_0,\Phi_0}(\eps)/\Stab_{A_0,\Phi_0}$.

\begin{claim}\label{claim:Stab}
There is a positive constant $\de=\de(A_0,\Phi_0,k)$ with the
following significance.
If $(A_i,\Phi_i)\in\bB_{A_0,\Phi_0}(\eps)$ for $i=1,2$ and 
there is a gauge transformation $u\in B_{\id}(\de)$ such that
$u(A_1,\Phi_1)=(A_2,\Phi_2)$, then $u\in\Stab_{A_0,\Phi_0}$.
\end{claim}

\begin{pf}
For small enough $\de$, Claim \ref{claim:GaugeGroupSlice} implies that
$u\in B_{\id}(\de)$ can be written uniquely as $u=u_1\gamma$ with 
$u_1\in  B_{\id}^\perp(\de)$ and $\gamma\in\Stab_{A_0,\Phi_0}$. Thus, 
$(A_2,\Phi_2) = u_1\ga(A_1,\Phi_1)\in \bB_{A_0,\Phi_0}(\eps)$. But the
neighborhood $\bB_{A_0,\Phi_0}(\eps)$ is
$\Stab_{A_0,\Phi_0}$-invariant, so we also have $\ga(A_1,\Phi_1)
\in \bB_{A_0,\Phi_0}(\eps)$. Therefore
$u_1=\id_E$ by the uniqueness assertion of Claim
\ref{claim:GaugeGroupSlice} and so $u=\ga\in\Stab_{A_0,\Phi_0}$.
This completes the proof of the claim.
\end{pf}

Claim \ref{claim:Stab} shows that
$\bB_{A_0,\Phi_0}(\eps)/\Stab_{A_0,\Phi_0}$ injects into the quotient
$\sC_{W,E}$ modulo the assumption that the gauge transformations
are close to $\id_E$. It remains to show that if
$(A_i,\Phi_i)\in\bB_{A_0,\Phi_0}(\eps)$ for $i=1,2$ and there is a gauge
transformation $u\in \ssG_E$ such that $u(A_1,\Phi_1)=(A_2,\Phi_2)$,
then $u$ is necessarily close to $\id_E$ and hence that
$\bB_{A_0,\Phi_0}(\eps)/\Stab_{A_0,\Phi_0}$ injects into the quotient
$\sC_{W,E}$.

\begin{claim}\label{claim:injective}
For small enough $\eps(A_0,\Phi_0)$, the projection 
$\pi:\bB_{A_0,\Phi_0}(\eps)/\Stab_{A_0,\Phi_0}\to\sC_{W,E}$ is
injective.
\end{claim}

\begin{pf}
Suppose $(A_i,\Phi_i)\in\bB_{A_0,\Phi_0}(\eps)$ for $i=1,2$ and that
$[A_1,\Phi_1]=[A_2,\Phi_2]\in \sC_{W,E}$, so that
$u(A_1,\Phi_1)=(A_2,\Phi_2)$ for some $u\in \ssG_E$. 

Since $u(A_0) = A_0 - (d_{A_0}u)u^{-1}$, we see that
$u\in\Stab_{A_0,\Phi_0}$ if 
and only $d_{A_0,\Phi_0}^0u = (d_{A_0}u,-u\Phi_0) = (0,-\Phi_0)$ or,
equivalently, if 
and only if $d_{A_0,\Phi_0}^0(u-\id_E) = 0$, since
$\id_E\in\Stab_{A_0,\Phi_0}$. Here, we view $u\in
L^2_{k+1}(\gl(E))$ via the isometric embedding
$\ssG_E\subset L^2_{k+1}(\gl(E))$ and write
$$
u-\id_E = u_0 - \ga,
$$
where $u_0 \in (\Ker d_{A_0,\Phi_0}^0)^\perp$ and $\ga \in \Ker
d_{A_0,\Phi_0}^0$. Our first task is to estimate
$\|u-\id_E\|_{L^2_{3,A_0}}$.  

By assumption,
$d_{A_0,\Phi_0}^0\ga = (d_{A_0}\ga,\ga\Phi) = 0$, so $d_{A_0}\ga = 0$ and
$\ga\Phi_0=0$. Since $u(A_1)\equiv A_1 - (d_{A_1}u)u^{-1}=A_2$, we have
$A_2u = A_1u - d_{A_1}u = A_1u - d_{A_0}u - [A_1-A_0,u]$, and therefore,
using $d_{A_0}\id_E = 0 = d_{A_0}\ga$, we have 
$$
d_{A_0}u_0 = d_{A_0}u = u(A_1-A_0) - (A_2-A_0)u.
$$
As $d_{A_0}^*(A_i-A_0) = (\cdot\Phi_0)^*(\Phi_i-\Phi_0)$, we obtain
\begin{align*}
d_{A_0}^*d_{A_0}u_0
&= -*(d_{A_0}u\wedge *(A_1-A_0)) + ud_{A_0}^*(A_1-A_0) \\
&\qquad - (d_{A_0}^*(A_2-A_0))u + *(*(A_2-A_0)\wedge d_{A_0}u) \\
&= -*(d_{A_0}u_0\wedge *(A_1-A_0)) + u(\cdot\Phi_0)^*(\Phi_1-\Phi_0) \\
&\qquad - ((\cdot\Phi_0)^*(\Phi_2-\Phi_0))u + *(*(A_2-A_0)\wedge d_{A_0}u_0).
\end{align*}
We define the Laplacian $\De^0_{A_0,\Phi_0}$ by setting
$$
\De^0_{A_0,\Phi_0} =  d^{0,*}_{A_0,\Phi_0}d^0_{A_0,\Phi_0}
= d_{A_0}^*d_{A_0} + (\cdot\Phi_0)^*\Phi_0.
$$
Now $u_0 = u-\id_E + \ga$ and $d_{A_0}\id_E = 0 = d_{A_0}\ga$ and
$\ga\Phi_0=0$ so, using $u\Phi_1=\Phi_2$, 
\begin{align*}
\De^0_{A_0,\Phi_0}u_0 &= \De^0_{A_0,\Phi_0}(u-\id_E+\ga) \\
&= d_{A_0}^*d_{A_0}(u-\id_E+\ga) + (\cdot\Phi_0)^*(u-\id_E)\Phi_0 \\
&= d_{A_0}^*d_{A_0}u + (\cdot\Phi_0)^*(u\Phi_0-u\Phi_1+u\Phi_1-\Phi_0) \\
&= d_{A_0}^*d_{A_0}u_0 
+ (\cdot\Phi_0)^*\left(u(\Phi_0-\Phi_1)+(\Phi_2-\Phi_0)\right).
\end{align*}
Our assumption that $(A_i,\Phi_i)\in\bB_{A_0,\Phi_0}(\eps)$ and the
embeddings $L^2_k\subset L^2_2$ imply that for $i=1,2$ we have
$$
\|(A_i,\Phi_i)-(A_0,\Phi_0)\|_{L^2_{2,A_0}}
\le \|(A_i,\Phi_i)-(A_0,\Phi_0)\|_{L^2_{k,A_0}} < \eps,
$$
Since $u_0\in (\Ker d^0_{A_0,\Phi_0})^\perp$, the standard elliptic estimate
for the Laplacian $\De^0_{A_0,\Phi_0}$, the fact that $|u|=1$ as
$u\in\ssG_E$, the Sobolev embedding $L^2_2\subset L^q$, $2\le q<\8$,
and multiplication $L^2_2\times L^2_2\to L^2_1$, and
our expression for $d_{A_0}^*d_{A_0}u = d_{A_0}^*d_{A_0}u_0$ combine to give
\begin{align*}
\|u_0\|_{L^2_{3,A_0}}
&\le C\|\De^0_{A_0,\Phi_0}u_0\|_{L^2_{1,A_0}} \\
&\le C\|d_{A_0}^*d_{A_0}u_0\|_{L^2_{1,A_0}} 
+ C\|(\cdot\Phi_0)^*(u(\Phi_0-\Phi_1)+(\Phi_2-\Phi_0))\|_{L^2_{1,A_0}} \\
&\le C\|d_{A_0}u_0\|_{L^2_{2,A_0}}
(\|A_1-A_0\|_{L^2_{2,A_0}}+\|A_2-A_0\|_{L^2_{2,A_0}}) \\
&\quad + C\|u\|_{L^2_{2,A_0}}\|\Phi_0\|_{L^2_{2,A_0}}
(\|\Phi_1-\Phi_0\|_{L^2_{2,A_0}} + \|\Phi_2-\Phi_0\|_{L^2_{2,A_0}}) \\
&\le C\|u_0\|_{L^2_{3,A_0}}\eps
+ C\|u\|_{L^2_{2,A_0}}\|\Phi_0\|_{L^2_{2,A_0}}\eps \\
&\le C\|u_0\|_{L^2_{3,A_0}}\eps + C\|\Phi_0\|_{L^2_{2,A_0}}\eps,
\end{align*}
for some constant $C=C(A_0,\Phi_0)$. For small enough $\eps$ we can
rearrange the preceding inequality and use the Sobolev embedding
$L^2_3\subset C^0$ to give 
$$
\|u_0\|_{L^2_{3,A_0}} \le C\|\Phi_0\|_{L^2_{2,A_0}}\eps
\quad\text{and}\quad 
\|u_0\|_{C^0} \le C\|\Phi_0\|_{L^2_{2,A_0}}\eps.
$$
Therefore, as $|u|=1$, we have the pointwise bounds
$$
1-|u_0|\le |\id_E-\ga| \le 1+|u_0|,
$$
and thus, for $\eps$ small, $|\id_E-\ga|>0$; the pointwise norm of
$z = |\id_E-\ga|$ is constant since
$d_{A_0}(\id_E-\ga)=0$. Consequently, $\ga_0 \equiv z^{-1}(\id_E-\ga)$ lies
in $\ssG_E$: clearly, $d_{A_0}\ga_0=0$, so $\ga_0\in\Stab_A$. If
$\Phi_0\equiv 0$, we trivially have $\ga_0\in\Stab_{\Phi_0}$; if
$\Phi_0\not\equiv 0$, then the equalities $|\ga_0\Phi_0|=|\Phi_0|$ (as
$\ga_0$ is a unitary gauge transformation) and $\ga_0\Phi_0 =
z^{-1}(\id_E-\ga)\Phi_0 = z^{-1}\Phi_0$, so $|\ga_0\Phi_0|=z^{-1}|\Phi_0|$, 
imply that $z=1$ (it is enough to have equality of $L^2$ norms here as $z$ is
constant). Hence, we also have $\ga_0\in\Stab_{\Phi_0}$ and in particular,
$\ga_0 \in \Stab_{A_0,\Phi_0} = \Stab_{A_0}\cap \Stab_{\Phi_0}$.

We now write $\ga_0^{-1}u = \ga_0^{-1}u_0 + \ga_0^{-1}(\id_E - \ga)
= \ga_0^{-1}u_0 + z\,\id_E$, so that
$$
\ga_0^{-1}u - \id_E = (z-1)\id_E + \ga_0^{-1}u_0.
$$
Clearly, $(z-1)\id_E\in \Ker d_{A_0,\Phi_0}^0$ by the remarks of the last
paragraph, while $\ga_0^{-1}u_0\in (\Ker d_{A_0,\Phi_0}^0)^\perp$ since
$u_0\in (\Ker d_{A_0,\Phi_0}^0)^\perp$ and $\ga_0^{-1}\in
\Stab_{A_0,\Phi_0}$. Similarly, $\ga_0^{-1}u(A_1,\Phi_1) =
\ga_0^{-1}(A_2,\Phi_2)\in \bB_{A_0,\Phi_0}(\eps)$, as the neighborhood
$\bB_{A_0,\Phi_0}(\eps)$ is $\Stab_{A_0,\Phi_0}$-invariant. 
Thus, replacing $u$ by
$\ga_0^{-1}u$ and $(A_2,\Phi_2)$ by
$(\tA_2,\tPhi_2)\equiv\ga_0^{-1}(A_2,\Phi_2)$ in our 
$L^2_{3,A_0}$ estimate for $u_0$ yields
$$
\|\ga_0^{-1}u_0\|_{L^2_{3,A_0}} \le C\|\Phi_0\|_{L^2_{2,A_0}}\eps,
$$
while $|z-1| \le \|u_0\|_{C^0} \le C\|\Phi_0\|_{L^2_{2,A_0}}\eps$, so we
find that
$$
\|\ga_0^{-1}u - \id_E\|_{L^2_{3,A_0}}
\le C\|\Phi_0\|_{L^2_{2,A_0}}\eps.
$$
If $\eps$ is small enough that $C\|\Phi_0\|_{L^2_{2,A_0}}\eps<\de$,
where $\de = \de(A_0,\Phi_0)$ is the constant of Claim \ref{claim:Stab}
with $k=2$, then $\ga_0^{-1}u\in \Stab_{A_0,\Phi_0}$ and so $u$
lies in $\Stab_{A_0,\Phi_0}$, as desired.
\end{pf}

\begin{claim}\label{claim:Hausdorff}
The quotient space $\sC_{W,E}$ is Hausdorff.
\end{claim}

\begin{pf}
Let $\Ga$ be the subspace
$\{\{(A,\Phi),u(A,\Phi)\}:(A,\Phi)\in \tsC_{W,E}\text{
and }u\in \ssG_E\}$ of $\tsC_{W,E}\times
\tsC_{W,E}$. If $\{(A_\al,\Phi_\al),u_\al(A_\al,\Phi_\al)\}$ is a
sequence in $\Ga$ which converges in $L^2_k$ to a point 
$\{(A_\8,\Phi_\8),(B_\8,\Psi_\8)\}$, then Lemma \ref{lem:GaugeSequence}
implies that there is a subsequence $\{\al'\}\subset\{\al\}$ such that
$\{u_\al\}$ converges in $L^2_{k+1}$ 
to $u_\8\in\ssG_E$ and $u_\8(A_\8)=B_\8$. But
then $u_{\al'}\Phi_{\al'}$ converges in $L^2_k$ to $u_\8\Phi_\8$ and so
$\Psi_\8 = u_\8\Phi_\8$. Thus, $(B_\8,\Psi_\8)=u_\8(A_\8,\Phi_\8)$ for some
$u_\8\in \ssG_E$, so $\Ga$ is closed and the quotient 
$\tsC_{W,E}/\ssG_E$ is Hausdorff.
\end{pf}

Claim \ref{claim:Hausdorff} gives Assertion (1) of the proposition and
Assertions (2), (3), and (4) now follow from the preceding arguments and
Claim \ref{claim:injective}. This completes the proof of the proposition.
\end{pf}

\begin{rmk}
(1) Alternatively, one can show that the quotient $L^2_k$ topology on
$\sC_{W,E}$ is metrizable via the $L^2$ metric (exactly as in
\cite[Lemma 4.2.4]{DK}) and thus $\sC_{W,E}$ is Hausdorff. 

(2) As Mrowka pointed out to us, one can sharpen the assertions of
Proposition \ref{prop:Slice}, at least for the quotient space $\sB_E$: one
finds that charts are provided by $L^4$-balls in $\Ker
d_{A_0}^*$ rather than the much smaller $L^2_{k,A_0}$-balls usually
employed; see \cite{FeehanSlice}.
\end{rmk}

It is convenient to extract the following global, Coulomb gauge-fixing
result (analogous to Proposition 2.3.4 in \cite{DK}) which we established
in the course of proving Proposition \ref{prop:Slice}:

\begin{lem}\label{lem:CoulombPair}
Let $X$ be a closed, oriented, Riemannian \spinc four-manifold and let $E$
be a Hermitian vector bundle over $X$. Suppose that $k\ge 2$ is an integer
and that $(A_0,\Phi_0)\in\tsC_{W,E}$.
Then there is a positive constant $\eps=\eps(A_0,\Phi_0,k)$ such that
for any $(A,\Phi)\in\tsC_{W,E}$ with
$$
\|(A,\Phi)-(A_0,\Phi_0)\|_{L^2_{k,A_0}(X)}<\eps,
$$
there is a gauge transformation $u\in \ssG_E$, unique up to an
element of $\Stab_{A_0,\Phi_0}$, such that
$$
d_{A_0,\Phi_0}^{0,*}(u(A,\Phi)-(A_0,\Phi_0)) = 0.
$$
\end{lem}

\subsection{Connections on SO(3) bundles and groups of gauge
transformations} 
\label{subsec:GaugeGroup}  
For some local patching arguments over simply-connected open regions
$Y\subset X$ in \S \ref{sec:regularity} and \S \ref{sec:compact} it is very
useful to be able to lift gauge transformations in $\sG_{\su(E)}(Y)$ to
gauge transformations in $\sG_E(Y)$. The following result tells us that
this is always possible when $Y$ is simply connected; it is an extension of
Theorem IV.3.1 in \cite{MorganGTNotes} from $\SU(2)$ to $\U(2)$ bundles.

\begin{prop}
\label{prop:GaugeGroupSequence}
Let $E$ be a Hermitian two-plane bundle over
a connected manifold $X$. Then there is an exact sequence
$$
1\rightarrow
\{\pm \id_E\}
\rightarrow
\sG_E
\rightarrow
\sG_{\su(E)}
\rightarrow
H^1(X;\ZZ/2\ZZ).
$$
\end{prop}

\begin{proof}
Let $P$ be the $\U(2)$ principal bundle underlying $E$ so $P^{\ad}=P/S^1_Z$
is the $\SO(3)$ principal bundle underlying $\su(E)$.  We can view elements
of $\sG_E$ as maps $u:P\rightarrow \SU(2)$ satisfying $u(pg)=g^{-1}u(p)g$
for $p\in P$ and $g\in\U(2)$, and similarly for elements of $\sG_{\su(E)}$.
The map $\sG_E\rightarrow \sG_{\su(E)}$ is then given by $u\rightarrow \Ad u$.
Because $X$ is connected, the kernel of this map is $\{\pm\id_E\}$, giving
the exactness of all but the last terms in the sequence.

Let $\nu_2\in H^1(\SO(3);\ZZ/2\ZZ)$ be the non-trivial cohomology class given
by the double cover $\SU(2)\rightarrow \SO(3)$.  Thus a map $u:V\rightarrow
\SO(3)$ lifts to $\SU(2)$ if and only if $u^*\nu_2=0$.  We would like to
define the map $\alpha:\sG_{\su(E)}\rightarrow H^1(X;\ZZ/2\ZZ)$ by
$\alpha(u)=u^*\nu_2$.  It is not immediately clear that $\alpha(u)$
actually lies in $H^1(X;\ZZ/2\ZZ)$.  The homotopy exact sequence
$$
\pi_1(\SO(3))\rightarrow
\pi_1(P^{\ad})\rightarrow
\pi_1(X)\rightarrow 1
$$
gives an exact sequence in cohomology
$H^1(\cdot\ ;\ZZ/2\ZZ)=\Hom(\pi_1(\cdot),\ZZ/2\ZZ)$,
$$
1\rightarrow H^1(X;\ZZ/2\ZZ)
\rightarrow H^1(P^{\ad};\ZZ/2\ZZ)
\rightarrow H^1(\SO(3);\ZZ/2\ZZ).
$$
Thus $u^*\nu_2$ is pulled back from a unique element
of $H^1(X;\ZZ/2\ZZ)$ if and only if $i^*_xu^*\nu_2=0$,
where $i_x:\SO(3)\rightarrow P^{\ad}$ is the inclusion
of a fiber.  

Because $u^*\nu_2$ depends only on
the homotopy class of $u$, the map $\alpha$
is constant on connected components. 
Let $\sG^x_{\su(E)}$ be the subgroup of elements of
$\sG_{\su(E)}$ which are the identity over $x\in X$.
The exact
sequence of groups
$$
1\rightarrow \sG_{\su(E)}^x\rightarrow \sG_{\su(E)}\rightarrow \SO(3)
\rightarrow 1
$$
and the corresponding exact sequence of homotopy
groups ($\pi_0$ to be specific) shows that
the inclusion $\sG^x_{\su(E)}$ induces a surjection of
connected components.  Now if $u\in\sG^x_{\su(E)}$
then $i_x^*u^*\nu_2=0$, and so $i_xu^*\nu_2=0$
for all $u\in\sG_{\su(E)}$.  This implies that the map
$\alpha$ takes values in $H^1(X;\ZZ/2\ZZ)$.

By the defining property of $\nu_2$, we have $\alpha(u)=0$
if and only if there is a lift of $u$ to
$\SU(2)$.  Because $\{\pm\id_E\}$ is central,
any $\Ad$-equivariant map $\tilde u:P\rightarrow \SU(2)$ descends
to $P^{\ad}\rightarrow \SU(2)$.
Thus $u\in\sG_{\su(E)}$ is induced by some $\tilde u\in\sG_E$
if and only if $\alpha(u)=0$.
\end{proof}

\begin{rmk}
It can also be shown that the map $\sG_{\su(E)}\rightarrow H^1(X;\ZZ/2\ZZ)$
is surjective.  We omit the proof here in the interests of brevity.
\end{rmk}

\subsection{Quadratic forms for coupled spinor bundles}
\label{subsec:QuadForm}
Give $\gl(E)$ (respectively, $\gl(W^\pm)$) the 
Hermitian fiber inner product and norm
$$
\langle M,N\rangle := \half\tr(M^\dagger N) \quad\text{and}\quad
|M|^2 = \langle M,M\rangle,
$$
for $M,N\in\Om^0(\gl(E))$ (respectively, $\gl(W^\pm)$), so that $|\id_E| =
1$, etc. If $\Phi\in\Om^0(W^+\otimes E)$, then
$\Phi^*\in\Om^0((W^+\otimes E)^*) 
= \Om^0((W^+)^*\otimes E^*)$ is defined by $\Phi^*(\Psi) :=
\langle\Psi,\Phi\rangle$ for $\Psi\in\Om^0(W^+\otimes E)$. 

With respect to the induced Hermitian fiber inner products on $\gl(W^+)$ and
$\gl(E)$, we have fiber-wise orthogonal direct sums, $\gl(W^+) =
\CC\,\id_{W^+}\oplus\fsl(W^+)$, and similarly for $\gl(E)$ and thus
$\gl(W^+)\otimes_\CC\gl(E)$. Let $\pi_{W^+}$ (respectively, $\pi_E$) be the
fiber-wise orthogonal projection from $\gl(W^+)$ 
(respectively, $\gl(E)$) onto the subspace $\fsl(W^+)$
(respectively, $\fsl(E)$), so that
$$
\pi_{W^+}M := M - \half\tr M\,\id_{W^+} \quad\text{and}\quad
\pi_E N := N - \half\tr N\,\id_E
$$
are the projection onto the traceless components of $M\in\gl(W^+)$ and
$N\in\gl(E)$, and 
$$
\pi_{W^+}^\perp M := \half\tr M\,\id_{W^+} \quad\text{and}\quad
\pi_E^\perp N := \half\tr N\,\id_E
$$
are the projections onto the trace components of $M\in\gl(W^+)$ and
$N\in\gl(E)$. The orthogonal projection from $\gl(W^+)\otimes_\CC\gl(E)$ 
onto $\fsl(W^+)\otimes_\CC\fsl(E)$ 
is obtained by writing $\pi_{W^+} = \id_{W^+} - \pi_{W^+}^\perp$ and 
$\pi_E = \id_E - \pi_E^\perp$, so that
\begin{equation}
\pi_{W^+}\otimes\pi_E
= \id_{W^+}\otimes\id_E 
 - \left(\pi_{W^+}\otimes \pi_E^\perp
+ \pi_{W^+}^\perp\otimes \pi_E\right) 
+ \pi_{W^+}^\perp\otimes \pi_E^\perp.
\label{eq:doubletrace}
\end{equation}
Consequently, if $\phi\in\Om^0(W^+)$ and $\Phi\in\Om^0(W^+\otimes E)$, we
may define
\begin{align}
(\phi\otimes\phi^*)_0 
&:= \pi_{W^+}(\phi\otimes\phi^*), \\
(\Phi\otimes\Phi^*)_{00} 
&:= (\pi_{W^+}\otimes\pi_E)(\Phi\otimes\Phi^*).\notag
\end{align}
Note that 
$$
\gl(W^+)\otimes_\CC\gl(E)
=
\fu(W^+)\otimes_\RR\fu(E) + i\fu(W^+)\otimes_\RR\fu(E).
$$
Plainly, $(\phi\otimes\phi^*)_0$ is a section of $\fsl(W^+)\cap\, i\fu(W^+) =
i\su(W^+)$.  The following identities are now well-known \cite{KrM},
\cite[Chapter 8]{Salamon}:

\begin{lem}\label{lem:sigmaphiphi}
For every $\phi\in\Om^0(W^+)$ or $\Om^0(E)$, $T\in\Om^0(\fsl(W^+))$, and
$\eta\in \Om^2(X,\CC)$, the following hold. The endomorphism
$(\phi\otimes\phi^*)_0$ is a section of $\Om^0(i\su(W^+))$ and
satisfies:
\begin{align}
|(\phi\otimes\phi^*)_0|^2 &= \quarter |\phi|^4, \tag{1}\\
|\rho(\eta)|^2 &= 2|\eta^+|^2\quad\text{and}\quad
|\rho^{-1}T|^2 = \half |T|^2. \tag{2}
\end{align}
\end{lem}

Similarly, $\Phi\otimes\Phi^*$ is a section of $i\fu(W^+\otimes_\CC E) =
\fu(W^+)\otimes_\RR\fu(E)$ and
$(\Phi\otimes\Phi^*)_{00}$ is a section of $\su(W^+)\otimes_\RR\su(E)$.  We
then have the following analogue of Lemma \ref{lem:sigmaphiphi}:

\begin{lem}\label{lem:UnpertSigmaPhiPhi}
For every $\Phi\in\Om^0(W^+\otimes E)$ and $M\in\Om^0(\gl(W^+\otimes E))$,
the following identities hold:
\begin{align}
\langle \Phi\otimes\Phi^*, M\rangle &= \half\langle M\Phi,\Phi\rangle, 
\tag{1}\\
{\textstyle{\left({\frac{5}{4}}-{\frac{1}{\sqrt{2}}}\right)}}
|\Phi|^4 &\le \langle(\Phi\otimes\Phi^*)_{00}\Phi,\Phi\rangle
\le {\textstyle{\frac{1}{\sqrt{2}}}}|\Phi|^4, \tag{2}\\
{\textstyle\left({\frac{5}{8}}-{\frac{1}{2\sqrt{2}}}\right)}|\Phi|^4 
&\le |(\Phi\otimes\Phi^*)_{00}|^2
\le {\textstyle{\frac{1}{2\sqrt{2}}}}|\Phi|^4. \tag{3}
\end{align}
\end{lem}

\begin{pf}
Using $\Phi=c^a\Phi^a$, so
$M\Phi = c^a M\Phi^a$ and $\Phi^*\Phi^a = 
\langle\Phi^a,\Phi\rangle = \barc^a$, we have
\begin{align*}
\langle \Phi\otimes \Phi^*, M\rangle
&= \half\tr ((\Phi\otimes \Phi^*)^\dagger M)
= \half\langle (\Phi\otimes \Phi^*)M\Phi^a,\Phi^a\rangle \\
&= \half\langle M\Phi^a,(\Phi\otimes \Phi^*)\Phi^a\rangle 
= \half\langle M\Phi^a,\Phi(\Phi^*\Phi^a)\rangle \\
&= \half\langle M\Phi^a,\barc^a\Phi\rangle 
= \half\langle c^aM\Phi^a,\Phi\rangle 
= \half\langle M\Phi,\Phi\rangle,
\end{align*} 
which gives (1). Next, we see that
\begin{align*}
|\Phi\otimes\Phi^*|^2
&= \langle\Phi\otimes\Phi^*,\Phi\otimes\Phi^*\rangle
= \half\langle(\Phi\otimes\Phi^*)\Phi,\Phi\rangle 
= \half|\Phi|^2\langle\Phi,\Phi\rangle
= \half|\Phi|^4,
\end{align*} 
and so the upper bound in (2) is obtained by
$$
|\langle(\Phi\otimes\Phi^*)_{00}\Phi,\Phi\rangle|
\le |(\Phi\otimes\Phi^*)_{00}||\Phi|^2 \le |\Phi\otimes\Phi^*||\Phi|^2
\le (1/\sqrt{2})|\Phi|^4 < \threequarter|\Phi|^4,
$$
where we
use the fact that $(\ \cdot\ )_{00} = \pi_{W^+}\otimes\pi_E$ is an
orthogonal projection to obtain the second inequality. 

To obtain the
lower bound in (2), suppose that
$\Phi=\phi^a\otimes\xi^a$ (implied summation), where
$\{\phi^a\}\subset \Om^0(W^+)$ and $\{\xi^a\}\subset \Om^0(E)$, then
$\Phi\otimes\Phi^* = (\phi^a\otimes\phi^{b,*})\otimes 
(\xi^a\otimes\xi^{b,*})$ and
\begin{align*}
(\pi_{W^+}^\perp\otimes\pi_E^\perp)(\Phi\otimes\Phi^*)
&= \pi_{W^+}^\perp(\phi^a\otimes\phi^{b,*})
\otimes\pi_E^\perp(\xi^a\otimes\xi^{b,*}) \\
&= \quarter\tr(\phi^a\otimes\phi^{b,*})
\tr(\xi^a\otimes\xi^{b,*})\,\id_{W^+\otimes E} \\
&= \quarter\langle\phi^a,\phi^b\rangle \langle\xi^a,\xi^b\rangle
\,\id_{W^+\otimes E} \\
&= \quarter\langle\phi^a\otimes\xi^a,\phi^b\otimes\xi^b\rangle
\,\id_{W^+\otimes E} 
= \quarter\langle\Phi,\Phi\rangle\,\id_{W^+\otimes E}.
\end{align*}
Using the last identity together with
$(\Phi\otimes\Phi^*)\Phi = \Phi\langle\Phi,\Phi\rangle
= \Phi|\Phi|^2$ and \eqref{eq:doubletrace} gives
\begin{align*}
\langle(\Phi\otimes\Phi^*)_{00}\Phi,\Phi\rangle
&= \langle (\Phi\otimes\Phi^*)\Phi,\Phi\rangle 
- \langle (\pi_{W^+}^\perp\otimes\pi_E + \pi_{W^+}\otimes\pi_E^\perp)
(\Phi\otimes\Phi^*)\Phi,\Phi\rangle \\
&\qquad + \langle \quarter|\Phi|^2\Phi,\Phi\rangle \\
&= \fivequarter|\Phi|^4  
- \langle (\pi_{W^+}^\perp\otimes\pi_E + \pi_{W^+}\otimes\pi_E^\perp)
(\Phi\otimes\Phi^*)\Phi,\Phi\rangle. 
\end{align*} 
But $\pi_{W^+}^\perp\otimes\pi_E + \pi_{W^+}\otimes\pi_E^\perp$ is just the
orthogonal projection onto to the middle two factors of the orthogonal
decomposition of $\gl(W^+)\otimes_\CC\gl(E)$, so
\begin{align*}
&\left|\langle (\pi_{W^+}^\perp\otimes\pi_E + \pi_{W^+}\otimes\pi_E^\perp)
(\Phi\otimes\Phi^*)\Phi,\Phi\rangle\right| \\
&\qquad \le \left|(\pi_{W^+}^\perp\otimes\pi_E 
+ \pi_{W^+}\otimes\pi_E^\perp)(\Phi\otimes\Phi^*)\right||\Phi|^2 \\
&\qquad \le |\Phi\otimes\Phi^*||\Phi|^2 \le (1/\sqrt{2})|\Phi|^4
< \threequarter|\Phi|^4,
\end{align*} 
and therefore
$$
|\langle(\Phi\otimes\Phi^*)_{00}\Phi,\Phi\rangle|
\ge \fivequarter|\Phi|^4 - (1/\sqrt{2})|\Phi|^4 > \half|\Phi|^4.
$$
Consequently, using $(\Phi\otimes\Phi^*)_{00}^\dagger
= (\Phi\otimes\Phi^*)_{00}$, we have
\begin{align*}
\half\langle(\Phi\otimes\Phi^*)_{00}\Phi,\Phi\rangle
&= \langle(\Phi\otimes\Phi^*)_{00},\Phi\otimes\Phi^*\rangle 
= \langle(\Phi\otimes\Phi^*)_{00},(\Phi\otimes\Phi^*)_{00}\rangle \\
&= |(\Phi\otimes\Phi^*)_{00}|^2,
\end{align*} 
which gives (3).
\end{pf}

The fact that $|(\Phi\otimes\Phi^*)_0|^2 \ge \quarter|\Phi|^4$, for
$\Phi\in\Om^0(W^+\otimes E)$, is used in \S \ref{sec:compact} to show
that the moduli space of $\PU(2)$ monopoles has an Uhlenbeck
compactification. (The analogous equality, $|(\phi\otimes\phi^*)_0|^2 =
\quarter|\phi|^4$ for $\phi\in\Om^0(W^+)$, is used in \cite{KrM,Witten}
to show that the moduli space of Seiberg-Witten monopoles is compact.)

We adapt some terminology from the proofs of Proposition 3.4 in \cite{FU}
and Lemma 4.3.25 in \cite{DK} which we will need for our proof of
transversality for the $\PU(2)$ monopole moduli space in \S
\ref{sec:transv}:

\begin{defn}\label{defn:RankforTransv}
The {\em rank\/} of a section $\Phi\in\Omega^0(W^\pm\otimes E)$ {\em at a
point\/} $x\in X$ is the rank of $\Phi|_x$ considered as a complex linear
map $(W^\pm)^*|_x\to E|_x$ (that is, its rank as a complex two-by-two
matrix).  We say that $\Phi$ is {\em rank one\/} on a subset $U\subset X$
if it has rank less than or equal to one at all points in $U$.  Similarly,
for $v\in\Omega^0(\La^+(T^*X)\otimes\su(E))$, the {\em
rank\/} of $v$ {\em at a point\/} $x\in X$ is the rank of $v$ considered as
a real linear map from $\Lambda^+(TX)|_x\to \su(E)|_x$ (that is, its rank
as a real three-by-three matrix). The {\em rank\/} of $v$ on a set
$U\subset X$ is the maximum rank of $v|_x$ over all points $x\in U$.
Finally, for $M\in \Om^0(\su(W^+)\otimes\su(E))$, the {\em rank\/} of $M$
{\em at a point\/} $x\in X$ is the rank of $M$ considered as a real linear
map from $\su(W^+)^*|_x\to \su(E)|_x$ (again, its rank as a real
three-by-three matrix).  The {\em rank\/} of $M$ on a set $U\subset X$ is
the maximum rank of $M|_x$ over all points $x\in U$.
\end{defn}

For the proof of Lemma
\ref{lem:RankOfTau} below we shall need a simple linear algebra identity. 
Suppose that $\Phi=\phi\otimes\xi$, for some $\phi\in\Om^0(W^+)$ and
$\xi\in\Om^0(E)$, and that that $\Psi = \psi\otimes\zeta$, with
$\phi\in\Om^0(W^+)$ and $\xi\in\Om^0(E)$. Then
$$
\Phi\otimes\Psi^* 
= (\phi\otimes\xi)\otimes(\psi\otimes\zeta)^*
= (\phi\otimes\psi^*)\otimes(\xi\otimes\zeta^*) 
\in \gl(W^+)\otimes_\CC\gl(E).
$$
Writing $B\in \gl(W^+)$ (or $\gl(E)$) as $B=B_h+B_s$, where
$B_h=\half(B+B^\dagger)$ is Hermitian and $B_s = \half(B-B^\dagger)$ is
skew-Hermitian, we have
\begin{align*}
\Phi\otimes\Psi^* 
&= \quarter(\phi\otimes\psi^* + \psi\otimes\phi^*)
\otimes(\xi\otimes\zeta^* + \zeta\otimes\xi^*) \\
&\quad + \quarter(\phi\otimes\psi^* - \psi\otimes\phi^*)
\otimes(\xi\otimes\zeta^* - \zeta\otimes\xi^*) \\
&\quad + \quarter(\phi\otimes\psi^* + \psi\otimes\phi^*)
\otimes(\xi\otimes\zeta^* - \zeta\otimes\xi^*) \\
&\quad + \quarter(\phi\otimes\psi^* - \psi\otimes\phi^*)
\otimes(\xi\otimes\zeta^* + \zeta\otimes\xi^*),
\end{align*}
and similarly for $\Psi\otimes\Phi^*$. Therefore,
\begin{align}
\Phi\otimes\Psi^* + \Psi\otimes\Phi^*
&= \half(\phi\otimes\psi^* + \psi\otimes\phi^*)
\otimes(\xi\otimes\zeta^* + \zeta\otimes\xi^*) \notag\\
&\quad + \half(\phi\otimes\psi^* - \psi\otimes\phi^*)
\otimes(\xi\otimes\zeta^* - \zeta\otimes\xi^*) \notag\\
\label{eq:OffDiag}
&= -\half i(\phi\otimes\psi^* + \psi\otimes\phi^*)
\otimes i(\xi\otimes\zeta^* + \zeta\otimes\xi^*) \\
&\quad + \half (\phi\otimes\psi^* - \psi\otimes\phi^*)
\otimes (\xi\otimes\zeta^* - \zeta\otimes\xi^*).
\notag
\end{align}
We shall need the following elementary observation in our proof of
transversality for the moduli space of $\PU(2)$ monopoles in \S
\ref{sec:transv}. Note that as $\rho: \Lambda^+\to \su(W^+)$ is an
isomorphism, the rank of a section $v$ of $\La^+\otimes\su(E)$ is equal to
the rank of the section $\rho(v)$ of $\su(W^+)\otimes\su(E)$.

\begin{lem}
\label{lem:RankOfTau}
If $\Phi \in C^0(W^+\otimes E)$ and $x\in X$, then the following hold: 
\begin{enumerate}
\item $\Rank_\RR(\Phi(x)\otimes\Phi^*(x))_{00} = 1$ 
if and only if $\Rank_\CC\Phi(x) = 1$,
\item $\Rank_\RR(\Phi(x)\otimes\Phi^*(x))_{00} = 3$ 
if and only if $\Rank_\CC\Phi(x) = 2$.
\end{enumerate}
\end{lem}

\begin{pf}
For convenience we write $\Phi$, $W^+$, and
$E$ for $\Phi|_x$, $W^+|_x$, and $E|_x$.  If $\Phi\in\Hom(W^+,E)$ has 
complex rank one, we can write $\Phi=\phi\otimes \xi$, for $\phi\in W^+$ and
$\xi\in E$.  Then $(\Phi\otimes\Phi^*)_{00}
=-i(\phi\otimes\phi^*)_0\otimes i(\xi\otimes \xi^*)_0$ and so
$(\Phi\otimes\Phi^*)_{00}\in\Hom(\su(W^+),\su(E))$ has real rank one. 

Conversely, suppose $(\Phi\otimes\Phi^*)_{00}$ has real rank one. Let
$\{\xi_1,\xi_2\}$ be an orthonormal basis for $E$ and write $\Phi =
\phi_1\otimes \xi_1+ \phi_2\otimes \xi_2$: if $\phi_1,\phi_2\in E$ are
linearly dependent then $\Phi$ has complex rank one and we are done, so
suppose that $\phi_1,\phi_2$ are linearly independent. Using the standard
basis $\{\si_1,\si_2,\si_3\}$ of Pauli matrices \eqref{eq:Pauli} for
$\su(2)$, namely 
\begin{equation}
\si_1 = \left(\begin{matrix} i & 0 \\ 0 & -i \end{matrix}\right), \quad
\si_2 = \left(\begin{matrix} 0 & 1 \\ -1 & 0 \end{matrix}\right), \quad
\si_3 = \left(\begin{matrix} 0 & i \\ i & 0 \end{matrix}\right),
\label{eq:Pauli}
\end{equation}
we obtain bases for the real vector spaces $\su(W^+)$ and $\su(E)$:
\begin{align*}
\si^e_1&=i(\xi_1\otimes \xi_1^*-\xi_2\otimes \xi_2^*), \quad
\si^e_2=\xi_1\otimes \xi_2^*-\xi_2\otimes \xi_1^*, \quad  
\si^e_3=i(\xi_1\otimes \xi_2^*+\xi_2\otimes \xi_1^*), \\
\si^w_1&=i(\phi_1\otimes \phi_1^*-\phi_2\otimes \phi_2^*), \quad
\si^w_2=\phi_1\otimes \phi_2^*-\phi_2\otimes \phi_1^*, \quad  
\si^w_3=i(\phi_1\otimes \phi_2^*+\phi_2\otimes \phi_1^*).
\end{align*}
Writing $\Phi = \Phi_1+\Phi_2$, we have
$(\Phi_a\otimes\Phi_a^*)_{00} = 
-i(\phi_a\otimes\phi_a^*)_0\otimes i(\xi_a\otimes\xi_a^*)_0$ 
for $a=1,2$ and
$$
(\Phi\otimes\Phi^*)_{00} = (\Phi_1\otimes\Phi_1^*)_{00} + 
(\Phi_1\otimes\Phi_2^*
+ \Phi_2\otimes\Phi_1^*)_{00} + (\Phi_2\otimes\Phi_2^*)_{00}.
$$
So, by the identity \eqref{eq:OffDiag} and noting that 
$i(\xi_1\otimes\xi_1^*)_0 = \half i(\xi_1\otimes\xi_1^*
- \xi_2\otimes\xi_2^*) = \half \si^e_1 = -i(\xi_2\otimes\xi_2^*)_0$,
and $i(\phi_1\otimes\phi_1^*)_0 = \half i(\phi_1\otimes\phi_1^*
- \phi_2\otimes\phi_2^*) = \half \si^w_1 = -i(\phi_2\otimes\phi_2^*)_0$
we see that
\begin{align*}
(\Phi\otimes\Phi^*)_{00}
&= -i(\phi_1\otimes\phi_1^*)_0\otimes i(\xi_1\otimes\xi_1^*)_0 
-i(\phi_2\otimes\phi_2^*)_0\otimes i(\xi_2\otimes\xi_2^*)_0 \\
&\quad -\half i(\phi_1\otimes\phi_2^* + \phi_2\otimes\phi_1^*)
\otimes i(\xi_1\otimes\xi_2^* + \xi_2\otimes\xi_1^*) \\
&\quad + \half (\phi_1\otimes\phi_2^* - \phi_2\otimes\phi_1^*)
\otimes (\xi_1\otimes\xi_2^* - \xi_2\otimes\xi_1^*) \\
&= -\quarter \si^w_1\otimes\si^e_1 - \quarter \si^w_1\otimes\si^e_1
-\half\si^w_3\otimes\si^e_3 + \half \si^w_2\otimes\si^e_2 \\
&= -\half \si^w_1\otimes\si^e_1
+ \half \si^w_2\otimes\si^e_2 - \half\si^w_3\otimes\si^e_3.
\end{align*}
Thus, $(\Phi\otimes\Phi^*)_{00}$ would have real rank three, a
contradiction, and so $\Phi$ must have complex rank one, which proves
Assertion 1. 

The preceding argument also shows that if $\Phi$ has complex rank two, then 
$(\Phi\otimes\Phi^*)_{00}$ has real rank three. Conversely, if
$(\Phi\otimes\Phi^*)_{00}$ has real rank three, then $\Phi$ cannot have
complex rank one by Assertion 1, so $\Phi$ has complex rank two.
\end{pf}

\subsection{The PU(2) monopole equations and holonomy perturbations}
\label{subsec:MonopoleEqns}
In \S \ref{subsubsec:UnpertMonopoleEqns} we describe the $\PU(2)$
monopole equations in their unperturbed form, following \cite{PTCambridge,
PTLocal}. In \S \ref{subsubsec:HolonomyPerturbations} we introduce the
holonomy perturbations and the perturbed $\PU(2)$
monopole equations, deferring a detailed discussion of most of the
technical regularity issues concerning holonomy perturbations to the Appendix.
In Donaldson's application to the 
extended anti-self-dual equation, some important
features ensure that the requisite analysis is relatively tractable: (i)
reducible connections can be excluded from the compactification of the
extended moduli spaces \cite[p. 283]{DonPoly}, (ii) the cohomology groups
for the elliptic complex of his extended equations have simple weak
semi-continuity properties with respect to Uhlenbeck limits
\cite[Proposition 4.33]{DonPoly}, and (iii) the zero locus being perturbed
is cut out of a finite-dimensional manifold \cite[p. 281, Lemma 4.35, \&
Corollary 4.38]{DonPoly}. For the development of Donaldson's method for
$\PU(2)$ monopoles described here, none of these simplifying features hold
and so the corresponding transversality argument is rather
complicated. Indeed, one can see from Proposition 7.1.32 in
\cite{DK} that because of the Dirac operator, the behavior of the cokernels
of the linearization of the $\PU(2)$ monopole equations can be quite
involved under Uhlenbeck limits. The holonomy perturbations considered by
Donaldson in \cite{DonPoly} are inhomogeneous, as he uses the perturbations
to kill the cokernels $\Coker d_A^+$ directly.  In contrast, the
perturbations we consider in \eqref{eq:IntroPT} are homogeneous and we
shall argue indirectly in \S \ref{sec:transv}
that the cokernels of the linearization vanish away
from the reducible and zero-section solutions.

\subsubsection{The unperturbed $\PU(2)$ monopole equations}
\label{subsubsec:UnpertMonopoleEqns}
Recall that we consider Hermitian two-plane bundles $E$ over $X$ whose
determinant line bundles $\det E$ are isomorphic to a fixed Hermitian line
bundle over $X$ endowed with a fixed $C^\8$, unitary connection $A_e$.  Let
$(\rho,W^+,W^-)$ be a \spinc structure on $X$, where $\rho:T^*X\to\End W$
is the Clifford map, and the Hermitian four-plane bundle $W=W^+\oplus W^-$
is endowed with a $C^\8$ \spinc connection.  Given a connection $A$ on $E$
with curvature $F_A\in L^2_{k-1}(\La^2\otimes\fu(E))$, then $(F_A^+)_0 \in
L^2_{k-1}(\La^+\otimes\su(E))$ denotes the traceless part of its self-dual
component. Equivalently, if $A$ is a connection on $\su(E)$ with curvature
$F_A\in L^2_{k-1}(\La^2\otimes\so(\su(E)))$, then $F_A^+ \in
L^2_{k-1}(\La^+\otimes\su(E))$ is its self-dual component, viewed as a
section of $\La^+\otimes\su(E)$ via the implicit isomorphism
$\ad:\su(E)\to\so(\su(E))$.  Let $D_A:L^2_k(W^+\otimes E)\to
L^2_{k-1}(W^-\otimes E)$ be the corresponding Dirac operator.

For an $L^2_k$ section $\Phi$ of $W^+\otimes E$, let $\Phi^*$ be its
pointwise Hermitian dual and let $(\Phi\otimes\Phi^*)_{00}$ be the
component of the Hermitian endomorphism $\Phi\otimes\Phi^*$ of $W^+\otimes
E$ which lies in $\su(W^+)\otimes\su(E)$. The Clifford map $\rho$
defines an isomorphism $\rho:\La^+\to\su(W^+)$ and thus an isomorphism
$\rho=\rho\otimes\id_{\su(E)}$ of $\La^+\otimes\su(E)$ with
$\su(W^+)\otimes\su(E)$. Then
\begin{align}
F_A^+ - \rho^{-1}(\Phi\otimes\Phi^*)_{00} &= 0, 
\label{eq:UnpertPT}\\
D_A\Phi &= 0, \notag
\end{align}
are the unperturbed equations considered in
\cite{OTVortex,OTQuaternion,PTCambridge,PTLocal} (with slightly differing
trace conditions) for a pair $(A,\Phi)$ consisting of a connection $A$ on
$\su(E)$ and a section $\Phi$ of $W^+\otimes E$. Equivalently, given a pair
$(A,\Phi)$ with $A$ a fixed-determinant connection on $E$, the equations
\eqref{eq:UnpertPT} take the same form except that $F_A^+$ is replaced
by $(F_A^+)_0$.  

\subsubsection{The perturbed $\PU(2)$ monopole equations}
\label{subsubsec:HolonomyPerturbations}
We next introduce perturbations of the $\PU(2)$ monopole equations 
\eqref{eq:UnpertPT} which will
enable us to prove the transversality result described in Theorem
\ref{thm:Transversality}. The perturbations in question
make use of holonomy and their construction is partly modelled on related
perturbations introduced by Donaldson, Floer and Taubes
\cite{BraamDon,DonOrient,Floer,DonPoly,TauCasson}. All of these
constructions require 
that a unitary connection $A$ on a Hermitian two-plane bundle $E$ over $X$
be regular enough that parallel translation along a $C^\8$ path in $X$ is
well-defined. For example, Floer employs a configuration space of $L^4_1$
connections modulo $L^4_2$ gauge-transformations over a three-manifold, so
these connections are at least continuous (as $L^4_1\subset C^0$ in
dimension three) and the standard theory of ordinary differential equations
applies to define parallel translation without further qualification
\cite{Floer}. In dimension four, parallel
translation can be defined without difficulty using $L^2_k$ connections
with $k\ge 3$. With a little more care one can see that parallel translation
can be defined by $L^2_2$ connections, even though $L^2_2\not\subset C^0$
\cite[\S 3.2]{MorganGTNotes}. These regularity issues become rather more
intractable, as one can see from the Sobolev restriction theorem, for
connections which are only $L^2_1$ \cite[Theorem V.5.4]{Adams}.  As
explained in the Appendix --- to which we refer the reader for a discussion of
the more technical points --- we shall ultimately restrict our attention to
configuration spaces of $L^2_k$ connections with $k\ge 3$.  The
perturbations will, by definition, be zero on a neighborhood of a point in
$X$ where the curvature is large and so our regularity theory for $L^2_1$
monopoles will apply near points where curvature has bubbled off in order
to prove removability of singularities (see \S
\ref{subsubsec:DefnPertUhlCompact}).

We follow standard convention by saying that a connection $A$ on a $G$
bundle $E$ over a connected manifold $Y$ is `irreducible' if its stabilizer
$\Stab_A$ is trivial, that is, the center of the Lie group $G$
\cite[p. 133]{DK}, rather than saying (more correctly)
that its holonomy group $\Hol_A(y_0)$
is not a proper Lie subgroup of $\Aut(E_{y_0}) \simeq G$ (where $y_0\in Y$
is any basepoint). However, the holonomy will be our primary concern in
this section, so some care is required as the two notions do not coincide
in general. Recall from \cite[Lemma 4.2.8]{DK} that $\Stab_A$ is isomorphic
to the centralizer of $\Hol_A(y_0)$ in $G$.  If $Y$ is simply connected,
then $\Hol_A(y_0)$ is a connected Lie subgroup of $G$ \cite[Theorem
II.4.2]{KobNom}. Thus, if $G = \SU(2)$ or $\SO(3)$ and $Y$ is simply
connected, then $A$ has trivial stabilizer in $G$ if and only if
$\Hol_A(y_0) \simeq G$
\cite[p. 133]{DK}. Indeed, if $G = \SO(3)$ and $\Stab_A =
\{\id\}$, then $H := \Hol_A(y_0) = \SO(3)$; otherwise, we would
have
\begin{itemize}
\item $H = \SO(2)$, with centralizer $Z(H) = \SO(2)$ (by \cite[Theorem
IV.2.3(ii)]{BtD}), contradicting $\Stab_A =\{\id\}$, or
\item $H = \{\id\}$, with centralizer $Z(H) = \SO(3)$, again
contradicting $\Stab_A = \{\id\}$.
\end{itemize}
The same argument holds for $G = \SU(2)$, but not for higher-dimensional
Lie groups. For example, if $G = \U(2) = \SU(2)\times_{\{\pm\id\}}S^1$, we
cannot exclude the possibility that a $\U(2)$ connection $A$ with $\Stab_A
= S_Z^1$ reduces to an $\SU(2)$ connection, that is, $A$ has holonomy
$\Hol_A(y_0) = \SU(2)$.

With the preceding comments in mind, the first ingredient in our
construction of the holonomy perturbations is a local section of $\su(E)$
given by the holonomy of an $\SO(3)$ connection $A$ on $\su(E)$. Let
$\gamma\subset X$ be a $C^\8$ loop based at a point $x_0\in X$ and let
$$
h_{\ga,x_0}(A)\in\SO(\su(E))|_{x_0}
$$ 
be the holonomy of the connection $A$ around the loop $\gamma$. The
exponential map $\exp:\so(3)\to\SO(3)$ gives a diffeomorphism from a ball $2B$
around the origin in $\so(3)$ to a ball around the identity in $\SO(3)$. Let
$\psi:\RR\to[0,1]$ be a $C^\8$ cutoff function such that
$\psi(|\zeta|)=1$ for $\zeta\in \half B$,
$\psi(|\zeta|)>0$ for $\zeta\in B$, and $\psi(|\zeta|)=0$ for $\zeta\in
\so(3)-B$. Then
\begin{equation}
\fh_{\gamma,x_0}(A) 
:= 
\psi\left(|\exp^{-1}(h_{\ga,x_0}(A))|\right)
\cdot\ad^{-1}(\exp^{-1}(h_{\ga,x_0}(A)))
\end{equation}
defines a gauge-equivariant map $\sA_E(X)\to\su(E)|_{x_0}$, where
$\ad:\su(E) \to \so(\su(E))$ is the standard isomorphism.

\begin{lem}
\label{lem:DenseHolonomy}
Let $U\subset X$ be a simply connected open subset and let $x_0$ be a point
in $U$.  If $A|_U$ is an irreducible $\SO(3)$ connection on $\su(E)|_U$ then
there are loops $\{\gamma_l\}_{l=1}^3
\subset U$, depending on $A$, such that the set
$\{h_{\gamma_l,x_0}(A)\}_{l=1}^3\subset \SO(\su(E))|_{x_0}$ lies in the
open subset of $\SO(\su(E))|_{x_0}$ given by the image under $\exp$ of the ball
$B\subset\so(\su(E))|_{x_0}$ around the identity. The set
$\{\fh_{\gamma_l,x_0}(A)\}_{l=1}^3$ is then a basis for $\su(E)|_{x_0}$.
\end{lem}

\begin{pf}
Since $A|_U$ is an irreducible $\SO(3)$ connection over a simply connected
manifold $U$, the holonomy group $\Hol_{A|U}(x_0)$ is equal to
$\SO(\su(E))|_{x_0}$ by the remarks preceding the
statement of the lemma. Hence, there are three loops
$\gamma_1,\gamma_2,\gamma_3$ such that the holonomies $h_{\gamma_l,x_0}(A)\in
B\subset \SO(\su(E))|_{x_0}$ give a basis
$\{\fh_{\gamma_l,x_0}(A)\}_{l=1}^3$ for $\su(E)|_{x_0}$.
\end{pf}

For a $C^\8$ connection $A$ we may extend $\fh_{\gamma,x_0}(A)$ to a $C^\8$
section $\hat\fh_{\gamma}(A)$ of $\su(E)$ by radial parallel translation,
with respect to $A$ over a small ball $B(x_0,2R_0)$ and then multiplying by
a $C^\8$ cutoff function $\varphi$ on $X$ which is positive on $B(x_0,R_0)$
and identically zero on $X-B(x_0,R_0)$. Thus if the set
$\{\fh_{\gamma_l,x_0}(A)\}_{l=1}^3$ spans $\su(E)|_{x_0}$, then
the set $\{\varphi\hat\fh_{\gamma_l}(A)|_y\}_{l=1}^4$ will span
$\su(E)|_y$ for $y\in B(x_0,R_0)$. The constant $R_0$ is chosen so that
$4R_0$ is smaller than the injectivity radius of $(X,g)$.

In general, for an $L^2_k$ connection $A$ with $k\ge 2$, the section
$\hat\fh_{\gamma}(A)$ will not be in $L^2_{k+1}$ and so we use the Neumann
heat operator, for fixed small $t>0$,
$$
K_t(A|_{B(x_0,2R_0)}):
L^2(B(x_0,2R_0),\su(E)) \to L^2_{k+1}(B(x_0,2R_0),\su(E)),
$$ 
as discussed in Appendix
\ref{app:Parallel} to construct an $L^2_{k+1}$ section
\begin{equation}
\fh_{\gamma}(A) := K_t(A|_{B(x_0,2R_0)})\hat\fh_{\gamma}(A)
\label{eq:HeatKernelSmoothedHolonomy}
\end{equation}
of $\su(E)$ over $B(x_0,2R_0)$ which converges to $\hat\fh_{\gamma}(A)$ in
$C^0(B(x_0,R_0))$ as $t\to 0$. Therefore, for small enough
$t=t(A)$, the set $\{\varphi\fh_{\gamma_l}(A)|_y\}_{l=1}^3$,
will span $\su(E)|_y$ for all $y\in B(x_0,R_0)$ just as before.

\begin{lem}
\label{lem:SpanFiber}
Let $k\ge 2$ be an integer and let $A_0$ be an $L^2_k$ unitary connection
on $E|_{B(x_0,2R_0)}$. Let $\ga_1,\ga_2,\ga_3$ be loops in $B(x_0,2R_0)$
based at $x_0$ such that $\{\fh_{\gamma_l,x_0}(A_0)\}_{l=1}^3$ spans
$\su(E)|_{x_0}$.  Then there is a positive constant
$\eps(A_0,\{\ga_l\})$ such that the following holds. If $A$ is an
$L^2_k$ unitary connection on $E|_{B(x_0,2R_0)}$ such that
$$
\|A-A_0\|_{L^2_{k,A_0}(B(x_0,2R_0))} < \eps,
$$
then the set $\{\fh_{\gamma_l,x_0}(A)\}_{l=1}^3$ spans $\su(E)|_{x_0}$.
\end{lem}

\begin{pf}
According to Lemma \ref{lem:HolonomyPointRegularity} the holonomy maps
$h_{\gamma,x_ 0}:\sA_E(X)\to \U(E)|_{x_0}$ are continuous and so
the maps $\fh_{\gamma_l,x_0}:\sA_E(X)\to \su(E)|_{x_0}$ are continuous on
small open neighborhoods of $A_0\in\sA_E(X)$.
Hence, for small enough $\eps$, the set
$\{\fh_{\gamma_l,x_0}(A)\}_{l=1}^3$ spans $\su(E)|_{x_0}$.
\end{pf}

We now specify the loops to be used in this construction, essentially
following the argument used in the proof of Lemma 2.5 in \cite{DonOrient}.
Let $\{B(x_j,4R_0)\}_{j=1}^{N_b}$, be a disjoint collection of open balls,
where $N_b$ is a fixed integer to be specified later (see \S
\ref{subsubsec:DefnPertUhlCompact}). According to the slice results of
\cite[p. 192]{DK}, \cite[Proposition 2.1]{TauL2}
for manifolds with boundary, the quotient space $\sB^*_E(B(x_j,2R_0))$ of
irreducible $L^2_{k}$ connections on $E|_{B(x_j,2R_0)}$ is a $C^\8$
manifold modelled on a separable Hilbert space. The quotient $L^2_{k}$
topology on $\sB^*_E(B(x_j,2R_0))$ is clearly metrizable and so Stone's
Theorem implies that $\sB^*_E(B(x_j,2R_0))$ is paracompact and thus admits
partitions of unity \cite[Corollary II.3.8]{Lang}.

For each $j=1,\dots,N_b$ and each point $[A_0]$ in
$\sB_E^*(B(x_j,2R_0))$, we can find loops
$\{\ga_{j,l,A_0}\}_{l=1}^3$, contained in $B(x_j,2R_0)$ and based at
$x_j$ such that $\{\fh_{\gamma_{j,l,A_0}}(A_0)\}_{l=1}^3$, spans
$\su(E)|_{x_j}$. For each such point $[A_0]$, Lemma
\ref{lem:SpanFiber} implies that there is an $L^2_{k,A_0}$ ball 
$$
B_{[A_0]}(\eps_{A_0}) := \{[A]\in \sB_E^*(B(x_j,2R_0)):
\dist_{L^2_{k,A_0}}([A],[A_0]) < \eps_{A_0}\}
$$
such that for all $[A]\in B_{[A_0]}(\eps_{A_0})$, the set
$\{\fh_{\gamma_{j,l,A_0}}(A)\}_{l=1}^3$ spans $\su(E)|_{x_j}$.  These
balls give an open cover of $\sB^*_E(B(x_j,2R_0))$ and hence
there is a locally finite refinement of this open cover,
$\{U_{j,\alpha}\}_{\alpha=1}^\infty$, and a positive partition 
$\chi_{j,\alpha}$ (see Proposition \ref{prop:PositivePartition})
subordinate to $\{U_{j,\alpha}\}_{\alpha=1}^\infty$ in the sense that
$$
\sum_{\alpha}\chi_{j,\alpha}[A]>0, \qquad
[A]\in \sB^*_E(B(x_j,2R_0)).
$$
Hence, for each $U_{j,\alpha}$, we obtain loops
$\{\gamma_{j,l,\alpha}\}_{l=1}^3\subset B(x_j,2R_0)$ such that for
all $[A]\in U_{j,\alpha}$, the sections $\fh_{\gamma_{j,l,\alpha}}(A)$
span $\su(E)|_{x_j}$. 

Let $\beta$ be a smooth cutoff function on $\RR$ such that $\beta(t) = 1$
for $t\le \half$ and $\beta(t)=0$ for $t\ge 1$, with $\beta(t)>0$ for
$t<1$. Then the $C^\8$ gauge-invariant maps
$\sA_E(X)\to \RR$, $A\mapsto \beta_j[A]$ given by
\begin{equation}
\beta_j[A]
:=
\beta\left(\frac{1}{\eps_0^2}\int_{B(x_j,4R_0)}
\beta\left(\frac{\dist_g(\cdot,x_j)}{4R_0}\right)|F_A|^2\,dV\right)
\label{eq:ConnEnergyCutoff}
\end{equation}
are zero when the energy of the connection $A$ is greater than or equal to
$\half\eps_0^2$ over a ball $B(x_j,2R_0)$. Here, $\eps_0$ is the constant of
Corollary \ref{cor:PTLocalCpt}.
Finally, we define $C^\8$ cutoff functions on $X$ by setting
\begin{equation}
\varphi_j(x) := \beta\left(\frac{\dist_g(x,x_j)}{R_0}\right), \qquad x\in X,
\end{equation}
so that $\varphi_j$ is positive on the ball $B(x_j,R_0)$ and zero on its
complement in $X$.

We can now define a gauge equivariant $C^\8$ map
$$
\sA_E(X)\to L^2_{k+1}(X,\su(E)), 
\qquad A\mapsto \fm_{j,l,\alpha}(A)
$$
by setting
\begin{equation}
\fm_{j,l,\alpha}(A)
:= \beta_j[A]\chi_{j,\alpha}[A|_{B(x_j,2R_0)}]
\varphi_j\fh_{\gamma_{j,l,\alpha}}(A). 
\label{eq:HolonomySection} 
\end{equation}
Thus at each point $A\in \sA_E(X)$ only a finite number of the
$\fm_{j,l,\alpha}(A)$ are non-zero and each map $\fm_{j,l,\alpha}$ is
$C^\8$ with uniformly bounded derivatives of all orders on
$\sA_E(X)$ (see Appendix \ref{app:HolonomyConvergence}).  

To define the perturbation of the Dirac operator in \eqref{eq:UnpertPT}, we
need elements of $\Omega^0(\Hom_\CC(W^+,W^-))$ and these are conveniently
provided by complex one-forms using the Clifford isomorphism
$$
\rho:\La^1_\CC\to \Hom_\CC(W^+,W^-),
$$
where $\La^1_\CC = \La^1\otimes_\RR\CC$.
These one-forms are parametrized by the following
Banach space: Let $\AAA$ be the index set $\{(j,l,\alpha):1\le j\le N_b,
1\le l\le 3, \alpha\in\NN\}$, let $\de = (\de_\al)_{\al=1}^\8 \in
\ell^1(\RR)$ be the sequence of positive weights described in Appendix
\ref{app:HolonomyConvergence}, let $r\ge k+1$, and let
\begin{equation}
\sP_\vartheta^r := \ell^1_{\de}(\AAA,C^r(X,\La^1_\CC))
\end{equation}
be the set of sequences $\vecvartheta := (\vartheta_{j,l,\alpha})$ in
$C^r(X,\La^1_\CC)$ such that
$$
\|\vecvartheta\|_{\ell^1_\delta(C^r(X))}
:= 
\sum_{j,l,\alpha}\delta_\alpha^{-1}
\|\vartheta_{j,l,\alpha}\|_{C^r(X)} < \8.
$$
Then $\sP_\vartheta^r$ is a Banach space with respect to
the above norm \cite[\S 1.7]{KadisonRingrose}.

\begin{rmk}
\begin{enumerate}
\item Note that we measure reducibility of connections and allow the loops
$\ga_{j,l,\alpha}$ to be contained in the {\em larger\/} balls
$B(x_j,2R_0)$ while the perturbations are supported on the {\em
smaller\/} balls $B(x_j,R_0)$.  Thus, if $A\in\sA_E(X)$ is a
connection such that $A|_{B(x_j,2R_0)}$ is reducible for all
$j\in\{1,\dots,N_b\}$ such that $\beta_j[A]>0$ and $\Phi\not\equiv 0$, 
then $[A,\Phi]$ cannot
be a point in $M_{W,E}^{*,0}$: our unique continuation result for $\PU(2)$
monopoles (see Theorem \ref{thm:LocalToGlobalReducible}) which are
reducible on an open subset of $X$ containing all balls
$\barB(x_j,R_0)$ with $\beta_j[A]>0$ would imply that $A$ is reducible
on all of $X$ and so $(A,\Phi)$ would be a reducible $\PU(2)$
monopole. Consequently, if $[A,\Phi]\in M_{W,E}^{*,0}$ then $A$ must be
irreducible on {\em some} ball $B(x_j,2R_0)$; otherwise, $A$ would be
reducible on all balls $B(x_j,2R_0)$ with $\beta_j[A]>0$ and so
reducible on $X$.
\item Note that the sections $\fm_{j,l,\alpha}$ are defined on the entire
space $\sA_E(X)$: They are zero for connections
$A\in\sA_E(X)$ which are reducible when restricted to $B(x_j,2R_0)$.
\item Our energy bound \eqref{eq:MonopoleEnergyBound}
for solutions $(A,\Phi)$ to the perturbed $\PU(2)$
monopole equations ensures that there is always at least one ball
$B(x_j,R_0)$ where the sections $\{\fm_{j,l,\alpha}(A)\}_{l=1}^3$
span $\su(E)|_{B(x_j,R_0)}$: see \S \ref{subsubsec:DefnPertUhlCompact}.
\end{enumerate}
\end{rmk}

Because only a finite number of the sections $\fm_{j,l,\alpha}(A)$
are non-zero at each $A\in\sA_E(X)$ the sum
\begin{equation}
\vecvartheta\cdot\vec\fm(A) 
:= \sum_{j,l,\alpha}\rho(\vartheta_{j,l,\alpha})
\otimes_\CC \fm_{j,l,\alpha}(A)
\end{equation}
gives a well-defined $L^2_{k+1}$ section of
$\Hom_\CC(W^+,W^-)\otimes_\CC\fsl(E)\simeq \La^1_\CC\otimes_\CC\fsl(E)$. 
Note that although each $\fm_{j,l,\alpha}(A)$ is a section of $\su(E)$,
the tensor product is over the complex numbers, where $\su(E)$ is contained
in the complex space $\fsl(E)=\su(E)\otimes_\RR \CC$. 
(Since $\fsl(E) = \su(E)\oplus i\su(E)$, then $\Hom_\CC(W^+,W^-)
\otimes_\CC\su(E) = \Hom_\CC(W^+,W^-)\otimes_\CC\gl(E)$
and similarly for $\gl(E)$.) The map
$$
\sA_E(X)\to L^2_{k+1}(X,\Hom_\CC(W^+,W^-)\otimes_\CC\fsl(E)), 
\qquad A\mapsto \vecvartheta\cdot\vecfm(A)
$$
is $C^\8$ and gauge equivariant and so defines a $C^\8$
section of the vector bundle
$$
\sA^*_E(X)\times_{\sG_E}
L^2_{k+1}(X,\Hom_\CC(W^+,W^-)\otimes_\CC\fsl(E))\to
\sB^*_E(X).
$$
By construction (see Appendix \ref{app:HolonomyConvergence})
the sum $\vecvartheta\cdot\vecfm(A)$ satisfies a $C^0$ estimate of
the form
\begin{equation}
\sup_{A\in\sA_E(X)}
\|\vecvartheta\cdot\vecfm(A)\|_{L^2_{k+1,A}(X)} 
\le C\|\vecvartheta\|_{\ell^1_\delta(C^r(X))} 
\quad\text{and}\quad
\|\vecvartheta\|_{\ell^1_\delta(C^r(X))}\le \eps_\vartheta,
\label{eq:BasicEstPert2}
\end{equation}
where $C=C(g,k)$ and
$\eps_\vartheta$ is a positive constant which we are free to specify.

We shall also need to construct a gauge equivariant $C^\8$ map from
$\sA^*_E(X)$ to $L^2_{k+1}(\gl(\La^+)\otimes_\RR\so(\su(E)))$. This map
will define a perturbation of the quadratic form in \eqref{eq:UnpertPT}
using the representation $\ad: \su(E)\rightarrow \so(\su(E))$ and is
parametrized by the Banach space
\begin{equation}
\sP_\tau^r
:=
\ell^1_{\delta}(\AAA,C^r(X,\gl(\La^+)))
\end{equation}
of sequences $\vectau := (\tau_{j,l,\alpha})$ in $C^r(X,\gl(\Lambda^+))$
such that
$$
\|\vectau\|_{\ell^1_\delta(C^r(X))}
:=
\sum_{j,l,\alpha}\delta_\alpha^{-1}\|\tau_{j,l,\alpha}\|_{C^r(X)}
<
\infty.
$$
Then the sum
\begin{equation}
\vectau\cdot\vecfm(A)
:= 
\sum_{j,l,\alpha}\tau_{j,l,\alpha}\otimes_\RR
\ad(\fm_{j,l,\alpha}(A))
\end{equation}
is pointwise finite and gives a well-defined $L^2_{k+1}$ section of 
$\gl(\Lambda^+)\otimes_\RR\so(\su(E))$. The map 
$$
\sA_E(X)\to L^2_{k+1}(X,\gl(\La^+)\otimes_\RR\so(\su(E))), 
\qquad A\mapsto \vectau\cdot\vecfm(A)
$$
is $C^\8$ and gauge equivariant and so defines a $C^\8$
section of the vector bundle
$$
\sA^*_E(X)\times_{\sG_E}
L^2_{k+1}(X,\gl(\La^+)\otimes_\RR\so(\su(E)))\to 
\sB^*_E(X).
$$
By construction (see Appendix \ref{app:HolonomyConvergence})
the sum $\vectau\cdot\vecfm(A)$ satisfies $C^0$ estimates
\begin{equation}
\sup_{A\in\sA_E(X)}
\|\vectau\cdot\vecfm(A)\|_{L^2_{k+1,A}(X)} 
\le C\|\vectau\|_{\ell^1_\delta(C^r(X))} 
\quad\text{and}\quad
\|\vectau\|_{\ell^1_\delta(C^r(X))}\le \eps_\tau,
\label{eq:BasicEstPert1} 
\end{equation}
where $C=C(g,k)$ and
$\eps_\tau$ is a positive constant which we are free to specify.

Our perturbed $\PU(2)$ monopole equations for a pair $(A,\Phi)$ on
$(\su(E),W^+\otimes E)$ then take the form
\begin{align}
\fs_1(A,\Phi) 
&:=
F^+_A
- \left(\id+\tau_0\otimes\id_{\su(E)}+\vectau\cdot\vec\fm(A)\right)
\rho^{-1}(\Phi\otimes\Phi^*)_{00}=0,
\label{eq:PT}\\
\fs_2(A,\Phi) 
&:=
D_A\Phi + \rho(\vartheta_0)\Phi + \vecvartheta\cdot\vec\fm(A)\Phi = 0, 
\notag
\end{align}
where $\tau_0\in C^r(X,\gl(\La^+))$ and $\vartheta_0 \in C^r(X,\La^1_\CC)$
are additional perturbation parameters.  For brevity, we shall
often denote $\vectau_A := \vectau\cdot\vec\fm(A)$ and $\vecvartheta_A :=
\vecvartheta\cdot\vec\fm(A)$.

\begin{rmk}
\begin{enumerate}
\item In \cite{FL2} we consider the question of transversality of the spaces of
reducible solutions or $\U(1)$ monopoles, which are identified with moduli
spaces of Seiberg-Witten monopoles. In particular, we show that the
moduli spaces of $\U(1)$ monopoles are cut out transversely for generic
$\tau_0\in\Om^0(\gl(\La^+))$.
\item If the section $\Phi$ is identically zero then the $\PU(2)$
monopole equations reduce to those for a projectively anti-self-dual
unitary connection $A$ on $E$. If the connection $A$ is reducible then all
of the sections $\fm_{j,l,\alpha}(A)$ are identically zero on $X$; the
perturbations $\vectau\cdot\vecfm(A)$ and $\vecvartheta\cdot\vecfm(A)$ in
\eqref{eq:PT} are then zero.
\item We emphasize that the
holonomy sections in the sequence $\vecfm(A)$ are always $L^2_{k+1}$ by
construction for $L^2_k$ connections $A$ (with $k\ge 2$)
and that the sums defining $\vectau\cdot\vecfm(A)$ and
$\vecvartheta\cdot\vecfm(A)$ are always finite for each $A$, although the
number of terms may tend to infinity as $A$ approaches a reducible
connection.
\end{enumerate}
\end{rmk}

The proof of the existence of an Uhlenbeck compactification for
$M_{W,E}$ in \S \ref{sec:compact}
requires the perturbations $\vectau\cdot\vecfm(A)$ and
$\vecvartheta\cdot\vecfm(A)$ to satisfy the following estimates in order to
obtain universal {\em a priori\/} $L^2_{1,A}$ bounds for $\Phi$ and $L^2$
bounds for $F_A$ (see Lemmas \ref{lem:L21AaprioriEstAPhi} and 
\ref{lem:L2aprioriEstFA}):
\begin{align}
\sup_{A\in\sA_E}\|\vectau\cdot\vecfm(A)\|_{L^\8(X)} &\le 1, 
\label{eq:CompactEstPertL21} \\
\sup_{A\in\sA_E}
\|\vecvartheta\cdot\vecfm(A)\|_{L^4(X)} &\le 1. \notag
\end{align}
To obtain the more delicate universal {\em a priori\/} $L^\8$
bounds for $\Phi$ and $F_A^+$ (see Lemmas \ref{lem:SigmaPhiPhi}
and \ref{lem:C0EstFAPhi}) required by our Uhlenbeck compactness argument in
\S \ref{sec:compact}, the
perturbations $\tau_0$, $\vectau\cdot\vecfm(A)$, and
$\vecvartheta\cdot\vecfm(A)$ must satisfy the following stronger estimates:
\begin{align}
\|\tau_0\|_{L^\8(X)}+ \sup_{A\in\sA_E}\|\vectau\cdot\vecfm(A)\|_{L^\8(X)} 
&\le \sixtyfourth, \label{eq:CompactEstPertC0}\\
\|\vartheta_0\|_{L^\8_1(X)} +
\sup_{A\in\sA_E}\|\vecvartheta\cdot\vecfm(A)\|_{L^\8_{1,A}(X)} 
&\le 1. \notag
\end{align}
We take a constant $1$ on the right-hand sides of
\eqref{eq:CompactEstPertL21} and the second inequality in
\eqref{eq:CompactEstPertC0}
for notational convenience only: these bounds do not need to be `small'.
There are continuous Sobolev embeddings $L^\8(X,\RR)\subset L^2_4(X,\RR)$
and $L^\8_1(X,\RR)\subset L^2_4(X,\RR)$ over a four-manifold $X$.  The
estimates \eqref{eq:CompactEstPertC0} then follow from the estimates
\eqref{eq:BasicEstPert2} and \eqref{eq:BasicEstPert1} when $k\ge 3$. {\em
Therefore, to obtain the inequalities \eqref{eq:CompactEstPertC0}, we
require that $k\ge 3$} in our configuration spaces of $L^2_k$ connections
$\sB_E$ and pairs $\sC_{W,E}$: see \S
\ref{app:HolonomyConvergence}. The bounds in \eqref{eq:CompactEstPertC0}
then follow for small enough choices of $\eps_\tau$ and $\eps_\vartheta$ in
the inequalities \eqref{eq:BasicEstPert2} and \eqref{eq:BasicEstPert1}. 

We shall need a slight generalization of Lemma \ref{lem:UnpertSigmaPhiPhi}
which applies to the perturbed quadratic form
$\tau\rho^{-1}(\Phi\otimes\Phi^*)_{00}$ when
$\tau\ne\id_{\La^+}$. Without loss of generality, we may restrict our
attention to automorphisms $\tau_\rho := \rho\tau\rho^{-1}$ of
$\su(W^+)$ such that $|\tau_\rho-\id_{\su(W^+)}| =
|\rho\tau\rho^{-1}-\id| < \eps_\tau$ and thus
$$
\langle\Phi_1,\Phi_2\rangle-\eps_\tau|\Phi_1||\Phi_2|
\le \langle\tau(\Phi_1),\Phi_2\rangle 
\le \langle\Phi_1,\Phi_2\rangle + \eps_\tau|\Phi_1||\Phi_2|
$$
for some universal positive constant $\eps_\tau\ll 1$ and all
$\Phi_1,\Phi_2\in\Om^0(W^+\otimes E)$.
Then Lemma \ref{lem:UnpertSigmaPhiPhi} yields the following bounds for the
perturbed quadratic form:

\begin{lem}\label{lem:SigmaPhiPhi}
There is universal positive constant $\eps_\tau\ll 1$ such that for all
$\|\tau-\id_{\La^+}\|_{C^0(X)}<\eps_\tau$ and
$\Phi\in\Om^0(W^+\otimes E)$ the following inequalities hold:
\begin{align}
\half|\Phi|^4 &\le \langle\tau_\rho(\Phi\otimes\Phi^*)_{00}\Phi,\Phi\rangle
\le \threequarter|\Phi|^4, \tag{1}\\
\quarter|\Phi|^4 &\le |\tau_\rho(\Phi\otimes\Phi^*)_{00}|^2
\le {\textstyle{\frac{3}{8}}}|\Phi|^4. \tag{2}
\end{align}
\end{lem}

\begin{rmk}
The constant $\eps_\tau=1/64$ will suffice. The constraint
$\|\tau-\id\|_{C^0(X)}<\eps_\tau$ is only used in \S \ref{subsec:BW} and
consequently in \S \ref{subsec:SeqCompact}, where we establish the
`universal' {\em a priori\/} bounds for $\PU(2)$ monopoles and prove the
existence of an Uhlenbeck compactification, respectively.
\end{rmk}

\subsection{The moduli space and the elliptic deformation complex}
\label{subsec:DeformationComplex}
We define the moduli space of $\PU(2)$ monopoles and compute the index of
its elliptic deformation complex.

For any integer $k\ge 2$, the 
$\PU(2)$ monopole equations \eqref{eq:PT} define a $C^\8$
map $\fs:=(\fs_1,\fs_2)$ of Hilbert manifolds,
\begin{equation}
\fs:\tsC_{W,E}\rightarrow
L^2_{k-1}(\La^+\otimes\su(E))\oplus L^2_{k-1}(W^-\otimes E),
\label{eq:DefnSMap}
\end{equation}
given by
$$
(A,\Phi)\mapsto 
\left(\begin{matrix}\fs_1(A,\Phi) \\ \fs_2(A,\Phi)\end{matrix}\right)
:= 
\left(\begin{matrix}F_A^+-
(\id+\tau_0\otimes\id_{\su(E)}
+\vectau\cdot\vecfm(A))\rho^{-1}(\Phi\otimes\Phi^*)_{00} \\
D_A\Phi + + \rho(\vartheta_0)\Phi + \vecvartheta\cdot\vecfm(A)\Phi
\end{matrix}\right).
$$
We can then define the moduli space of {\em $\PU(2)$ monopoles\/} 
by setting 
\begin{equation}
M_{W,E} := \{(A,\Phi)\in \tsC_{W,E}:
\fs(A,\Phi)=0\}/\ssG_E
\label{eq:ModuliSpace}
\end{equation}
and denote the moduli space of $\PU(2)$ monopoles which are neither reducible
nor zero-section pairs by
$$
M_{W,E}^{*,0} := M_{W,E}\cap\sC^{*,0}_{W,E}.
$$
We let $M^{\asd}_E := \{[A]:F_A^+=0\}\subset \sB_E$ be the moduli
space of anti-self-dual connections on $\su(E)$.

The $\SO(3)$ bundle $\su(E)$ has first Pontrjagin number 
and second Stiefel-Whitney class given by
\begin{equation}
\kappa := -\quarter p_1(\su(E)) = c_2(E)-\quarter c_1(E)^2
\quad\text{and}\quad
w_2(\su(E)) = c_1(E)\pmod{2},
\label{eq:PontrjaginNumber}
\end{equation}
and, for a connection $A$ on $\su(E)$, 
we have the following Chern-Weil integral identity,
\begin{equation}
-\quarter p_1(\su(E))
=
\frac{1}{8\pi^2}\int_X \left(|F_A^-|^2-|F_A^+|^2\right)\,dV 
\label{eq:CWIntegral}
\end{equation}
where we view $F_A$ as a section of $\La^2\otimes\su(E)$ via the
isomorphism $\ad:\su(E)\simeq\so(\su(E))$ (see \cite[\S 2.1.4]{DK}).

Since the map $\fs$ is $\ssG_E$-equivariant 
it defines a section of the Hilbert vector bundle 
$\fV$ over $\sC_{W,E}^{*,0}$ with total space
\begin{equation}
\label{eq:BanachVectorBundle}
\fV := \tsC_{W,E}^{*,0}\times_{\ssG_E}
\left(L^2_{k-1}(\La^+\otimes\su (E))\oplus 
L^2_{k-1}(W^-\otimes E)\right). 
\end{equation}
While the equation $\fs[A,\Phi]=0$ is, of course, defined on 
$\sC_{W,E}$,
$\fV$ does not extend from $\sC_{W,E}^{*,0}$
to a vector bundle over $\sC_{W,E}$. The
moduli space $M_{W,E}^{*,0}\subset \sC_{W,E}^{*,0}$ is then the
zero locus of the section $\fs$ of Hilbert vector bundle $\fV$ over 
the Hilbert manifold $\sC_{W,E}^{*,0}$: it 
will be a {\em regular submanifold\/} if $\fs$ vanishes
transversely, that is, if the differential
\begin{equation}
(D\fs)_{A,\Phi}:T_{[A,\Phi]}\sC_{W,E}^{*,0}\to 
L^2_{k-1}(\La^+\otimes\su (E))\oplus L^2_{k-1}(W^-\otimes E)
\label{eq:DifferentialSMap}
\end{equation}  
is surjective at all points $[A,\Phi]$ in $\fs^{-1}(0)\cap
\sC_{W,E}^{*,0}$; recall from \S \ref{subsec:ConfigSpace}
that the tangent space $T_{[A,\Phi]}\sC_{W,E}^{*,0}$ is
canonically identified with $\Ker d_{A,\Phi}^{0,*}\subset 
\tsC_{W,E}^{*,0}$.

Suppose $(A,\Phi)$ is a pair in $\sA_E\times\Om^0(W^+\otimes E)$.  
Recall from Proposition \ref{prop:GaugeGroup} that
the differential at the identity $\id_E\in \ssG_E$, 
of the map $\ssG_E\to \sA_E\times\Om^0(W^+\otimes E)$
given by $u\mapsto u(A,\Phi)$ is 
$$
\Om^0(\su(E))\oplus i\RR_Z
\to 
\Om^1(\su(E)) \oplus \Om^0(W^+\otimes E),
\quad
\zeta\mapsto -d_{A,\Phi}^0\zeta = (-d_A\zeta,\zeta\Phi).
$$
Similarly, the differential 
$(a,\phi)\mapsto d_{A,\Phi}^1(a,\phi) := (D\fs)_{A,\Phi}(a,\phi)$,
$$
\Om^1(\su(E)) \oplus \Om^0(W^+\otimes E)
\to
\Om^+(\su(E)) \oplus \Om^0(W^-\otimes E),
$$
of the map $\fs$ at the point $(A,\Phi)$ is given by
\begin{equation}
d_{A,\Phi}^1(a,\phi) :=   
\left(\begin{matrix} 
\begin{split}
d_A^+a - \vectau\cdot(\delta\vecfm)(a)(\Phi\otimes\Phi^*)_{00}
-&(\id+\tau_0\otimes\id_{\su(E)} \\
+\vectau\cdot\vecfm(A))
\rho^{-1}&(\Phi\otimes\phi^* + \phi\otimes\Phi^*)_{00} 
\end{split}
\\
D_A\phi + \rho(\vartheta_0)\phi + \vecvartheta\cdot\vecfm(A)\phi
+ \rho(a)\Phi+\vecvartheta\cdot(\delta\vecfm)(a)\Phi\end{matrix}\right),
\label{eq:DefnLinearization}
\end{equation}
where $\delta\vecfm = \delta\vecfm/\delta A$.
(In the sequels we also find it convenient to use $L_{A,\Phi}$
to denote the linearization of the map $\fs$ at the point $(A,\Phi)$.) 
The differential of the composition $\ssG_E\mapsto
\Om^+(\su(E))\oplus\Om^0(W^-\otimes E)$ given by
$$
u\mapsto u(A,\Phi)\mapsto \fs(u(A,\Phi)) = u(\fs(A,\Phi))
= (u\fs_1(A,\Phi)u^{-1},u\fs_2(A,\Phi))
$$
is then 
$$
\Om^0(\su(E))\oplus i\RR_Z 
\to 
\begin{matrix}\Om^+(\su(E)) \\
\oplus \\
\Om^0(W^-\otimes E)
\end{matrix},
\quad
\zeta\mapsto d_{A,\Phi}^1\circ d_{A,\Phi}^0\zeta
= \left([\zeta,\fs_1(A,\Phi)],\zeta\fs_2(A,\Phi)\right).
$$
Therefore, $d_{A,\Phi}^1\circ d_{A,\Phi}^0 = 0$ if and only if
$\fs(A,\Phi)=0$, that is, if  and only if 
$(A,\Phi)$ is a $\PU(2)$ monopole. 
Consequently, the sequence 
\begin{equation}
\begin{CD}
\Om^0(\su(E))\oplus i\RR_Z
@>{d_{A,\Phi}^0}>> 
\begin{matrix}
\Om^1(\su(E))\\
\oplus \\
\Om^0(W^+\otimes E)
\end{matrix}
@>{d_{A,\Phi}^1}>> 
\begin{matrix}\Om^+(\su(E)) \\
\oplus \\
\Om^0(W^-\otimes E)
\end{matrix}
\end{CD} 
\label{eq:DefComplex}
\end{equation}
is a complex if and only if $(A,\Phi)$ is a $\PU(2)$ monopole. The $L^2$
adjoint of $-d_{A,\Phi}^0$ is given by
\begin{equation}
-d_{A,\Phi}^{0,*}(a,\phi) = -d_A^*a + (\cdot\Phi)^*\phi.
\label{eq:CoulombDefn}
\end{equation}
The operator
\begin{equation}
\sD_{A,\Phi} := d_{A,\Phi}^{0,*} + d^1_{A,\Phi}: 
\begin{matrix}
\Om^1(\su(E))\\
\oplus \\
\Om^0(W^+\otimes E)
\end{matrix}
\longrightarrow
\begin{matrix}
\Om^0(\su(E))\oplus i\RR_Z \\
\oplus \\
\Om^+(\su(E)) \\
\oplus \\
\Om^0(W^-\otimes E) \\
\end{matrix}
\label{eq:RolledOperator}
\end{equation}
is elliptic (thus Fredholm)
and so \eqref{eq:DefComplex} is an {\em elliptic
deformation complex\/} for the $\PU(2)$ monopole equations
\eqref{eq:PT}, with cohomology groups
$$
H_{A,\Phi}^0 := \Ker d_{A,\Phi}^0, \quad
H_{A,\Phi}^1 := \Ker d_{A,\Phi}^1/\Imag d_{A,\Phi}^0, \quad\text{and}\quad
H_{A,\Phi}^2 := \Coker d_{A,\Phi}^1,
$$
analogous to the usual elliptic
deformation complex \cite[Eq. (4.2.26)]{DK}
for the anti-self-dual equation, $F_A^+=0$.
(An elliptic deformation complex for stable pairs for holomorphic
bundles is given by Bradlow and Daskalopoulos in \cite[\S 2]{BradlowDask}.)

Recall that $H_{A,\Phi}^0 = \Ker d_{A,\Phi}^0$ is just the Lie algebra
of the stabilizer $\Stab_{A,\Phi}$ of the point
$[A,\Phi]\in M_{W,E}$ and $H^1_{A,\Phi}$ is
the Zariski or formal tangent space. 
Thus, for any point $[A,\Phi]\in M_{W,E}^{*,0}$ we have
$H_{A,\Phi}^0 = 0$. If $H_{A,\Phi}^2=0$, then $\Coker (D\fs)_{A,\Phi}=0$ and so
$[A,\Phi]$ is a regular point of the zero locus of the section $\fs$ of
$\fV$ over $\sC_{W,E}^{*,0}$. So, if $H_{A,\Phi}^0 = 0$ and 
$H_{A,\Phi}^2 = 0$, then $[A,\Phi]$ is a smooth point of $M_{W,E}$
with tangent space
$\Ker\sD_{A,\Phi} = \Ker(d^{0,*}_{A,\Phi} + d^1_{A,\Phi})
= H_{A,\Phi}^1$,
as we see from \eqref{eq:DifferentialSMap}. Provided the zero set
$\fs^{-1}(0)$ is regular, then $M_{W,E}^{*,0}$ will be a smooth
manifold of dimension $-\Ind \sD_{A,\Phi}$. 

The perturbation terms in \eqref{eq:PT} define gauge equivariant maps
$$
\tsC_{W,E}^{*,0}
\to 
L^2_{k-1}(\La^+\otimes\su(E))\oplus L^2_{k-1}(W^-\otimes E),
$$
given by
$$
(A,\Phi)\mapsto 
\left(\begin{matrix}
\vectau\cdot\vecfm(A)\rho^{-1}(\Phi\otimes\Phi^*)_{00} \\
\vecvartheta\cdot\vecfm(A)\Phi
\end{matrix}\right).
$$
For $k\ge 2$, the Sobolev multiplication theorem implies that
$\vecvartheta\cdot\vecfm(A)\Phi$ is in $L^2_k(W^-\otimes E)$, while
$\vectau\cdot\vecfm(A)\rho^{-1}(\Phi\otimes\Phi^*)_{00}$ is in
$L^2_k(\La^+\otimes\su(E))$ when $k\ge 3$ and in
$L^p_2(\La^+\otimes\su(E))$, $1\le p <2$ when $k=2$. By the Rellich
embedding theorem, the inclusions $L^2_k\subset L^2_{k-1}$ and
$L^p_2\subset L^2_1$, $p>1$, are compact. In particular, it follows that
the linearization of the perturbed $\PU(2)$ monopole equations
\eqref{eq:DefnLinearization} differs from the
linearization of the unperturbed equations \eqref{eq:UnpertPT} by a compact
operator \cite[Theorem VI.2]{Adams}.

\begin{prop}\label{prop:DimOfModuliSpace}
If the map $\fs$ of \eqref{eq:DefnSMap} vanishes transversely
for the parameters $(\tau_0,\vartheta_0,\vectau,\vecvartheta)$ then
$M^{*,0}_{W,E}(\tau_0,\vartheta_0,\vectau,\vecvartheta)$ 
is a smooth manifold of dimension
\begin{align*} 
\dim M^{*,0}_{W,E}
&= -2p_1(\su(E))-3(1-b^1(X)+b^+(X))\\
&\quad + \textstyle{\frac{1}{2}}p_1(\su(E))
+\half(c_1(W^+)+c_1(E))^2 -\half\sigma(X) - 1\\
&=  \dim M^{\asd}_E+2\Ind_\CC D_A-1,
\end{align*}
where $\sigma(X)$ is the signature of $X$.
\end{prop}

\begin{pf}
Since $H^0_{A,\Phi}=0$ and $H^2_{A,\Phi}=0$ (by hypothesis)
at any point $[A,\Phi]$ in 
$M^{*,0}_{W,E}$, we have
$\dim M^{*,0}_{W,E} = -\Ind\sD_{A,\Phi}$.
By our regularity result, Proposition \ref{prop:GlobalReg}, any point
$[A,\Phi]$ in $M^{*,0}_{W,E}$ has a smooth representative
$(A,\Phi)$. Therefore, from the expressions for $d_{A,\Phi}^{0,*}$ in
\eqref{eq:CoulombDefn} and for $d_{A,\Phi}^1$ in
\eqref{eq:DefnLinearization} and the Sobolev multiplication and embedding
theorems, we find that the operator
$\sD_{A,\Phi}=d_{A,\Phi}^{0,*}+d_{A,\Phi}^1:L^2_k\to L^2_{k-1}$ differs
from 
$$ 
\left(\begin{matrix}d_A^* + d_A^+& 0 \\ 0 & D_A\end{matrix}\right):
\begin{matrix}
L^2_k(\La^1\otimes\su(E))
\\ \oplus \\ L^2_k(W^+\otimes E)
\end{matrix}
\rightarrow 
\begin{matrix}
L^2_{k-1}(\su(E))\oplus i\RR_Z\\ 
\oplus \\ 
L^2_{k-1}(\La^+\otimes\su(E))\\ 
\oplus \\
L^2_{k-1}(W^-\otimes E)
\end{matrix}
$$
by a compact operator and so has the same (real) index:
$$
\Ind\sD_{A,\Phi} = 
\Ind(d_A^* + d_A^+) + \Ind_\RR D_A - 1.
$$
We recall from \cite[Eq. (4.2.22)]{DK} that the operator
$$
d_A^*+ d_A^+:
L^2_k(\La^1\otimes\su(E))\rightarrow 
L^2_{k-1}(\su(E))\oplus L^2_{k-1}(\La^+\otimes\su(E))
$$
has index
$$
\Ind(d_A^* + d_A^+) = -2p_1(\su(E))-3(1-b^1(X)+b^+(X)).
$$
The complex index of the Dirac operator
$D_A$ is given by
\begin{align*}
\Ind_\CC D_A&=  \langle (\hat A(X))\ch(E)e^{\half c_1(W^+)},[X]\rangle 
\\
&=  \left\langle(1-\textstyle{\frac{1}{24}}p_1(X))(2+c_1(E)\right. \\
&\qquad \left. +\half c_1(E)^2-c_2(E))
(1+\half c_1(W^+)+\textstyle{\frac{1}{8}}c_1(W^+)^2),[X]\right\rangle
\\
&=  \langle -\textstyle{\frac{1}{12}}p_1(X)+\quarter c_1(W^+)^2+
  \half c_1(E)^2-c_2(E)+\half c_1(E)c_1(W^+),[X]\rangle
\\
&=  -\quarter\sigma(X)+\quarter c_1(W^+)^2+\half c_1(E)^2
  -c_2(E)+\half c_1(E)c_1(W^+) 
\\
&=  (\quarter c_1(E)^2-c_2(E))+\quarter((c_1(W^+)+c_1(E))^2-\sigma(X)) 
\\
&=  \quarter p_1(\su(E))+\quarter((c_1(W^+)+c_1(E))^2-\sigma(X)).
\end{align*}
Adding these indices and noting that 
$\Ind_\RR D_A=2\Ind_\CC D_A$
gives our dimension formula.
\end{pf}


\section{Regularity}
\label{sec:regularity}
Our goal in this section is to establish the basic regularity results for
solutions to the $\PU(2)$ monopole equations \eqref{eq:PT}. In \S
\ref{subsec:L2_1InhomoGlobal} 
we show that global $L^2_1$ solutions to the first-order elliptic system
comprising \eqref{eq:PT} and the Coulomb gauge equation 
(see \S \ref{subsec:ConfigSpace}) are
necessarily $C^\8$, while in \S \ref{subsec:Lp_kInhomoGlobal} we show
that global $L^2_k$ solutions to \eqref{eq:PT} are equivalent via an
$L^2_{k+1}$ determinant-one, unitary automorphism of $E$
--- provided $k\ge 2$ --- to $C^\8$ solutions to
\eqref{eq:PT} in Coulomb gauge relative to some $C^\8$ reference pair. 
In \S \ref{subsec:L2_1InhomoLocal} and \S \ref{subsec:Uhlenbeck} we
establish local versions of the results of 
\S \ref{subsec:L2_1InhomoGlobal} and \S \ref{subsec:Lp_kInhomoGlobal}. The
regularity results and estimates of this section will be needed repeatedly
throughout our work, so we state them here in sufficient generality to
cover all of these applications. In the present article, we require the
regularity results for our proof of transversality for the moduli space of
$\PU(2)$ monopoles (see \S \ref{sec:transv}) and they form the
cornerstone of our proofs of removable singularities and existence of an
Uhlenbeck compactification (see \S \ref{sec:compact}). 
Furthermore, in sequels to this article \cite{FL3,FL4}, the
regularity results of this section are used to show that
$L^2_1$ gluing solutions to \eqref{eq:PT} are necessarily
$C^\8$ and to analyse the Uhlenbeck-boundary behavior of the gluing and
obstruction maps parametrizing the moduli space ends. 

In order to simultaneously address all of the intended applications,  
the equations we find it convenient to consider here 
are a quasi-linear, inhomogeneous elliptic system consisting of a 
generalization of the equations \eqref{eq:PT} and Coulomb gauge equation
for a pair $(A,\Phi) \in \sA_E(X)\times\Om^0(W^+\otimes E)$.
Specifically, we allow inhomogeneous, right-hand
terms: the need for this generalization
arises in our proofs of removable singularities and of
regularity for gluing solutions to 
\eqref{eq:PT} and in analysing the Uhlenbeck-boundary behavior of gluing
maps.  Suppose that $(A_0,\Phi_0)$ is a fixed $C^\8$ reference pair in
$\sA_E(X)\times\Om^0(X,W^+\otimes E)$: writing $(A,\Phi) =
(A_0,\Phi_0)+(a,\phi)$, combining \eqref{eq:PT} with the Coulomb gauge
equation, and allowing inhomogenous terms, we obtain an elliptic system of
equations for a pair $(a,\phi)$ in $\Om^1(X,\su(E))\oplus\Om^0(X,W^+\otimes
E)$,
\begin{align}
d_{A_0}^*a - (\cdot\Phi_0)^*\phi &= \zeta, 
\label{eq:Coulomb+PT}\\
\fs(A_0+a,\Phi_0+\phi) &= (v_0,\psi_0), \notag
\end{align}
Considering $A$ to be a connection on $\su(E)$ and using the isomorphism
$\ad:\su(E)\to\so(\su(E))$ to view $F_A$ as a section of $\La^2\otimes\su(E)$,
we write \eqref{eq:Coulomb+PT} as
\begin{align*}
d_{A_0}^*a - (\cdot\Phi_0)^*\phi &= \zeta, \\
F_{A_0+a}^+ 
- \left(\id+\tau_0\otimes\id_{\su(E)}+\vectau\cdot\vecfm(A)\right)
\rho^{-1}((\Phi_0+\phi)\otimes(\Phi_0+\phi)^*)_{00} &=  v_0, \\ 
(D_{A_0+a}+\rho(\vartheta_0)
+\vecvartheta\cdot\vecfm(A))(\Phi_0+\phi) &= \psi_0,
\end{align*}
where $(v_0,\psi_0)\in\Om^+(X,\su(E))\oplus\Om^0(X,W^-\otimes E)$. 
Recalling that $d_{A_0,\Phi_0}^0$ and $d_{A_0,\Phi_0}^1$ are the
differential operators in the elliptic deformation complex
\eqref{eq:DefComplex} for the $\PU(2)$ monopole equations \eqref{eq:PT},
the above system may be rewritten in the form
\begin{align*}
d_{A_0,\Phi_0}^{0,*}(a,\phi) &= \zeta, \\
d_{A_0,\Phi_0}^1(a,\phi) + \{(a,\phi),(a,\phi)\} &= 
-\fs(A_0,\Phi_0) + (v_0,\psi_0) =: (v,\psi),
\end{align*}
where $(v,\psi) \in \Om^+(X,\su(E))\oplus \Om^0(X,W^-\otimes E)$ and
the differentials $d_{A_0,\Phi_0}^{0,*}$ and $d_{A_0,\Phi_0}^1$ are given by 
\eqref{eq:CoulombDefn} and \eqref{eq:DefnLinearization}.
It will be convenient to view the
quadratic term $\{(a,\phi),(a,\phi)\}$ as being defined via the following
bilinear form, 
$$
\{(a,\phi),(b,\varphi)\}
:= \left(\begin{matrix}
(a\wedge b)^+ 
- \left(\id + \tau_0\otimes\id_{\su(E)}+\vectau\cdot\vecfm(A)\right)
\rho^{-1}(\phi\otimes\varphi^*)_{00} \\
\rho(a)\varphi
\end{matrix}\right), 
$$
where $(b,\varphi)\in \Om^1(X,\su(E))\oplus\Om^0(X,W^+\otimes E)$; we will
further abbreviate $\{(a,\phi),(a,\phi)\}$ by $q(a,\phi)$ when convenient.
Our elliptic system \eqref{eq:Coulomb+PT} then takes the simple shape 
\begin{equation}
\sD_{A_0,\Phi_0}(a,\phi) + \{(a,\phi),(a,\phi)\} = (\zeta,v,\psi),
\label{eq:PTEllReg}
\end{equation}
recalling from \eqref{eq:RolledOperator} that $\sD_{A_0,\Phi_0} =
d_{A_0,\Phi_0}^{0,*} + d_{A_0,\Phi_0}^1$.  This is the form of the
(inhomogeneous) Coulomb gauge and $\PU(2)$ monopole equations we will
use for the majority of the basic regularity arguments. 

Some of the regularity results and estimates of this section generalize
corresponding results for the first-order anti-self-dual equation
\cite{DK,FU} and, to a certain extent, those of Uhlenbeck \cite{UhlRem} and
Parker \cite{Parker} for the second-order (coupled) Yang-Mills equations.
As usual for a quasi-linear, first-order elliptic system with a quadratic
non-linearity, over an $n$-dimensional manifold, 
the main difficulty is to prove regularity for $L^{n/2}_1$
solutions (that is, at the critical Sobolev exponent). Once the solutions
are known to be in $L^\8$, then standard linear elliptic regularity
theory can applied \cite{Hormander,Morrey}. It is worth noting at the
outset, though, that despite an extensive literature on quasi-linear
elliptic systems due to Ladyzhenskaya-Ural'tseva, Morrey and others
\cite{LU,Morrey} these classical results do not meet the usual
demands of gauge theory and this explains why we develop the precise
results we require here. The {\em constants\/}
appearing in our estimates generally depend on the underlying Riemannian
metric on $X$, the fixed \spinc connection on $W$, the fixed unitary
connection on $\det E$, and the perturbations $(\tau_0,\vectau,\vecvartheta)$: 
this dependence is not explicitly noted, as these ingredients in the
$\PU(2)$ monopole equations \eqref{eq:PT} are fixed once and for all.

\subsection{Regularity for $L^2_1$ solutions to the inhomogeneous
Coulomb gauge and PU(2) monopole equations}
\label{subsec:L2_1InhomoGlobal}
We show in this subsection that an $L^2_1$ solution $(a,\phi)$ to the
$\PU(2)$ monopole and Coulomb-gauge equations \eqref{eq:PTEllReg}, with an
$L^2_k$ inhomogeneous term (with $k \ge 1$) is in $L^2_{k+1}$. Thus, if the
inhomogeneous term is in $C^\8$ then $(a,\phi)$ is in $C^\8$. In passing
from an $L^2_1$ to an $L^2_2$ solution we need only consider the case where
$(\vectau,\vecvartheta) = 0$, while no restriction is placed on the
perturbation $(\vectau,\vecvartheta)\in\sP$ given an $L^2_2$ solution
$(a,\phi)$.

Our regularity result contains, as special cases, 
Proposition 4.4.13 in \cite{DK} and Theorem 8.8 in \cite{FU},
for anti-self-dual connections, and Theorem 8.11 in
\cite{Salamon} for Seiberg-Witten monopoles. The proof we give below for
$\PU(2)$ monopoles is rather different. We provide a fairly detailed
argument here since regularity of $L^2_2$ solutions to an elliptic equation
with a quadratic non-linear term does not quite follow directly from
standard theory for non-linear elliptic systems (for example, \cite[Theorem
6.8.1]{Morrey}). 

The two principal steps are, first, to get $L^p_1$ regularity of an $L^2_1$
solution $(a,\phi)$ when $2<p<4$ and $(a,\phi)$ is sufficiently $L^4$-small
and, second, to apply elliptic boostrapping and the Sobolev multiplication
and embedding theorems to get $C^\8$ regularity of an $L^p_1$ solution
$(a,\phi)$ when $2<p<\8$. We will use these sharp regularity results and
estimates repeatedly throughout our work, so we give the argument in some
detail here. The main ingredient in the first step is supplied by
Proposition \ref{prop:L2_1InhomoReg} and uses a Fredholm alternative
argument to pass from $L^2_1$ to the slightly stronger $L^p_1$ regularity
\cite[Theorem 5.3]{GT}. Although we will often be able to simply assume
that the inhomogenous term is in $C^\8$, rather than just in $L^p$, we will
later need these intermediate regularity results in our development of the
gluing theory for $\PU(2)$ monopoles \cite{FL3,FL4} to show that $L^2_1$
solutions $(a,\phi)$ to the system \eqref{eq:PTEllReg} with a certain
inhomogeneous term $(\zeta,v,\psi)$, where the latter is initially only
known to be in $L^2_1$ or $L^2_2$, are actually in $C^\8$ (together, of
course, with $(\zeta,v,\psi)$).

\begin{lem}\label{lem:OneSidedInverse}
Let $\fB_1$, $\fB_2$ be Banach spaces and let $T\in\Hom(\fB_1,\fB_2)$ have
a right (left) inverse $P$. If $T'\in\Hom(\fB_1,\fB_2)$ satisfies
$\|T'-T\|<\|P\|^{-1}$, then $T'$ also has a right (left) inverse.
\end{lem}

\begin{pf}
If $P\in\Hom(\fB_2,\fB_1)$ is a right inverse for $T$, so $TP=\id_1$, then
$\|(T'-T)P\|\le \|T'-T\|\|P\| < 1$ and $\id_1+(T'-T)P$ is an invertible element
of the Banach algebra $\End(\fB_1)$. Define 
$P'=P(1+(T'-T)P)^{-1}$, so $T'P'=\id_1$ and $P'$ is a right inverse for
$T'$. Similarly for left inverses.
\end{pf}

The preceding elementary consequence of the usual characterization of
invertible elements of a Banach algebra will be frequently invoked in this
section and in particular, in the proof of the proposition below.

\begin{prop}\label{prop:L2_1InhomoReg}
Let $X$ be a closed, oriented, Riemannian four-manifold with metric $g$,
\spinc structure $(\rho,W)$, and let $E$ be a Hermitian two-plane bundle
over $X$.  Let $(A_0,\Phi_0)$ be a $C^\8$ pair on the $C^\8$ bundles
$(\su(E),W^+\otimes E)$ over $X$ and let $2\le p < 4$.  Then there are
positive constants $\eps=\eps(A_0,\Phi_0,p)$ and $C=C(A_0,\Phi_0,p)$ with
the following significance. Suppose that $(a,\phi)\in
L^2_1(X,\La^1\otimes\su(E))\oplus L^2_1(X,W^+\otimes E)$ is an $L^2_1$
solution on $(\su(E),W^+\otimes E)$ to the elliptic system
\eqref{eq:PTEllReg} over $X$ with $(\vectau,\vecvartheta) = 0$, where
$(\zeta,v,\psi)$ is in $L^p$.  If $\|(a,\phi)\|_{L^4(X)}<\eps$ then
$(a,\phi)$ is in $L^p_1$ and
$$
\|(a,\phi)\|_{L^p_{1,A_0}(X)} \le C\left(\|(\zeta,v,\psi)\|_{L^p(X)}
+ \|(a,\phi)\|_{L^2(X)}\right).
$$
\end{prop}

\begin{pf}
The operator $\sD_{A_0,\Phi_0}$ is
Fredholm (since it is elliptic with $C^\8$ coefficients and $X$ is closed),
and so has a finite-dimensional kernel and cokernel. In particular, 
$\Ker\sD_{A_0,\Phi_0}|_{L^2_1}\subset C^\8$. Let $\Pi_1$
be the $L^2$-orthogonal projection onto 
$\Ker\sD_{A_0,\Phi_0}|_{L^2_1}$ and let $\Pi_2$
be the $L^2$-orthogonal projection onto 
$\Ker\sD_{A_0,\Phi_0}^*|_{L^2} = (\Imag\sD_{A_0,\Phi_0}|_{L^2_1})^\perp$, 
the $L^2$-orthogonal complement of the
image of $\Imag\sD_{A_0,\Phi_0}|_{L^2_1}$. 
The $L^2$-adjoint $\sD_{A_0,\Phi_0}^*$ is again a
first-order linear elliptic operator with $C^\8$ coefficients and so 
$\Ker\sD_{A_0,\Phi_0}^*|_{L^2}\subset C^\8$. 
We may then rewrite our
quasi-linear elliptic equation \eqref{eq:PTEllReg} in the form
\begin{align*}
&\Pi_2^\perp\sD_{A_0,\Phi_0}\Pi_1^\perp(a,\phi) 
+ \Pi_2^\perp \{(a,\phi),\Pi_1^\perp(a,\phi)\} \notag\\
&= -\Pi_2\left(\sD_{A_0,\Phi_0}(a,\phi) 
+ \{(a,\phi),\Pi_1^\perp(a,\phi)\}\right) \\
&\qquad - \{(a,\phi),\Pi_1(a,\phi)\} 
+(\zeta,v,\psi) =: \Upsilon. \notag 
\end{align*}
where $\Pi_i^\perp = \id - \Pi_i$ for $i=1,2$. Thus, we need to consider
the existence and uniqueness problem for solutions $(b,\varphi)$ to  
\begin{equation}
\Pi_2^\perp\sD_{A_0,\Phi_0}(b,\varphi) + \Pi_2^\perp\{(a,\phi),(b,\varphi)\}
= \Upsilon, \label{eq:InvertiblePTEllReg}
\end{equation}
where $(b,\varphi) \in L^p_1\cap (\Ker\sD_{A_0,\Phi_0})^\perp$. Plainly, the
operator
$$
\Pi_2^\perp\sD_{A_0,\Phi_0}: L^p_1\cap(\Ker\sD_{A_0,\Phi_0})^\perp 
\to L^p_1\cap(\Ker\sD_{A_0,\Phi_0}^*)^\perp 
$$
is (left and right) invertible and hence this
will also be true for any operator from $L^p_1\cap(\Ker\sD_{A_0,\Phi_0})^\perp$
to $L^p\cap(\Ker\sD_{A_0,\Phi_0}^*)^\perp$, such as 
$\Pi_2^\perp\sD_{A_0,\Phi_0} + \Pi_2^\perp\{(a,\phi),\ \cdot\ \}$,
which is sufficiently
close in the $\Hom(L^p_{1,A_0}, L^p)$ operator norm
by Lemma \ref{lem:OneSidedInverse}.
Since $\sD_{A_0,\Phi_0}$ and $\sD_{A_0,\Phi_0}^*$ are 
first-order linear elliptic operators with $C^\8$ coefficients and 
$\sD_{A_0,\Phi_0}\Pi_1=0$ and
$\sD_{A_0,\Phi_0}^*\Pi_2=0$, standard elliptic theory implies that
\begin{align}
\|\Pi_2(\xi,w,\varrho)\|_{L^r_{k,A_0}} 
&\le C\|(\xi,w,\varrho)\|_2, \qquad (\xi,w,\varrho)\in L^r_k, 
\label{eq:KernelEst} \\
\|\Pi_1(b,\varphi)\|_{L^r_{k,A_0}} 
&\le C\|(b,\varphi)\|_2, \qquad (b,\varphi)\in L^r_k, \notag
\end{align}
for some constant $C= C(A_0,\Phi_0,k,r)$, whenever $k\ge 1$ and $rk \ge
4/3$ or $k=0$ and $r\ge 2$.

Since $(a,\phi)$ is in $L^2_1$, the terms $\sD_{A_0,\Phi_0}(a,\phi)$ and
$\{(a,\phi),(a,\phi)\}$ and $\{(a,\phi),\Pi_1^\perp(a,\phi)\}$ are in
$L^2$, while the term $\{(a,\phi),\Pi_1(a,\phi)\}$ is in $L^2_1$. The terms
$\Pi_2\sD_{A_0,\Phi_0}(a,\phi)$ and $\Pi_2\{(a,\phi),\Pi_1^\perp(a,\phi)\}$
are each in $C^\8$, while $(\zeta,v,\psi)$ is in $L^p$, and so
the right-hand side $\Upsilon$ of \eqref{eq:InvertiblePTEllReg} is
in $L^p\cap(\Ker\sD_{A_0,\Phi_0}^*)^\perp$. 

Let $q=4p/(4-p)$, so $4<q<\8$; 
there is a continuous multiplication map
$L^4\times L^q\to L^p$ and an embedding $L^p_1\subset L^q$.
So, as $\Pi_2^\perp = \id - \Pi_2$, 
\begin{align*}
\|\Pi_2^\perp\{(a,\phi),(b,\varphi)\}\|_{L^p} 
&\le \|\{(a,\phi),(b,\varphi)\}\|_{L^p} 
+ C\|\{(a,\phi),(b,\varphi)\}\|_{L^2} \\
&\le C\|(a,\phi)\|_{L^4}\|(b,\varphi)\|_{L^p_{1,A_0}},
\end{align*}
for some positive constant $C=C(A_0,\Phi_0,p)$. 
By Lemma \ref{lem:OneSidedInverse}
there is a positive constant $\eps=\eps(A_0,\Phi_0,p)$ such that
if $\|(a,\phi)\|_{L^4}<\eps$,  
the linear operator $\Pi_2^\perp\sD_{A_0,\Phi_0} + \Pi_2^\perp\{(a,\phi),\
\cdot\ \}$ is (left and right) invertible as a map 
$L^p_1\cap (\Ker\sD_{A_0,\Phi_0})^\perp 
\to L^p\cap (\Ker\sD_{A_0,\Phi_0}^*)^\perp$. Since 
$\Pi_1^\perp(a,\phi)$ is a solution to 
\eqref{eq:InvertiblePTEllReg} when $p=2$, then it is the unique solution
for $p=2$. 

Let $(a_p,\phi_p)\in L^p_1\cap (\Ker\sD_{A_0,\Phi_0})^\perp$ be the $L^p_1$
solution to \eqref{eq:InvertiblePTEllReg}. Then $(a_p,\phi_p)$
is also an $L^2_1$ solution and by the uniqueness
assertion we must have that $(a_p,\phi_p) =
\Pi_1^\perp(a,\phi)$. Thus, $\Pi_1^\perp(a,\phi)$ is in $L^p_1$ and so
$(a,\phi)$ is in $L^p_1$, since $\Pi_1(a,\phi)$ is in $C^\8$.
Finally, the standard estimate for $\sD_{A_0,\Phi_0}$, the estimate
\eqref{eq:KernelEst}, and equation \eqref{eq:PTEllReg} yield
\begin{align*}
\|(a,\phi)\|_{L^p_{1,A_0}} 
&\le C\left(\|\sD_{A_0,\Phi_0}(a,\phi)\|_{L^p} 
+ \|\Pi_1(a,\phi)\|_{L^2}\right) \\
&\le C\left(\|\{(a,\phi),(a,\phi)\}\|_{L^p} 
+ \|(\zeta,v,\psi)\|_{L^p} 
+ \|\Pi_1(a,\phi)\|_{L^2}\right),
\end{align*}
and so the desired bound for $(a,\phi)$ follows by the Sobolev
multiplication $L^4\times L^q\to L^p$, the embedding $L^p_1\subset L^q$, and
rearrangement. 
\end{pf}

Proposition \ref{prop:L2_1InhomoReg} will have two main applications: the
first is in our proof of removable singularities for $\PU(2)$ monopoles
and the second is in our proof of $C^\8$ regularity for $L^2_1$
solutions to the $\PU(2)$ monopole equation obtained by gluing \cite{FL3,FL4}. 

As is well-known from standard gauge theory, it is not possible to
construct a well-defined quotient space using $L^2_1$ pairs modulo $L^2_2$
gauge transformations. We construct a quotient using $L^2_k$ pairs modulo
$L^2_{k+1}$ gauge transformations, with $k\ge 2$.
We first establish a regularity result for the inhomogeneous
$\PU(2)$ monopole and Coulomb gauge equations, while in \S
\ref{subsec:Lp_kInhomoGlobal} we show that any $\PU(2)$ monopole in $L^2_k$ is
$L^2_{k+1}$-gauge equivalent to a $\PU(2)$ monopole in $C^\8$. 

\begin{prop}\label{prop:Lp_1InhomoReg}
Continue the notation of Proposition \ref{prop:L2_1InhomoReg}.
Let $k\ge 1$ be an integer and let $2<p<\8$. Let
$(A_0,\Phi_0)$ be a $C^\8$ pair on the $C^\8$ bundles $(\su(E),W^+\otimes E)$
over $X$. Suppose that either
\begin{itemize}
\item $(a,\phi)\in L^p_1(X,\La^1\otimes\su(E))\oplus L^p_1(X,W^+\otimes
E)$, with $(\vectau,\vecvartheta) = 0$, or
\item $(a,\phi)\in L^2_2(X,\La^1\otimes\su(E))\oplus L^2_2(X,W^+\otimes E)$
\end{itemize}
is a solution on $(\su(E),W^+\otimes E)$ to the elliptic system
\eqref{eq:PTEllReg} over $X$, where $(\zeta,v,\psi)$ is in $L^2_k$. Then
$(a,\phi)$ is in $L^2_{k+1}$ and there is a universal polynomial
$Q_k(x,y)$, with positive real coefficients, depending at most on
$(A_0,\Phi_0),k$, such that $Q_k(0,0)=0$ and
$$
\|(a,\phi)\|_{L^2_{k+1,A_0}(X)} 
\le Q_k\left(\|(\zeta,v,\psi)\|_{L^2_{k,A_0}(X)},
\|(a,\phi)\|_{L^p_{1,A_0}(X)}\right).
$$
In particular, if $(\zeta,v,\psi)$ is in $C^\8$ then $(a,\phi)$ is in
$C^\8$ and if $(\zeta,v,\psi)=0$, then
$$
\|(a,\phi)\|_{L^2_{k+1,A_0}(X)} \le C\|(a,\phi)\|_{L^p_{1,A_0}(X)}.
$$
\end{prop}

\begin{pf}
We consider the cases $k=1$, $k=2$, and $k\ge 3$ separately. We may further
assume without loss of generality that $2<p<4$. 

\begin{case}[{\boldmath $k=1$}]
Let $p_0 = p$ and
$q_0=q=4p/(4-p)$, 
so that $1/p=1/4+1/q$ and $4<q<\8$. We have a continuous
Sobolev multiplication map $L^4\times L^q\to L^p$ and an embedding
$L^p_1\subset L^q$. Since
\begin{equation}
\sD_{A_0,\Phi_0}(a,\phi) = - \{(a,\phi),(a,\phi)\}
+ (\zeta,v,\psi), \label{eq:Bootstrap}
\end{equation}
then elliptic regularity for $\sD_{A_0,\Phi_0}$ implies that $(a,\phi)\in
L^{p_1}_1$, where $p_1=2p/(4-p)$ and $q_1 = 4p_1/(4-p_1)$, so $2<p_1<4$ and
$4<q_1<\8$. Here, we used   
the Sobolev multiplication $L^4\times
L^{q_1}\to L^{p_1}$ to get $\{(a,\phi),(a,\phi)\}$ in $L^{p_1}$ and
and the embedding $L^2_1\subset L^{p_1}$ to get $(\zeta,v,\psi)$ in $L^{p_1}$.
Let $p = 2+\eps$ and $\de = \eps/(2-\eps)$, and note that
$$
p_1 = \frac{2p}{4-p} = \frac{2p}{2-\eps} = (1+\de)p > 2+\eps.
$$   
If $p_j < 4$, we inductively define $p_{j+1} = 2p_j/(4-p_j)$ and
$q_{j+1} = 4p_j/(4-p_j)$ for $j\ge 0$. Therefore, we have $p_j > 2 +\eps$
and thus  
$$
p_{j+1} 
= 
\frac{2p_j}{4-p_j} 
> 
\frac{2p_j}{2-\eps} = (1+\de)p_j > p_j > 2 +\eps, 
$$
and so $p_{j+1}>(1+\de)^{j+1}p$ for $j\ge 0$.

By repeating the above regularity argument when $2 < p_{j+1} < 4$,
using 
\eqref{eq:Bootstrap} at each stage, we see that $(a,\phi)\in L^{p_j}_1$ for
$j\ge 0$. We continue the induction until for 
large enough $j\ge 0$, we find that $p' :=
p_{j+1}\ge 8/3$ and $q' := q_{j+1}\ge 8$, so $(a,\phi)\in
L^{8/3}_1\subset L^8$. Therefore, with $(\zeta,v,\psi)\in
L^2_1\subset L^4$ and using the Sobolev multiplication $L^8\times L^8\to L^4$
to get $\{(a,\phi),(a,\phi)\}$ in $L^4$, 
equation \eqref{eq:Bootstrap} gives $(a,\phi)\in L^4_1$. The
Sobolev multiplication $L^4_1\times L^4_1\to L^2_1$ implies that the
quadratic term $\{(a,\phi),(a,\phi)\}$ is in $L^2_1$ and so 
\eqref{eq:Bootstrap} gives $(a,\phi)\in L^2_2$, as required.
\end{case}

\begin{case}[{\boldmath $k=2$}]
{}From the case $k=1$ we continue the induction until $p' 
> 8/3$ and $q' > 4$. So, with $(\zeta,v,\psi)$ now in
$L^2_2\subset L^{p'}_1$, 
equation \eqref{eq:Bootstrap} gives $(a,\phi)\in L^{p'}_2\subset C^0$. The
Sobolev multiplication $L^{p'}_2\times L^2_2\to L^2_2$ implies that the
quadratic term $\{(a,\phi),(a,\phi)\}$ is then in $L^2_2$ and so 
\eqref{eq:Bootstrap} gives $(a,\phi)\in L^2_4$, as required.
\end{case}

\begin{case}[{\boldmath $k\ge 3$}]
There is continuous Sobolev multiplication 
map $L^2_k\times L^2_k\to L^2_k$ and so the
quadratic term $\{(a,\phi),(a,\phi)\}$
is in $L^2_k$, since $(a,\phi)$ is in $L^2_k$. Therefore, 
\eqref{eq:Bootstrap} gives $(a,\phi)\in L^2_{k+1}$.
\end{case}

This completes the proof of the proposition. 
\end{pf}

By combining Propositions \ref{prop:L2_1InhomoReg} and \ref{prop:Lp_1InhomoReg}
we obtain the desired regularity result for
$L^2_1$ solutions to the inhomogeneous Coulomb gauge and $\PU(2)$
monopole equations:

\begin{cor}
\label{cor:L2_1InhomoReg}
Continue the notation of Proposition \ref{prop:L2_1InhomoReg}.
Let $(A_0,\Phi_0)$ be a $C^\8$ pair on the $C^\8$
bundles $(\su(E),W^+\otimes E)$ over $X$.  Then there is a positive constant
$\eps=\eps(A_0,\Phi_0)$ such that the following hold. 
Suppose that either
\begin{itemize}
\item $(a,\phi)\in 
L^2_1(X,\La^1\otimes\su(E))\oplus L^2_1(X,W^+\otimes E)$, with
$(\vectau,\vecvartheta) = 0$, or
\item $(a,\phi)\in L^2_2(X,\La^1\otimes\su(E))\oplus L^2_2(X,W^+\otimes E)$
\end{itemize}
is a solution on $(\su(E),W^+\otimes E)$ to the elliptic system
\eqref{eq:PTEllReg} over $X$, where $(\zeta,v,\psi)$ is in $L^2_k$ and 
$\|(a,\phi)\|_{L^4(X)}<\eps$ and $k\ge 0$ is an integer.
Then $(a,\phi)$ is in $L^2_{k+1}$ and 
there is a universal polynomial $Q_k(x,y)$, with positive real coefficients,
depending at most on $(A_0,\Phi_0),k$, such that $Q_k(0,0)=0$ and
$$
\|(a,\phi)\|_{L^2_{k+1,A_0}(X)} 
\le Q_k\left(\|(\zeta,v,\psi)\|_{L^2_{k,A_0}(X)},
\|(a,\phi)\|_{L^2(X)}\right).
$$
In particular, if $(\zeta,v,\psi)$ is in $C^\8$ then $(a,\phi)$ is in
$C^\8$ and if $(\zeta,v,\psi)=0$, then
$$
\|(a,\phi)\|_{L^2_{k+1,A_0}(X)} \le C\|(a,\phi)\|_{L^2(X)}.
$$
\end{cor}

\begin{rmk}
A similar sharp regularity result for solutions to the Coulomb and
anti-self-dual equations (that is, $(d_\Ga^*+d_\Ga^+)a + (a\wedge a)^+ =
v$) on the product bundle 
over $S^4$ is given by Proposition 4.4.13 in \cite{DK}. The reader is
forwarned that Corollary \ref{cor:L2_1InhomoReg} does not apply directly to
show that $L^2_1$ solutions to the $\PU(2)$ monopole equation in
\cite{FL3} are $C^\8$ due to the unfavorable dependence of the
constant $\eps(A',\Phi',p)$ on the approximate $\PU(2)$ monopole
$(A',\Phi')$ when $p>2$. 
However, a local version of this result, namely Corollary
\ref{cor:L2_1InhomoRegLocal} below, is applicable in this
situation. The point is explained further in \cite{FL3}. 
\end{rmk} 

\subsection{Regularity of $L^2_k$ solutions to the PU(2) monopole equations}
\label{subsec:Lp_kInhomoGlobal}
We explain in this subsection why an $L^2_k$ monopole $(A,\Phi)$ on the
$C^\8$ bundles $(\su(E),W^+\otimes E)$ over $X$ is gauge equivalent, via an
$L^2_{k+1}$ determinant-one, unitary automorphism $u$ of the bundle $E$, to
a $C^\8$ solution $u(A,\Phi)$ on $(\su(E),W^+\otimes E)$ over $X$ when $k\ge
2$. This regularity result contains, as special cases, Theorem 8.8 in
\cite{FU}, for anti-self-dual connections, and Theorem 8.11 in
\cite{Salamon} for Seiberg-Witten monopoles.   

We will need the following observation concerning the symmetry of the
Coulomb gauge equation for pairs; the corresponding fact for connections is
explained in \cite[p. 56]{DK}.

\begin{lem}\label{lem:CoulombSymmetry}
Let $(A,\Phi)$, $(A_0,\Phi_0)$ be $L^2_k$ pairs on the bundles
$(\su(E),W^+\otimes E)$ over $X$. If $(A_0,\Phi_0)$ is in Coulomb gauge
relative to $(A,\Phi)$, so $d_{A,\Phi}^{0,*}((A_0,\Phi_0)-(A,\Phi))=0$ then
$(A,\Phi)$ is in Coulomb gauge relative to $(A_0,\Phi_0)$, so
$d_{A_0,\Phi_0}^{0,*}((A,\Phi)-(A_0,\Phi_0))=0$.
\end{lem}

\begin{pf}
The equation $d_{A,\Phi}^{0,*}((A_0,\Phi_0)-(A,\Phi))=0$ is the Euler-Lagrange
equation for the functional $\sG_E^{2,k+1}\ni u\mapsto  
\|u(A_0,\Phi_0)-(A,\Phi)\|_{L^2}^2$, while the equation 
$d_{A_0,\Phi_0}^{0,*}((A,\Phi)-(A_0,\Phi_0))=0$ is the Euler-Lagrange
equation for the functional $\sG_E^{2,k+1}\ni v\mapsto  
\|v(A,\Phi)-(A_0,\Phi_0)\|_{L^2}^2$.

But for any $u\in\sG_E^{2,k+1}$ we have $\|u(A_0,\Phi_0)-(A,\Phi)\|_{L^2} =
\|u^{-1}(A,\Phi)-(A_0,\Phi_0)\|_{L^2}$ and so if the functional 
$u\mapsto \|u(A_0,\Phi_0)-(A,\Phi)\|_{L^2}^2$ has a critical point at
$u=\id_E$, then functional
$u^{-1}\mapsto\|u^{-1}(A,\Phi)-(A_0,\Phi_0)\|_{L^2}^2$ 
will also have a critical point at $u=\id_E$.
\end{pf}

\begin{prop}\label{prop:GlobalReg}
Let $X$ be a closed, oriented, Riemannian four-manifold with \spinc bundle
$W$ and let $E$ be a Hermitian two-plane bundle over $X$.  Let $k\ge 2$ be
an integer and suppose that $(A,\Phi)$ is an $L^2_k$ solution to
\eqref{eq:PT} on the $C^\8$ bundles $(\su(E),W^+\otimes E)$ over $X$. Then
there is a $L^2_{k+1}$ determinant-one, unitary automorphism $u$ of the
bundle $E$ over $X$ such that $u(A,\Phi)$ is $C^\8$ over $X$.
\end{prop}

\begin{pf}
It suffices, of course, to show that there is an $L^2_{k+1}$ gauge
transformation $u$ of $E$ over such that $u(A,\Phi)$ is in $L^2_{k+1}$.
The $C^\8$ pairs $(A_0,\Phi_0)$ on $(\su(E),W^+\otimes E)$ form a dense
subspace of the space of $L^2_k$ pairs and so, given $\eps=\eps(A,\Phi)>0$,
there is a $C^\8$ pair $(A_0,\Phi_0)$ such that
$$
\|(A,\Phi) - (A_0,\Phi_0)\|_{L^2_{k,A_0}} < \eps.
$$
For small enough $\eps$, Proposition \ref{prop:Slice} gives an $L^2_{k+1}$
determinant-one, unitary automorphism $u$ of the bundle $E$ such that
$u^{-1}(A_0,\Phi_0)$ is in Coulomb gauge relative to $(A,\Phi)$:
$$
d_{A,\Phi}^{0,*}(u^{-1}(A_0,\Phi_0)-(A,\Phi)) = 0.
$$
Now $u(d_{A,\Phi}^{0,*}(u^{-1}(A_0,\Phi_0)-(A,\Phi))) 
= d_{u(A,\Phi)}^{0,*}((A_0,\Phi_0)-u(A,\Phi))$ and so
$$
d_{u(A,\Phi)}^{0,*}((A_0,\Phi_0)-u(A,\Phi)) = 0.
$$
Therefore, $(A_0,\Phi_0)$ is in Coulomb gauge relative to $u(A,\Phi)$ and 
Lemma \ref{lem:CoulombSymmetry} implies that $u(A,\Phi)$ is in Coulomb
gauge relative to $(A_0,\Phi_0)$,
$$
d_{A_0,\Phi_0}^{0,*}(u(A,\Phi)-(A_0,\Phi_0)) = 0.
$$
Let $(a,\phi) = u(A,\Phi)-(A_0,\Phi_0)$, so that $(a,\phi) \in
L^2_k(X,\La^1\otimes\su(E))\oplus L^2_k(X,W^+\otimes E)$; the Coulomb gauge
condition then takes the simpler form $d_{A_0,\Phi_0}^{0,*}(a,\phi) = 0$.
Since $u(A,\Phi) = (A_0,\Phi_0) + (a,\phi)$ is an $L^2_k$ monopole,
then $\fs((A_0,\Phi_0) + (a,\phi)) = 0$ and so $(a,\phi)$ is an $L^2_k$
solution to the quasi-linear elliptic equation 
$$
\sD_{A_0,\Phi_0}(a,\phi) + \{(a,\phi),(a,\phi)\} = - \fs(A_0,\Phi_0),
$$
with $C^\8$ data $-\fs(A_0,\Phi_0)$. The conclusion now follows from
Proposition~\ref{prop:Lp_1InhomoReg} with $\zeta = 0$ and 
$(v,\psi)=- \fs(A_0,\Phi_0)$.
\end{pf}

One of the convenient practical consequences of Proposition
\ref{prop:GlobalReg} is that we can always work, modulo global $L^2_{k+1}$
gauge transformations, with $C^\8$ rather than $L^2_k$
monopoles.  For the remainder of this article, therefore, 
the term `$\PU(2)$ monopole' is generally 
reserved for $C^\8$ solutions to the
perturbed $\PU(2)$ monopole equations \eqref{eq:PT}. In a similar vein, we
generally reserve the terms `gauge transformation' or `bundle map'
for gauge transformations or bundle maps which are in $C^\8$.

Given Proposition \ref{prop:GlobalReg}, we then
have the following analogue of Proposition 4.2.16 in \cite{DK}, the
corresponding result for the moduli space of anti-self-dual
connections. The proof is standard and so is left to the reader.

\begin{cor}
Continue the hypotheses of Proposition \ref{prop:GlobalReg}.
Then for any $k\ge 2$ the natural inclusion of topological spaces
$M_{W,E}^k \hookrightarrow M_{W,E}^{k+1}$ is a homeomorphism.
\end{cor}

Thus, the topology of the moduli space $M_{W,E}^k$ 
of $L^2_k$ monopoles is independent of the Sobolev spaces used
in its construction for $k\ge 2$ and so we simply denote the moduli space
by $M_{W,E}$.   

\subsection{Local regularity and interior estimates for $L^2_1$ solutions
to the inhomogeneous Coulomb gauge and PU(2) monopole equations}
\label{subsec:L2_1InhomoLocal} 
In this section we specialize the results of \S 
\ref{subsec:L2_1InhomoGlobal} to the case where the reference pair is a
trivial $\PU(2)$ monopole, so $(A_0,\Phi_0)=(\Ga,0)$ on the bundles
$(\su(E),W^+\otimes E)$, over an open subset $\Om\subset X$, where $\Ga$ is
a flat connection.

We continue to assume that $X$ is a closed, oriented four-manifold with
metric $g$, \spinc bundle $W$, and Hermitian two-plane bundle $E$ extending
those on $\Om\subset X$.  We use the inhomogenous estimates and regularity
results in our proof of removable singularities for $\PU(2)$ monopoles in
\S \ref{subsec:RemSing} and in our development of the gluing theory for
$\PU(2)$ monopoles in sequels \cite{FL3,FL4} to the present article ---
especially to show that a global $L^2_1$ gluing solution $(a,\phi)$ is
actually $C^\8$. We use the homogeneous estimates and regularity results in
\S \ref{sec:compact}, in our proof of the existence of an Uhlenbeck
compactification for the moduli space of $\PU(2)$ monopoles.

We have the following local versions of Propositions
\ref{prop:L2_1InhomoReg} and \ref{prop:Lp_1InhomoReg} and
Corollary \ref{cor:L2_1InhomoReg}:

\begin{prop}
\label{prop:L2_1InhomoRegLocal}
Continue the notation of the preceding paragraph. Let $\Om'\Subset\Om$ be a
precompact open subset and let $2\le p < 4$.  Then there are positive
constants $\eps=\eps(\Om,p)$ and
$C=C(\Om',\Om,p)$ with the following
significance. Suppose that $(a,\phi)$ is an $L^2_1(\Om)$ solution to the
elliptic system \eqref{eq:PTEllReg} over $\Om$ for $(\vectau,\vecvartheta)
= 0$, with $(A_0,\Phi_0)=(\Ga,0)$ and where
$(\zeta,v,\psi)$ is in $L^p(\Om)$. If $\|(a,\phi)\|_{L^4(\Om)}<\eps$ then
$(a,\phi)$ is in $L^p_1(\Om')$ and
$$
\|(a,\phi)\|_{L^p_{1,\Ga}(\Om')} 
\le C\left(\|(\zeta,v,\psi)\|_{L^p(\Om)} 
+ \|(a,\phi)\|_{L^2(\Om)}\right).
$$
\end{prop}

\begin{pf}
Choose an open subset $\Om''$ such that $\Om'\Subset\Om''\Subset\Om$.
Let $\chi$ be $C^\8$ cutoff function such that $\supp\chi\subset\Om''$ and
$\chi=1$ on $\Om'$. 
Let $\be$ be a cutoff function such that $\be=1$ on
$\supp\chi$ and $\supp\be\subset \Om''$.
Since $(a,\phi)$ is a solution to 
\eqref{eq:PTEllReg} with right-hand side $(\zeta,v,\psi)$ over $\Om$, then
$\chi(a,\phi)$ solves
\begin{align*}
&\sD_{\Ga,0}\chi(a,\phi) + \{(a,\phi),\chi(a,\phi)\} \\
&\qquad = \chi\sD_{\Ga,0}(a,\phi) + d\chi\otimes (a,\phi) 
+ \chi\{(a,\phi),(a,\phi)\} \\
&\qquad = \chi(\zeta,v,\psi) + d\chi\otimes (a,\phi) 
=: (\zeta',v',\psi').  
\end{align*}
Thus, $\chi(a,\phi)$ is an $L^2_1$ solution 
over $X$ to the linear elliptic system over $X$,
$$
\sD_{\Ga,0}(b,\varphi) + \{\be(a,\phi),(b,\varphi)\}
= (\zeta',v',\psi'),
$$
with $L^2_1$ coefficient $\be(a,\phi)$ and $L^2_1$ right-hand side
$(\zeta',v',\psi')$. Since $\|\be(a,\phi)\|_{L^4(X)} \le 
\|(a,\phi)\|_{L^4(\Om)}$, the proof of Proposition \ref{prop:L2_1InhomoReg}
implies that $\chi(a,\phi)$ is in $L^p_1$ over $X$ if
$\eps=\eps(g,\Om,p)$ is sufficiently small. Thus, $(a,\phi)$ is in $L^p_1$ over
$\Om'$ with 
$$
\|\chi(a,\phi)\|_{L^p_{1,\Ga}(X)} 
\le C\left(\|(\zeta',v',\psi')\|_{L^p(X)}
+ \|\chi(a,\phi)\|_{L^2(X)}\right), 
$$
and therefore, for $2\le p < 4$, we have
$$
\|(a,\phi)\|_{L^p_{1,\Ga}(\Om')} 
\le C\left(\|(\zeta,v,\psi)\|_{L^p(\Om'')}
+ \|(a,\phi)\|_{L^p(\Om'')}\right).
$$
The preceding bound and the Sobolev embedding $L^2_1\subset L^p$ give 
the estimate
$$
\|(a,\phi)\|_{L^p(\Om'')} 
\le c\|(a,\phi)\|_{L^2_{1,\Ga}(\Om'')} 
\le C\left(\|(\zeta,v,\psi)\|_{L^2(\Om)}
+ \|(a,\phi)\|_{L^2(\Om)}\right).
$$
Combining these inequalities then gives the required
$L^p_1$ estimate for $(a,\phi)$ over $\Om'$.
\end{pf}

\begin{prop}
\label{prop:Lp_1InhomoRegLocal}
Continue the notation of Proposition \ref{prop:L2_1InhomoRegLocal}.
Let $k\ge 1$ be an integer, and let $2<p<\8$. Suppose that either
\begin{itemize}
\item $(a,\phi)$ is an $L^p_1(\Om)$, when $(\vectau,\vecvartheta) = 0$, or
\item $(a,\phi)$ is an $L^2_2(\Om)$
\end{itemize}
solution to the elliptic system
\eqref{eq:PTEllReg} over $\Om$ with $(A_0,\Phi_0)=(\Ga,0)$, 
where $(\zeta,v,\psi)$ is in $L^2_k(\Om)$. Then $(a,\phi)$ is in
$L^2_{k+1}(\Om')$ and there is a universal polynomial $Q_k(x,y)$, with
positive real coefficients, depending at most on
$k$, $\Om'$, $\Om$, such that
$Q_k(0,0)=0$ and
$$
\|(a,\phi)\|_{L^2_{k+1,\Ga}(\Om')} 
\le Q_k\left(\|(\zeta,v,\psi)\|_{L^2_{k,\Ga}(\Om)},
\|(a,\phi)\|_{L^p_{1,\Ga}(\Om)}\right).
$$
If $(\zeta,v,\psi)$ is in $C^\8(\Om)$ then $(a,\phi)$ is in $C^\8(\Om')$
and if $(\zeta,v,\psi)=0$, then
$$
\|(a,\phi)\|_{L^2_{k+1,\Ga}(\Om')} 
\le C\|(a,\phi)\|_{L^p_{1,\Ga}(\Om)}.
$$
\end{prop}

\begin{pf}
Again, choose an open subset $\Om''$ such that $\Om'\Subset\Om''\Subset\Om$,
let $\chi$ be $C^\8$ cutoff function such that $\supp\chi\subset\Om$ and
$\chi=1$ on $\Om''$ and let $\be$ be a cutoff function such that $\be=1$ on
$\supp\chi$ and $\supp\be\subset\Om$.
Since $(a,\phi)$ is a solution to 
\eqref{eq:PTEllReg} with right-hand side $(\zeta,v,\psi)$ over $\Om$, then
$\chi(a,\phi)$ solves
\begin{align*}
&\sD_{\Ga,0}\chi(a,\phi) + \{\chi(a,\phi),\chi(a,\phi)\} \\
&\qquad = \chi\sD_{\Ga,0}(a,\phi) + d\chi\otimes (a,\phi) 
+ \chi\{(a,\phi),(a,\phi)\}  + \chi(\chi-1)\{(a,\phi),(a,\phi)\} \\
&\qquad = \chi(\zeta,v,\psi) + d\chi\otimes (a,\phi)
+ \chi(\chi-1)\{(a,\phi),(a,\phi)\} \\
&\qquad =: (\zeta',v',\psi').  
\end{align*}
Note that $\chi(\zeta,v,\psi)$ is in $L^2_k(X)$, while $d\chi\otimes
(a,\phi)$ is in $L^p_1(X)$ and the Sobolev multiplication $L^q\times L^q\to
L^{p_1}$ implies that $\chi(\chi-1)\{(a,\phi),(a,\phi)\}$ is in
$L^{p_1}(X)$, where $p_1=q/2=2p/(4-p) > p$. Thus, $(\zeta',v',\psi')$ is in
$L^{p_1}(X)$.  

The proof of the case $k=1$ in Proposition
\ref{prop:Lp_1InhomoReg} then implies that $\chi(a,\phi)$ is an
$L^{p_1}_1(X)$ so $(a,\phi)$ is in $L^{p_1}_1(\Om'')$. We now repeat this
process for each of the remaining steps in the proof of Proposition
\ref{prop:Lp_1InhomoReg}, at each stage
on successively smaller open subsets $\Om'''$
such that $\Om'\Subset\Om'''\Subset\Om''$, 
until we obtain $(a,\phi)$ in
$L^2_{k+1}(\Om')$ and the desired $L^2_{k+1}(\Om')$ estimate.
\end{pf}

\begin{cor}\label{cor:L2_1InhomoRegLocal}
Continue the notation of Proposition \ref{prop:Lp_1InhomoRegLocal}.  Then
there is a positive constant $\eps=\eps(\Om)$ with the following
significance. Suppose that either
\begin{itemize}
\item $(a,\phi)$ is an $L^2_1(\Om)$, when $(\vectau,\vecvartheta) = 0$, or
\item $(a,\phi)$ is an $L^2_2(\Om)$
\end{itemize}
solution to the elliptic system \eqref{eq:PTEllReg} over $\Om$ with
$(A_0,\Phi_0)=(\Ga,0)$, where $(\zeta,v,\psi)$ is in $L^2_k(\Om)$ and
$\|(a,\phi)\|_{L^4(\Om)}<\eps$.  Then $(a,\phi)$ is in $L^2_{k+1}(\Om')$
and there is a universal polynomial $Q_k(x,y)$, with positive real
coefficients, depending at most on $k$,
$\Om'$, $\Om$, such that $Q_k(0,0)=0$ and
$$
\|(a,\phi)\|_{L^2_{k+1,\Ga}(\Om')} 
\le Q_k\left(\|(\zeta,v,\psi)\|_{L^2_{k,\Ga}(\Om)},
\|(a,\phi)\|_{L^2(\Om)}\right).
$$
If $(\zeta,v,\psi)$ is in 
$C^\8(\Om)$ then $(a,\phi)$ is in $C^\8(\Om')$ and if
$(\zeta,v,\psi)=0$, then
$$
\|(a,\phi)\|_{L^2_{k+1,\Ga}(\Om')} \le C\|(a,\phi)\|_{L^2(\Om)}.
$$
\end{cor}

Corollary \ref{cor:L2_1InhomoRegLocal} thus yields a sharp local
elliptic regularity result for $\PU(2)$ monopoles $(A,\Phi)$ in $L^2_1$ 
which are given to us in Coulomb gauge relative to $(\Ga,0)$: this
regularity result is the key ingredient in our proof (given in \S
\ref{subsec:RemSing}) of removable point singularities for $\PU(2)$
monopoles.

\begin{prop}
\label{prop:L2_1CoulMonoRegLocal}
Continue the notation of Corollary \ref{cor:L2_1InhomoRegLocal}.  Then
there is a positive constant
$\eps=\eps(\Om)$ and, if $k\ge 1$ is an
integer, there is a positive constant
$C=C(\Om',\Om,k)$ with the following
significance. Suppose that either
\begin{itemize}
\item $(A,\Phi)$ is an $L^2_1$, when $(\vectau,\vecvartheta) = 0$, or
\item $(A,\Phi)$ is an $L^2_2$
\end{itemize}
solution to the $\PU(2)$ monopole equations \eqref{eq:PT} over $\Om$,
which is in Coulomb gauge over $\Om$ relative to $(\Ga,0)$, so
$d_\Ga^*(A-\Ga) = 0$, 
and obeys $\|(A-\Ga,\Phi)\|_{L^4(\Om)}<\eps$.
Then $(A-\Ga,\Phi)$ is in $C^\8(\Om')$ and for any $k\ge 1$,
$$
\|(A-\Ga,\Phi)\|_{L^2_{k,\Ga}(\Om')} 
\le C\|(A-\Ga,\Phi)\|_{L^2(\Om)}.
$$
\end{prop}

\begin{pf}
Corollary \ref{cor:L2_1InhomoRegLocal} applies to the $L^2_1(\Om)$ pair
$(a,\phi) = (A-\Ga,\Phi)$ and yields the required regularity and
estimates for $(A-\Ga,\Phi)$ with $(\zeta,v,\psi) = 0$ in
\eqref{eq:PTEllReg}.
\end{pf}

\subsection{Estimates for PU(2) monopoles in a good local
gauge}\label{subsec:Uhlenbeck} 
It remains to combine the local regularity results and estimates of \S
\ref{subsec:L2_1InhomoLocal}, for $\PU(2)$ monopoles $(A,\Phi)$ where
the connection $A$ is assumed to be in Coulomb gauge relative to the
product $\SO(3)$ connection $\Ga$, with Uhlenbeck's local, Coulomb
gauge-fixing theorem. We then obtain regularity results and estimates for
$\PU(2)$ monopoles $(A,\Phi)$ with small curvature $F_A$, parallel to those
of Theorem 2.3.8 and Proposition 4.4.10 in \cite{DK} for anti-self-dual
connections.

In order to apply Corollary \ref{cor:L2_1InhomoRegLocal} we need
Uhlenbeck's Coulomb gauge-fixing result
\cite[Theorem 2.1 \& Corollary 2.2]{UhlLp}). Let $B$ (respectively,
$\barB$) be the open (respectively, closed) unit ball centered at
the origin in $\RR^4$ and let $G$ be a compact
Lie group. In order to provide universal constants we assume $\RR^4$ has
its standard metric, though the results of this subsection naturally hold
for any $C^\8$ Riemannian metric, with comparable constants for metrics
which are suitably close. 

\begin{thm}\label{thm:CoulombBallGauge}
There are positive constants $c$ and $\eps$ with the following
significance. If $2\le p < 4$ is a constant and $A\in
L^p_1(B,\La^1\otimes\fg)\cap L^p_1(\rd B,\La^1\otimes\fg)$ is a
connection matrix whose curvature satisfies $\|F_A\|_{L^p(B)} < \eps$,
then there is an gauge transformation 
$u\in L^p_2(B, G)\cap  L^p_2(\rd B, G)$ such that $u(A)
:= uAu^{-1} - (du)u^{-1}$ satisfies
\begin{align}
d^*u(A) &= 0\text{ on }B, \tag{1}\\
{\textstyle{\frac{\rd}{\rd r}}}\lrcorner u(A) &= 0\text{ on }\rd B, \tag{2}\\
\|u(A)\|_{L^p_1(B)} &\le c\|F_A\|_{L^p(B)}.\tag{3}
\end{align}
If $A$ is in $L^p_k(B)$, for $k\ge 2$, then $u$ is in $L^p_{k+1}(B)$. The gauge
transformation $u$ is unique up to multiplication by a constant element of $G$.
\end{thm}

\begin{rmk}
If $G$ is abelian then the requirement that $\|F_A\|_{L^p(B)} < \eps$ can be
omitted. 
\end{rmk} 

It is often useful to rephrase Theorem \ref{thm:CoulombBallGauge} in two
other slightly different ways. Suppose $A$ is an $L^2_k$ connection on a
$C^\8$ principal $G$ bundle $P$ over $B$ with $k\ge 2$ and $\|F_A\|_{L^2(B)} <
\eps$. Then the assertions of Theorem \ref{thm:CoulombBallGauge} are 
equivalent to each of the following:
\begin{itemize}
\item There is an $L^2_{k+1}$ trivialization $\tau:P\to B\times G$ such
that (i) $d_\Gamma^*(\tau(A)-\Gamma)=0$, where $\Gamma$ is the product
connection on $B\times G$, (ii) $\frac{\rd}{\rd r}
\lrcorner (\tau(A)-\Gamma) = 0$, and (iii)
$\|(\tau(A)-\Gamma)\|_{L^2_1(B)} \le c\|F_A\|_{L^2(B)}$.

\item There is an $L^2_{k+1}$ flat connection $\Gamma$ on $P$ such
that (i) $d_\Gamma^*(A-\Gamma)=0$, (ii) $\frac{\rd}{\rd r}
\lrcorner (A-\Gamma) = 0$, and (iii)
$\|(A-\Gamma)\|_{L^2_1(B)} \le c\|F_A\|_{L^2(B)}$, and
an $L^2_{k+1}$ trivialization $P|_B\simeq B\times G$
taking $\Gamma$ to the product connection.
\end{itemize}

We can now combine Theorem \ref{thm:CoulombBallGauge} with Proposition
\ref{prop:L2_1CoulMonoRegLocal} to give the following analogue of Theorem
2.3.8 in \cite{DK} --- the interior estimate for anti-self-dual connections
with $L^2$-small curvature.  

\begin{cor}
\label{cor:PTLocalReg}
Let $B\subset\RR^4$ be the open unit ball with center at the origin with
spinc structure $(\rho,W)$, let $U\Subset B$ be an open subset, and let
$\Ga$ be the product connection on $B\times\SO(3)$.  Then there is a
positive constant $\eps$ and if $\ell\ge 1$ is an integer, there is a positive
constant $C(\ell,U)$ with the following significance.  Suppose that either
\begin{itemize}
\item $(A,\Phi)$ is an $L^2_1$, when $(\vectau,\vecvartheta) = 0$, or
\item $(A,\Phi)$ is an $L^2_k$, $k\ge \max\{2,\ell\}$,
\end{itemize}
solution to the $\PU(2)$ monopole equations \eqref{eq:PT} over $B$ and that
the curvature of the $\SO(3)$
connection matrix $A$ obeys $\|F_A\|_{L^2(B)} < \eps$.
Then there is an $L^2_{k+1}$ 
gauge transformation $u:B\to\SU(2)$ such
that $(u(A)-\Ga,u\Phi)$ is in $C^\8(B)$ with 
$d^*(u(A)-\Ga)=0$ over $B$ and
$$
\|(u(A)-\Ga,u\Phi)\|_{L^2_{\ell,\Gamma}(U)} 
\le C\|F_A\|_{L^2(B)}^{1/2}.
$$
\end{cor}

\begin{pf}
Let $\eps_1$ be the constant in Theorem \ref{thm:CoulombBallGauge} and note
that for $\eps\le\eps_1$, Theorem
\ref{thm:CoulombBallGauge} (taking $G=\SO(3)$) and the Sobolev embedding
$L^2_1(B)\subset L^4(B)$ imply that there is an $L^2_{k+1}$ gauge
transformation $u:B\to\SU(2)$ (by lifting the $\SO(3)$ gauge
transformation) such that
\begin{align*}
d_\Ga^*(u(A)-\Ga) &= 0, 
\\
\|u(A)-\Ga\|_{L^4(B)} &\le c_1\|F_A\|_{L^2(B)}< c_1\eps.
\end{align*}
On the other hand, the quadratic equation for $\Phi$ in
\eqref{eq:PT} and Lemma \ref{lem:SigmaPhiPhi} give the $L^4$ and
$L^2$ estimates,
\begin{align*}
\|u\Phi\|_{L^4(B)} &= \|\Phi\|_{L^4(B)}\le 2\|F_A^+\|_{L^2(B)}^{1/2}, 
\\
\|u\Phi\|_{L^2(B)} &\le 2\|F_A\|_{L^1(B)}^{1/2}.
\end{align*}
Let $\eps_2$ be the constant in Proposition
\ref{prop:L2_1CoulMonoRegLocal}.  Hence, if $c_1\eps\le\eps_2$ and
$4\sqrt{\eps}\le\eps_2$, then Proposition \ref{prop:L2_1CoulMonoRegLocal}
implies that $(u(A)-\Ga,u\Phi)$ obeys
$$
\|(u(A)-\Ga,u\Phi)\|_{L^2_{\ell,\Gamma}(U)} 
\le C\|(u(A)-\Ga,u\Phi)\|_{L^2(B)},
$$
since $d_\Ga^*(u(A)-\Ga) = 0$ and $\|(u(A)-\Ga,u\Phi)\|_{L^4(B)} <
\eps_2$. The desired estimate follows by combining these inequalities for
small enough $\eps \le 1$.
\end{pf}

Again, it is often useful to rephrase Corollary \ref{cor:PTLocalReg} in the two
other slightly different ways. Suppose $\ell\ge 1$ and that $(A,\Phi)$ is a
$\PU(2)$ monopole in $L^2_k$ on $(\su(E),W^+\otimes E)$ over the unit ball
$B\subset\RR^4$ with $k\ge \max\{2,\ell\}$, 
$\|F_A\|_{L^2(B)}<\eps$ and $U\Subset B$. Then the
assertions of Corollary \ref{cor:PTLocalReg} are  
equivalent to each of the following: 
\begin{itemize}
\item There is a $C^\8$ trivialization $\tau:E|_B\to
B\times \CC^2$ and a $L^2_{k+1}$
determinant-one, unitary bundle automorphism $u$ of
$E|_B$ such that, with respect to the product connection $\Gamma$ on
$B\times \su(2)$, we have
(i) $d_\Gamma^*(\tau u(A)-\Gamma)=0$, and (ii)
$\|(\tau u(A)-\Gamma,\tau u\Phi)\|_{L^2_{\ell,\Gamma}(U)} 
\le C\|F_A\|_{L^2(B)}^{1/2}$.

\item There is an $L^2_{k+1}$ flat connection $\Gamma$ on $\su(E)|_B$ such
that (i) $d_\Gamma^*(A-\Gamma)=0$, and (ii)
$\|(A-\Gamma,\Phi)\|_{L^2_{\ell,\Gamma}(U)} \le c\|F_A\|_{L^2(B)}^{1/2}$,
and an $L^2_{k+1}$ trivialization $\su(E)|_B\simeq B\times \su(2)$
taking $\Gamma$ to the product connection.
\end{itemize}

Corollary \ref{cor:PTLocalReg}
immediately yields the following local compactness
result for $\PU(2)$ monopoles analogous to the local compactness
result for anti-self-dual connections in \cite[Corollary 2.3.9]{DK}. 

\begin{cor}\label{cor:PTLocalCpt}
Let $B\subset\RR^4$ be the open
unit ball and \spinc structure $(\rho,W)$.
Then there is a positive constant
$\eps_0(g,A_{\det W^+},A_{\det E})$ with the following significance. 
Let $U\Subset B$ be an open subset
and let $k\ge 2$ be an integer.
If $(A_\al,\Phi_\al)$ is a sequence of 
$\PU(2)$ monopoles in $L^2_k$ on $(\su(E),W^+\otimes E)$ over $B$ such that 
$$
\|F_{A_\al}\|_{L^2(B)} < \eps_0,
$$
then the following holds. There are a subsequence $\{\al'\}\subset\{\al\}$,
a sequence of determinant-one, unitary
$L^2_{k+1}$ automorphisms $\{u_{\al'}\}$ of $E|_B$
and a sequence of gauge equivalent pairs 
$(\tA_{\al'},\tPhi_{\al'}) := 
u_{\al'}(A_{\al'},\Phi_{\al'})$ which converge in $L^2_{k,\loc}$ on $U$ to a
$\PU(2)$ monopole $(\tA,\tPhi)$ over $U$.
\end{cor}

We will also need interior estimates for $\PU(2)$ monopoles in a
good local gauge over more general simply-connected
regions than the open balls
considered in Corollary \ref{cor:PTLocalReg}. Specifically, recall
that a domain $\Om\subset X$ is {\em strongly simply-connected\/} if it has
an open covering by balls $D_1,\dots,D_m$ (not necessarily geodesic) such
that for $1\le r\le m$ the intersection $D_r\cap (D_1\cup\cdots\cup
D_{r-1})$ is connected. We recall (see \cite[Proposition 2.2.3]{DK} or
\cite[Proposition I.2.6]{Kobayashi}):

\begin{prop}
\label{prop:FlatProduct}
If $\Gamma$ is a $C^\8$ flat connection on a principal $G$ bundle
$P$ over a simply-connected manifold $\Om$, then there is a $C^\8$
isomorphism $P\simeq \Om\times G$ taking $\Gamma$ to
the product connection on $\Om\times G$. 
\end{prop}

More generally, if $A$ is $C^\8$ connection on a $G$ bundle $P$ over 
a simply-connected manifold-with-boundary $\barOm=\Om\cup\rd\Om$
with $L^p$-small curvature (with $p>2$), then Uhlenbeck's theorem implies
that $A$ is $L^p_2$-gauge equivalent to a connection which is $L^p_1$-close
to an $L^p_1$ flat connection on $P$ (see \cite[Corollary 4.3]{UhlChern} or
\cite[p. 163]{DK}). The following {\em a priori\/}
interior estimate is a straightforward generalization of
\cite[Proposition 4.4.10]{DK}: the method of proof is identical to that
described in \cite[pp. 161--162]{DK} and uses a patching argument for gauge
transformations employed by Uhlenbeck in the proof of Theorem 3.6 in
\cite{UhlLp}. The required bound for the connection in terms of its
curvature is obtained by covering the given open region with balls and
applying the estimate of Corollary \ref{cor:PTLocalReg} in place of Theorem
2.3.8 in \cite{DK}.

We recall from Proposition \ref{prop:GaugeGroupSequence} that any automorphism
of the $\SO(3)$ bundle $\su(E)|_\Om$ over a simply-connected open subset
$\Om\subset X$ lifts to a determinant-one, unitary bundle automorphism of
$E|_\Om$. The method of \cite[pp. 161--162]{DK} then yields:

\begin{prop}
\label{prop:MonopoleGoodGaugeIntEst}
Let $X$ be a closed, oriented, Riemannian four-manifold with \spinc
structure $(\rho,W)$ and let $\Om\subset X$ be a strongly simply-connected
open subset. Then there is a positive constant $\eps(\Om)$ with the
following significance. For $\Om'\Subset \Om$ a precompact open subset and
an integer $\ell\ge 1$, there is a constant $C(\ell,\Om',\Om)$ such that
the following holds. Suppose $(A,\Phi)$ is a $\PU(2)$ monopole in $L^2_k$ on
$(\su(E),W^+\otimes E)$ over $\Om$ with $k\ge \max\{2,\ell\}$ such that
$$
\|F_A\|_{L^2(\Om)} < \eps.
$$
Then there is an $L^2_{k+1}$ flat connection $\Gamma$ on $\su(E)|_{\Om'}$ such
that 
$$
\|(A-\Ga,\Phi)\|_{L^2_{\ell,\Ga}(\Om')} \le C\|F_A\|_{L^2(\Om)}^{1/2},
$$
and an $L^2_{k+1}$ trivialization $\su(E)|_{\Om'}\simeq \Om'\times \su(2)$
taking $\Gamma$ to the product connection.
\end{prop}


\section{Uhlenbeck compactification for the moduli space of $\PU(2)$
monopoles} 
\label{sec:compact}
Our goal in this section is to prove the existence of an Uhlenbeck-type
compactification of the
moduli space of $\PU(2)$ monopoles analogous to the Uhlenbeck
compactification of the moduli space of anti-self-dual connections
\cite{DK}. A compactification of this type, for the unperturbed moduli space
of $\PU(2)$ monopoles, was first outlined by Pidstrigach and Tyurin
\cite{PTCambridge,PTLocal}. By construction our holonomy
perturbations are continuous with respect Uhlenbeck limits, as outlined in 
\S \ref{subsubsec:DefnPertUhlCompact}.

In \S \ref{subsec:BW} we establish the Bochner-Weitzenb\"ock formulas
for the coupled Dirac operators $D_A$ and
$D_A^*$ and the {\em a priori\/} bounds satisfied by the
section $\Phi$ and the curvature $F_A$, if $(A,\Phi)$ is a $\PU(2)$
monopole on $(\su(E),W^+\otimes E)$; these generalize the well-known {\em a
priori\/} bounds for Seiberg-Witten monopoles
\cite{KrM,MorganSWNotes,Salamon}. In our proof of removable singularities
we need to take into account the variation of \spinc structures with the
Riemannian metric and the rescaling behavior of the $\PU(2)$ equations:
this variation is discussed in \S \ref{subsec:ScaleInv}.

In \S \ref{subsec:RemSing} we prove the removability of point
singularities for $\PU(2)$ monopoles. T. Parker established the
removability of point singularities for solutions to the second-order
coupled Yang-Mills equations, namely, $d_A^*F_A = q(\Phi)$ and $D_A\Phi =
0$, for a unitary connection $A$ on $E$ and a section $\Phi$ of $W\otimes
E$, where $q$ is a certain quadratic form
\cite{Parker}. This generalizes the corresponding results of Uhlenbeck in
the case $\Phi=0$, where the above system then reduces to the second-order
Yang-Mills equation \cite{FU,UhlRem}; proofs of removable singularities for
anti-self-dual connections are given in \cite{DK} (see also \cite{UhlChern}).  
D. Salamon has given a proof of removable singularities for Seiberg-Witten
monopoles \cite[Chapter 9]{Salamon}. The arguments of Uhlenbeck and Parker rely
on pointwise curvature and energy decay estimates; Salamon's argument
relies on energy decay estimates and elliptic regularity results for
Seiberg-Witten monopoles. The proof we give for $\PU(2)$ monopoles is
rather different and instead relies heavily on our $C^\8$ regularity result for
Coulomb-gauge $\PU(2)$ monopoles in $L^2_1$ (Proposition
\ref{prop:L2_1CoulMonoRegLocal}): this is similar to the strategy
used by Donaldson and Kronheimer in \cite{DK}.

The technical ingredients we need for patching sequences of local gauge
transformations are described in \S \ref{subsec:PatchGauge} and the
Uhlenbeck closure $\barM_{W,E}$ is defined in \S
\ref{subsec:DefnUhlCompact}, by analogy with the corresponding definition
for the moduli space of anti-self-dual connections in \cite[\S 4.4]{DK}.
In \S \ref{subsubsec:DefnPertUhlCompact} we describe how the holonomy
perturbations extend continously with respect to Uhlenbeck limits and
induce holonomy perturbations on all lower-level moduli spaces of $\PU(2)$
monopoles. We then define the Uhlenbeck closure for the moduli space of
perturbed $\PU(2)$ monopoles.

In \S \ref{subsec:SeqCompact} we prove Theorem \ref{thm:Compactness}, which
asserts that the Uhlenbeck closure $\barM_{W,E}$ is 
compact. The main analytical ingredients in the proof comprise the
regularity results and estimates of \S \ref{sec:regularity}; the {\em
a priori\/} bounds of \S \ref{subsec:BW} provide a `universal energy
bound' for a $\PU(2)$ monopole $(A,\Phi)$ and this bound plays the same
role here in establishing the existence of an Uhlenbeck compactification 
as the usual topological bound for the energy of an
anti-self-dual connection in \cite{DK}.
The scale invariance of the $\PU(2)$ monopole equation
described in \S \ref{subsec:ScaleInv} is used here in much the same
way that the conformal invariance of the anti-self-dual equation is
exploited in \cite{DK}.  

The parameters $\tau_0,\vartheta_0,\vectau,\vecvartheta$ are chosen so that
the perturbation estimates \eqref{eq:CompactEstPertC0} are satisfied.
These bounds follow from Proposition \ref{prop:SmoothAcrossReducibles},
with $k\ge 2$ for the first inequality in \eqref{eq:CompactEstPertC0}
and $k\ge 3$ for the second, if
\begin{equation}
\|\vecvartheta\|_{\ell^1_\delta(C^r(X))} < \eps_\vartheta
\quad\text{and}\quad
\|\vectau\|_{\ell^1_\delta(C^r(X))} < \eps_\tau,
\label{eq:CompactConstraintOnVecVarthetaTau}
\end{equation}
with suitable positive constants $\eps_\tau$ and $\eps_\vartheta$.  These
constraints are needed in the proofs of the {\em a priori\/} estimates in
Lemmas \ref{lem:L21AaprioriEstAPhi}, \ref{lem:L2aprioriEstFA}, and
\ref{lem:C0EstFAPhi} of \S \ref{subsec:BW} and hence in the proof of
Theorem \ref{thm:SeqCompact} in \S \ref{subsec:SeqCompact}. {\em For the
remainder of the article we therefore require that $k\ge 3$.\/}

\subsection{Bochner formulas and
a priori estimates for $\PU(2)$ monopoles}
\label{subsec:BW} 
In this section we apply the Bochner-Weitzenb\"ock formula
for $D_A^*D_A$ together with Kato's inequality and the maximum principle to
derive {\em a priori\/} estimates for solutions $(A,\Phi)$ to the $\PU(2)$
monopole equations \eqref{eq:PT}. These estimates then lead to a vanishing
result generalizing that of \cite[pp. 781--782]{Witten}. 

We first have a generalization of the usual Bochner-Weitzenb\"ock
identity for Seiberg-Witten monopoles
\cite{KrM,MorganSWNotes,Salamon,TauSymp,Witten}.

\begin{lem}\label{lem:BWDirac}
Let $X$ be an oriented, Riemannian, four-manifold with \spinc structure
$(\rho,W)$, and let $E$ be a
Hermitian bundle over $X$. If $A$ denotes a $\U(2)$ connection on $E$,
$\cov_A:\Om^0(W^\pm\otimes E)\to\Om^1(W^\pm\otimes E)$ are the
covariant derivatives defined by $A$, and $D_A:\Om^0(W^+\otimes
E)\to\Om^0(W^-\otimes E)$ the Dirac operator, then
\begin{align}
D_A^*D_A &= \cov_A^*\cov_A + \quarter R
+ \rho_+(F_A^+) + \half\rho_+(F^+_{A_L}), \tag{1a}\\
D_AD_A^* &= \cov_A^*\cov_A + \quarter R
+ \rho_-(F_A^-) + \half\rho_-(F^+_{A_L}), \tag{2a}
\end{align}
where $R$ is the scalar curvature of the Riemannian metric and 
$A_L = A_{\det
W^+}$ is the induced connection on $L = \det W^+\simeq\det W^-$. If $A$ denotes
an $\SO(3)$ connection on $\su(E)$, then
\begin{align}
D_A^*D_A &= \cov_A^*\cov_A + \quarter R
+ \rho_+(F_A^+) + \half\rho_+(F^+_{A_L} + F_{A_e}^+), \tag{1b}\\
D_AD_A^* &= \cov_A^*\cov_A + \quarter R
+ \rho_-(F_A^-) + \half\rho_-(F^-_{A_L} + F_{A_e}^+), \tag{2b}
\end{align}
where $A_e = A_{\det E}$ is the fixed connection on $\det E$ and the
identification $\ad:\su(E)\simeq \so(\su(E))$ is implicit.
\end{lem}

\begin{pf}
We just consider the first pair of identities, as the second follow
immediately from the fact that $F(A_E) = F(A_{\su(E)})+\half F(A_{\det
E})\otimes \id_E$.  Over any sufficiently small local coordinate
neighborhood $U$ in $X$ we have a local spin structure and a Hermitian spin
bundle $S$ such that $W|_U=S\otimes L^{1/2}$, where $L^{1/2}$ is a
Hermitian line bundle over $U$ such that $(L^{1/2})^{\otimes 2}=L|_U$ and
with induced unitary connection $\half A_L$.  Applying the
Bochner-Weitzenb\"ock identity of \cite[Theorem II.8.17]{LM} to the
Hermitian bundle $S\otimes L^{1/2}\otimes E$ over $U$ with unitary
connection $\cov_A$ given by the tensor product of $\cov^S$, 
$\cov_{2^{-1}A_L}^{L^{1/2}}$, and $\cov_A^E$  gives
$$
D_A^2 = \cov_A^*\cov_A + \quarter R + \rho(F_{L^{1/2}\otimes E}),
$$
where $F_{L^{1/2}\otimes E}$ is the curvature of
the tensor product connection $\cov^{L^{1/2}}_{2^{-1} A_L}\otimes\id_E +
\id_{L^{1/2}}\otimes\cov^E_A$ on $L^{1/2}\otimes E$. Since 
$$
F_{L^{1/2}\otimes E} = F_{2^{-1}A_L}\otimes\id_E + \id_{L^{1/2}}\otimes F_A
= \half F_{A_L}\otimes\id_E + \id_{L^{1/2}}\otimes F_A,
$$
then
$$
\rho(F_{L^{1/2}\otimes E}) 
= \half\rho(F_{A_L})\otimes\id_E + \id_{L^{1/2}}\otimes \rho(F_A),
$$
and hence on $\Om^0(U,W\otimes E)$, we have
$$
D_A^2 = \cov_A^*\cov_A + \quarter R + \half\rho(F_{A_L}) + \rho(F_A),
$$
which is plainly independent of the local splitting $W=S\otimes L^{1/2}$
and so gives an identity on $\Om^0(X,W\otimes E)$.
{}From the decomposition $\rho=\rho_+\oplus\rho_-$ we have
$$
\rho^\pm(F_A)\Phi^\pm = \rho^\pm(F_A^\pm)\Phi^\pm \quad\text{and}\quad
\rho^\pm(F_{A_L})\Phi^\pm = \rho^\pm(F_{A_L}^\pm)\Phi^\pm, 
$$
for any $\Phi^\pm \in \Om^0(X,W^\pm\otimes E)$, 
and so the result follows.
\end{pf}

As in the case of the abelian monopole
equations, the equations \eqref{eq:PT} and Lemma \ref{lem:BWDirac} combine to
give {\em a priori\/} estimates for solutions $(A,\Phi)$ which play
essential role in the proof of existence of the Uhlenbeck compactification.

\begin{lem}\label{lem:L21AaprioriEstAPhi}
Let $X$ be a closed, oriented, Riemannian four-manifold with \spinc
structure $(\rho,W)$, and Hermitian line bundle $\det E$ with fixed unitary
connection. Then there is a positive constant $K_1$, depending only on the
$L^2$ norms of the scalar curvature $R$, the curvatures $F(A_{\det E})$ and
$F(A_{\det W^+})$ such that the following holds.  If $(A,\Phi)$ is an
$L^2_1$ solution on $(\su(E),W^+\otimes E)$ to the $\PU(2)$ monopole
equations \eqref{eq:PT} over $X$, then 
$$
\|\Phi\|_{L^4(X)} \le K_1
\quad\text{and}\quad
\|\cov_A\Phi\|_{L^2(X)} \le K_1.
$$
\end{lem}

\begin{pf}
We may assume that $(A,\Phi)$ is a $C^\8$ pair. 
Then Lemma \ref{lem:BWDirac} gives 
$$
(D_A^*D_A\Phi,\Phi) = (\cov_A^*\cov_A\Phi,\Phi) + \quarter(R\Phi,\Phi)
+ (\rho(F_A^+)\Phi,\Phi) + \half(\rho(F^+_{A_L}+F_{A_e}^+)\Phi,\Phi),
$$
while the first equation in \eqref{eq:PT} gives
$$
\rho(F_A^+) = \rho\vectau_A\rho^{-1}(\Phi\otimes\Phi^*)_{00}.
$$
So, using the second equation in \eqref{eq:PT} and integration by parts, we
obtain
\begin{align*}
 \|\vecvartheta_A\Phi\|_{L^2}^2 &= \|D_A\Phi\|_{L^2}^2 \\
&= \|\cov_A\Phi\|_{L^2}^2  
+ \quarter(R\Phi,\Phi)
+ (\rho\vectau_A\rho^{-1}(\Phi\otimes\Phi^*)_{00}\Phi,\Phi) \\
&\quad + \half(\rho(F^+_{A_L}+F_{A_e}^+)\Phi,\Phi).
\end{align*}
Consequently, Lemmas \ref{lem:sigmaphiphi} and \ref{lem:SigmaPhiPhi},
H\"older's inequality, and the estimate for $\vectau_A$ 
in \eqref{eq:CompactEstPertC0} imply that
$$
\half\|\Phi\|_{L^4}^4 + \|\cov_A\Phi\|_{L^2}^2 
\le \left(\quarter\|R\|_{L^2} + \half\|F^+_{A_L}\|_{L^2} 
+ \half\|F_{A_e}^+\|_{L^2}
+ \|\vecvartheta_A\|_{L^4}^2\right)\|\Phi\|_{L^4}^2.
$$
Thus, if $\Phi\not\equiv 0$, the above inequality and the 
estimate for $\vecvartheta_A$ in \eqref{eq:CompactEstPertC0}
gives the required $L^4$ estimate for
$\Phi$. Then the above inequality and the $L^4$ estimate for
$\Phi$ gives the $L^2$ bound for
$\cov_A\Phi$. The bounds hold trivially if $\Phi\equiv 0$.
\end{pf}

The preceding {\em a priori} $L^4$ estimate on the section $\Phi$ yields
{\em a priori} $L^2$ bounds on the components of the curvature, $F_A^+$ and
$F_A^-$, if $(A,\Phi)$ is a $\PU(2)$ monopole:

\begin{lem}\label{lem:L2aprioriEstFA} 
Continue the notation and hypotheses of Lemma \ref{lem:L21AaprioriEstAPhi}.
Then there is a constant a positive constant $K_2^+$, depending only on the
$L^2$ norms of the scalar curvature $R$, the curvatures $F(A_{\det E})$ and
$F(A_{\det W^+})$ and a positive constant $K_2^-$, depending only on
$K_2^+$ and $p_1(\su(E))$, such that the following holds. If $(A,\Phi)$ is
an $L^2_1$ solution on $(\su(E),W^+\otimes E)$ to the $\PU(2)$ monopole
equations \eqref{eq:PT}, then the following hold:
$$
\|F_A^+\|_{L^2} \le K_2^+ 
\quad\text{and}\quad
\|F_A^-\|_{L^2} \le K_2^-.
$$
\end{lem}

\begin{pf}
Lemma \ref{lem:SigmaPhiPhi} and the first $\PU(2)$ monopole equation in
\eqref{eq:PT} imply that
$$
\|F_A^+\|_{L^2} \le \threequarter\|\Phi\|_{L^4}^2
$$ 
and so the $L^4$ bound for $\Phi$ in Lemma \ref{lem:L21AaprioriEstAPhi}
gives the first estimate. The Chern-Weil integral identity
\eqref{eq:CWIntegral} implies that
$$
\int_X|F_A^-|^2 dV \le \int_X|F_A^+|^2 dV 
+ 8\pi^2\left(c_2(E) - \quarter c_1(E)^2\right),
$$
and so the second estimate follows from the first. 
\end{pf}

Lemmas \ref{lem:L21AaprioriEstAPhi} and
\ref{lem:L2aprioriEstFA} yield an {\em a priori\/} bound $K$ on the
`energy' of a $\PU(2)$ monopole $(A,\Phi)$ in terms of the scalar
curvature $R$, the connections on $\det W^+$ and $\det E$ and $p_1(\su(E))$:
\begin{equation}
\int_X\left(|F_A|^2 + |\Phi|^4 + |\cov_A\Phi|^2\right)\,dV
\le K, \label{eq:MonopoleEnergyBound}
\end{equation}
We use Lemmas \ref{lem:L21AaprioriEstAPhi} and \ref{lem:L2aprioriEstFA} to
provide the `energy bound' assumed in our proof of removable singularities
(Theorem \ref{thm:PTRemovSing}) and we use Lemma \ref{lem:L2aprioriEstFA}
to give a lower bound on the second Chern class of an ideal monopole
appearing in the Uhlenbeck compactification of the moduli space of
$\PU(2)$ monopoles (see \S \ref{subsec:DefnUhlCompact} and the
conclusion of the proof of Theorem \ref{thm:SeqCompact} in \S
\ref{subsec:SeqCompact}). 

As in the case of the $\U(1)$ monopole equations
\cite{KrM,Witten}, the Bochner-Weitzenb\"ock
identity and the maximum principle yield {\em a priori} $C^0$ estimates
for $\Phi$ and $F_A^+$ when $(A,\Phi)$ is a $\PU(2)$ monopole.

\begin{lem}\label{lem:C0EstFAPhi}
Continue the notation and hypotheses of Lemma \ref{lem:L21AaprioriEstAPhi}.
Then there is a non-negative constant $K_3$, depending only on the $C^0$
norms of the scalar curvature $R$ and the curvatures $F(A_{\det E})$ and
$F(A_{\det W^+})$, such that the following holds.  If $(A,\Phi)$ is a $C^1$
solution on $(\su(E),W^+\otimes E)$ to the $\PU(2)$ monopole equations
\eqref{eq:PT} then the following hold:
$$
\|\Phi\|_{C^0(X)}^2 \le K_3
\quad\text{and}\quad
\|F_A^+\|_{C^0(X)} \le K_3.
$$
\end{lem}

\begin{pf}
We may assume that the pair $(A,\Phi)$ is $C^\8$.  The
analogue of Equation (6.18) of \cite{FU} for {\em Hermitian\/} bundles (see
\cite[\S 8.3]{Salamon}) reads
$$
\half\De|\Phi|^2 + |\cov_A\Phi|^2 
= \Real\langle\cov_A^*\cov_A\Phi,\Phi\rangle,
$$
where $\De = d^*d$ on $\Om^0(X,\RR)$.  Since $\Phi$ is continuous,
$|\Phi|$ achieves its maximum at some point $x_0\in X$, so
$\De|\Phi|^2(x_0) = -\cov_{e_\mu}\cov_{e_\mu}|\Phi|^2(x_0)\ge 0$ and thus
at $x_0$:
$$
\Real\langle\cov_A^*\cov_A\Phi,\Phi\rangle\ge |\cov_A\Phi|^2.
$$
Let $\{e_\mu\}$ be a local oriented, orthonormal frame for the tangent bundle
$TX$ such that $(\cov_{e_\mu}e^\nu)_{x_0}=0$, where $\{e^\mu\}$ is the dual
coframe for the tangent bundle $T^*X$. 
Since $D_A\Phi=-\vecvartheta_A\Phi$ by the second equation in
\eqref{eq:PT}, at $x_0$ we have 
\begin{align*}
\langle D_A^*D_A\Phi,\Phi\rangle
&= -\langle D_A^*(\vecvartheta_A\Phi),\Phi\rangle \\
&= -\langle \rho(e^\mu)(\cov_{e_\mu}\vecvartheta_A)\Phi,\Phi\rangle 
-\langle \rho(e^\mu)\vecvartheta_A\cov_{e_\mu}^A\Phi,\Phi\rangle,
\end{align*}
and so (using the inequality $ab\le \quarter a^2 + b^2$) we get the
following bound at $x_0$:
\begin{align*}
|\langle D_A^*D_A\Phi,\Phi\rangle|
&\le 4|\cov_A\vecvartheta_A||\Phi|^2 + 4|\vecvartheta_A||\cov_A\Phi||\Phi| \\
&\le 4|\cov_A\vecvartheta_A||\Phi|^2 
+ 4\left(|\vecvartheta_A|^2|\Phi|^2 + \quarter|\cov_A\Phi|^2\right) \\
&= 4\left(|\cov_A\vecvartheta_A| + |\vecvartheta_A|^2\right)|\Phi|^2 
+ |\cov_A\Phi|^2.
\end{align*}
Recall that the first $\PU(2)$ monopole equation in \eqref{eq:PT} gives
$$
\rho(F_A^+) = \rho\vectau_A\rho^{-1}(\Phi\otimes\Phi^*)_{00}.
$$
The endomorphism $\rho\vectau_A\rho^{-1}(\Phi\otimes\Phi^*)_{00}$ of
$W^+\otimes E$ lies in $\su(W^+)\otimes\su(E)$ and in particular is
Hermitian, so
$$
\Real\langle\rho\vectau_A\rho^{-1}
(\Phi\otimes\Phi^*)_{00}\Phi,\Phi\rangle
= \langle\rho\vectau_A\rho^{-1}(\Phi\otimes\Phi^*)_{00}\Phi,\Phi\rangle,
$$
and similarly for the endomorphism $\rho(F_A^+)$.
We now combine Lemmas \ref{lem:SigmaPhiPhi} and \ref{lem:BWDirac} and the
first equation in \eqref{eq:PT} to get an estimate
for $\Phi$ at the point $x_0$:
\begin{align*}
\half|\Phi|^4
&\le \langle\rho\vectau_A\rho^{-1}
(\Phi\otimes\Phi^*)_{00}\Phi,\Phi\rangle \\
&= \langle\rho(F_A^+)\Phi,\Phi\rangle \\
&= \Real\langle D_A^*D_A\Phi,\Phi\rangle 
- \Real\langle \cov_A^*\cov_A\Phi,\Phi\rangle 
- \quarter\langle R\Phi,\Phi\rangle \\
&\quad - \half\langle\rho(F^+_{A_L}+F_{A_e}^+)\Phi,\Phi\rangle.
\end{align*}
We now combine the last inequality with our estimates for 
$\langle D_A^*D_A\Phi,\Phi\rangle$ and 
$\langle \cov_A^*\cov_A\Phi,\Phi\rangle$ 
to obtain the following bound for $\Phi$ at $x_0$:
\begin{align*}
\half|\Phi|^4
&\le 4\left(|\cov_A\vecvartheta_A| + |\vecvartheta_A|^2\right)|\Phi|^2 
+ |\cov_A\Phi|^2\\
&\quad - |\cov_A\Phi|^2 - \quarter R|\Phi|^2
+ \half\left(|F^+_{A_L}| + |F_{A_e}^+|\right)|\Phi|^2 \\
&\le 4\left(|\cov_A\vecvartheta_A| + |\vecvartheta_A|^2\right)|\Phi|^2 
- \quarter R|\Phi|^2
+ \half\left(|F^+_{A_L}| + |F_{A_e}^+|\right)|\Phi|^2,
\end{align*}
Either $\Phi(x_0)=0$, and so $\Phi$ is
identically zero, or at the point $x_0$ the preceding inequality implies that
$$
|\Phi|^2 \le 8\left(|\cov_A\vecvartheta_A| + |\vecvartheta_A|^2\right) 
- \half \inf_X R + |F^+_{A_L}| + |F_{A_e}^+|,
$$
which gives the first desired estimate. 
Then the second desired estimate follows from Lemma \ref{lem:SigmaPhiPhi},
the first equation in \eqref{eq:PT}, and the estimate for $\vecvartheta_A$
in \eqref{eq:CompactEstPertC0}.
\end{pf}

We use Lemma \ref{lem:C0EstFAPhi} in \S \ref{subsec:SeqCompact} to
show that the curvature of a $\PU(2)$ monopole connection $A$
concentrates at points with integer multiplicities given by the second
Chern classes of limiting (ideal) anti-self-dual connections over $S^4$
(see Lemma \ref{lem:CharClassLimit}).
These $C^0$ estimates yield the following analogue of Witten's
vanishing theorem \cite[\S 3]{Witten} for $\U(1)$ monopoles over
four-manifolds with non-negative scalar curvature.

\begin{cor}\label{cor:Vanishing}
Continue the notation and hypotheses of Lemma \ref{lem:C0EstFAPhi} and
suppose $K_3\le 0$. If $(A,\Phi)$ is a $C^1$
solution on $(\su(E),W^+\otimes E)$ to the $\PU(2)$ monopole equations
\eqref{eq:PT} then $\Phi \equiv 0$ and $F_A^+ \equiv 0$.
\end{cor}

\begin{rmk}
The proof of Lemma \ref{lem:C0EstFAPhi} shows that 
$$
K_3 = \max\left\{0,-\half\inf_X R + 8\|\vecvartheta_A\|_{L^\8_{1,A}(X)}
+ \|F^+_{A_L}\|_{C^0(X)} + \|F_{A_e}^+\|_{C^0(X)}\right\}.
$$
In particular, we see that if $X=S^4$ has its round metric of scalar
curvature $R=1$, standard \spinc structure with $c_1(W^+)=0$ and $F_{A_L} =
0$, the Hermitian bundle $E$ has $c_1(E)=0$ and $F_{A_e}=0$, and we have
$\vecvartheta = 0$, then $\Phi\equiv 0$ and $A$ is an anti-self-dual $\SO(3)$
connection.  
\end{rmk}

\subsection{Scale invariance of the PU(2) monopole equations}
\label{subsec:ScaleInv}
In this section we describe the behavior of the $\PU(2)$ monopole under
rescaling of the Riemannian metric. As is well-known, the anti-self-dual
equation is conformally invariant. Although the $\PU(2)$ monopole equations
are not conformally invariant they are, like the Seiberg-Witten equations,
invariant under constant rescalings of the metric in a sense we describe
below.

Suppose that $\la>0$ is a constant and that the Riemannian metric $g$ on
$T^*X$ is replaced by $\lambda^2g$. Since the Clifford map $\rho$ is
compatible with $g$, the Clifford map 
$$
\rho_{\lambda^2g} \equiv \lambda\rho_g:T^*X \to \End(W)
$$ 
is compatible with $\lambda^2g$, as for any $\alpha\in\Omega^1(X,\RR)$ it
satisfies 
$$
\rho_{\lambda^2g}(\alpha)^\dagger\rho_{\lambda^2g}(\alpha)
=
\lambda^2\rho_{g}(\alpha)^\dagger\rho_{g}(\alpha)
=
\lambda^2g(\alpha,\alpha)\id_W.
$$
It extends in the usual way to a
linear map $\rho_{\lambda^2g}:\Lambda^\bullet(T^*X)\otimes\CC\to\End(W)$. The
Levi-Civita connection on $T^*X$ for the metric $g$ coincides with the 
Levi-Civita connection on $T^*X$ for the rescaled metric $\lambda^2g$, so
the $\SO(4)$ connection on $T^*X$ induced by the unitary connection on $W$
and the Clifford map $\rho_{\lambda^2g}:T^*X \to \End(W)$ is still
torsion free. 

\begin{lem}\label{lem:RescalePT}
If $(A,\Phi)$ is a solution to the $\PU(2)$ monopole equations
\eqref{eq:PT} for the metric $g$ on $T^*X$, then $(A,\lambda\Phi)$ is a
solution to the $\PU(2)$ monopole equations for the rescaled metric
$\lambda^2 g$ on $T^*X$, where $\lambda$ is a positive constant.
\end{lem}

\begin{proof}
The projection $P^+(\lambda^2g)=\half(1+*_{\lambda^2g})$ 
from $\Lambda^2(T^*X)$ to $\Lambda^+(T^*X)$ is given by
$P^+(\lambda^2g) = P^+(g)$, while the induced map
$\rho_{\lambda^2g}:\Lambda^2(T^*X)\to\End(W)$ is given by
$\rho_{\lambda^2g} = \lambda^2\rho_g$. Therefore,
\begin{align*}
\rho_{\lambda^2g}(\vartheta)
&= 
\lambda\rho_g(\vartheta), \quad \vartheta\in\Om^1(X,\CC),
\\
D_A^{\lambda^2g} &= \lambda D_A^g,
\\
\rho_{\lambda^2g}(P^+(\lambda^2g)F_A)
&= \lambda^2\rho_{g}(P^+(g)F_A).
\\
\end{align*}
Consequently, we see from \eqref{eq:PT} that if $(A,\Phi)$ is a solution
for the metric $g$ then $(A,\lambda\Phi)$ is a solution for the metric
$\lambda^2 g$. 
\end{proof}

\begin{rmk} 
By adapting the proof of Theorem II.5.24 in \cite{LM} we see that if $g$ is
replaced by the conformally equivalent metric $h^{-2}g$, then
$D_A^{h^{-2}g}=h^{5/2}D_A^gh^{-3/2}$. 
Thus, while the
proof of Lemma \ref{lem:RescalePT} adapts without change to show that the
first equation in \eqref{eq:PT} is invariant under the transformation
$(A,\Phi)\mapsto (A,h\Phi)$ when $g\mapsto h^{-2}g$, the second
equation (when $\vartheta=0$) is invariant under the transformation
$(A,\Phi)\mapsto (A,h^{3/2}\Phi)$. It is this incompatibility which
prevents the $\PU(2)$ monopole equations from being conformally
invariant, although they are scale invariant in the sense described above.
\end{rmk}

The proof of the existence of an Uhlenbeck compactification (Theorem
\ref{thm:Compactness}) in \S \ref{sec:compact} relies on local
regularity and removable singularity results 
for solutions to \eqref{eq:PT} over the unit ball $B$ in $\RR^4$. The
requirements that the $L^2$ norm of the curvature $F_A$ and the $L^2_1$
norm of section $\Phi$ be sufficiently small  
are met via a rescaling argument. 

\subsection{Removable singularities}
\label{subsec:RemSing}
Given the sharp local elliptic regularity result of Proposition
\ref{prop:L2_1CoulMonoRegLocal} for 
$\PU(2)$ monopoles in Coulomb gauge in $L^2_1$, we can now establish a
removable singularities theorem for $\PU(2)$ monopoles analogous to Theorem
4.1 in
\cite{UhlRem}, in the case of the Yang-Mills equations, and Theorem 8.1 in
\cite{Parker}, in the case of the coupled Yang-Mills equations. Our method
is modelled on the proof of Theorem 4.4.12 in \cite{DK} --- the removable
singularities result for the anti-self-dual equation --- which uses a
local elliptic regularity result for $L^2_1$ solutions to the 
Coulomb gauge and inhomogeneous
anti-self-dual equation (namely, Proposition 4.4.13 in \cite{DK}) and which
in turn has its antecedent in the proof of Theorem 4.5 in
\cite{UhlChern}. This, of course, is not the only possible approach:
Uhlenbeck's original argument \cite{UhlRem} employed a differential
inequality to obtain a pointwise decay estimate for solutions near the
singular point and this was the method generalized by Parker to the case of
certain coupled Yang-Mills equations; see also \cite[\S 9.2]{Salamon} for
a proof of removable singularities for Seiberg-Witten monopoles which also
uses differential inequalities.

It will be convenient to define the following annuli in $X$, given a point
$x_0\in X$ and a positive constant $r$:
\begin{align*}
\Om(x_0;r) &:=\{x\in X:\eighth r <\dist_g(x,x_0)< r\}, \\
\Om'(x_0;r) &:=\{x\in X:\quarter r <\dist_g(x,x_0)< \half r\} 
\Subset \Om(x_0,r).
\end{align*}
If $X=\RR^4$ and $x_0=0$ and $r=1$, we denote $\Om=\Om(1)$ and
$\Om'=\Om'(1)$. We will need the following special case of Proposition
\ref{prop:MonopoleGoodGaugeIntEst}. 

\begin{lem}\label{lem:MonopoleAnnulusEst}
Let $\RR^4$ have a $C^\8$ Riemannian metric $g$, let $\Om\subset\RR^4$ be
an open subset with \spinc structure $(\rho,W)$, let $E$ be a Hermitian
two-plane product bundle over $\Om$, and let $\Om'\Subset\Om$ be an open
subset.  Then there are positive constants $C,\eps$ such that the following
holds.  If $(A,\Phi)$ is a $\PU(2)$ monopole in $C^\8$ on
$(\su(E),W^+\otimes E)$ over $\Om$ with
$$
\|F_A\|_{L^2(\Om)} < \eps,
$$
then there is a $C^\8$ flat connection $\Gamma$ on $\su(E)|_{\Om'}$ such that
$$
\|A-\Ga\|_{L^4(\Om')} + \|\cov_\Ga(A-\Ga)\|_{L^2(\Om')}
\le c\|F_A\|_{L^2(\Om)},
$$
and a $C^\8$ trivialization $\su(E)|_{\Om'}\simeq \Om'\times \su(2)$
taking $\Gamma$ to the product connection.
\end{lem}

The proof of Theorem \ref{thm:PTRemovSing} relies on a cutting off
procedure to `smooth out' the $\PU(2)$ monopole near the singular point,
using a family of cutoff functions which we now define.  
Let $\chi:\RR\to[0,1]$ be a bump
function such that $\chi(t)=1$ if $t\ge \half$ and $\chi(t)=0$ if $t\le
\quarter$. For any $r\in(0,\varrho)$, where $\varrho$ is the injectivity
radius of $(X,g)$, define a $C^\8$ cutoff function on $X$ by setting
$\chi_r(\cdot)=\chi(\dist_g(\cdot,x_0)/r)$. Thus, we have $\chi_r = 0$ on
the ball $B(x_0,\quarter r)$, while $\chi_r = 1$ on
$X-B(x_0,\half r)$ and so $d\chi_r$ is supported in $\Om'(x_0,r)$.

\begin{thm}\label{thm:PTRemovSing}
Let $B\subset\RR^4$ be a geodesic ball with $C^\8$ metric $g$ and center at
the origin, \spinc structure $(\rho,W)$ over $B$, and Hermitian two-plane
bundle $E$ over $B\less\{0\}$.  Suppose $(A,\Phi)$ is a $C^\8$ solution to
the $\PU(2)$ monopole equations \eqref{eq:PT} on $(\su(E),W^+\otimes E)$
with $(\vectau,\vecvartheta)=0$
over the punctured ball $B\less\{0\}$ with {\em finite energy\/},
$$
\int_{B\less\{0\}}\left(|F_A|^2+|\Phi|^4+|\cov_A\Phi|^2\right)\,dV <\8,
$$
Then there are a Hermitian two-plane bundle $\tE$ over $B$ with $\det\tE =
\det E$, a $C^\8$ $\PU(2)$ monopole $(\tA,\tPhi)$ on
$(\su(\tE),W^+\otimes\tE)$ over the ball $B$, and a $C^\8$, determinant-one
unitary bundle isomorphism $u:E|_{B\less\{0\}}\to
\tE|_{B\less\{0\}}$ such that
$$
u(A,\Phi) = (\tA,\tPhi) \quad\text{over}\quad B\less\{0\}.
$$
\end{thm}

\begin{rmk}
We restrict our attention to the case of $(\vectau,\vecvartheta)=0$ in the
$\PU(2)$ monopole equations \eqref{eq:PT} since the holonomy perturbations
are undefined for the $L^2_1$ connections which arise in the proof of
Theorem \ref{thm:PTRemovSing}. However, there is no loss of generality in
making this restriction as the holonomy perturbations vanish near points
where curvature has bubbled off.
\end{rmk}

\begin{pf}
We may suppose without loss of generality that the ball $B$ has radius less
than or equal to one. Since $\|F_A\|_{L^2(B)}<\8$
then $\|F_A\|_{L^2(B_{r_0})}$ tends to zero as $r_0$ tends to
zero.  Hence, for small enough $r_0\in (0,1]$, we may suppose
that $\|F_A\|_{L^2(B_{r_0})}<\eps$, where $\eps$ is the constant of Lemma
\ref{lem:MonopoleAnnulusEst}. Then, for any $r\in(0,r_0)$, Lemma
\ref{lem:MonopoleAnnulusEst} provides a $C^\8$ flat connection $\Gamma_r'$ on
$\su(E)|_{\Om'(r)}$ such that
\begin{equation}
\|A-\Gamma_r'\|_{L^4(\Om'(r))}\le c\|F_A\|_{L^2(\Om(r))},
\label{eq:AuL4Bound}
\end{equation}
where the constant $c$ is {\em independent of $r\in (0,r_0)$\/}. To see
that the constant $c$ is indeed scale invariant, note that by Lemma
\ref{lem:RescalePT} the pair $(A,r^{-1}\Phi)$ is a $\PU(2)$ monopole over
$B$ with respect to the rescaled metric $g_r := r^{-2}g$, so $\Om_g'(r) =
\Om_{g_r}'(1)$ and $\Om_g(r) = \Om_{g_r}(1)$. We then apply Lemma
\ref{lem:MonopoleAnnulusEst} to the annuli $\Om_{g_r}'(1)\Subset
\Om_{g_r}(1)$ and observe that the $L^4$ norm on one-forms and the $L^2$
norm on two-forms are scale invariant.

Lemma \ref{lem:MonopoleAnnulusEst} also provides a $C^\8$ trivialization
$\su(E)|_{\Om'(r)}\simeq \Om'(r)\times\su(2)$ taking $\Gamma_r'$ to the
product connection on $\Om'(r)\times\su(2)$. We can then define a smooth
$\SO(3)$ bundle $V_r$ over $B$ by setting
$$
V_r := 
\begin{cases}
B(\half r)\times\su(2) &\text{ over }B(\half r),
\\
\su(E) &\text{ over }B-B(\quarter r),
\end{cases}
$$
recalling that $\Om'(r) = B(\half r)-\barB(\quarter r)$.  Let $\Gamma_r'$
denote the $C^\8$ flat connection on $V_r$ over the ball $B(\half r)$,
extending $\Gamma_r'$ on $\su(E)|_{\Om'(r)}$ via the product connection on
$B(\half r)\times\su(2)$, and let $E_r$ be the smooth $\U(2)$ bundle over
$B$ with $\det E_r = \det E$ and $\su(E_r) = V_r$ over $B$.

Now let $(A_r,\Phi_r)$ be the $C^\8$ pair on $B$ defined by 
$$
(A_r,\Phi_r) :=
\begin{cases}
(\Gamma_r'+\chi_r(A-\Gamma_r'),\chi_r\Phi) &\text{over }B(\half r),\\
(A,\Phi)          &\text{over }B\less B(\quarter r),
\end{cases}
$$
where we note that $\chi_r=0$ on $B(\quarter r)$ and $\chi_r=1$ on $B-
B(\half r)$.  To estimate the $L^2$ norm of $F_{A_r}$, note that over
$B-B(\half r)$ we have $A_r=A$ and $F_{A_r} = F_A$, while over $B(\half r)$
we have
$$
F_{A_r} = \chi_rF_A + d\chi_r\wedge (A-\Gamma_r')
+ (\chi_r^2-\chi_r)(A-\Gamma_r')\wedge (A-\Gamma_r').
$$
Hence, by \eqref{eq:AuL4Bound}, we have
\begin{align*}
\|F_{A_r}\|_{L^2(\Om'(r))} &\le \|F_A\|_{L^2(\Om'(r))}
+ \|d\chi_r\|_{L^4(\Om'(r))}\|A-\Gamma_r'\|_{L^4(\Om'(r))}
+ \|A-\Gamma_r'\|_{L^4(\Om'(r))}^2.
\end{align*}
Therefore, since $\|d\chi_r\|_{L^4(\Om'(r))}\le c_0$, for a constant $c_0$
independent of $r\in(0,\8)$, there is a positive constant $c$, independent
of $r\in(0,r_0)$, such that
\begin{equation}
\|F_{A_r}\|_{L^2(B_{r_0})} 
\le c\left(\|F_A\|_{L^2(B_{r_0})} + \|F_A\|_{L^2(B_{r_0})}^2\right)
\le c\|F_A\|_{L^2(B_{r_0})}.
\label{eq:RemSingCurvEst}
\end{equation}
Thus, $\|F_{A_r}\|_{L^2(B_{r_0})}$ tends to zero as $r_0\to 0$, uniformly
with respect to $r\in(0,r_0)$.  Fix $r_0$ small enough so that
$\|F_{A_r}\|_{L^2(B_{r_0})}<\eps_1$ for all $r\in(0,r_0)$, where $\eps_1$
is the constant of Theorem \ref{thm:CoulombBallGauge}.  
Hence, there are a family of $C^\8$ flat
connections $\Gamma_r$, $r\in(0,r_0)$, on the $\SO(3)$ bundles $\su(E_r)$
over $B_{r_0}$ such that 
\begin{equation}
d_{\Ga_r}^*(A_r-\Gamma_r)=0\quad\text{and}\quad
\|A_r-\Gamma_r\|_{L^2_{1,\Gamma_r}(B_{r_0})}\le c\|F_{A_r}\|_{L^2(B_{r_0})},
\label{eq:AtaurL2_1Bound}
\end{equation}
for a positive constant $c$ independent of $r\in(0,r_0)$ and a family of
$C^\8$ bundle isomorphisms $\tau_r:\su(E_r)|_{B_{r_0}}\simeq
B_{r_0}\times\su(2)$ inducing $C^\8$ bundle isomorphisms $\tau_r:E_r\simeq
B_{r_0}\times\CC^2$, via a choice of fixed $C^\8$ trivialization $\det
E_r|_{B_{r_0}} = \det E|_{B_{r_0}} \simeq B_{r_0}\times \CC$.

Since $\|\Phi\|_{L^2_{1,A}(B)}<\8$, then $\|\Phi\|_{L^2_{1,A}(B_{r_0})}$
and $\|\Phi\|_{L^4(B_{r_0})}$ tend to zero as $r_0$ as tends to zero.  As
$\Phi_r=\chi_r\Phi$ over $B(\half r)$ and using \eqref{eq:AuL4Bound}, we
have
\begin{align*}
\|\Phi_r\|_{L^2_{1,\Ga_r}(B_{r_0})} 
&\le \|\Phi_r\|_{L^2(B_{r_0})} 
+ \|\cov_{\Ga_r}\Phi_r\|_{L^2(B_{r_0})} \\
&\le \|\Phi_r\|_{L^2(B_{r_0})} 
+ \|(A_r-\Gamma_r)\cdot\Phi_r\|_{L^2(B_{r_0})}
+ \|\cov_{A_r}\Phi_r\|_{L^2(B_{r_0})} \\
&\le c\left(1+\|(A_r-\Gamma_r)\|_{L^4(B_{r_0})}\right)
\|\Phi_r\|_{L^4(B_{r_0})} 
+ \|\cov_{A_r}\Phi_r\|_{L^2(B_{r_0})} \\
&\le c\left(1+\|F_{A_r}\|_{L^2(B_{r_0})}\right)
\|\Phi\|_{L^4(B_{r_0})} 
+ \|\cov_{A_r}\Phi_r\|_{L^2(B_{r_0})} \\
&\le c\|\Phi\|_{L^4(B_{r_0})} 
+ \|\cov_{A_r}(\chi_r\Phi)\|_{L^2(\Om(\quarter r,\half r))} 
+ \|\cov_A\Phi\|_{L^2(\Om(\half r,r_0))} \\
&\le c\|\Phi\|_{L^2_{1,A}(B_{r_0})} + \|d\chi_r\cdot\Phi\|_{L^2(\Om'(r))}
+ \|\chi_r\cov_{A_r}\Phi\|_{L^2(\Om'(r))} \\
&\le c\|\Phi\|_{L^2_{1,A}(B_{r_0})} 
+ \|d\chi_r\|_{L^4(B_{r_0})}\|\Phi\|_{L^4(B_{r_0})} 
+ \|\cov_A\Phi\|_{L^2(\Om'(r))} \\
&\quad + \|(\chi_r -1)(A-\Gamma_r')\cdot\Phi\|_{L^2(\Om'(r))} \\
&\le c\|\Phi\|_{L^2_{1,A}(B_{r_0})} 
+ \|A-\Gamma_r'\|_{L^4(\Om'(r))}\|\Phi\|_{L^4(\Om'(r))},
\end{align*}
using $A_r = \Gamma_r'+\chi_r(A-\Gamma_r') = A+(\chi_r -1)(A-\Gamma_r')$
over $\Om'(r)$, and so
\begin{equation}
\|\Phi_r\|_{L^2_{1,\Ga_r}(B_{r_0})}\le c\|\Phi\|_{L^2_{1,A}(B_{r_0})}, 
\label{eq:PhitaurL2_1Bound}
\end{equation}
for some constant $c$ independent of $r\in(0,r_0)$. Therefore,
$\|\Phi_r\|_{L^2_{1,\Ga_r}(B_{r_0})}$ tends to zero as $r_0\to 0$,
uniformly with respect to $r\in(0,r_0)$.

Hence, the estimates \eqref{eq:RemSingCurvEst},
\eqref{eq:AtaurL2_1Bound}, and \eqref{eq:PhitaurL2_1Bound}
combine to give a uniform bound,
\begin{equation}
\|(A_r-\Gamma_r,\Phi_r)\|_{L^2_{1,\Ga_r}(B_{r_0})} 
\le c\left(\|F_A\|_{L^2(B_{r_0})} + \|\Phi\|_{L^2_{1,A}(B_{r_0})}\right),
\label{eq:APhitaurL2_1Bound}
\end{equation}
for some constant $c$ independent of $r\in(0,r_0)$. Note that
\begin{align*}
&L^2_{1,\Gamma_r}(B_{r_0},\La^1\otimes\su(E_r))
\oplus L^2_{1,\Gamma_r}(B_{r_0},W^+\otimes E_r) \\
&\quad =
L^2_{1,\Gamma}(B_{r_0},\La^1\otimes\su(2))
\oplus L^2_{1,\Gamma}(B_{r_0},W^+\otimes\CC^2),
\end{align*}
via the $C^\8$ isomorphisms $\tau_r:E_r|_{B_{r_0}}\to B_{r_0}\times \CC^2$,
with
\begin{align*}
\|(A_r-\Gamma_r,\Phi_r)\|_{L^2_{1,\Ga_r}(B_{r_0})}
&=
\|(\tau_r(A_r)-\Gamma,\tau_r(\Phi_r))\|_{L^2_{1,\Ga}(B_{r_0})}, 
\\
d_\Ga^*(\tau_r(A_r)-\Gamma) &= 0.
\end{align*}
By the weak compactness of the unit ball in 
the Hilbert space $L^2_1(B_{r_0})$, there is a
sequence $r_\al\to 0$ such that the pairs
$(A_{r_\al}-\Gamma_{r_\al},\Phi_{r_\al})$ converge weakly in $L^2_1(B_{r_0})$
to a limit $(\tA,\tPhi)$ in $L^2_1(B_{r_0})$. For brevity, we denote
$(A_r^\tau,\Phi_r^\tau) := (\tau_r(A_r),\tau_r(\Phi_r))$. 

\begin{claim}\label{claim:LimitMonopole}
Continue the above notation. Then the following hold:
\begin{enumerate}
\item After passing to a subsequence, the pairs
$(A^\tau_{r_\al},\Phi^\tau_{r_\al})$ converge in $C^\8$
over compact subsets of $B_{r_0}\less\{0\}$ to $(\tA^\tau,\tPhi^\tau)$ and so
$(\tA^\tau,\tPhi^\tau)$ is in $C^\8$ over $B_{r_0}\less\{0\}$; 
\item The pair $(\tA^\tau,\tPhi^\tau)$ is an $L^2_1$ solution over
$B_{r_0}$ to the elliptic system
$$
\fs(\tA^\tau,\tPhi^\tau) = 0 \quad\text{and}\quad d_\Ga^*(\tA^\tau-\Ga) = 0.
$$
\end{enumerate}
\end{claim}

\begin{pf}
For any open subset $U\Subset B_{r_0}\less\{0\}$, the ball $B(\half
r)$ containing the support of the cutoff function $\chi_r$ does not
meet $U$ when $r$ is sufficiently small and so
$(A^\tau_r,\Phi^\tau_r)$ is a $\PU(2)$ monopole in $C^\8$ over $U$ in
Coulomb gauge with respect to $(\Ga,0)$ with uniformly $L^2$-small
curvature. Hence, choosing $U'\Subset B_{r_0}\less\{0\}$ so that
$U'\Supset U$, Proposition \ref{prop:L2_1CoulMonoRegLocal} implies that for
any integer $k\ge 1$ and $r$ small enough that
$B(\half r)\cap U' = \emptyset$, we have uniform bounds
\begin{align*}
\|(A^\tau_r-\Gamma,\Phi^\tau_r)\|_{L^2_{k,\Gamma}(U)} 
&\le C\left(\|(A^\tau_r-\Gamma,\Phi^\tau_r)\|_{L^2(U')} 
+ \|F(A^\tau_r)\|_{L^2(U')}\right) \\
&\le C\left(\|\Phi\|_{L^4(B)} + \|F_A\|_{L^2(B)}\right) < \8,
\end{align*}
for some constant $C(g,k,U)$ independent of $r\in(0,r_0)$.  Therefore, by
passing to a subsequence, the pairs $(A^\tau_{r_\al},\Phi^\tau_{r_\al})$
converge in $C^\8$ over compact subsets of $B_{r_0}\less\{0\}$ to a $C^\8$
pair $(\tA^\tau,\tPhi^\tau)$ over $B_{r_0}\less\{0\}$ as $r_\al\to 0$. But
for any $U\Subset B_{r_0}\less\{0\}$ and small enough $r\in
(0,r_0)$, we have $\fs(A_r^\tau,\Phi_r^\tau) = \fs(A_r,\Phi_r) = 0$ over $U$,
and so
$$
\fs(\tA^\tau,\tPhi^\tau) 
= \lim_{\al\to\8}\fs(A^\tau_{r_\al},\Phi^\tau_{r_\al})
= 0\quad\text{over }U.
$$
Hence, $\fs(\tA^\tau,\tPhi^\tau)=0$ over $B_{r_0}\less\{0\}$ and so
$(\tA^\tau,\tPhi^\tau)$ is a 
$\PU(2)$ monopole in $C^\8$ over $B_{r_0}\less\{0\}$.
This proves Assertion (1) of the claim.

Since $(\tA^\tau,\tPhi^\tau)$ is a $C^\8$ monopole over
$B_{r_0}\less\{0\}$, then $\fs(\tA^\tau,\tPhi^\tau)=0\text{ a.e.}$ 
over $B_{r_0}$ and so
$(\tA^\tau,\tPhi^\tau)$ is an $L^2_1$ monopole over $B_{r_0}$.
Let $W^{2,1}_0(B_{r_0})\subset L^2_1(B_{r_0})$ be the closure in
$L^2_1(B_{r_0})$ of the pairs $C^\8_0(B_{r_0},\La^1\otimes\su(E))\oplus
C^\8_0(B_{r_0},W^+\otimes E)$ with compact support in the open ball
$B_{r_0}$. Then, for any $(b,\varphi)\in W^{2,1}_0(B_{r_0})$ we have
\begin{align*}
(b,d_\Ga^*(\tA^\tau-\Gamma))_{L^2(B_{r_0})} 
&= (d_\Ga b,\tA^\tau-\Gamma)_{L^2(B_{r_0})}
= \lim_{\al\to\8}(d_\Ga b,(A^\tau_{r_\al}-\Gamma))_{L^2(B_{r_0})} \\
&= \lim_{\al\to\8}(b,d_\Ga^*(A^\tau_{r_\al}-\Gamma))_{L^2(B_{r_0})} = 0,
\end{align*}
and so $d_\Ga^*(\tA^\tau-\Gamma)=0$, as required. This completes the proof
of the claim.
\end{pf}

By Claim \ref{claim:LimitMonopole}, the pair
$(\tA^\tau,\tPhi^\tau)$ is an $L^2_1$ monopole over $B_{r_0}$
in Coulomb gauge relative to $(\Ga,0)$. {}From the estimate
\eqref{eq:APhitaurL2_1Bound}
and the Sobolev embedding $L^2_1(B_{r_0})\subset L^4(B_{r_0})$, we can
ensure that for sufficiently small $r_0$,
$$
\|(A^\tau_r-\Gamma,\Phi^\tau_r)\|_{L^4(B_{r_0})}  < \eps_2,
$$
for all $r\in(0,r_0)$ and so $\|(\tA^\tau,\tPhi^\tau)\|_{L^4(B_{r_0})}
< \eps_2$, where $\eps_2$ is the constant of Proposition
\ref{prop:L2_1CoulMonoRegLocal} and therefore, $(\tA^\tau,\tPhi^\tau)$ is a
$C^\8$ monopole over $B_{r_0}$. (As usual this means, more
precisely, that there is a $C^\8$ monopole over $B_{r_0}$ which
coincides with $(\tA^\tau,\tPhi^\tau)$ over $B_{r_0}$ except over a subset
of measure zero; since $(\tA^\tau,\tPhi^\tau)$ is already $C^\8$ on the
punctured ball $B_{r_0}\less\{0\}$ then this $C^\8$ monopole is equal to
$(\tA^\tau,\tPhi^\tau)$ except possibly at the origin.)

Finally, the $C^\8$ bundle isomorphisms
$\tau_r:\su(E_r)|_{B_{r_0}}\simeq B_{r_0}\times\su(2)$ may be viewed as 
$\SU(2)$ automorphisms acting on the $\SO(3)$ bundle
$B_{r_0}\times\su(2)$ by initially choosing a fixed $C^\8$ trivialization
$$
E|_{B\less\{0\}} \simeq B\less\{0\} \times \CC^2, \quad\text{with}\quad
\su(E)|_{B\less\{0\}} \simeq B\less\{0\} \times \su(2).
$$
Then Lemma \ref{lem:GaugeSequence} implies that after passing to a
subsequence, the sequence of $\SU(2)$ automorphisms $\tau_{r_\al}$ converges
in $C^\8$ over compact subsets of $B_{r_0}\less\{0\}$ to a $C^\8$ limit
$\si$ over $B_{r_0}\less\{0\}$. Then
$(\tA^\tau,\tPhi^\tau)=(A^\si,\Phi^\si)$ over the punctured ball
$B_{r_0}\less\{0\}$, while $(\tA^\tau,\tPhi^\tau)$ is smooth over
$B_{r_0}$. Thus, the finite-energy $C^\8$  monopole $(A,\Phi)$ 
over the punctured ball $B_{r_0}\less\{0\}$ is equivalent via a $C^\8$,
determinant-one, unitary bundle isomorphism to
$(\tA^\tau,\tPhi^\tau)$ over $B_{r_0}\less\{0\}$, where
$(\tA^\tau,\tPhi^\tau)$ is a $C^\8$ monopole over
$B_{r_0}$. This completes the proof of the theorem.
\end{pf}

\begin{rmk}
The proof of Theorem \ref{thm:PTRemovSing} does not imply, of course, that the
section $\tPhi^\tau$ is zero at the center of the ball $B$. 
Even though the $C^\8$ sections  
$\Phi^\tau_{r_\al}$ are zero at the center, the subsequence only converges in
$L^2_1(B_{r_0})$ over $B_{r_0}$
to a limit $\tPhi^\tau$. Similarly, while Lemma
\ref{lem:C0EstFAPhi} provides a
uniform $C^0$ bound for the sections $\Phi^\tau_{r_\al}$ over $B_{r_0}$, we
would need a uniform, $C^{0,\nu}$ H\"older bound, for some $\nu\in (0,1)$,
in order to extract a convergent subsequence.
\end{rmk}

\subsection{Patching arguments}
\label{subsec:PatchGauge}
The standard proof of the compactness theorem for the moduli space of
anti-self-dual connections employs a patching argument for gauge
transformations to obtain $C^\8$ convergence (modulo gauge transformations)
on compact subsets of $X\less\{x_1,\dots,x_m\}$ for a sequence of
anti-self-dual connections $A_\al$ on a Hermitian
bundle $E$ over $X$. The gauge
transformations that require patching are obtained by repeated application
of Corollary 2.3.9 in \cite{DK} to geodesic balls where the $L^2$ norm of
the curvature $F_{A_\al}$ is less than $\eps$: since the $L^2$ norm of the
curvature is scale invariant, these possibly small balls may be rescaled to
standard size with metrics which are approximately Euclidean as in
\cite[Corollary 2.3.9]{DK}. 

Throughout this subsection, $(A_\al,\Phi_\al)$ will denote a sequence of
$C^\8$ pairs (not necessarily $\PU(2)$ monopoles)
on $(\su(E),W^+\otimes E)$ over $\Omega$ and $u_\al$ will denote a
sequence of $C^\8$ determinant-one, unitary automorphisms of a Hermitian
bundle $E$ (that is, gauge transformations in $\sG_E$ and not $\ssG_E$), where
$\Omega$ is an oriented, Riemannian four-manifold with
\spinc structure $(\rho,W)$. Convergence will mean convergence in $C^\8$
on compact subsets which, as usual, can be replaced by $L^2_{k+1,\loc}$
convergence of $L^2_k$ pairs $(A_\al,\Phi_\al)$ provided $k\ge 2$.

The following four patching results follow almost
immediately from the proofs of Lemmas 4.4.5--4.4.7 and Corollary 4.4.8 in
\cite{DK} (where the sequence of connections $A_\al$ is not assumed to be
anti-self-dual). Their proofs are omitted and instead we refer the reader
to \cite{DK} or \cite{UhlLp} for a detailed account; patching arguments of
this type are used by Uhlenbeck in her proof of Theorem 3.6 \cite{UhlLp},
where the connections (not necessarily anti-self-dual or Yang-Mills)
are just assumed to be in $L^p_1$ and the gauge
transformations are in $L^p_2$ with $p>2$.

\begin{lem}\label{lem:DKLemma4.4.5}
Suppose that $(A_\al,\Phi_\al)$ is a sequence of pairs on $(\su(E),W^+\otimes
E)$ over a base manifold $\Om$ (possibly non-compact), and let
$\tOm\Subset\Om$ be an interior domain. Suppose that there are gauge
transformations $u_\al\in\sG_E$ and  $\tu_\al\in\sG_{E|\tOm}$ such that 
$u_\al(A_\al,\Phi_\al)$ converges over $\Om$ and $\tu_\al(A_\al,\Phi_\al)$
converges over $\tOm$. Then for any compact set $K\Subset\tOm$ there are a
subsequence $\{\al'\}\subset \{\al\}$ and gauge transformations
$w_{\al'}\in\sG_E$ such that $w_{\al'}=\tu_{\al'}$ on a neighborhood of
$K$ and the pairs $w_{\al'}(A_{\al'},\Phi_{\al'})$ converge over $\Om$.
\end{lem}

We have the following two extensions of this result.

\begin{lem}\label{lem:DKLemma4.4.6}
Let $\Om$ be exhausted by an increasing sequence of precompact
open subsets $U_1\Subset U_2\Subset\cdots\Subset\Om$ with $\cup_n U_n =
\Om$. Suppose $(A_\al,\Phi_\al)$ is a sequence of pairs on $(\su(E),W^+\otimes
E)$ over $\Om$ and that for each $n$ there are a subsequence
$\{\al'\}\subset\{\al\}$ and gauge transformations
$u_{\al'}\in\sG_{E|U_n}$ such that $u_{\al'}(A_{\al'},\Phi_{\al'})$
converges over $U_n$. Then there are a subsequence
$\{\al''\}\subset\{\al\}$ and gauge transformations $u_{\al''}\in\sG_E$
such that $u_{\al''}(A_{\al''},\Phi_{\al''})$ converges over $\Om$.
\end{lem}

\begin{lem}\label{lem:DKLemma4.4.7}
Suppose that $\Om = \Om_1\cup\Om_2$ and that $(A_\al,\Phi_\al)$ is a
sequence of pairs on $(\su(E),W^+\otimes E)$ over $\Om$. If there are sequences
of gauge transformations $v_\al\in\sG_{E|\Om_1}$ and $w_\al\in
\sG_{E|\Om_2}$ such that $v_\al(A_\al,\Phi_\al)$ converges over $\Om_1$ and
$w_\al(A_\al,\Phi_\al)$ converges over $\Om_2$, then there are a 
subsequence $\{\al'\}\subset \{\al\}$ and gauge
transformations $u_{\al'}\in\sG_E$ such that 
the pairs $u_{\al'}(A_{\al'},\Phi_{\al'})$ converge
over $\Om$.
\end{lem}
 
Lemmas \ref{lem:DKLemma4.4.5}, \ref{lem:DKLemma4.4.6}, 
and \ref{lem:DKLemma4.4.7}
combine to yield the following analogue of Corollary 4.4.8 in \cite{DK}.

\begin{cor}\label{cor:DKCorollary4.4.8}
Suppose that $(A_\al,\Phi_\al)$ is a sequence of pairs on $(\su(E),W^+\otimes
E)$ over $\Om$ such that the following holds. For each point $x\in\Om$
there are a neighborhood $D$ of $x$, a subsequence $\{\al'\}\subset
\{\al\}$, and gauge transformations $v_{\al'}\in\sG_{E|_D}$ such that
the pairs $v_{\al'}(A_{\al'},\Phi_{\al'})$ converge over $D$. Then there are
a single subsequence $\{\al''\}\subset
\{\al\}$, and gauge transformations $u_{\al''}\in\sG_E$ such that
the pairs $u_{\al''}(A_{\al''},\Phi_{\al''})$ converge over $\Om$. 
\end{cor}

We now assume that $(A_\al,\Phi_\al)$ is a sequence of $\PU(2)$
monopoles on $(\su(E),W^+\otimes E)$ over $\Om$ and obtain the required
convergence from our local elliptic estimates for 
$\PU(2)$ monopoles and Uhlenbeck's
gauge-fixing theorem. The following result is the analogue of Proposition
4.4.9 in \cite{DK} which applies to a sequence of anti-self-dual connections.

\begin{prop}\label{prop:LocalUhlenbeck}
Let $Y$ be an oriented four-manifold with Riemannian metric $g$ and \spinc
structure $(\rho,W)$. Suppose that
$(A_\al,\Phi_\al)$ is a sequence of $\PU(2)$ monopoles, on the bundles
$(\su(E),W^+\otimes E)$ over $Y$, with the following property. For each point
$y\in Y$ there is a geodesic ball $B_g(y,r_y)$ with center $y$, radius
$r_y$, and index $\al_y$ such that
$$
\|F_{A_\al}\|_{L^2(B_g(y,r_y))} < \eps_0, \qquad \al\ge \al_y,
$$
where $\eps_0(g,A_{\det W^+},A_{\det E})$ is the constant of
Corollary \ref{cor:PTLocalCpt}. Then there is a subsequence
$\{\al''\}\subset\{\al\}$ and a sequence of $C^\8$ gauge transformations
$u_{\al''}\in\sG_E$ such that $u_{\al''}(A_{\al''},\Phi_{\al''})$ converges
in $C^\8$ on compact subsets over $Y$.
\end{prop}

\begin{pf}
Fix a point $y\in Y$.
If $B_g(y,r_y)$ is a geodesic ball with center $y$ and $g$-radius $r_y$ then
$B_{g_r}(y,1)$ is a geodesic ball with center $y$ and $g_r$-radius one, where
$g_r = r_y^{-2}g$. Then $(A_\al,r_y^{-1}\Phi_\al)$ is sequence of
$\PU(2)$ monopoles over $B_{g_r}(y,1)$ such that
$$
\|F_{A_\al}\|_{L^2(B_{g_r}(y,1))} < \eps_0,\qquad \al\ge \al_y.
$$
Corollary \ref{cor:PTLocalCpt} implies that there are a subsequence
$\{\al'\}\subset\{\al\}$ and a sequence of gauge transformations
$\{u_{\al'}\}$ over $B_{g_r}(y,1)$ such that the sequence
$u_{\al'}(A_{\al'},r_y^{-1}\Phi_{\al'})$ converges over
$B_{g_r}(y,\half)$. 

Therefore, for each point $y\in Y$, we have a sequence of gauge
transformations $\{u_{\al'}\}$ over $B_g(y,r_y)$ such that the sequence
$u_{\al'}(A_{\al'},\Phi_{\al'})$ converges over $B_g(y,\half r_y)$.
The conclusion now follows from Corollary \ref{cor:DKCorollary4.4.8}. 
\end{pf}

\subsection{Definition of the Uhlenbeck closure}
\label{subsec:DefnUhlCompact}
The definition of the Uhlenbeck closure of the moduli space of solutions to
the perturbed $\PU(2)$ monopole equations \eqref{eq:PT} is slightly more
involved than that of the unperturbed $\PU(2)$ monopole equations
\eqref{eq:UnpertPT}. For this reason it is convenient to define the
Uhlenbeck closure in the unperturbed case before considering the general case.

\subsubsection{Definition of the Uhlenbeck closure for the moduli space of
unperturbed $\PU(2)$ monopoles}
\label{subsubsec:DefnUnpertUhlCompact}
Let $M_{W,E}$ (temporarily) denote the moduli
space of gauge-equivalence classes of solutions $(A,\Phi)$ on
$(\su(E),W^+\otimes E)$ to the unperturbed $\PU(2)$ monopole equations
\eqref{eq:UnpertPT}. We define the Uhlenbeck
closure $\barM_{W,E}$ of the moduli space $M_{W,E}$
by analogy with the definition of the Uhlenbeck closure of
the moduli space of anti-self-dual connections \cite[\S 4.4]{DK}.
The {\em moduli set\/}
$IM_{W,E}$ of unperturbed {\em ideal $\PU(2)$ monopoles\/} on
$(\su(E),W^+\otimes E)$ is given by
$$
IM_{W,E} := \bigcup_{\ell=0}^N M_{W,E_{-\ell}}\times\Sym^\ell(X),
$$
where $N\ge N_p$ and $N_p$ is the constant defined in
\eqref{eq:MaxBubblePoints}, and $E_{-\ell}$ denotes a Hermitian two-plane
bundle with $\det E_{-\ell} = \det E$ and $c_2(E_{-\ell}) = c_2(E)-\ell$,
for $\ell \ge 0$. 

\begin{defn}\label{defn:UhlenbeckTop}
Suppose $[A_\al,\Phi_\al,\by_\al]$ is a sequence of points in $IM_{W,E}$
and $[A_0,\Phi_0,\bx]$ is a point in $IM_{W,E}$, where $(A_\al,\Phi_\al)$
and $(A_0,\Phi_0)$ are monopoles on $(\su(E_\al),W^+\otimes
E_\al)$ and $(\su(E_0),W^+\otimes E_0)$ over $X$, respectively, with
$\det E_\al = \det E_0 = \det E$ and $c_2(E_\al), c_2(E_0)\le
c_2(E)$.  Then the sequence of points $[A_\al,\Phi_\al,\by_\al]$ {\em
converges\/} to $[A_0,\Phi_0,\bx]$ (or, the sequence of triples
$(A_\al,\Phi_\al,\by_\al)$ {\em converges weakly\/} to $(A_0,\Phi_0,\bx)$)
if the following hold:
\begin{itemize}
\item There is a sequence of
$L^2_{k+1,\loc}$ determinant-one, unitary bundle isomorphisms 
$u_\al:E_\al|_{X\less\bx}\to
E_0|_{X\less\bx}$ such that the sequence of $\PU(2)$ monopoles
$u_\al(A_\al,\Phi_\al)$ converges to
$(A_0,\Phi_0)$ in $L^2_{k,\loc}$ over $X\less\bx$;
\item The sequence 
$|F_{A_\al}|^2+8\pi^2\sum_{y\in\by_\al}\de(y)$ converges
in the weak-* topology on measures to $|F_{A_0}|^2 +
8\pi^2\sum_{x\in\bx}\de(x)$.
\item We have $c_2(E) = c_2(E_0) + |\bx|$.
\end{itemize}
\end{defn}

Give $IM_{W,E}$ the Uhlenbeck topology specified by Definition
\ref{defn:UhlenbeckTop} and let
$\barM_{W,E}\subset IM_{W,E}$ be the closure of
$M_{W,E}$ in $IM_{W,E}$. The topological space
$IM_{W,E}$ is second-countable and Hausdorff.

\subsubsection{Definition of the Uhlenbeck closure for the moduli space of
perturbed $\PU(2)$ monopoles}
\label{subsubsec:DefnPertUhlCompact}
The basic idea underlying the choice of holonomy perturbations described in \S
\ref{subsubsec:HolonomyPerturbations} is a generalization of an earlier
construction due to S. K. Donaldson for the moduli spaces of solutions to
the `extended anti-self-dual equations', to which the Freed-Uhlenbeck
generic metrics theorem does not apply \cite[\S IV(v),
pp. 282--286]{DonPoly}.  As in the case of the moduli space of
anti-self-dual connections we shall see in \S \ref{sec:compact} that there
is an upper bound $M$ (which is determined by $g$, $A_{\det W^+}$, $A_{\det
E}$, and $p_1(\su(E))$) on the
total energy $\|F_A\|_{L^2(X)}^2$ for any solution $(A,\Phi)$ to the
perturbed $\PU(2)$ monopole equations \eqref{eq:PT} and so an upper bound
$2M/\eps_0^2$ on the number of disjoint balls $B(x_j,4R_0)$ with energy
greater than or equal to $\half\eps_0^2$. Hence, if $N_b\ge 2M/\eps_0^2 +
1$, at least one ball $B(x_j,4R_0)$ in the collection
$\{B(x_j,4R_0)\}_{j=1}^{N_b}$ has energy less than $\half\eps_0^2$.

Suppose $(A_\al,\Phi_\al)$ is a sequence of $\PU(2)$ monopoles which
converges to an ideal $\PU(2)$ monopole $(A_0,\Phi_0,\bx)$ in
$\sC_{W,E_{-\ell}}\times\Sym^\ell(X)$. If a point $x\in\bx$ lies in a
ball $B(x_j,2R_0)$, the corresponding sections $\fm_{j,l,\be}(A_\al)$
supported on $\barB(x_j,R_0)$ converge to zero (by construction) for all
$l,\be$. Thus, the solution $(A_0,\Phi_0)$ will satisfy a version of the
$\PU(2)$ monopole equations \eqref{eq:PT} with the perturbations supported
on $\barB(x_j,R_0)$ omitted.  (In the situation considered by Donaldson,
ideal extended anti-self-dual connections also satisfy a family of
equations \cite[Eq. (4.37)]{DonPoly} which depend on the bubble points in
$X$.) Therefore, the ideal limit $[A_0,\Phi_0,\bx]$ is a point in the fiber
$\bM_{W,E_{-\ell}}|_{\bx}$ 
over a point $\bx$ in the base $\Sym^\ell(X)$.  Here, $M_{W,E_{-\ell}}|_{\bx}$
is simply the moduli space of solutions to the perturbed monopole equations
\eqref{eq:PT} with the connection energy cutoff functions $\beta_j[A]$ of
\eqref{eq:ConnEnergyCutoff} (used in the definition of the perturbing
sections $\fm_{j,l,\beta}(A)$ of
\eqref{eq:HolonomySection}) replaced by cutoff functions
\begin{equation}
\beta_j[A_0,\bx]
:=
\beta\left(\frac{1}{\eps_0^2}\int_{B(x_j,4R_0)}
\beta\left(\frac{\dist_g(\cdot,x_j)}{4R_0}\right)
\left(|F_{A_0}|^2+8\pi^2\sum_{x\in\bx}\de(x)\right)\,dV\right)
\label{eq:IdealConnEnergyCutoff},
\end{equation}
where $\eps_0$ is the constant of Corollary
\ref{cor:PTLocalCpt}. Then $\bM_{W,E_{-\ell}}$ is the moduli space of 
triples $(A_0,\Phi_0,\bx)$
solving the `lower-level' $\PU(2)$ monopole equations
\begin{align}
F_A^+
- \left(\id + \tau_0\otimes\id_{\su(E)}+\vectau\cdot\vecfm(A,\bx)\right)
\rho^{-1}(\Phi\otimes\Phi^*)_{00} &= 0, 
\label{eq:PTLowerStratum}\\
D_A\Phi + \rho(\vartheta_0)\Phi +
\vecvartheta\cdot \vecfm(A,\bx)\Phi &= 0. \notag
\end{align}
In the more familiar case of the Uhlenbeck compactification of the moduli
space of solutions to the unperturbed $\PU(2)$ monopole equations
\eqref{eq:UnpertPT}, the spaces $\bM_{W,E_{-\ell}}$  
above would be replaced by the products
$M_{W,E_{-\ell}}\times \Sym^\ell(X)$.

The definition of the Uhlenbeck topology for the moduli space of
solutions to the perturbed $\PU(2)$ monopole equations \eqref{eq:PT} is almost
identical to that of the unperturbed case. The only difference is in the
definition of the set of ideal solutions to \eqref{eq:PT}.
In the presence of holonomy perturbations, the Uhlenbeck closure 
$\barM_{W,E}$ is therefore defined to be the closure of
$M_{W,E}$ in
$$
IM_{W,E} := \bigcup_{\ell=0}^N\bM_{W,E_{-\ell}}
\subset \bigcup_{\ell=0}^N \sC_{W,E_{-\ell}}\times\Sym^\ell(X),
$$
where $\bM_{W,E_{-0}} := M_{W,E}$, while $N\ge N_p$ and $N_p$
is the constant defined in \eqref{eq:MaxBubblePoints}.

\subsection{Sequential compactness}
\label{subsec:SeqCompact}
In this section we apply our elliptic regularity and removable singularity
results to prove our main compactness result, namely Theorem
\ref{thm:Compactness}, which asserts the existence of an Uhlenbeck
compactification for the moduli space of $\PU(2)$ monopoles, analogous to
that given by Theorem 4.4.3 \cite{DK} in the case of the moduli space of
anti-self-dual connections.

As in \cite{DK}, the proof of Theorem \ref{thm:Compactness} follows by an
entirely routine argument (which we leave to the reader) from the special
case below which is an analogue of similar compactness results for
anti-self-dual connections; see, for example, \cite[Theorem 4.4.4]{DK},
\cite[Theorem 3.2]{FrM}, \cite[Chapter 8]{FU}, \cite[Theorem
3.1]{Sedlacek}, \cite[Proposition 4.4]{TauPath}, and
\cite[Proposition 5.1]{TauFrame}.

\begin{thm}\label{thm:SeqCompact}
Let $X$ be a closed, oriented, smooth four-manifold with $C^\8$ Riemannian
metric, \spinc structure $(\rho,W)$ with 
\spinc connection, and a Hermitian two-plane
bundle $E$ with unitary connection on $\det E$.  Then there is a positive
integer $N_p$, depending at most on the curvatures of the fixed connections
on $W$ and $\det E$ together with $c_2(E)$, such that for all $N\ge N_p$
the following holds: Any infinite sequence in $M_{W,E}$ has a weakly
convergent subsequence, with limit point in $\cup_{\ell=0}^N
\bM_{W,E_{-\ell}}$.
\end{thm}

\begin{pf}
The basic argument follows that of \cite[pp. 163--165]{DK} and
\cite[Proposition 4.4]{TauPath} for the moduli space of anti-self-dual
connections.  Let $[A_\al,\Phi_\al]$ be a sequence of points in $M_{W,E}$
and let $(A_\al,\Phi_\al)$ be a corresponding sequence of $\PU(2)$
monopoles in $C^\8$ on $(E,W^+\otimes E)$.  By passing to a subsequence we
can assume that the sequence of positive measures $\mu_\al :=
|F_{A_\al}|^2$ on $X$ converges to a measure $\mu_\8$ on $X$ in the weak-*
topology on measures, so
$$
\lim_{\al\to\8}\int_X |F_{A_\al}|^2\,dV
=
\int_X\mu_\8 =: M_\8\le K,
$$
where $K < \8$ is the constant in our universal energy bound 
\eqref{eq:MonopoleEnergyBound} for
a $\PU(2)$ monopole over $X$.  Hence there are at most $M_\8/\eps_0^2$
distinct points in $X$, labelled $\{x_1,\dots,x_m\}$, which do not lie in a
geodesic ball $B(x,r)$ of $\mu_\8$-measure less than $\eps_0^2$, where
$\eps_0$ is the constant of Proposition
\ref{prop:LocalUhlenbeck} and which appears in
\eqref{eq:ConnEnergyCutoff} and \eqref{eq:IdealConnEnergyCutoff}. 
Thus, for any $r>0$, we have
$$
\lim_{\al\to\8}\int_{B(x_i,r)} |F_{A_\al}|^2\,dV 
= \int_{B(x_i,r)}\mu_\8 \ge \eps_0^2,
\qquad i = 1,\dots,m,
$$
and so we may define real numbers $\ka_i\ge \eps_0^2/8\pi^2$ by setting
$$
\ka_i := \lim_{r\to 0}\lim_{\al\to\8}
{\frac{1}{8\pi^2}}\int_{B(x_i,r)} |F_{A_\al}|^2\,dV 
= \lim_{r\to 0}{\frac{1}{8\pi^2}}\int_{B(x_i,r)}\mu_\8.
$$
We may suppose, without loss of generality, that $m\ge 1$. 
If a point $x_i$ lies in a ball $\barB(x_j,2R_0)$ then the holonomy
perturbation sections are zero over $B(x_j,R_0)$ since $8\pi^2\ka_i
\ge \half\eps_0^2$.  
Hence, the points $\{x_1,\dots,x_m\}$ are contained in a
large open subset of $X$ where $(\vectau,\vecvartheta) = 0$.

By passing to a subsequence, Proposition \ref{prop:LocalUhlenbeck} supplies
determinant-one, unitary
gauge transformations $u_\al$ over $X\less\{x_1,\dots,x_m\}$ such that the
sequence $u_\al(A_\al,\Phi_\al)$ converges over $X\less\{x_1,\dots,x_m\}$
to a pair $(A_0,\Phi_0)$ on $(\su(E),W^+\otimes E)|_{X\less\{x_1,\dots,x_m\}}$,
such that the triple $(A_0,\Phi_0,\bx)$ solves the lower-level
$\PU(2)$ monopole equations \eqref{eq:PTLowerStratum}.  Plainly,
$$
\int\limits_{X\less\{x_1,\dots,x_m\}}\left(|F_{A_0}|^2
+|\Phi_0|^4+|\cov_{A_0}\Phi_0|^2\right)\,dV 
\le K < \8.
$$
By the removability of point singularities for finite-energy $\PU(2)$
monopoles (Theorem \ref{thm:PTRemovSing}), there are a Hermitian 
two-plane bundle $E_0$ with $\det E_0 = \det E$ over $X$, 
a $\PU(2)$ monopole $(\tA_0,\tPhi_0)$ on
$(\su(E_0),W^+\otimes E_0)$ and a determinant-one, unitary bundle isomorphism
$u_0$ from $E|_{X\less\{x_1,\dots,x_m\}}$ to
$E_0|_{X\less\{x_1,\dots,x_m\}}$ such that 
$u_0(A_0,\Phi_0)=(\tA_0,\tPhi_0)$ over $X\less\{x_1,\dots,x_m\}$.  

The limiting measure $\mu_\8$ has the form
$$
\mu_\8 = |F_{A_0}|^2 + 8\pi^2\sum_{i=1}^m\ka_i\de_{x_i},
$$
where the $\de_{x_i}$ have unit mass concentrated at the points $x_i$. It
remains to show that the $\ka_i$ are positive integers. For this purpose we
use an argument similar to that used to prove Lemma 3.8 in \cite{FrM} (due
to Taubes) which fits better with our later development of the gluing
theory for $\PU(2)$ monopoles, though one could also use the Chern-Simons
functional for this purpose as in \cite[p. 164]{DK}.  The proof of Theorem
3.2 in \cite{FrM} is a modification of an earlier compactness result,
Proposition 4.4 in
\cite{TauPath}, for connections with $L^2$ bounded curvature but which are
only approximately anti-self-dual in a suitable sense.
Theorem \ref{thm:SeqCompact} follows easily from the next lemma:

\begin{lem}\label{lem:CharClassLimit}
The bundle $E_0$ has Chern classes $c_1(E_0) = c_1(E)$ and
$c_2(E_0)=c_2(E)-\ell$, where $\ell=\sum_{i=1}^m\ka_i$ and the constants
$\ka_i$ are positive integers for $i=1,\dots,m$.
\end{lem}

\begin{pf}
The equality of first Chern classes follows from the preceding
paragraph. The proof that each $\ka_i$ is an integer requires a brief
digression in order to discuss the limiting behavior of the connections
$A_\al$ near the points $x_i\in X$.

Fix an index $i\in\{1,\dots,m\}$, let $\varrho$ be the injectivity radius
of $(X,g)$, and fix a constant $\de\in(0,\quarter\varrho)$.  Choose an
orthogonal frame for $\su(E_0)|_{x_i}$ and use parallel translation via the
connection $A_0$ along radial geodesics from $x_i\in X$ to trivialize
$\su(E_0)$ over the ball $B(x_i,\varrho)$ and let
$w_{0,i}:\su(E_0)|_{B(x_i,\varrho)} \to B(x_i,\varrho)\times \su(2)$ be the
resulting smooth bundle map. We have $\|F_{A_0}\|_{L^\8(X)} \le C$, for
some positive constant $C$, and so $\|F_{A_0}\|_{L^2(B(x_i,\de))} \le
C\de^2$. Thus, we may suppose that $\de$ is fixed small enough that Theorem
\ref{thm:CoulombBallGauge} provides an $\SU(2)$ gauge transformation
$v_{0,i}$ of $B(x_i,\de)\times \su(2)$ such that
$$
\|a_{0,i}\|_{L^4(B(x_i,\de))} + \|\cov_\Ga a_{0,i}\|_{L^2(B(x_i,\de))}
\le c\|F_{A_0}\|_{L^2(B(x_i,2\de))},
$$
where $a_{0,i} := u_{0,i}(A_0)-\Ga\in\Om^1(B(x_i,\de),\su(2))$ 
and $u_{0,i} := v_{0,i}\circ w_{0,i}$ and $\Ga$ is 
the product connection on $B(x_i,\de)\times \su(2)$.
The sequence of connections $u_\al(A_\al)$ converges in $C^\8$ on compact
subsets of the punctured balls $B(x_i,\de)\less\{x_i\}$ to the $C^\8$
connection $A_0$ on $\su(E_0)|_{B(x_i,\de)\less\{x_i\}}$; therefore, the
sequence of connections $u_{0,i}u_\al(A_\al)$ converges in $C^\8$ on compact
subsets of the punctured balls $B(x_i,\de)\less\{x_i\}$ to the $C^\8$
connection $u_{0,i}(A_0)$ on $B(x_i,\de)\less\{x_i\}\times\su(2)$.

Write $u_{0,i}u_\al(A_\al) = \Ga + a_{i,\al}$ over $B(x_i,\de)\less\{x_i\}$,
where $a_{i,\al}\in\Om^1(B(x_i,\de)\less\{x_i\},\su(2))$. Let $\Om(x_i;\half
r,2r)$ denote the open annulus $\{x\in X:\half r < \dist_g(x,x_i) < 2r\}$
in $X$.  Then, for any $r\in (0,\half\de)$, there is an index
$\al_0(r)$ such that
$$
\int\limits_{\Om(x_i,\half r,2r)}\left(|a_{i,\al}-a_{0,i}|^4 
+ |\cov_\Ga(a_{i,\al}-a_{0,i})|^2\right)\,dV < r^4, 
\qquad\al\ge\al_0.
$$
Since $\|F_{A_0}\|_{L^\8(X)} \le C$, we have
\begin{equation}
\int_{B(x_i,2r)}\left(|a_{0,i}|^4 
+ |\cov_\Ga a_{0,i}|^2\right)\,dV 
\le c\int_{B(x_i,2r)}|F_{A_0}|^2\,dV \le Cr^4,
\label{eq:a0Est}
\end{equation}
and therefore
\begin{equation}
\int\limits_{\Om(x_i,\half r,2r)}\left(|a_{i,\al}|^4 
+ |\cov_\Ga a_{i,\al}|^2\right)\,dV \le Cr^4, 
\qquad\al\ge\al_0.
\label{eq:aiSeqEst}
\end{equation}
Let $\chi:\RR\to[0,1]$ be a bump function such that $\chi(t)=0$ for $t\le
\half$ and $\chi(t)=1$ for $t\ge 2$. Define a cutoff function $\chi_{i,r}:X\to
[0,1]$ by setting $\chi_{i,r}(x)=1-\chi(\dist_g(x,x_i)/r)$ so that
$\chi_{i,r}=1$ on $B(x_i,\half r)$ and $\chi_{i,r}=0$ on $X-B(x_i,2r)$.
Fix a Riemannian metric $g_i$ on $S^4$ 
which coincides with $g$ on $B(x_i,\de)=B(n,\de)$ (after identifying the point
$x_i\in X$ with the north pole $n\in S^4$) and extending $g_i$
outside $B(x_i,2\de)=B(n,2\de)$ to a smooth metric on $S^4$. 
Define a sequence of $\SO(3)$ bundles $V_{i,r,\al}$ over $S^4$ by setting 
$$
V_{i,r,\al} := \begin{cases} \su(E)&\text{over }B(n,2r), \\
S^4\less\{n\}\times \su(2) &\text{over }S^4\less\{n\},
\end{cases}
$$
where the identification of the $\SO(3)$ bundles $\su(E)$ and
$S^4\less\{n\}\times \su(2)$ over the annulus $B(n,2r)\less\{n\}=
B(x_i,2r)\less\{x_i\}$ is induced from the the $\SO(3)$ bundle isomorphism
$u_{0,i}\circ u_\al:\su(E)|_{B(n,2r)\less\{n\}}\to B(n,2r)\less\{n\}\times
\su(2)$.  We cut off the
sequence of connections $A_\al$ on $\su(E)$ over the annulus $\Om(x_i;\half
r,2r)$ and thus obtain a sequence of $C^\8$ connections $A_{i,r,\al}$ on
the sequence of $\SO(3)$ bundles $V_{i,r,\al}$ over $S^4$ by setting
$$
A_{i,r,\al} 
:= 
\begin{cases}
A_\al&\text{on }\su(E)|_{B(n,\half r)}, \\
\Ga + \chi_{i,r}a_{i,\al} 
&\text{on }S^4\less\{n\}\times \su(2).
\end{cases}
$$
Recall from Lemma \ref{lem:C0EstFAPhi} that there is a 
constant $C$ independent of $\alpha$ such that
$\|F^+_{A_\al}\|_{L^\8(X)} \le C$ and so
\begin{equation}
\|F^+_{A_\al}\|_{L^2(B(x_i,2r))} \le Cr^2 \quad\text{for all $\al$}.
\label{eq:F+SeqEst}
\end{equation}
Since
$$
F^+_{A_{i,r,\al}}
= \chi_{i,r}F^+_{A_\al} + (d\chi_{i,r}\wedge a_{i,\al})^+ 
 + \chi_{i,r}(\chi_{i,r}-1)(a_{i,\al}\wedge a_{i,\al})^+,
$$
the estimates \eqref{eq:a0Est}, \eqref{eq:aiSeqEst}, and \eqref{eq:F+SeqEst}
imply that
\begin{align*}
\|F^+_{A_{i,r,\al}}\|_{L^2(S^4)}
&\le \|F^+_{A_\al}\|_{L^2(B(x_i,2r))} 
+ \sqrt{2}\|d\chi_{i,r}\|_{L^4(X)}\|a_{i,\al}\|_{L^4(\Om(x_i;\half r,2r))} \\
&\qquad + \sqrt{2}\|a_{i,\al}\|_{L^4(\Om(x_i;\half r,2r))}^2 \\
&\le C(r+r^2), \qquad\al\ge \al_0.
\end{align*}
Therefore, 
\begin{equation}
\lim_{r\to 0}\lim_{\al\to\8}\|F^+_{A_{i,r,\al}}\|_{L^2(S^4)} = 0.
\label{eq:SDBubbleCurvLimit}
\end{equation}
Similarly, as
$$
F_{A_{i,r,\al}}
= \chi_{i,r}F_{A_\al} + d\chi_{i,r}\wedge a_{i,\al}
 + \chi_{i,r}(\chi_{i,r}-1)a_{i,\al}\wedge a_{i,\al},
$$
the estimates \eqref{eq:a0Est} and \eqref{eq:aiSeqEst} yield
\begin{align*}
&\|F_{A_{i,r,\al}}
-F_{A_\al}\|_{L^2(\Om(x_i;\half r,2r))} \\
&\le \|\cov_\Ga a_{i,\al}\|_{L^2(\Om(x_i;\half r,2r))}
+ c\|a_{i,\al}\|_{L^4(\Om(x_i;\half r,2r))}
+ \|a_{i,\al}\|_{L^4(\Om(x_i;\half r,2r))}^2 \\
&< C(r + r^2), \qquad \al\ge \al_0.
\end{align*}
Therefore, 
\begin{equation}
\lim_{r\to 0}\lim_{\al\to\8}
\|F_{A_{i,r,\al}}
-F_{A_\al}\|_{L^2(\Om(x_i;\half r,2r))} = 0.
\label{eq:BubbleCurvDiffLimit}
\end{equation}
We can now complete the argument that the $\ka_i$ are integers:

\begin{claim}\label{claim:BubbleCharge}
There is an $\SU(2)$ bundle $E_i$ over $S^4$ such $c_2(E_i) = \ka_i$ for
each $i=1,\dots,m$.
\end{claim}

\begin{pf}
Fix an index $i\in\{1,\dots,m\}$. Over $S^4$, the $\SO(3)$ bundles
$V_{i,r,\al}$ lift to $\SU(2)$ bundles $E_{i,r,\al}$ with
$V_{i,r,\al} = \su(E_{i,r,\al})$ and $p_1(V_{i,r,\al}) = -4c_2(E_{i,r,\al})$.
The second Chern classes of the $\SU(2)$ bundles $E_{i,r,\al}$ are
given by
$$
c_2(E_{i,r,\al}) = {\frac{1}{8\pi^2}}\int_{S^4}
\left(|F^-_{A_{i,r,\al}}|^2 - |F^+_{A_{i,r,\al}}|^2\right)\,dV,
$$
recalling that the isomorphisms
$\ad:\su(E_{i,r,\al})\to\so(\su(E_{i,r,\al}))$ are implicit and that we view
$F_{A_{i,r,\al}}$ as sections of $\La^2\otimes\su(E_{i,r,\al})$.
Therefore, by \eqref{eq:BubbleCurvDiffLimit} and \eqref{eq:SDBubbleCurvLimit}
and the fact that $F_{A_{i,r,\al}}
 = F_{A_\al}$ on $B(x_i,\half r)$, we have
\begin{align*}
\lim_{r\to 0}\lim_{\al\to\8}c_2(E_{i,r,\al}) 
&= \lim_{r\to 0}\lim_{\al\to\8}{\frac{1}{8\pi^2}}\int_{S^4}
|F^-_{A_{i,r,\al}}|^2\,dV \\
&= \lim_{r\to 0}\lim_{\al\to\8}{\frac{1}{8\pi^2}}\int_{S^4}
\left(|F^-_{A_{i,r,\al}}|^2+|F^+_{A_{i,r,\al}}|^2\right)\,dV \\
&= \lim_{r\to 0}\lim_{\al\to\8}{\frac{1}{8\pi^2}}\int_{B(x_i,2r)}
|F_{A_{i,r,\al}}
|^2\,dV \\
&= \lim_{r\to 0}\lim_{\al\to\8}{\frac{1}{8\pi^2}}\int_{B(x_i,2r)}
|F_{A_\al}|^2\,dV = \ka_i,
\end{align*}
where the final equality follows by definition of $\ka_i$.  Thus, for small
enough $r$ and large enough $\al$, we have $c_2(E_{i,r,\al})=c_2(E_i)$ for
some fixed $\SU(2)$ bundle $E_i$ over $S^4$ and so $\ka_i=c_2(E_i)$,
completing the proof of the claim.
\end{pf}

By Claim \ref{claim:BubbleCharge} the $\ka_i$ are positive integers for
$i=1,\dots,m$.
We can now compute the second Chern class of the limit bundle $E_0$.
The Chern-Weil identity \eqref{eq:CWIntegral} implies that, for all $\al$,
$$
- \quarter p_1(\su(E))
=
c_2(E) - \quarter c_1(E)^2
= {\frac{1}{8\pi^2}}\int_X
\left(|F^-_{A_\al}|^2 - |F^+_{A_\al}|^2\right)\,dV. 
$$
Therefore, by \eqref{eq:F+SeqEst} we have
\begin{align*}
c_2(E) - \quarter c_1(E)^2
&= \lim_{r\to 0}\lim_{\al\to\8}{\frac{1}{8\pi^2}}\int_X
\left(|F^-_{A_\al}|^2 - |F^+_{A_\al}|^2\right)\,dV \\
&= \lim_{r\to 0}\lim_{\al\to\8}{\frac{1}{8\pi^2}}
\int_{X\less\cup_{i=1}^mB(x_i,2r)}
\left(|F^-_{A_\al}|^2 - |F^+_{A_\al}|^2\right)\,dV \\
&\quad + \sum_{i=1}^m\lim_{r\to 0}\lim_{\al\to\8}{\frac{1}{8\pi^2}}
\int_{B(x_i,2r)}
\left(|F^-_{A_\al}|^2 - |F^+_{A_\al}|^2\right)\,dV \\
&= {\frac{1}{8\pi^2}}\int_X\left(|F^-_{A_0}|^2 - |F^+_{A_0}|^2\right)\,dV 
+ \sum_{i=1}^m\ka_i \\
&= c_2(E_0)- \quarter c_1(E_0)^2 + \sum_{i=1}^m\ka_i.
\end{align*}
Now $c_1(E) = c_1(E_0)$ and thus $c_2(E_0)=c_2(E)-\sum_{i=1}^m\ka_i$. This
completes the proof of Lemma \ref{lem:CharClassLimit}. 
\end{pf}

Therefore, after passing to a subsequence, the sequence of points
$[A_\al,\Phi_\al]$ in $M_{W,E}$ converges to an ideal monopole
$[A_0,\Phi_0,\bx]$ in $M_{L,E_0}\times\Sym^\ell(X)$, for some integer
$\ell\ge 0$, and Lemma
\ref{lem:CharClassLimit} implies the Chern classes of the limit bundle
$E_0$ are given by $c_1(E_0) = c_1(E)$ and $c_2(E_0) = c_2(E) - \ell$.

It remains to give an upper bound for the integer $\ell = c_2(E) -
c_2(E_0)$. The Chern-Weil identity \eqref{eq:CWIntegral} implies that
$$
\ell = c_2(E) - c_2(E_0) 
= \frac{1}{8\pi^2}\int_X\left(|F_A^-|^2-|F_A^+|^2\right)\,dV
- \frac{1}{8\pi^2}\int_X\left(|F_{A_0}^-|^2-|F_{A_0}^+|^2\right)\,dV
$$
and thus, by Lemma \ref{lem:L2aprioriEstFA}, we have
\begin{equation}
\ell \le \frac{1}{8\pi^2}\int_X|F_A^-|^2\,dV 
+ \frac{1}{8\pi^2}\int_X|F_{A_0}^+|^2\,dV 
\le N_p,
\label{eq:MaxBubblePoints}
\end{equation}
for some positive integer $N_p=N_p(c_1(E),c_2(E),g,F(A_{\det W}),F(A_{\det
E}))$.  This completes the proof of Theorem \ref{thm:SeqCompact}.
\end{pf}
 
\begin{rmk}
The compactness result in this section for the moduli space of $\PU(2)$
monopoles has an antecedent in \cite[Theorem 3.2]{FrM} 
(due to Taubes) in the following
sense. The latter theorem provides a weak compactness result for
connections $A$ satisfying the `infinite-dimensional part' of the
anti-self-dual equation, namely
$\Pi_{A;\mu}F^+_A = 0$,
where $\mu\notin\Spec d_A^+d_A^*$ and $\Pi_{A;\mu}$
is the $L^2$-orthogonal projection onto the eigenvectors of $d_A^+d_A^*$
with eigenvalue less than $\mu$, together with the curvature bounds
$\|F_A\|_{L^2(X)} + \|d_AF_A^+\|_{L^2(X)} \le C$ 
for some constant $C$ independent of $A$. The analogous point here is that
although $\int_X|F_A|^2\,dV$ is not a topological invariant unless
$F_A^+=0$ or $F_A^-=0$, just as in the case of the $\PU(2)$ monopoles,
it is enough for the purposes of obtaining a weak compactness result to
have uniform bounds on the $L^2$ norm of $F_A$ together with an $L^p$
bound on $F_A^+$ for some $p > 2$.
\end{rmk}

There is one further compactness result we will need, analogous to
Uhlenbeck's original compactness theorem for connections (not necessarily
satisfying any elliptic equation) with $L^p$ bounds
on curvature with $p>2$ \cite[Theorem 1.5]{UhlLp}.

\begin{prop}\label{prop:TopStratumCpt}
Let $p>2$ and $K>0$ be constants.  If $[A_\al,\Phi_\al]$ is
an infinite sequence in $M_{W,E}$ satisfying 
$$
\|F_{A_\al}\|_{L^p(X)} \le K,
$$
then there is a subsequence
$\{\al'\}\subset\{\al\}$ such that the sequence $[A_{\al'},\Phi_{\al'}]$
converges to a point $[A_\8,\Phi_\8]\in M_{W,E}$.
\end{prop}

\begin{pf}
(1) Let $(A_\al,\Phi_\al)$ be a corresponding sequence of $C^\8$ pairs.
Since $p>2$, H\"older's inequality implies that for any geodesic ball
$B(x,r)\subset X$ we have
$$
\|F_{A_\al}\|_{L^2(B(x,r)}
\le cr^{2-(4/p)}\|F_{A_\al}\|_{L^p(B(x,r))} 
\le cKr^{2-(4/p)}.
$$
Hence, for small enough $r$, Proposition
\ref{prop:LocalUhlenbeck} applies and there is a subsequence
$\{\al'\}\subset\{\al\}$ and a sequence of $C^\8$ gauge transformations
$u_{\al'}$ such that the sequence $u_{\al'}(A_{\al'},\Phi_{\al'})$
converges in $C^\8$ to a limit $(A_\8,\Phi_\8)$ over all of $X$, with
no exceptional points. 
\end{pf} 

As we shall see in \S \ref{subsec:SmoothPert}, Proposition
\ref{prop:TopStratumCpt} allows us to work with perturbation parameters
$(\tau_0,\vectau,\vecvartheta)$ and a metric $g$ which are $C^\8$ rather
than just $C^r$, as required by the application of the Sard-Smale theorem
in our proof of transversality in
\S \ref{sec:transv}. 


\section{Transversality}
\label{sec:transv}
In this section we show that for generic perturbation
parameters $(\vartheta_0,\tau_0,\vectau,\vecvartheta)$
the moduli space of solutions to the
perturbed $\PU(2)$ monopole equations
\eqref{eq:PT} is a smooth manifold away from the zero-section and reducible
pairs. 

The outline of the proof is of the now standard form introduced in
\cite{DonApplic} and \cite{FU}. In \S \ref{subsec:ParamModuli}, we
define a parametrized moduli space and explain why transversality for the
moduli space (Corollary \ref{cor:GenericProjection}), follows from
transversality for the parametrized moduli space (Theorem
\ref{thm:SmoothParamModuliSpace}) via the Sard-Smale theorem.  In \S
\ref{subsec:SmoothParamModuli} we show that the parametrized moduli space
is a smooth Banach manifold (Theorem \ref{thm:SmoothParamModuliSpace}).
The proof of Theorem \ref{thm:SmoothParamModuliSpace} relies on the fact
that a $\PU(2)$ monopole which is reducible on an admissible open
subset of the manifold $X$ is reducible on the entire manifold (Theorem
\ref{thm:LocalToGlobalReducible}) and this is proved in \S
\ref{subsec:LocalToGlobal}.

The proof of Theorem \ref{thm:Transversality} does not apply to $\PU(2)$
monopoles which are zero-sections or which are reducible.  We describe the
cokernels of $D\fs$ evaluated at these pairs in the sequel \cite{FL2} to
the present article.

\subsection{The parametrized moduli space}
\label{subsec:ParamModuli}
It is convenient to first consider the question of transversality for the
top stratum $M_{W,E}^{*,0}$ of the Uhlenbeck compactification
$\barM_{W,E}$ and then consider the very slight modification required to obtain
simultaneous transversality for all the lower-level moduli spaces
$\bM_{W,E_{-\ell}}^{*,0}|_\Sigma \subset
\sC_{W,E_{-\ell}}^{*,0}\times\Sigma$, for smooth strata
$\Sigma\subset\Sym^{\ell}(X)$.  

\subsubsection{Transversality for the top-level moduli space}
\label{subsubsec:TopStratumTransv}
The condition in Proposition \ref{prop:DimOfModuliSpace} that the section
$\fs$, vanish transversely is, of course, not necessarily true for all
parameters $(\tau_0,\vartheta_0,\vectau,\vecvartheta)$, on which $\fs$ depends.  As in
the cases of the moduli spaces of anti-self-dual connections
\cite{DonApplic,DK,FU} and Seiberg-Witten monopoles \cite{KrM,Witten}, 
we first show that the family of moduli spaces parametrized by the
perturbations $(\tau_0,\vartheta_0,\vectau,\vecvartheta)$ is smooth and
then apply the Sard-Smale theorem \cite{Smale} to conclude that for generic
perturbations $(\tau_0,\vartheta_0,\vectau,\vecvartheta)$ (that is, a
subset of perturbations which is the complement of some first-category
subset), the moduli space $\fs^{-1}(0) =
M_{W,E}^{*,0}(\tau_0,\vartheta_0,\vectau,\vecvartheta)$ is smooth.

Set $\sP^r_0 := C^r(X,\gl(\La^+))\oplus C^r(X,\Lambda^1\otimes\CC)$
and let $\sP^r := \sP^r_0\oplus\sP^r_\tau\oplus\sP^r_\vartheta$ denote
our Banach space of $C^r$ perturbation parameters and define a
$\ssG_E$-equivariant map
$$
\fS:= (\fS_1,\fS_2)
:\sP^r\times \tsC_{W,E}
\to L^2_{k-1}(\La^+\otimes\su(E))\oplus L^2_{k-1}(W^-\otimes E))
$$
by setting
\begin{align}\label{eq:ParamMonopole}
\fS(\tau_0,\vartheta_0,\vectau,\vecvartheta,A,\Phi) 
 &:=
\begin{pmatrix}
\fS_1(\tau_0,\vartheta_0,\vectau,\vecvartheta,A,\Phi) \\
\fS_2(\tau_0,\vartheta_0,\vectau,\vecvartheta,A,\Phi)
\end{pmatrix} \\
\noindent
& = 
\begin{pmatrix}
F_A^+
-\left(\id + \tau_0\otimes\id_{\su(E)}
+\vectau\cdot\vecfm(A)\right)\rho^{-1}(\Phi\otimes\Phi^*)_{00} \\
D_A\Phi+\rho(\vartheta_0)\Phi+\vecvartheta\cdot \vecfm(A)\Phi
\end{pmatrix}, \notag
\end{align}
where $(A,\Phi)$ is a pair on $(\su(E),W^+\otimes E)$ and the isomorphism
$\ad:\su(E)\simeq \so(\su(E))$ is implicit,
$\ssG_E$ acts trivially on the space of perturbations
$\sP^r$, and so $\fS^{-1}(0)/\ssG_E$ is a subset of
$\sP^r\times\sC_{W,E}$.  We let
$\fM_{W,E}$ denote the {\em parametrized moduli space\/}
$\fS^{-1}(0)/\ssG_E$ and let $\fM_{W,E}^{*,0}=\fM_{W,E}\cap
(\sP^r\times \sC_{W,E}^{*,0})$.

\begin{rmk}
While we assumed for convenience in \S \ref{sec:regularity} and \S
\ref{sec:compact} that the parameters 
$g,\tau_0,\vartheta_0,\vectau,\vecvartheta$ were 
$C^\8$, the only difference if the parameters are only assumed to be
$C^r$, for some finite $r$, is the slight increase in bookkeeping
required to keep track of the regularity of solutions to \eqref{eq:PT}
and other associated elliptic systems.
\end{rmk}

Just as in \S \ref{subsec:DeformationComplex}, 
the $\ssG_E$-equivariant map $\fS$
defines a section of a Banach vector bundle 
$\ubarfV$ over
$\sP^r\times \sC_{W,E}^{*,0}$ with total space
$$
\ubarfV := \sP^r\times\tsC_{W,E}^{*,0}\times_{\ssG_E}
\left(L^2_{k-1}(\La^+\otimes\su(E))\oplus L^2_{k-1}(W^-\otimes E)\right), 
$$
so $\fs :=
\fS(\tau_0,\vartheta_0,\vectau,\vecvartheta,\cdot)$ is a section over $\tsC_{W,E}^{*,0}$ of
the Banach vector bundle
$\fV:=\ubarfV|_{(\tau_0,\vartheta_0,\vectau,\vecvartheta)}$ in
\eqref{eq:BanachVectorBundle}.  In particular, the parametrized moduli
space $\fM_{W,E}^{*,0}$ is the zero set of the section $\fS$ of the vector
bundle $\ubarfV$ over $\sP^r\times \sC_{W,E}^{*,0}$.

\begin{thm}\label{thm:SmoothParamModuliSpace}
The zero set in $\sP^r\times \sC_{W,E}^{*,0}$ 
of the section $\fS$ is regular
and, in particular, the moduli space $\fM^{*,0}_{W,E}$ is a smooth Banach
submanifold of $\sP^r\times \sC_{W,E}^{*,0}$.
\end{thm}

To preserve continuity, we defer the proof of Theorem
\ref{thm:SmoothParamModuliSpace} to \S \ref{subsec:SmoothParamModuli}.
The differential 
$$
D\fS := (D\fS)_{[\tau_0,\vartheta_0,\vectau,\vecvartheta,A,\Phi]}
$$ 
of the section $\fS$ at a point
$[\tau_0,\vartheta_0,\vectau,\vecvartheta,A,\Phi]$ in 
$\sP^r\times\sC^{*,0}_{W,E}$ is given by
\begin{equation}
\label{eq:LinearDfS}
D\fS(\de\tau_0,\delta\vartheta_0,\delta \vectau, \delta\vecvartheta,a,\phi)
= \begin{pmatrix}
D\fS_1(\de\tau_0,\delta\vartheta_0,\delta \vectau,\delta\vecvartheta,a,\phi) \\
D\fS_2(\de\tau_0,\delta\vartheta_0,\delta \vectau,\delta\vecvartheta,a,\phi) 
\end{pmatrix}, 
\end{equation}
where $(a,\phi)\in 
\bK_{A,\Phi}\subset L^2_k(\La^1\otimes\su(E))\oplus L^2_k(W^+\otimes E)$
represents a vector in the tangent space
$(T\sC^{*,0}_{W,E})_{[A,\Phi]} = T_{[A,\Phi]}\sC^{*,0}_{W,E}$ and
$(\de\tau_0,\delta\vartheta_0,\delta\vectau,\delta\vecvartheta)\in\sP^r$. The
differentials in
\eqref{eq:LinearDfS} are given explicitly by
\begin{align}
D\fS_1\left(\de\tau_0,\de\vartheta_0,\de\vectau,\de\vecvartheta,a,\phi\right)
&= 
d^+_Aa -\delta\tau_0\rho^{-1}(\Phi\otimes\Phi^*)_{00} 
\notag\\
&\quad
-\left(\id + \tau_0\otimes\id_{\su(E)}
+\vectau\cdot\vecfm(A)\right)
\rho^{-1}(\phi\otimes \Phi^*+\Phi\otimes\phi^*)_{00} 
\notag\\
&\quad
-\sum_{j,l,\alpha}
(\delta\tau_{j,l,\alpha}\otimes\ad(\fm_{j,l,\alpha}(A))
\rho^{-1}(\Phi\otimes\Phi^*)_{00}
\notag\\
&\quad 
-\sum_{j,l,\alpha}\tau_{j,l,\alpha}\otimes\ad
\left({\frac {\delta \fm_{j,l,\alpha}}{\delta A}}a\right)
\rho^{-1}(\Phi\otimes\Phi^*)_{00}, 
\label{eq:LinearOfS1}\\
D\fS_2\left(\de\tau_0,\de\vartheta_0,\de\vectau,\de\vecvartheta,a,\phi\right) 
&=
D_A\phi + \vecvartheta\cdot\vecfm(A)\phi + \rho(a)\Phi 
\notag\\
&\quad
+ \rho(\delta\vartheta_0)\Phi
+\sum_{j,l,\alpha}
\rho(\delta\vartheta_{j,l,\alpha})\otimes \fm_{j,l,\alpha}(A)\Phi 
\notag\\
&\quad
+ \sum_{j,l,\alpha}
\rho(\vartheta_{j,l,\alpha})\otimes 
\left({\frac {\delta \fm_{j,l,\alpha}}{\delta A}}a\right)\Phi. 
\label{eq:LinearOfS2}
\end{align}
We note that from their definitions in \S
\ref{subsubsec:HolonomyPerturbations}
the perturbations (and their variations) are {\em zero order\/}, unlike the
first order perturbations considered in \cite{TelemanMonopole}.

Recall from the arguments of \S
\ref{subsec:ConfigSpace} that $D\fS(\cdot,d_{A,\Phi}^0\zeta)=0$ for all
$\zeta\in L^2_{k+1}(\su(E))\oplus i\RR_Z$ since $\fS$ is
$\ssG_E$-equivariant. By Proposition
\ref{prop:GlobalReg} we may assume, without loss of generality, that
the pair $(A,\Phi)$ in $\tsC_{W,E}^{*,0}$ is a $C^r$ representative for
the point $[A,\Phi]$ in the zero set 
$\fs^{-1}(0)\subset\sC_{W,E}^{*,0}$.  
Since the tangent space  
$(T\sC^{0,*}_{W,E})_{[A,\Phi]}$ may be
identified with $\bK_{A,\Phi}:=\Ker d_{A,\Phi}^{0,*}$ (see \S
\ref{subsec:ConfigSpace}), we have
$$
D\fS(0,0,0,0,a,\phi) = d_{A,\Phi}^1(a,\phi) 
= (d_{A,\Phi}^{0,*} + d_{A,\Phi}^1)(a,\phi),
$$
for $(a,\phi)\in \bK_{A,\Phi}$, so the differential 
$D\fS|_{\{0\}\times T\sC^{*,0}_{W,E}}$ is Fredholm, where 
$\{0\}\times T\sC^{0,*}_{W,E} = T(\{\tau_0,\vartheta_0,\vectau,\vecvartheta\}\times
\sC^{*,0}_{W,E})$.  
Thus, $\fS$ is a Fredholm section when restricted to the fixed parameter
fibers $\{\tau_0,\vartheta_0,\vectau,\vecvartheta\}\times\sC^{*,0}_{W,E}\subset
\sP^r\times\sC^{*,0}_{W,E}$. The Sard-Smale theorem (in the form of
Proposition 4.3.11 in
\cite{DK}) then implies that there is a first-category subset of the space
$\sP^r$ such that the zero sets in $\sC^{*,0}_{W,E}$ of the sections $\fs =
\fS(\tau_0,\vartheta_0,\vectau,\vecvartheta,\cdot)$ are regular for all $C^r$
perturbations $(\tau_0,\vartheta_0,\vectau,\vecvartheta)$ in the complement
of this subset. Now
$$
M^{*,0}_{W,E}(\tau_0,\vartheta_0,\vectau,\vecvartheta)
=  
\fs^{-1}(0)\cap\sC^{*,0}_{W,E}
$$
and so for generic parameters $(\tau_0,\vartheta_0,\vectau,\vecvartheta)$, 
the moduli space
$M^{*,0}_{W,E}(\tau_0,\vartheta_0,\vectau,\vecvartheta)$ 
is a smooth manifold of the expected dimension. In summary, we have: 

\begin{cor}\label{cor:GenericProjection}
There is a first-category subset of the space $\sP^r$, 
such that for all $C^r$ perturbations 
$(\tau_0,\vartheta_0,\vectau,\vecvartheta)$ in the complement of this subset,
the zero locus of the section $\fs$ is
regular and so the moduli space 
$M^{*,0}_{W,E}(\tau_0,\vartheta_0,\vectau,\vecvartheta)
= \fs^{-1}(0)\cap\sC^{*,0}_{W,E}$ 
is a smooth submanifold of $\tsC^{*,0}_{W,E}$ with the expected dimension. 
\end{cor}

We recall that a subset $S$ of a topological
space $\sP$ is a set of the {\em first category}
if its complement $\sP-S$ is a countable intersection of dense 
open sets or, equivalently, if $S$ is a countable union of closed subsets
of $\sP$ with empty interior; if $\sP$ is a complete metric space, then
Baire's theorem implies that $\sP-S$ is dense in $\sP$ \cite{Rudin}. In our
applications, $\sP$ will either be a Banach or Fr\'echet space (with
a complete metric), so $\sP-S$ will always be dense if $S$ is a
first-category subset.

\subsubsection{Reduction to the case of $C^\8$ parameters}
\label{subsec:SmoothPert}
The restriction to $C^r$ parameters
$(\tau_0,\vartheta_0,\vectau,\vecvartheta)$, necessary to apply the
Sard-Smale theorem in \S \ref{subsubsec:TopStratumTransv}, proves
inconvenient in practice. We shall see that these restrictions can now be
removed, so we need only use $C^\8$ parameters
$(\tau_0,\vartheta_0,\vectau,\vecvartheta)$. Although we did not need the
metric $g$ to be generic in order for our transversality proof to work, we
will nonetheless require the metric $g$ to be generic in the sequels to the
present article and so, {\em a priori\/}, $g$ would also be restricted to a
certain Banach manifold of metrics on $X$.  An argument almost identical to
the one we describe here can be used to show that one need only
consider generic $C^\8$ metrics in those applications.

There is an argument due to Taubes --- for the moduli space of
Seiberg-Witten monopoles --- which reduces the case of transversality for
$C^\8$ parameters to the case of H\"older or Sobolev parameters
\cite[\S 9.4]{Salamon}. (A related result for generic metrics due to Freed
and Uhlenbeck appears as Proposition 3.20 in \cite{FU}, although it is only
stated for the moduli space of anti-self-dual $\SU(2)$ connections with
second Chern class one over a simply-connected, negative definite
four-manifold.)

We adapt Taubes argument here to the case of the moduli spaces of $\PU(2)$
monopoles. We define
\begin{align*}
\sP :&= \sP_0\oplus \sP_\tau\oplus\sP_\vartheta,\\
\sP_0 :&= \Omega^0(X,\gl(\La^+))\oplus\Omega^1(X,\CC), \\
\sP_\tau &= \ell^1_{\delta}(\AAA,\Omega^0(X,\gl(\La^+)))
:= \bigcap_{r=0}^\8\ell^1_{\delta}(\AAA,C^r(X,\gl(\La^+))),\\
\sP_\vartheta &= \ell^1_{\de}(\AAA,\Omega^1(X,\CC))
:= \bigcap_{r=0}^\8\ell^1_{\de}(\AAA,C^r(X,\La^1\otimes\CC)),\notag
\end{align*}
and let $\sP^r$ be the Banach space
\begin{align*}
\sP^r :&= \sP^r_0\oplus \sP^r_\tau\oplus \sP^r_\vartheta \\
&= C^r(X,\gl(\La^+))\oplus C^r(X,\Lambda^1\otimes\CC)
\\
&\qquad\oplus \ell^1_{\delta}(\AAA,C^r(X,\gl(\La^+)))
\oplus\ell^1_{\de}(\AAA,C^r(X,\La^1\otimes\CC)).
\end{align*}
Define metrics $d_0$, $d_\tau$ and $d_\vartheta$ on $\Met(X)$, $\sP_0$,
$\sP_\tau$ and $\sP_\vartheta$ by setting
\begin{align*}
d_0(\tau_{01},\vartheta_{01};\tau_{02},\vartheta_{02}) 
&:= \sum_{r=0}^\8 
\frac{2^{-r}(\|\tau_{01}-\tau_{02}\|_{C^r}
+\|\vartheta_{01}-\vartheta_{02}\|_{C^r})}
{1+\|\tau_{01}-\tau_{02}\|_{C^r}+\|\vartheta_{01}-\vartheta_{02}\|_{C^r}},
\\
d_\tau(\vectau_1,\vectau_2) 
&:= \sum_{r=0}^\8 
\frac{2^{-r}\|\vectau_1-\vectau_2\|_{\ell^1_{\delta}(C^r)}}
{1+\|\vectau_1-\vectau_2\|_{\ell^1_{\delta}(C^r)}}, \\
d_\vartheta(\vecvartheta_1,\vecvartheta_2) 
&:= \sum_{r=0}^\8 
\frac{2^{-r}\|\vecvartheta_1-\vecvartheta_2\|_{\ell^1_{\delta}(C^r)}}
{1+\|\vecvartheta_1-\vecvartheta_2\|_{\ell^1_{\delta}(C^r)}}, 
\end{align*}
and observe that $\sP_0$, $\sP_\tau$ and $\sP_\vartheta$ are
complete metric spaces, with the above metrics inducing the $C^\8$
topologies.  Thus, $\sP$ is a complete metric space with respect to the
product metric $d = d_0\times d_\tau\times d_\vartheta$, which
induces the $C^\8$ topology on $\sP$. The $C^r$ topology of $\sP^r$ is
induced by the product metric $d^r := d_0^r\times
d_\tau^r\times d_\vartheta^r$.

Let $\sP_{\reg}\subset\sP$ be the subspace of parameters for
which the zero set of $\fs_\bp$ is regular and note that
$$
\sP_{\reg} = \bigcap_{n\ge 1}\sP_{n,\reg},
$$
where $\sP_{n,\reg}\subset\sP_{\reg}$ is the subspace of $C^\8$
parameters $\bp$ such that the differential $D\fs_\bp$
is surjective for all pairs $(A,\Phi)$ in
$\fs_\bp^{-1}(0)\cap\tsC_{W,E}^{*,0}$, satisfying
\begin{equation}
\|F_A\|_{L^p} \le n \quad\text{and}\quad 
\nu_0(A,\Phi;\bp) \ge 1/n, \label{eq:PairBound}
\end{equation}
for some $p>2$,
where $\nu_0(A,\Phi;\bp)$ is the least eigenvalue of the Laplacian
$\De_{A,\Phi;\bp}^0 := d_{A,\Phi}^{0,*}d_{A,\Phi}^0$ computed with
respect to the metric $g$. Define, in the analogous way, the subspaces
$\sP^r_{\reg}$ and $\sP^r_{n,\reg}$ of $\sP^r$.
 
The uniform upper $L^p$ bound on $F_A$ precludes
bubbling, while the uniform lower bound on $\nu_0(A,\Phi)$ keeps $(A,\Phi)$
bounded away from the reducible or zero-section pairs. Let $\nu_2(A,\Phi;\bp)$
be the least eigenvalue of the Laplacian 
$\De_{A,\Phi;\bp}^2 := d_{A,\Phi}^1d_{A,\Phi}^{1,*}$ computed with respect
to the parameters $\bp=(\tau_0,\vartheta_0,\vectau,\vecvartheta)$. Then 
$(D\fs_\bp)_{A,\Phi}$ is surjective if and only if $\nu_2(A,\Phi;\bp)>0$.

\begin{claim}\label{claim:RegOpen}
The subspace $\sP_{n,\reg}\subset\sP$ is open in the $C^\8$ topology and
$\sP^r_{n,\reg}\subset\sP^r$ is open in the $C^r$ topology.
\end{claim}

\begin{pf}
Let $\{\bp_\al\}_{\al=1}^\8\subset
\sP\less\sP_{n,\reg}$ be a sequence of parameters
and suppose that $\bp_\al$ converges to $\bp\in\sP$ in the
$C^\8$ topology. Then there is a sequence of solutions $(A_\al,\Phi_\al)$
to \eqref{eq:PT} in $\tsC_{W,E}^{*,0}$, with parameters $\bp_\al$,
which satisfy the bounds in \eqref{eq:PairBound} and for which
$\nu_2(A_\al,\Phi_\al;\bp_\al)=0$. Proposition \ref{prop:TopStratumCpt}
and the $L^p$ bounds in 
\eqref{eq:PairBound} imply that, after passing to a subsequence, there is
a sequence of gauge transformations $\{u_\al\}\subset\sG_E$
such that $u_\al(A_\al,\Phi_\al)$ converges (strongly) in $L^2_k$ to a
solution $(A,\Phi)$ in $\tsC_{W,E}^{*,0}$ to \eqref{eq:PT} which
satisfies the curvature bound in \eqref{eq:PairBound}. Standard
perturbation theory implies that the eigenvalues
$\nu_i(u_\al(A_\al,\Phi_\al);\bp_\al)$ converge to 
$\nu_i(A,\Phi;\bp)$ for $i=0,2$ \cite{Kato}, so the triple $(A,\Phi;\bp)$
satisfies the eigenvalue bound in \eqref{eq:PairBound}. The eigenvalue
$\nu_2(A,\Phi;\bp)$ must be zero --- otherwise the eigenvalues
$\nu_2(u_\al(A_\al,\Phi_\al);\bp_\al)$ would be positive for large enough
$\al$. Hence, $\bp\notin\sP_{n,\reg}$ and so $\sP_{n,\reg}$ is open. The
proof that $\sP^r_{n,\reg}$ is an open subset of $\sP^r$ is identical.
\end{pf}

\begin{claim}\label{claim:Dense}
The subspace $\sP_{n,\reg}\subset\sP$ is dense in the $C^\8$ topology.
\end{claim}

\begin{pf}
By Corollary \ref{cor:GenericProjection}, the space $\sP^r_{\reg}$ is the
complement in $\sP^r$ of a first-category subset and so is dense by Baire's
theorem; clearly, $\sP^r_{n,\reg}$ is also dense in $\sP^r$, since
$\sP^r_{\reg}\subset\sP^r_{n,\reg}$.  Let $\bp\in \sP$ be a $C^\8$
parameter, and let $\{\bp_\al\}\subset\sP^r_{n,\reg}$ be a sequence of
$C^r$ parameters such that $d^r(\bp,\bp_\al)<2^{-\al-1}$, so $\bp_\al$
converges in $C^r$ to $\bp$.  Since $\sP^r_{n,\reg}\subset\sP^r$ is open by
Claim
\ref{claim:RegOpen} and the $C^\8$ parameters $\sP$ are dense in
$\sP^r_{n,\reg}$ we may choose, for each
$\al$, a $C^\8$ parameter $\bp_\al'\in
\sP^r_{n,\reg}$ such that $d^r(\bp_\al,\bp_\al')< 2^{-\al-1}$. Since
$\bp_\al'$ 
is $C^\8$, then $\{\bp_\al'\}\subset\sP_{n,\reg}$ and by construction we
have $d^r(\bp,\bp_\al')<2^{-\al}$ and so the sequence $\{\bp_\al'\}$
converges in $C^r$ to $\bp\in\sP$.

Therefore, for each $r$, we obtain a sequence $\{\bp_\al'(r)\}\subset 
\sP_{n,\reg}$ which converges in $C^r$ to $\bp\in\sP_{\reg}$. But then
the diagonal sequence $\{\bp_\al'(\al)\}\subset 
\sP_{n,\reg}$ converges to $\bp$ in $C^r$ for each $r$ 
(it satisfies $d^r(\bp,\bp_\al'(\al))<2^{-\al}$ for all $\al\ge r$)
and so the sequence converges in $C^\8$ to $\bp\in \sP$, as required.
\end{pf}

{}From Claims \ref{claim:RegOpen} and \ref{claim:Dense}
we conclude that $\sP_{\reg}$ is a countable intersection of dense, open
subsets of $\sP$ and hence is the complement of a first-category subset (in
particular, the subset $\sP_{\reg}\subset\sP$ is dense by Baire's theorem).
Hence, the space $\sP_{\reg}$ of $C^\8$ parameters
$(\tau_0,\vartheta_0,\vectau,\vecvartheta)$ such that the moduli spaces
$M^{*,0}_{W,E}(\tau_0,\vartheta_0,\vectau,\vecvartheta)$ are regular is 
the complement of a first-category subset
of $\sP$. {}From this and Corollary \ref{cor:GenericProjection} we
conclude: 

\begin{cor}\label{cor:GenericSmoothAnalyticProjection}
Let $X$ be a closed, oriented, smooth four-manifold with $C^\8$
metric $g$. There is a first-category subset of the space $\sP$, 
such that for all $C^\8$ perturbations 
$(\tau_0,\vartheta_0,\vectau,\vecvartheta)$ in the complement of this subset,
the zero locus of the section $\fs$ is
regular and so the moduli space
$M^{*,0}_{W,E}(\tau_0,\vartheta_0,\vectau,\vecvartheta) =
\fs^{-1}(0)$ is a smooth submanifold of
$\sC^{*,0}_{W,E}$ with the expected dimension.
\end{cor}

\begin{rmk}
The same argument shows that the standard Freed-Uhlenbeck
generic theorems (specifically, Corollaries 4.3.15, 4.3.18, and 4.3.19 in
\cite{DK} and the refinement Lemma 2.4 in \cite{KMStructure} in the 
non-simply connected case) for the moduli spaces of anti-self-dual
connections on an $\SU(2)$ or $\SO(3)$ bundle over $X$ continue to hold for
the complement of a first-category subset of $C^\8$ metrics, rather than
just $C^r$ metrics.
\end{rmk}

\subsubsection{Simultaneous transversality for the top and lower-level
moduli spaces} 
\label{subsubsec:LowerStratumTransv}
Let $\Sigma$ be a smooth stratum of $\Sym^{\ell}(X)$.  Recall from \S
\ref{subsubsec:DefnPertUhlCompact} that a universal choice of
sufficiently large $N_b$
gaurantees that if $[A,\Phi,\bx]$ is any point in $\barM_{W,E}$ and
$A$ is irreducible, then $A$ has at least one ball $B(x_j,R_0)$ which
supports holonomy perturbations.

The $\PU(2)$ monopole equations cutting out the locus
$\bM_{W,E_{-\ell}}^{*,0}|_\Sigma \subset
\sC_{W,E_{-\ell}}^{*,0}\times \Sigma$ from
the Uhlenbeck compactification $\barM_{W,E}$ are equations for
triples $(A,\Phi,\bx)\in \tsC_{W,E_{-\ell}}^{*,0}\times \Sigma$.
We can again define a $\ssG_E$-equivariant $C^\8$ map
$$
\fS:\sP^r\times \tsC_{W,E_{-\ell}}^{*,0}\times \Sigma
\to L^2_{k-1}(\La^+\otimes\su(E))\oplus L^2_{k-1}(W^-\otimes E)
$$
by setting
$$
\fS(\tau_0,\vartheta_0,\vectau,\vecvartheta,A,\Phi,\bx) 
:= 
\begin{pmatrix}
F_A^+
-\left(\id + \tau_0\otimes\id_{\su(E)}+\vectau\cdot\vecfm(A,\bx)\right)
\rho^{-1}(\Phi\otimes\Phi^*)_{00} \\
D_A\Phi + \rho(\vartheta_0)\Phi
+\vecvartheta\cdot \vecfm(A,\bx)\Phi
\end{pmatrix}. 
$$
The proof of Corollary \ref{cor:GenericSmoothAnalyticProjection}
now shows that $\bM_{W,E_{-\ell}}^{*,0}|_\Sigma \subset
\sC_{W,E_{-\ell}}^{*,0}\times \Sigma$ is a smooth submanifold
of the expected dimension for generic parameters
$(\tau_0,\vartheta_0,\vectau,\vecvartheta)$, 
\begin{equation}
\label{eq:DimOfLowerStrata}
\dim \bM_{W,E_{-\ell}}^{*,0}|_\Sigma
=
\dim M_{W,E_{-\ell}}^{*,0}+\dim \Sigma.
\end{equation}
Furthermore, considering regular values of the projection maps onto the
second factors $\Sigma$, Sard's Theorem also shows that the fibers
$M_{W,E_{-\ell}}^{*,0}|_\bx$ are smooth manifolds of the expected dimension
for generic points $\bx\in\Sym^{\ell}(X)$.  Indeed, the
only tangent vectors in each stratum $\Sigma$ which might not appear in the
image of the projection are those arising from the radial vector on the
annuli $B(x_j,4R_0)\backslash\barB(x_j,2R_0)$.  This observation shows that
the projection from $\bM_{W,E_{-\ell}}^{*,0}$ to $\Sigma$ is transverse to
certain submanifolds of $\Sigma$ which will allow dimension counting
arguments in \cite{FL2}.  Issues related to dimension counting in the
presence of holonomy perturbations are also discussed by Donaldson in
\cite[pp. 282--287]{DonPoly}. We can now conclude the proof of our main
transversality result:

\begin{proof}[Proof of Theorem~\ref{thm:Transversality}, given 
Corollary \ref{cor:GenericSmoothAnalyticProjection}]
For the case $\ell=0$, the transversality assertion is given by Corollary
\ref{cor:GenericSmoothAnalyticProjection} and the dimension formula is
provided by Proposition \ref{prop:DimOfModuliSpace}. The case $\ell>0$ then
follows from the discussion in the preceding paragraphs.
\end{proof}

\subsection{Smoothness of the parametrized moduli space}
\label{subsec:SmoothParamModuli}
We prove Theorem \ref{thm:SmoothParamModuliSpace} by
showing that the $\ssG_E$-equivariant map
$\fS:\sP^r
\times\tsC_{W,E}^{*,0}\to L^2_{k-1}(\La^+\otimes\su(E))\oplus
L^2_{k-1}(W^+\otimes E)$ vanishes transversely and
so the parametrized moduli space
$\fM^{*,0}_{W,E}=\fS^{-1}(0)/\ssG_E$ is a smooth Banach
manifold. The broad strategy is reminiscent of that of 
\cite[Chapter 3]{FU} and \cite[\S 4.3.5]{DK},
where the analogous result is established for the moduli space of
anti-self-dual connections parametrized by the Banach space of $C^r$
(conformal classes of) metrics: our proof of surjectivity of the
differential $D\fS$ at a point $(\tau_0,\vartheta_0,\vectau,\vecvartheta,A,\Phi)$ in
the zero set $\fS^{-1}(0)$ ultimately relies on the fact that, for
$\Phi\not\equiv 0$, a monopole $(A,\Phi)$ which is reducible on an {\em
admissible\/} open subset of $X$ is necessarily reducible over all of $X$
(Theorem
\ref{thm:LocalToGlobalReducible}). If there are sections
$\fm_{j,l,\alpha}(A)$ which are non-zero on $B(x_j,R_0)$, we say that
the connection $A$ has {\em holonomy perturbations supported on\/}
$B(x_j,R_0)$; the set $\{\fm_{j,l,\alpha}(A)\}_{l=1}^3$ then spans
$\su(E)|_{B(x_j,R_0)}$ for at least one index $\alpha$.  Then an admissible
open set for the pair $(A,\Phi)$ is one containing $\barB(x_j,R_0)$ for all
$j$ such that $\beta_j[A]>0$.  (The supports of all the sections
$\fm_{j,l,\alpha}(A)$ are contained in $\cup_{j=1}^{N_b}\barB(x_j,R_0)$
and so any open subset of $X$ containing $\cup_{j=1}^{N_b}\barB(x_j,R_0)$
is admissible.) Recall from \S
\ref{subsubsec:DefnPertUhlCompact} that 
because of our choice of $N_b$, if
$[A,\Phi]$ is any point in $M_{W,E}^{*}$, then at least one
ball $B(x_j,R_0)$ supports holonomy perturbations for $A$.

Theorem \ref{thm:SmoothParamModuliSpace} is an almost immediate
consequence of

\begin{prop}
\label{prop:Surjectivity}
Suppose $(\tau_0,\vartheta_0,\vectau,\vecvartheta,A,\Phi)$ is a point in 
$\fS^{-1}(0)$ with $A$ irreducible and $\Phi\not\equiv 0$.
If $(v,\psi)$ is in the cokernel of the differential
$D\fS$ at the point $(\tau_0,\vartheta_0,\vectau,\vecvartheta,A,\Phi)$, then
$(v,\psi)|_{B(x_j,R_0)}\equiv 0$ for each
ball $B(x_j,R_0)$ supporting holonomy perturbations for $A$.
\end{prop}

We first observe that elements of the cokernel of $D\fS :=
D\fS_{(\tau_0,\vartheta_0,\vectau,\vecvartheta,A,\Phi)}$ have a
restricted form of the unique continuation property (sufficient for
our purposes) by Aronszajn's theorem \cite{Aron}:

\begin{lem}
\label{lem:NonVanishingCokernel}
If $(v,\psi)\in\Ker D\fS(D\fS)^*$ and $(v,\psi)|_U\equiv 0$ on some
non-empty open subset $U\subset X$ containing all balls $B(x_j,R_0)$
supporting holonomy perturbations for $A$, then $(v,\psi)\equiv 0$ on $X$.
\end{lem}

\begin{proof}
By hypothesis, the pair $(v,\psi)$ solves the second-order elliptic equation
$$
(D\fS)(D\fS)^*(v,\psi) = 0 \quad\text{on $X$},
$$
where the Laplacian
$D\fS = (D\fS_1,D\fS_2)$, given by equations (4.3) and (4.4), has
$C^{r-1}$ coefficients and $(v,\psi)$ is at least $C^{r+1}$. Also,
$(v,\psi) \equiv 0$ on the set of closed balls supporting holonomy
perturbations,
$$
\barB_I(A) := \bigcup_{j\in I(A)}\barB(x_j,R_0),
$$
where $I(A) := \{1\le j\le N_b: \beta_j[A]>0\text{ and }
A|_{B(x_j,2R_0)}\text{ is irreducible}\}$. Now, on the subset $X-B_{I(A)}$
where all of the holonomy perturbations and in particular their derivatives
with respect to $A$ vanish (see their definition in \S
\ref{subsubsec:HolonomyPerturbations}), the Laplacian $(D\fS)(D\fS)^*$ is a
purely differential operator. In particular, it extends to a differential
operator with $C^{r-1}$ coefficients over $X$, say $(D\fS^0)(D\fS^0)^*$,
given by the linearized $\PU(2)$ monopole equations
\eqref{eq:LinearOfS1} and \eqref{eq:LinearOfS2} with all terms involving
holonomy perturbations and their derivatives set equal to zero. On the
other hand, the pair $(v,\psi)$ also solves the resulting second-order
elliptic {\em differential\/} equation
$$
(D\fS^0)(D\fS^0)^*(v,\psi) = 0 \quad\text{on $X$},
$$
since $(v,\psi)=0$ on $\barB_I(A)$ and $D\fS^0 \neq D\fS$ only on $B_I(A)$,
while $D\fS^0 = D\fS$ on $X-B_{I(A)}$.

Without loss of generality, we may scale the Dirac equation in (2.21)
by $1/\sqrt{2}$. We then have
$$
D\fS^0(D\fS^0)^*
= 
\left(\begin{matrix}
d^+_Ad_A^{+,*} & 0 \\ 0 & \half D_AD_A^*
\end{matrix}\right)
+ \text{First-order differential terms},
$$
and so by the Bochner formulas of Lemma \ref{lem:BWDirac} and [30,
Eq. (6.26)], the Laplacian $D\fS^0(D\fS^0)^*$ is a second order elliptic
differential operator with scalar principal symbol (given by the metric
$\half g$ on $T^*X$). The desired conclusion then follows from Aronszajn's
unique continuation theorem \cite{Aron}.
\end{proof}

\begin{rmk}
\begin{enumerate}
\item Aronszajn's theorem does not apply without the given restriction on
the open set $U$ the proof of Lemma \ref{lem:NonVanishingCokernel}, as the
Laplacian $(D\fS)(D\fS)^*$ is not a purely differential operator over all
of $X$.  One can see from equations \eqref{eq:LinearOfS1} and
\eqref{eq:LinearOfS2} that the problem terms are those appearing in the
last line of each displayed equation: the operator
$\delta\fm_{j,\ell,\alpha}/\delta A$ acting on $a\in\Omega^1(\su(E))$ is an
integral operator, as is clear from the formula for the differential of the
holonomy with respect to the connection in \eqref{eq:DifferentialHolonomy}.
\item Lemma \ref{lem:NonVanishingCokernel} can also be proved without
using Aronszajn's theorem explicitly and instead applying the
Agmon-Nirenberg unique continuation theorem (Theorem
\ref{thm:UniqueFlow}) to the equation $(D\fS)^*(v,\psi) = 0$, and mimicking
the existing application in the proof of Theorem
\ref{thm:LocalToGlobalReducible}. Indeed, this second proof of Lemma
\ref{lem:NonVanishingCokernel} is virtually identical to the proof of
Theorem \ref{thm:LocalToGlobalReducible}.  We leave the details to the
interested reader, as the preceding use of Aronszajn's theorem appears
easier to us. We note that Aronszajn's theorem can be derived from that of
Agmon-Nirenberg (see \cite{Agmon}).
\end{enumerate}
\end{rmk}

If $(v,\psi)$ is an $L^2_{k-1}$ element of
the cokernel of $D\fS$,
then elliptic regularity for the Laplacian $D\fS(D\fS)^*$, with $C^{r-1}$
coefficients, implies that $(v,\psi)$ is in $C^{r+1}$. Our proof of
Theorem \ref{thm:SmoothParamModuliSpace} also relies on the following
unique `unique continuation' result for reducible
monopoles:

\begin{thm} 
\label{thm:LocalToGlobalReducible}
If $(A,\Phi)$ is a $C^r$ solution to the perturbed $\PU(2)$ monopole
equations \eqref{eq:PT} 
with $\Phi\not\equiv 0$ over a connected, oriented, smooth four-manifold
$X$ with $C^r$ Riemannian metric and $(A,\Phi)$ is reducible on a non-empty
open subset $U\subset X$
with $\barB(x_j,R_0)\subset U$ for all $j$ such that $\beta_j[A]>0$,
then $(A,\Phi)$ is reducible on $X$.
\end{thm}

The proof of Theorem \ref{thm:LocalToGlobalReducible} is
lengthy, so we defer it to \S \ref{subsec:LocalToGlobal}.

\begin{proof}[Proof of Theorem~\ref{thm:SmoothParamModuliSpace}, given 
Proposition~\ref{prop:Surjectivity} and
Theorem~\ref{thm:LocalToGlobalReducible}] Let
$(\tau_0,\vartheta_0,\vectau,\vecvartheta,A,\Phi)$ be a $C^r$
representative for a point in $\fM_{W,E}^{*,0}$, so that $A$ is
irreducible and $\Phi\not\equiv 0$, and suppose $(v,\psi)$ is in the
cokernel of $D\fS$. By definition of $N_b$ in \S
\ref{subsubsec:HolonomyPerturbations}, the set $J(A):= \{j:1\le j\le
N_b\text{ and }\beta_j[A]>0\}$ is non-empty.  Since $A$ is irreducible
on $X$, Theorem \ref{thm:LocalToGlobalReducible} implies that
$A|_{B(x_{j'},2R_0)}$ must be irreducible for some
$j'\in\{1,\dots,N_b\}$ such that $\beta_{j'}[A]>0$; otherwise,
$A|_{B(x_{j'},2R_0)}$ would be reducible for all $j$ such that
$\beta_j[A]>0$ and Theorem \ref{thm:LocalToGlobalReducible} would
imply that $A$ would be reducible over all of $X$, contradicting our
assumption that $A$ is irreducible.  But then Proposition
\ref{prop:Surjectivity} implies that $(v,\psi)|_{B(x_j,R_0)}
\equiv 0$ for all $j$ such that $\beta_j[A]>0$ and $A|_{B(x_j,2R_0)}$ is
irreducible and so $(v,\psi) \equiv 0$ on $X$ by Lemma
\ref{lem:NonVanishingCokernel}.
\end{proof}

The proof of Proposition \ref{prop:Surjectivity} occupies the remainder of
this subsection.  We first note that since $\Phi$ is in the kernel
$D_A+\vecvartheta\cdot\vecfm(A)$, it has the unique continuation
property by Aronszajn's Theorem \cite{Aron}:

\begin{lem}
\label{lem:PerturbedDiracUniqueContinuation}
If $(D_A+\rho(\vartheta_0)+
\vecvartheta\cdot\vecfm(A))\Phi=0$ and $\Phi|_U\equiv 0$ for
some non-empty open subset $U\subset X$, then $\Phi\equiv 0$.
\end{lem}

\begin{pf}
The perturbed Dirac operator
$D_A+\rho(\vartheta_0)+\vecvartheta\cdot\vecfm(A)$ differs from $D_A$
by a zeroth order term and so
$$
(D_A+\rho(\vartheta_0)+\vecvartheta\cdot\vecfm(A))^*
(D_A+\rho(\vartheta_0)+\vecvartheta\cdot\vecfm(A))
= D_A^*D_A + \text{First order terms}.
$$
The Bochner formula of Lemma \ref{lem:BWDirac} then implies that the above
Laplacian is a second order elliptic differential operator with scalar
principal symbol (given by the metric $g$ on $T^*X$). The conclusion now
follows from Aronszajn's unique continuation theorem \cite{Aron}.
\end{pf}

We shall use the following linear algebra result to show that $v \in
C^{r+1}(X,\La^+\otimes\su(E))$ vanishes on a ball:

\begin{lem}
\label{lem:PointwiseSurjectiveForQuadratic}
Let $M,N$ be elements of $(\Lambda^+\otimes\su(E))|_x$.  Suppose
$\fm_1,\fm_2,\fm_3$ span $\su(E)|_x$. If
\begin{equation}
\langle\tau_0M,N\rangle +
\sum_{l=1}^3\langle(\tau_l\otimes\ad(\fm_l))M,N\rangle = 0
\label{eq:LemOrthogHypothesis}
\end{equation}
for all $\tau_0,\tau_1,\tau_2,\tau_3
\in\gl(\Lambda^+)|_x$,
then either $M=0$ or $N=0$.
\end{lem}

\begin{pf}
If $\langle\tau_0M,N\rangle=0$ for all
$\tau_0\in\gl(\Lambda^+)|_x$, then by 
the proof of \cite[Lemma 3.7]{FU}, the images in $\su(E)|_x$ 
of $M, N \in \Hom(\Lambda^+|_x,\su(E)|_x)$ are orthogonal.  
(Although their lemma refers to a element of 
$\Hom(\Lambda^+|_x,\su(E)|_x)$ and an element of 
$\Hom(\Lambda^-|_x,\su(E)|_x)$, we can choose any isomorphism
between $\Lambda^+|_x$ and $\Lambda^-|_x$ to translate the result.)
We can therefore assume that
$M$ has rank one and $N$ has rank less than or equal to two. (If $M$ is
rank two and $N$ is rank one, we can reverse their
roles by using adjoints.)  If both $M\neq 0$ and $N\neq 0$ let
\begin{align*}
M & = u\otimes m, \\
N & = v_1\otimes n_1+ v_2\otimes n_2,
\end{align*}
where $u,v_1,v_2\in\Lambda^+|_x$ and $m,n_1,n_2\in\su(E)|_x$. Since
the images of $M$ and $N$ in $\su(E)|_x$ are orthogonal we have $\langle
m,n_1\rangle=0=\langle m,n_2\rangle$; without loss of generality we can
assume that $\langle n_1,n_2\rangle =0$.  
(If $N$ is rank one, $n_2$ can be any element of $\su(E)|_x$ completing
$m,n_1$ to an orthogonal basis of $\su(E)|_x$.)
Under the isomorphism $\su(E)|_x\simeq
\RR^3$ the adjoint representation is given by the cross-product.  We can
find $f_l\in\RR$ such that $n_2=\sum_l f_l \fm_l$, so
\begin{align}
\sum_{l=1}^3\langle\ad(f_l \fm_l)m,n_1\rangle 
& = \langle [n_2,m],n_1\rangle \neq 0, \label{eq:Bracket1}\\
\sum_{l=1}^3\langle\ad(f_l \fm_l)m,n_2\rangle 
& = \langle [n_2,m],n_2\rangle = 0.
\label{eq:Bracket2}
\end{align}
By assumption, $M\ne 0$ and $N\ne 0$, so $u\ne 0$ and either $v_1 \ne 0$ or
$v_2 \ne 0$; we may suppose without loss of generality that $v_1 \ne 0$.
Thus we can find $\tau\in\gl(\Lambda^+)|_x$ such that $\tau u=v_1$
and so choosing $\tau_l = f_l\tau$ for $l=1,2,3$ and $\de\tau_0=0$, we
have 
\begin{align*}
\sum_{l=1}^3\langle \left(\tau_l\otimes \ad(\fm_l)\right)M,N\rangle 
&= \sum_{l=1}^3\langle \tau\otimes\ad(f_l \fm_l)M,N\rangle 
= \langle(\tau\otimes \ad(n_2))M,N\rangle \\
&= \langle(\tau\otimes\ad(n_2))(u\otimes m), 
v_1\otimes n_1 + v_2\otimes n_2\rangle \\
&= \langle \tau u\otimes [n_2,m],v_1\otimes n_1 + v_2\otimes n_2\rangle \\
&= \langle \tau u,v_1\rangle \langle [n_2,m],n_1\rangle + 
\langle \tau u,v_2\rangle \langle [n_2,m],n_2\rangle \\
&= |v_1|^2\langle [n_2,m],n_1\rangle \ne 0 
\quad\text{by \eqref{eq:Bracket1} and \eqref{eq:Bracket2}},
\end{align*}
contradicting our hypothesis in \eqref{eq:LemOrthogHypothesis}. Hence,
either $M=0$ or $N=0$, as desired.
\end{pf}

\begin{rmk}
\label{rmk:HigherRank}
Lemma \ref{lem:PointwiseSurjectiveForQuadratic} does not hold
if the rank of $E$ is greater than two.  If $a,b\in\su(E)$
the above arguments would only allow one to conclude that
$$0=\langle [\fm,a],b\rangle =\langle \fm,[a,b]\rangle$$
for all $\fm\in\su(E)$ so $[a,b]=0$ and $a,b$ are simultaneously
diagonalizable.  The subspace of diagonal elements of $\su(n)$ has
dimension $n-1$, so this would not contradict the orthogonality
of $a,b$ if $n>2$.
\end{rmk}

\begin{lem}
\label{lem:DS1Surjectivity}
Continue the hypotheses of Proposition \ref{prop:Surjectivity} and
suppose $B(x_j,R_0)$ is a ball supporting holonomy
perturbations for $A$. Then $v\equiv 0$ on $B(x_j,R_0)$.
\end{lem}

\begin{pf}
By hypothesis, there are holonomy sections
$\fm_l \equiv \fm_{j,l,\alpha}(A)$, $l=1,2,3$, which span
$\su(E)|_y$, for any point $y\in B(x_j,R_0)$. Let
$\de\tau_l := \delta\tau_{j,l,\alpha} \in
\Om^0(\gl(\La^+))$, $l=1,2,3$, denote corresponding
coefficients, and let $\delta\vectau$ be a sequence with all
other coefficients zero.

By the hypothesis of Proposition \ref{prop:Surjectivity} we have
$(D\fS(\de\tau_0,\de\vartheta_0,\de\vectau,\de\vecvartheta,a,\phi),
(v,\psi))_{L^2}=0$ for all
$(\de\tau_0,\de\vartheta_0,\de\vectau,\de\vecvartheta,a,\phi)$ and so
\begin{align*}
0
& =
(D\fS(\delta\tau_0,0,\de\vectau,0,0,0),(v,\psi))_{L^2} \\
& =
(D\fS_1(\delta\tau_0,0,\de\vectau,0,0,0),v)_{L^2} \\
&=((\de\tau_0\otimes\id_{\su(E)}+\de\vectau\cdot\vecfm(A))
\rho^{-1}(\Phi\otimes\Phi^*)_{00}),v)_{L^2} \\
&= (\de\tau_0\rho^{-1}(\Phi\otimes\Phi^*)_{00},v)_{L^2}
+ \sum_{l=1}^3 
(\de\tau_l\otimes\ad(\fm_l)\rho^{-1}(\Phi\otimes\Phi^*)_{00},v)_{L^2}.
\end{align*}
Taking a sequence of $\de\tau_l$'s which approximate
$\de\tau_{l,y}\delta(\cdot,y)$, where $\delta(\cdot,y)$ is the Dirac delta
distribution supported at $y$ and $\de\tau_{l,y} \in \gl(\Lambda^+)|_y$, we
obtain the pointwise identity
$$
\langle\de\tau_{0,y}\rho^{-1}(\Phi\otimes\Phi^*)_{00}|_y,v|_y\rangle +
\sum_{l=1}^3 
\langle
\de\tau_{l,y}\otimes\ad(\fm_l)\rho^{-1}(\Phi\otimes\Phi^*)_{00}|_y,v|_y
\rangle
= 0,
$$
for all $\de\tau_{l,y} \in \gl(\Lambda^+)|_y$, $l=0,1,2,3$.  Lemma
\ref{lem:PointwiseSurjectiveForQuadratic} then implies that either 
$\rho^{-1}(\Phi\otimes\Phi^*)_{00}|_y = 0$, and thus $\Phi|_y = 0$ by Lemma
\ref{lem:UnpertSigmaPhiPhi}, or else $v|_y = 0$. If $v|_y
\ne 0$ then $\Phi$ would be zero on the nonempty open
subset $\{v\ne 0\}\cap B(x_j,R_0)$. But then  
Lemma \ref{lem:PerturbedDiracUniqueContinuation} would imply that
$\Phi\equiv 0$ on $X$, contradicting our assumption that $\Phi\not\equiv
0$. Thus, $v\equiv 0$ on $B(x_j,R_0)$, as desired. 
\end{pf}

The following similar argument shows that $\psi\equiv 0$ on the ball
$B(x_j,R_0)$.  Note that having only $v\equiv 0$ or $\psi\equiv 0$ on an
open set does not suffice to contradict Lemma
\ref{lem:NonVanishingCokernel}, as the non-vanishing result of
Lemma \ref{lem:NonVanishingCokernel} applies to the pair $(v,\psi)$.
We again begin with a linear algebra result:

\begin{lem}
\label{lem:PointwiseNonVanishDirac}
Let $S^+\in (W^+\otimes E)|_x$ and $S^-\in (W^-\otimes E)|_x$.
If $\fm_0,\dots, \fm_3$ span $\fu(E)|_x$ and
$$
\sum_{l=0}^3 \langle (\vartheta_l\otimes \fm_l) S^+,S^-\rangle = 0
$$
for all $\vartheta_0,\dots,\vartheta_3\in\Hom_\CC(W^+,W^-)|_x$, then
$S^+ =0$ or $S^- = 0$.
\end{lem}

\begin{pf}
Because $\{\fm_l\}_{l=0}^3$ spans $\fu(E)|_x$ and
$\gl(E)|_x=\fu(E)|_x\oplus\, i\fu(E)|_x$, we have $\gl(E)|_x =
\fu(E)|_x\otimes_\RR\CC$ and the set $\{\fm_l\}_{l=0}^3$ is a complex
basis for $\gl(E)|_x$.  Thus, any element of
$$
\Hom_\CC(W^+\otimes_\CC E,W^-\otimes_\CC E)|_x
\simeq \Hom_\CC(W^+,W^-)|_x\otimes_\CC \gl(E)|_x
$$ 
can be written as
$\sum_{l=0}^3\vartheta_l\otimes\fm_l$, for some
$\vartheta_l\in\Hom_\CC(W^+,W^-)|_x$, $l=0,\dots,3$. Thus, if $S^+
\ne 0$, the hypothesis implies that $S^- = 0$ and conversely, if $S^-
\ne 0$ then $S^+=0$.
\end{pf}

\begin{lem}
\label{lem:DS2Surjectivity}
Continue the hypotheses of Proposition \ref{prop:Surjectivity} and
suppose $B(x_j,R_0)$ is a ball supporting holonomy
perturbations for $A$. Then $\psi\equiv 0$ on $B(x_j,R_0)$.
\end{lem}

\begin{pf}
By hypothesis, there are holonomy sections $\fm_l :=
\fm_{j,l,\alpha}(A)$, $l=1,2,3$, which span $\su(E)|_y$, for any
point $y\in B(x_j,R_0)$, and so $\{\fm_l\}_{l=0}^3$ spans 
$\fu(E)|_y$, where $\fm_0 := i\cdot\id_E$. 
Let $\de\vartheta_l := \delta\vartheta_{j,l,\alpha} \in
\Om^1(X,\CC)$, $l=1,2,3$, denote corresponding
coefficients, and let $\delta\vecvartheta$ be a perturbation sequence with
all other coefficients zero.

By the hypothesis of Proposition \ref{prop:Surjectivity} we have
$(D\fS(\de\tau_0,\de\vartheta_0,\de\vectau,\de\vecvartheta,a,\phi),
(v,\psi))_{L^2}=0$
for all $(\de\tau_0,\de\vartheta_0,\de\vectau,\de\vecvartheta,a,\phi)$ and so
\begin{align*}
0 
&= 
(D\fS(0,\delta\vartheta_0,0,\delta\vecvartheta,0,0),(v,\psi))_{L^2} 
= (D\fS_2(0,\delta\vartheta_0,0,\delta\vecvartheta,0,0),\psi)_{L^2} \\
&=
(\rho(\delta\vartheta_0)\Phi+\de\vecvartheta\cdot\vecfm(A)\Phi,\psi)_{L^2} 
= \sum_{l=0}^3
(\rho(\de\vartheta_l)\otimes\fm_l\Phi,\psi)_{L^2}.
\end{align*}
Taking a sequence of $\de\vartheta_l$'s which approximate
$\de\vartheta_{l,y}\delta(\cdot,y)$, where $\delta(\cdot,y)$ is the Dirac
delta distribution supported at $y$ and $\de\vartheta_{l,y} \in
T^*X|_y\otimes\CC$, we obtain the pointwise identity
$$
\sum_{l=0}^3 
\langle\rho(\de\vartheta_{l,y})\otimes\fm_l\Phi|_y,\psi|_y\rangle = 0,
$$
for all $\de\vartheta_{l,y} \in T^*X|_y\otimes\CC$, $l=0,1,2,3$.  Lemma
\ref{lem:PointwiseNonVanishDirac} then implies that either 
$\Phi|_y = 0$ or $\psi|_y = 0$. If $\psi|_y
\ne 0$ then $\Phi$ would be zero on the nonempty open
subset $\{\psi\ne 0\}\cap B(x_j,R_0)$. But then Lemma
\ref{lem:PerturbedDiracUniqueContinuation} would imply that $\Phi\equiv 0$
on $X$, again contradicting our assumption that $\Phi\not\equiv 0$.
Thus, $\psi\equiv 0$ on $B(x_j,R_0)$, as desired.
\end{pf}

We can then conclude the proof of Proposition \ref{prop:Surjectivity}:

\begin{proof}[Proof of Proposition~\ref{prop:Surjectivity}]
If $(v,\psi)$ is in the cokernel of
$D\fS$ 
and $B(x_j,R_0)$ is a ball supporting holonomy perturbations for $A$, then
Lemmas \ref{lem:DS1Surjectivity} and \ref{lem:DS2Surjectivity} imply that
$(v,\psi)\equiv 0$ on $B(x_j,R_0)$.
\end{proof}

\subsection{Local reducibility implies global reducibility}
\label{subsec:LocalToGlobal}
The goal of this section is to prove Theorem
\ref{thm:LocalToGlobalReducible}. The argument has two main ingredients:
a local extension result for stabilizers of
pairs which are reducible on a ball and
a description of how these local stabilizers fit together to give a
stabilizer and thus a
reducible pair on the whole manifold.

\begin{rmk}\label{rmk:ExplainLocReduc} The fact that an anti-self-dual
connection which is reducible on an open subset is necessarily reducible on
all of $X$ is an essential part of Donaldson and Kronheimer's proof of
transversality for the moduli space of anti-self-dual connections in
\cite[\S 4.3]{DK}. The original argument of Freed and Uhlenbeck
\cite[pp. 57-58]{FU} constructs a parallel section
$\zeta$ of $\su(E)$ on the set $\{F_A\neq 0\}$.  Because the connection $A$
is Yang-Mills is anti-self-dual and therefore Yang-Mills, so
$d_A^*F_A=0=d_AF_A$, the set $\{F_A\neq 0\}$ is open, dense and connected.
The section $\zeta$ cannot develop any holonomy on $\{F_A=0\}$, so it
extends across all of $X$, showing that $A$ is globally reducible.  This
argument does not work in the case of $\PU(2)$ monopoles because the
connection $A$ is not necessarily Yang-Mills and our argument does not show
that the existence of a nonzero element $(v,\psi)$ in $\Coker D\fS$ implies
that the connection $A$ is reducible on a dense open subset of $X$.
\end{rmk}

We first state the local extension result for pair stabilizers and defer
its lengthy proof until after that of
Theorem \ref{thm:LocalToGlobalReducible}.  
Generalizations due to Taubes of the analogous result
for anti-self-dual connections, namely Lemma 4.3.21 in
\cite{DK}, appear as Theorems 4 and 5 in \cite{TauUnique}. As Taubes
points out in \cite[p. 35]{TauUnique}, unique continuation theorems for
solutions to the anti-self-dual equation do not seem to follow from
standard results for elliptic partial differential equations (such as those
of Aronszajn
\cite{Aron}) since the anti-self-dual equation does not linearize as an
elliptic equation for the connection. Since the $\PU(2)$ monopole equations
do not linearize as an elliptic system for pairs, the same remarks apply
here as well. Rather than rely on the Agmon-Nirenberg theorem for the
unique continuation property for a general class of ordinary differential
equations (Theorem \ref{thm:UniqueFlow}), Taubes proves the required unique
continuation property directly for the ordinary differential equation
induced by the anti-self-dual equation on a cylinder. (As Mrowka pointed
out to us, it should also be possible to deduce the unique continuation
results of \cite{TauUnique} by studying the anti-self-dual equation on a
ball and applying the Fredholm theory of
\cite{APS}.)  Recall that $B(x_0,r_0)\subset X$ denotes an open geodesic
ball with center at the point $x_0$ and radius $r_0$. Also, recall that if
$[A,\Phi]$ is a point in $M_{W,E}$ then Proposition
\ref{prop:GlobalReg} implies that it has a `smooth' (that is, $C^r$)
representative $(A,\Phi)$ solving \eqref{eq:PT}.

\begin{prop}
\label{prop:LocalUniqueExtension}
Let $X$ be an oriented, smooth four-manifold with $C^r$ Riemannian metric
$g$ and injectivity radius $\varrho=\varrho(x_0)$ at a point $x_0$. Suppose
that $0<r_0<r_1\le\half\varrho$. Let $(A,\Phi)$ be a $C^r$ pair solving the
$\PU(2)$ monopole equations \eqref{eq:PT} on $X$. If $u$ is a $C^{r+1}$
gauge transformation of $E|_{B(x_0,r_0)}$ satisfying $u(A,\Phi)=(A,\Phi)$
on $B(x_0,r_0)$, and if either $B(x_j,R_0)\cap B(x_0,r_1)=\emptyset$ or
$B(x_j,R_0)\subset B(x_0,r_0)$, for all balls $B(x_j,R_0)$ 
for which $\beta_j[A]>0$, then there is an extension of $u$ to a
$C^{r+1}$ gauge transformation $\hat u$ of $E|_{B(x_0,r_1)}$ with $\hat
u(A,\Phi)=(A,\Phi)$ on $B(x_0,r_1)$.
\end{prop}

\begin{rmk}
This extension result only holds on domains $B(x_0,r_1)\backslash
\barB(x_0,r_0)$ where the perturbations vanish because the perturbation
terms $\fm_{j,l,\alpha}(A)$ depend not just on the connection $A$ and
its derivatives at a point, but on the connection $A$ over open
neighborhoods in $X$. Although the unique continuation theorem of \cite{AN}
does allow certain integral terms, this still does not cover the 
perturbations we consider here with their non-local dependence on $A$.
\end{rmk}

We digress briefly to introduce some useful facts about stabilizer
subgroups of $\ssG_E$.

\begin{lem}
\label{lem:StabOfPhi}
If $u\in\Stab_\Phi$ for $\Phi\in C^0(X,W^+\otimes E)$ and
 $u\neq \id_E$ then $\Phi$ is rank one.
If $u\in S^1_Z$ and $u\neq\id_E$, then $\Phi\equiv 0$ on $X$.
\end{lem}

\begin{pf}
Because $u$ and $\Phi$ are continuous, the equality $u\Phi=\Phi$ holds at
each point $x\in X$ and so $u|_x\in \gl(E)|_x$
must be the identity on the image of $\Phi$ in $E|_x$. If $\Phi|_x$ is rank
two, then $u|_x$ must be the identity on $E|_x$, while if $u|_x\neq
\id_{E_x}$, then $\Phi|_x$ can be at most rank one.  If $u\in S^1_Z$, then
$u\Phi=e^{i\theta}\Phi=\Phi$ and so $\Phi=0$.
\end{pf}

Next we consider the stabilizers of reducible pairs. Recall from
\cite[Chapter II]{KobNom}, \cite[III.3.3]{MorganGTNotes} that the stabilizer
$\Stab_A\subset\ssG_E$ may be identified with a
subgroup of $\Aut(E|_x)$, for any point $x\in X$, by parallel
translation with respect to the connection $A$ and hence with a subgroup of
$\U(2)$ by choosing an orthonormal frame for $E|_x$; these subgroups
are again denoted by $\Stab_A$. 
If $E=L_1\oplus L_2$, let $S^1_{L_1}$ denote
the gauge transformations given by $e^{i\theta}\id_{L_1}\oplus \id_{L_2}$.

\begin{lem}
\label{lem:ReducibleStabilizers}
Let $(A,\Phi)$ be a $\PU(2)$ monopole in $C^r$ on $(E,W^+\otimes E)$.
Let $U\subset X$ be a connected open set,
with $U\cap B(x_j,R_0)$ empty 
or $B(x_j,2R_0)\subset U$ for all balls $B(x_j,R_0)$
such that $\beta_j[A]>0$.
If $\Phi\not\equiv 0$
and $A|_U$ is reducible, then $A|_U$ is reducible 
with respect
to a splitting $E|_U=L_1\oplus L_2$ where $\Phi$ is rank one on $U$, with
image contained either in $L_1$ or in $L_2$.  In the first case
$\Stab_{(A,\Phi)|_U}=S^1_{L_2}$, in the latter
$\Stab_{(A,\Phi)|_U}=S^1_{L_1}$.
\end{lem}

\begin{pf}
Because $A|_U$ is reducible, all the holonomy sections
$\fm_{j,l,\alpha}(A)$ vanish on $U$ if $B(x_j,2R_0)\subset U$.  If
$B(x_j,R_0)\cap U$ is empty the holonomy sections also vanish on $U$, as they
are supported on $B(x_j,R_0)$.  Since $\Phi\not\equiv 0$, we have
$(F_A^+)_0\not\equiv 0$ on $U$ by Lemma
\ref{lem:UnpertSigmaPhiPhi} and the equation
$(F_A^+)_0=(\id+\tau_0)\rho^{-1}(\Phi\otimes\Phi^*)_{00}$ of
\eqref{eq:PT}, so $A|_U$ is not projectively
flat. Therefore, $\Stab_{A|_U}\simeq T^2$ and 
$A|_U$ is reducible with respect
to a splitting $E|_U=L_1\oplus L_2$ by Lemma \ref{lem:ConnStab}.

Because the connection $A|_U$ is reducible with respect to the splitting
$E|_U=L_1\oplus L_2$, we can write $A|_U=A_1\oplus A_2$, where $A_1,A_2$ are
unitary connections on $L_1,L_2$, so $F_A|_U=F_{A_1}\oplus F_{A_2}$ and
$$
(F_A^+)_0 = \begin{pmatrix} \half(F_{A_1}^+ - F_{A_2}^+) & 0 \\
0 & -\half(F_{A_1}^+ - F_{A_2}^+) \end{pmatrix}
= \varpi\otimes\si_1,
$$
where $\varpi = -{\frac{i}{2}}(F_{A_1}^+ - F_{A_2}^+)\in \Om^+(U,\RR)$ and
$\si_1\in \su(2)$ is one of the Pauli matrices \eqref{eq:Pauli}. Hence,
$(F_A^+)_0$ is rank one on $U$ and the equation
$(F_A^+)_0=(\id+\tau_0)\rho^{-1}(\Phi\otimes\Phi^*)_{00}$ implies that
$(\Phi\otimes\Phi^*)_{00}$ is rank one on $U$;   
Lemma \ref{lem:RankOfTau} then
implies that $\Phi$ is also rank one on $U$ and so 
$\Phi|_U=\phi\otimes\xi$, for some $\phi\in C^r\Omega^0(U,W^+)$
and $\xi\in C^r\Omega^0(U,E)$, and
$$(\Phi\otimes\Phi^*)_{00}
=-i(\phi\otimes\phi^*)_0\otimes i(\xi\otimes\xi^*)_0.$$
Writing $\xi=\xi_1+\xi_2$ for $\xi_j\in C^r\Omega^0(U,L_j)$, we see that
$$(\xi\otimes \xi^*)_0 =
\begin{pmatrix} \half (|\xi_1|^2- |\xi_2|^2) & \xi_1\otimes \xi_2^* \\
                \xi_2\otimes \xi_1^*          & -\half (|\xi_1|^2- |\xi_2|^2)
\end{pmatrix}.$$
Since $(\id+\tau_0)\rho^{-1}(\Phi\otimes\Phi^*)_{00}
=-(\id+\tau_0)\rho^{-1}(i(\phi\otimes\phi^*)_0)\otimes
i(\xi\otimes\xi^*)_0$, we have
\begin{align*}
&(\id+\tau_0)\rho^{-1}(\Phi\otimes\Phi^*)_{00} \\
&\qquad = -(\id+\tau_0)
\rho^{-1}(i(\phi\otimes\phi^*)_0)\otimes 
\begin{pmatrix} {\frac{i}{2}}(|\xi_1|^2- |\xi_2|^2) & i(\xi_1\otimes \xi_2^*)\\
          i(\xi_2\otimes \xi_1^*)      & -{\frac{i}{2}}(|\xi_1|^2- |\xi_2|^2)
\end{pmatrix},
\end{align*}
and by comparison with our matrix expression for $(F_A^+)_0$ we see that
$\xi_1\otimes\xi_2^*=0$ on $U$ and thus at each point of $U$, either
$\xi_1=0$ or $\xi_2=0$.  Since $\nabla_A\xi=\nabla_{A_1}\xi_1\oplus
\nabla_{A_2}\xi_2$, and the perturbations
vanish on $U$,
the equation $(D_A+\vecvartheta\cdot\vecfm(A))\Phi=0$ reduces to
$(D_{A_1}+\rho(\vartheta_0))(\phi\otimes\xi_1)=0$ and 
$(D_{A_2}+\rho(\vartheta_0))(\phi\otimes\xi_2)=0$. The unique
continuation result for the perturbed Dirac operator, Lemma
\ref{lem:PerturbedDiracUniqueContinuation}, implies that if
$\phi\otimes\xi_1$ vanishes on an open subset of $U$, then
$\phi\otimes\xi_1\equiv 0$ on $U$, and similarly for $\phi\otimes\xi_2$.
If $\xi_1$ is non-zero at a point and
thus non-zero on an open neighborhood, $\xi_2$ vanishes on this open set
and by unique continuation $\xi_2\equiv 0$ on all of the connected set $U$.
Symmetrically, if $\xi_2$ is non-zero at a point, then $\xi_1\equiv 0$ on
$U$.  Thus, $\Phi=\phi\otimes\xi_1$ or $\Phi=\phi\otimes\xi_2$.

The stabilizer of $A|_U$ is $S^1_{L_1}\times S^1_{L_2}$.  If $\xi_1=0$, then
$\Stab_{\Phi|_U}=\Map(U,S^1_{L_1})$ while if $\xi_2=0$, then 
$\Stab_{\Phi|_U}=\Map(U,S^1_{L_2})$.
\end{pf}

We see that elements of the stabilizer of a pair cannot
exhibit holonomy, in the sense of the following lemma.

\begin{lem} 
\label{lem:UniqueStabilizer}
Suppose $(A,\Phi)$ is a $\PU(2)$ monopole with $\Phi\not\equiv 0$ and
that $U_1$, $U_2$ are connected open subsets of $X$. If there are gauge
transformations $u_i\in\ssG_{E|U_i}$, $i=1,2$, such that
$u_i\in\Stab_{A|_{U_i}}$, $u_i\in\Stab_{(A,\Phi)|_{U_1\cap U_2}}$, $u_i\neq
\id$, and there is a point $x\in U_1\cap U_2$ such that $u_1 = u_2$ on
$E|_x$, then $u_1=u_2$ on $E|_{U_1\cap U_2}$.
\end{lem}

\begin{pf}
Let $V\subset U_1\cap U_2$ be the dense open subset of points $\{\Phi\ne
0\}\cap U_1\cap U_2$.  Because there is a gauge transformation
$u_i\in\ssG_{E|U_i}$ with $u_i\Phi=\Phi$ over $U_i$, then $\Phi|_{U_1\cap
U_2}$ must be rank one by Lemma \ref{lem:StabOfPhi} and there is an
orthogonal decomposition $E|_V=\Imag\Phi\oplus (\Imag\Phi)^\perp$.  Since
$u_i\Phi=\Phi$ on $V$, both $u_i$ respect this decomposition and must be
the identity on $\Imag\Phi|_V$.  Thus, on $V$ we can write
$$
u_i=\begin{pmatrix} 1 & 0 \\ 0 & e^{i\theta_i}\end{pmatrix}
$$
with respect to this decomposition.  Now $\det u_i=e^{i\theta_i}$
and because $u_i\in \ssG_{E|U_i}$, the function
$\det u_i$ is constant on $U_i$ and $e^{i\theta_i}\in S^1$. Thus, if
$u_1=u_2$ on $E|_x$, then $u_1=u_2$ on all points in $V$ which can be
connected to $x$ by a path in $U_1$ and a path in $U_2$ (note that these
need not be the same paths).  Because $U_1$ and $U_2$ are connected,
$u_1=u_2$ over all of $V$.  Now $V$ is dense in $U_1\cap U_2$ and the
$u_i$ are continuous so $u_1=u_2$ over all of $U_1\cap U_2$.
\end{pf}

Theorem \ref{thm:LocalToGlobalReducible} now follows
from Proposition \ref{prop:LocalUniqueExtension}:

\begin{proof}[Proof of Theorem~\ref{thm:LocalToGlobalReducible}
given Proposition \ref{prop:LocalUniqueExtension}]
Let $(A,\Phi)$ be a
$C^r$ solution to the $\PU(2)$ monopole equations \eqref{eq:PT} with
$\Phi\not\equiv 0$ and $A$ reducible on a non-empty open set $U\subset X$.
Let $J(A)=\{j: \beta_j[A]\neq 0\}$ and let 
$\barB_J(A)=\cup_{j\in J(A)}\barB(x_j,R_0)$.
Let $U^\circ$ be a connected component of $U-\barB_J(A)$.
Lemma \ref{lem:ReducibleStabilizers} 
then implies that $\Stab_{(A,\Phi)|_{U^\circ}}\simeq S^1\neq S^1_Z$
and so there is a $C^{r+1}$ gauge transformation
$u$ of $E|_{U^\circ}$ such that $u(A,\Phi)=(A,\Phi)$ on $U^\circ$.  
Proposition \ref{prop:LocalUniqueExtension} allows
the extension of $u$ to open
subsets of $X\backslash \barB_J(A)$
containing $U^\circ$. 
To be explicit,
let $x\in U^\circ$ and $r_0=\dist_g(x,\rd U^\circ)$ and
choose $r_1>r_0$ such that
$r_1\leq\half\varrho$
and $r_1\le\min_{j\in J(A)}\dist_g(x,B(x_j,R_0))$.  By Proposition
\ref{prop:LocalUniqueExtension} there is an extension
of $u\in\Stab_{(A,\Phi)|_{B(x,r_0)}}$ to an element 
$\hat u\in\Stab_{(A,\Phi)|_{B(x,r_1)}}$. By Lemma
\ref{lem:UniqueStabilizer} we have $\hatu = u$ on 
$U^\circ \cap B(x_0,r_1)$ and not just $B(x_0,r_0)$. This gives an
extension of the stabilizer $u$ for $(A,\Phi)$ over $U^\circ$ to a
stabilizer $\hatu$ over the slightly larger open set $U^\circ\cup
B(x,r_1)$.  Since we do not assume that $X$ is simply connected, we must
check that the extension obtained by repeating this process yields a
single-valued gauge transformation over $X-\barB_J(A)$.

The consistency of two extensions follows from Lemma
\ref{lem:UniqueStabilizer}.  Let $u_i$, $i=1,2$, be two extensions of $u$ to
connected open sets $U_i^\circ$ containing $U^\circ$, so we have
$u_i\in\Stab_{(A,\Phi)|_{U_i^\circ}}$ with $u_i=u$ on $E|_{U^\circ}$.
Because $u_1=u=u_2$ on $E|_{U^\circ}$, Lemma \ref{lem:UniqueStabilizer}
implies that $u_1=u_2$ on $E|_{U_1^\circ\cap U_2^\circ}$.  Therefore, the
extensions of $u\in\Stab_{(A,\Phi)|_{U^\circ}}$ fit together to form a
global gauge transformation $\hatu\in\Stab_{A,\Phi}$ on $X-
\barB_J(A)$ such that $\hatu=u$ on $E|_{U^\circ}$.

Thus, given that $A$ is reducible on an open set $U$ containing 
$\barB_J(A)$, the above argument produces an element
$u\in\Stab_{(A,\Phi)|_{X-\barB_J(A)}}$ with
$u\notin S^1_Z$.  Then Lemma \ref{lem:ConnStab}
implies that $A$ is reducible on $X-\barB_J(A)$.
For $j=1,\dots,N_b$, let $V_j$
be an open, connected subset of $U\cap B(x_j,2R_0)$ containing
$\barB(x_j,R_0)$, such that
$V_j\cap(X-B_J(A))$ is connected
and set
$$
X_1 = X-\barB_J(A) 
\quad\text{and}\quad
X_j=(X-\barB_J(A))\cup\bigcup_{1\le k \le j-1}V_k, \qquad 2\le j\le N_b+1.
$$
Each subset $X_j$ and $V_j\cap X_j$ is connected and
$X_1\subset X_2\subset\cdots \subset X_{N_b+1} = X$. 
We extend $u\in\Stab_{(A,\Phi)|_{X-\barB_J(A)}}$ inductively over each $X_j$.
Plainly, we have $u\in \Stab_{A|_{V_j\cap X_j}}$ and 
$\Stab_{A|_{V_j}}\subset \Stab_{A|_{V_j\cap X_j}}$, so we first check that
$\Stab_{A|_{V_j}}=\Stab_{A|_{V_j\cap X_j}}$. 
By hypothesis, $A|_U$ is reducible and so for each subset $V_j\subset U$,
we have that $A|_{V_j}$ is reducible.  Because $\Phi|_{V_j\cap
X_j}\not\equiv 0$ and thus $F_A \not\equiv 0$ on $V_j\cap X_j$ by
\eqref{eq:PT} --- so $A|_{V_j\cap X_j}$ is not flat --- then Lemma
\ref{lem:ConnStab} implies that $\Stab_{A|_{V_j\cap X_j}}\simeq T^2$. The
same argument yields $\Stab_{A|_{V_j}}\simeq T^2$. (We
use the assumptions on the connectedness of $V_j$ and
$V_j\cap X_j$ here: If the sets
were not connected, the stabilizers would be $\oplus T^2$, a direct sum over
connected components.)  Thus $\Stab_{A|_{V_j}}=
\Stab_{A|_{V_j\cap X_j}}\simeq T^2$.

Hence, there is an element $u''\in \Stab_{A|_{V_j}}$ such that $u''=u'$ on
$V_j\cap X_j$, where $u'\in\Stab_{(A,\Phi)|_{X_j}}$ with $u'\ne \id$.
Together, $u'$ and $u''$ give an element $u\in\Stab_{A|_{X_{j+1}}}$ which
is not in $S^1_Z$. The connection $A$ is then reducible on $X_{j+1}$, which
implies that all the holonomy perturbations vanish on $X_{j+1}$.  
Lemma \ref{lem:ReducibleStabilizers} then shows that the pair $(A,\Phi)$ is
reducible on $X_{j+1}$ and we obtain a stabilizer $u\in
\Stab_{(A,\Phi)|_{X_{j+1}}}$, with $u\ne\id$. The construction of
$u\in\Stab_{A,\Phi}$, $u\ne\id$, is then completed by induction on $j$.
\end{proof}

\begin{rmk}
The analogue of Theorem \ref{thm:LocalToGlobalReducible} 
(namely, that local reducibility implies
global reducibility) does not hold for anti-self-dual connections without
further restrictions on the topology of $X$. For example, in \cite{FU} it
is assumed that the four-manifold $X$ is simply connected (see Lemma
4.3.21).  As described by Kronheimer and Mrowka
in \cite{KrM} one can have {\em locally reducible\/}
anti-self-dual connections (called `twisted reducibles' in \cite{KrM}); see
their Lemma 2.4 for a sharp version of the Freed-Uhlenbeck generic
metrics theorem (Corollaries 4.3.15, 4.3.18, and 4.3.19 in \cite{DK}) which
holds when the requirement that $X$ be simply connected is dropped. In our
case, we see from the proof of Theorem \ref{thm:LocalToGlobalReducible}
that globally irreducible,
locally reducible solutions $(A,\Phi)$ to \eqref{eq:PT} do not exist (at
least when $\Phi\not\equiv 0$) because the stabilizer
$u\in\Stab_{A,\Phi}$ must stabilize the section $\Phi$ and not just the
connection $A$. 
\end{rmk}

The proof of Proposition \ref{prop:LocalUniqueExtension} takes up the
remainder of this section. 

\subsubsection{The Agmon-Nirenberg unique continuation theorem}
As in the case of the anti-self-dual equation \cite[Lemma 4.3.21]{DK}, our
proof of the unique continuation property for $\PU(2)$ monopoles in 
radial gauge relies on the following special case of a more general result
due to Agmon and Nirenberg for an ordinary differential equation
on a Hilbert space \cite{AN}:

\begin{thm}\label{thm:UniqueFlow}
\cite[Theorem 2 (ii)]{AN}
Let $\fH$ be a Hilbert space and let $\sP:
\Dom(\sP(r))\subset\fH\to\fH$ be a family of symmetric linear operators for
$r\in [r_0,R)$.  Suppose that $\eta\in C^1([r_0,R),\fH)$ with
$\eta(r)\in \Dom(\sP(r))$ and $\sP\eta \in C^0([r_0,R),\fH)$ such that
\begin{align}\label{eq:FlowBound}
\left\|\frac{d\eta}{dr}-\sP(r)\eta(r)\right\|\leq c_1\left\| \eta(r)\right\|,
\end{align}
for some positive constant $c_1$ and
all $r\in [r_0,R)$.  If the function $r\mapsto (\eta(r),\sP(r)\eta(r))$
is differentiable for $r\in[r_0,R)$ and satisfies
\begin{align}\label{eq:OperatorBound}
\frac{d}{dr}\left(\eta,\sP \eta\right)
-2\Real\left({\frac{d\eta}{dr}},\sP\eta\right)
\geq
-c_2 \left\| \sP \eta\right\|\|\eta\| -c_3\|\eta\|^2,
\end{align}
for positive constants $c_2,c_3$ and every $r\in[r_0,R)$,
then the following holds: If $\eta(r_0)=0$ then $\eta(r)=0$ for all $r\in
[r_0,R)$.  
\end{thm}

Applications of \cite[Theorem 2]{AN} to the proof of unique continuation
results for first-order elliptic and parabolic
partial differential equations were considered by Agmon in \cite[Chapter
II]{Agmon}.  In our application $\sP(r)$ will be a family of first-order
partial differential operators which are self-adjoint over the closed
manifold $X$.  Theorem \ref{thm:UniqueFlow} has also been applied by
D. Salamon to prove the unique continuation property for harmonic spinors
\cite[Appendix E]{Salamon}. One of the difficulties in applying the
Agmon-Nirenberg theorem to the $\PU(2)$ monopole equations (in Gaussian
polar coordinates, $(r,\theta)$) is the requirement that the one-parameter
family of operators be {\em self-adjoint\/} with respect to a {\em fixed
inner product\/} on a {\em fixed vector space\/}. The additional
complication, not present in \cite[\S 4.3.4]{DK}, is that the induced
\spinc structures on geodesic spheres in $X$ vary with the induced family
of $r$-dependent metrics.

\begin{rmk}\label{rmk:ANComplexHilbert}
\cite[Remark, p. 209]{AN}
If $\Dom(\sP(r)) = D$ is independent of $r$ and $\sP(r)$, $r\in [r_0,R)$,
is a differentiable family of self-adjoint operators, then the left-hand
side of \eqref{eq:OperatorBound} simplifies, of course, to give the condition
\begin{equation}
\left(\eta,{\frac{d\sP}{dr}}\eta\right) \ge
-c_2 \left\| \sP\eta\right\|\|\eta\| -c_3\|\eta\|^2,
\label{eq:RealOperatorBound}
\end{equation} 
since, noting that $\sP$ is self-adjoint,
\begin{align*}
{\frac{d}{dr}}\left(\eta,\sP \eta\right)
-2\Real\left({\frac{d\eta}{dr}},\sP\eta\right)
&= \left({\frac{d\eta}{dr}},\sP\eta\right)
+ \left(\eta,{\frac{d\sP}{dr}}\eta\right)
+ \left(\eta,\sP{\frac{d\eta}{dr}}\right) \\
&\quad - \left({\frac{d\eta}{dr}},\sP\eta\right)
- \left(\sP\eta,{\frac{d\eta}{dr}}\right) \\
&= 
\left(\eta,{\frac{d\sP}{dr}}\eta\right).
\end{align*}
If $\sP(r)$, $r\in [r_0,R)$, is a differentiable family of self-adjoint
operators then it is easy to see that \eqref{eq:RealOperatorBound} follows
from the simpler condition
\begin{equation}
\left\|\frac{d\sP}{dr}\eta\right\| \le c_2\|\sP\eta\| + c_3\|\eta\|,
\label{eq:SimpleRealOperatorBound}
\end{equation}
provided $\sP\eta \in C^0([r_0,R),\fH)$.
\end{rmk}

\subsubsection{The $\PU(2)$ monopole equations in Gaussian polar
coordinates} 
Our first task is to write the pair $\PU(2)$ monopole equations
\eqref{eq:PT} as an ordinary differential equation with respect to 
Gaussian polar coordinates $(r,\theta)$ centered at a point $x_0\in X$. For
the analogous ordinary differential equation in the case of the
anti-self-dual equation, see
\cite{DK,MMR,TauUnique} and for the Seiberg-Witten equations see
\cite{KrM,Salamon}. 

Recall that $\varrho$ is the injectivity radius of $(X,g)$ at the point
$x_0$, so $\exp_{x_0}:B(0,\varrho)\subset (TX)_{x_0}\to
B(x_0,\varrho)\subset X$ is a diffeomorphism.  For each $\xi$ in the unit
sphere $S^3\subset (TX)_{x_0}$, let $\{e_i(\xi)\}_{i=1}^3$ be an
oriented, orthonormal basis for $(\RR\xi)^\perp = (TS^3)_\xi
\subset(TX)_{x_0}$; that is, let $\{e_i\}$ be an orthonormal frame for
$TS^3$.  Let $\ga_\xi(r)$ be the geodesic $\exp_{x_0}(r\xi)$,
$r\in[0,\varrho)$, so that $\ga_\xi(0)=x_0$ and $|\ga_\xi'(r)|=|\xi|$, and
let $\tau_\xi(r):(TX)_{x_0}\to (TX)_{\ga_\xi(r)}$ denote parallel
translation with respect to the metric's Levi-Civita connection along
$\ga_\xi(r)$. Let $e_i(r,\xi):=\tau_\xi(r)e_i(\xi)$ for $r\ge 0$, so that
$\{\ga_\xi'(r),e_i(r,\xi)\}$ is an orthonormal frame for
$(TX)_{\ga_\xi(r)}$ which is parallel along the radial geodesics
$\ga_\xi(r)$ and satisfies $\ga_\xi'(0)=\xi$ and $e_i(0,\xi)=e_i(\xi)$.
Denote the radial vector $\ga_\xi'(r)\in (TX)_{\ga_\xi(r)}$ by
${\frac{\rd}{\rd r}}:={\frac{\rd}{\rd r}}|_{\ga_\xi(r)}$ when no confusion
can arise. Thus, $\{{\frac{\rd}{\rd r}},e_i\}$ is an oriented, orthonormal
frame for $TX$ over $B(x_0,\varrho)\less\{0\}$, which is parallel along
radial geodesics; let $\{dr,e^i\}$ be the corresponding dual frame for
$T^*X$ over $B(x_0,\varrho)\less\{0\}$. 

With respect to the parametrization $S^3\times (0,\varrho) \simeq
B(x_0,\varrho)-\{x_0\}$, given by $(r,\xi)\mapsto \exp_{x_0}(r\xi)$, the
metric $g$ on $B(x_0,\varrho)-\{x_0\}$ pulls back to
$$
g = (dr)^2 + g_r,
$$
where $g_r$ is the metric on $S^3$ pulled back from the restriction
$g|_{S^3(x_0,r)}$ to the geodesic sphere $S^3(x_0,r) := \{x\in
X:\dist_g(x,x_0)=r\}$. Let $*_{g_r}$ denote the Hodge star operator for the
metric $g_r$ on $S^3$ and, for emphasis, we write $*_g$ for the Hodge star
operator for the metric $g$ on $X$.

Suppose that a pair 
$(A,\Phi)$ on $(\su(E),W^+\otimes E)$ is a
$C^r$ solution to the $\PU(2)$ monopole equations \eqref{eq:PT} over $X$,
\begin{align*}
\rho(F_A^+) &= \rho\tau\rho^{-1}(\Phi\otimes\Phi^*)_{00}, \\
D_{A,\vartheta_0}\Phi &= 0,
\end{align*}
where $\tau := \id_{\La^+}+\tau_0$ is an automorphism of $\La^+$ and
$D_{A,\vartheta_0}\Phi := D_A\Phi+\rho(\vartheta_0)\Phi$.  We have not
included the holonomy perturbations $\vectau\cdot\vecfm(A)$ and
$\vecvartheta\cdot\vecfm(A)$ because they vanish near $x_0$ by hypothesis. 

We obtain an isomorphism $E|_{B(x_0,\varrho)}\simeq E_0\times
B(x_0,\varrho)$ of complex two-plane bundles by choosing a unitary frame
for $E_0 := E|_{x_0}$ and using parallel translation via the $\U(2)$
connection on $E$ defined by $A$ and $A_{\det E}$ along radial geodesics
emanating from $x_0$. Let $A=B+Cdr$ denote the induced $\SO(3)$ connection
on the bundle $\su(E_0)\times S^3\times (0,\varrho)$ over $S^3\times
(0,\varrho)$ and note that $A$ is in radial gauge with respect to the point
$x_0$, so $C:= A({\frac{\rd}{\rd r}}) = 0$. We let $B=B(r)$,
$r\in(0,\varrho)$, denote the resulting one-parameter family of $\SO(3)$
connections on the bundle $\su(E_0)\times S^3$ over $S^3$. In exactly the
same way, we obtain an induced one-parameter family of $\U(2)$ connections
on the bundle $W_0^+\times S^3$ over $S^3$, $r\in(0,\varrho)$, induced by
the isomorphism $W^+|_{B(x_0,\varrho)-\{x_0\}} \simeq W_0^+\times
S^3\times(0,\varrho)$.

A section $\Phi$ of the bundle $W^+\otimes E$ over $B(x_0,\varrho)$ pulls
back, via the isomorphism $(W^+\otimes E)|_{B(x_0,\varrho)-\{x_0\}}\simeq
(W_0^+\otimes E_0)\times S^3\times (0,\varrho)$, to a one-parameter family of
sections $\Psi(r)$ of the bundle $W_0^+\otimes E_0\times S^3$ over
$S^3$. The automorphism $\tau$ of $\La^{+,g}$ pulls back to a one-parameter
family of automorphisms $\sigma(r)$ of $T^*S^3$, for $r\in(0,\varrho)$,
using the isomorphism
$$
(0,\varrho)\times T^*S^3 \to \La^{+,g}(T^*X),
\qquad (r,\al)\mapsto *_{g_r}\al + dr\wedge \al.
$$
The $g$-compatible Clifford map $\rho:T^*X\to\Hom(W^+,W^-)$ and the
isomorphism $W^+|_{B(x_0,\varrho)-\{x_0\}} \simeq W_0^+\times
S^3\times(0,\varrho)$ define a family of $g_r$-compatible Clifford maps
$\gamma(r):T^*S^3\to\End(W_0^+)$ by setting
$$
\gamma(r)
:=
\rho(dr)\rho(\ \cdot\ ).
$$
Indeed, to see this, observe that $g(dr,dr)=1$ and so for  
a family of one-forms $\alpha(r)$ on $S^3$, defined by the isomorphism
$B(x_0,\varrho)-\{x_0\} \simeq S^3\times (0,\varrho)$ and a one-form
$\alpha$ on $B(x_0,\varrho)$,  we have
\begin{align*}
\gamma_r(\alpha)^\dagger\gamma_r(\alpha)
&=
\rho(\alpha)^\dagger\rho(dr)^\dagger
\rho(dr)\rho(\alpha)
=
\rho(\alpha)^\dagger\rho(\alpha)
\\
&=
g_r(\alpha,\alpha)\,\id_{W^+_0},
\end{align*}
as required. The map $\gamma$ extends to a one-parameter family of Clifford
maps $\gamma(r):\Lambda^\bullet(T^*S^3)\otimes\CC \to \End(W_0^+)$ in the
usual way. For example, $\gamma_r(\alpha\wedge\beta) :=
\gamma_r(\alpha)\gamma_r(\beta)$, for $\alpha, \beta \in
\Omega^1(S^3)$, in which case we see that $\gamma_r(\alpha\wedge\beta) =
\rho(\alpha)\rho(\beta)$. 

With the above understood, we can proceed to rewrite the $\PU(2)$ monopole
equations \eqref{eq:PT} over the ball
$B(x_0,\varrho)$ as an ordinary differential equation
for a one-parameter family of pairs $(B(r),\Psi(r))$ on 
$(\su(E_0),W^+_0\otimes E_0)$ over $S^3$.
The curvature $F_A$ of the connection $A$
over $B(x_0,\varrho)$ is given by
$$
F_A = F_{B} - {\frac{dB}{dr}}\wedge dr.
$$
For any $\omega\in\Omega^2(X,\RR)$, we have
$\rho_+(\omega) = \rho_+(*_g\omega)$, since $\rho_+|_{\Lambda^-} =
0$. If the radial component of $\omega$ vanishes and
we consider $\omega|_{B(x_0,\varrho)}$ as a one-parameter family of
two-forms $\omega(r)$ on $S^3$, we see that $*_g\omega =
-(*_{g_r}\omega)\wedge dr$.  Combining these observations yields
\begin{align*}
\rho(F_A^+) 
&= 
\rho_+(F_A)
=  
\rho_+\left(F_{B} - \frac{dB}{dr}\wedge dr\right)
\\
&=
\rho_+\left(-(*_{g_r}F_{B})\wedge dr - \frac{dB}{dr}\wedge dr\right)
\\
&=
\rho_-(dr)\rho_+\left(*_{g_r}F_{B} + \frac{dB}{dr}\right),
\end{align*}
and therefore,
$$
\rho(F_A^+)
=
\gamma\left(*_{g_r}F_{B} + \frac{dB}{dr}\right).
$$
The section $\rho\tau\rho^{-1}(\Phi\otimes\Phi^*)_{00}$ 
of $\su(W^+)\otimes\su(E)$ over $B(x_0,\varrho)$ pulls back to the 
one-parameter family of sections
$\gamma\sigma\gamma^{-1}(\Psi\otimes\Psi^*)_{00}$ of 
$\su(W^+_0)\otimes\su(E_0)\times S^3$ over $S^3$ via the isomorphism
$W^+|_{B(x_0,\varrho)-\{x_0\}}\simeq
W_0^+\times S^3\times (0,\varrho)$ and similarly for $E$.

Let $\theta(r)$ be the induced one-parameter family of complex one-forms on
$S^3$ defined by $\vartheta_0$ on $B(x_0,\varrho)$ and the isomorphism
$B(x_0,\varrho)-\{x_0\} \simeq S^3\times (0,\varrho)$, so $\vartheta_0 =
-fdr + \theta$ on $S^3\times(0,\vartheta)$, where
$f(r)\in\Omega^0(S^3,\CC)$.  Given the preceding identifications, the Dirac
operator term in \eqref{eq:PT} can then be written over $S^3\times\{r\}$ as
\begin{align*}
D_{A,\vartheta_0}\Phi 
&= \rho(dr)\cov^{A}_{\frac{\rd}{\rd r}}\Psi
+ \sum_{i=1}^3\rho(e^i)\cov^{A}_{e_i}\Psi
+ \rho(\vartheta_0)\Psi
\\
&= \rho(dr){\frac{d\Psi}{dr}}
+ \sum_{i=1}^3\rho(e^i)\cov^{B}_{e_i}\Psi
+ \rho(\theta)\Psi - f\rho(dr)\Psi
\\
&= \rho(dr)\left({\frac{d\Psi}{dr}}
- \sum_{i=1}^3\rho(dr)\rho(e^i)\cov^{B}_{e_i}\Psi 
- \rho(dr)\rho(\theta)\Psi - f\Psi\right)\\
&= \rho(dr)\left({\frac{d\Psi}{dr}}
- \sum_{i=1}^3\gamma(e^i)\cov^{B}_{e_i}\Psi
- \gamma(\theta)\Psi - f\Psi\right),
\end{align*}
and therefore,
$$
D_{A,\vartheta_0}\Phi
= 
\rho(dr)\left({\frac{d\Psi}{dr}} - D_{B,\theta,f}\Psi\right).
$$
Hence, the $\PU(2)$ monopole equations \eqref{eq:PT} can be written as
\begin{align*}
\gamma\left({\frac{dB}{dr}} + *_{g_r}F_{B}\right) 
&= \gamma\sigma\gamma^{-1}(\Psi\otimes\Psi^*)_{00}, \\
\rho(dr){\frac{d\Psi}{dr}}
&= \rho(dr)D_{B,\theta,f}\Psi.
\end{align*}
We use $D_{B}:\Om^0(S^3,W_0^+\otimes E_0)\to
\Om^0(S^3,W_0^+\otimes E_0)$ to denote the one-parameter family of 
Dirac operators defined by the family of $g(r)$-compatible Clifford maps
$\gamma(r):T^*S^3\to\End W_0^+$, the family of $\U(2)$ connections on
$W_0^+\times S^3$, the family of $\SO(3)$ connections $B(r)$ on
$\su(E_0)\times S^3$ over $S^3$, and the family of determinant connections
on $\det E_0\times S^3$ over $S^3$.  
Since the Clifford map gives an isomorphism
$\gamma:T^*S^3\simeq\su(W_0^+)\times S^3$, the $\PU(2)$ monopole equations
then take the shape
\begin{align}
{\frac{dB}{dr}} + *_{g_r}F_{B}
&= \sigma\gamma^{-1}(\Psi\otimes\Psi)_{00},
\label{eq:PolarMonopole}\\
{\frac{d\Psi}{dr}} &= D_{B,\theta,f}\Psi,
\notag
\end{align}
for a one-parameter family $(B(r),\Psi(r))$, $r\in(0,\varrho)$, on
$(E_0,W_0^+\otimes E_0)$ over $S^3$. 

\subsubsection{The ordinary differential equation for the difference pair}
Let $(A,\Phi)$ be a $C^r$ solution to the $\PU(2)$ monopole equations
\eqref{eq:PT} over $X$. Suppose, as in the hypothesis of Proposition
\ref{prop:LocalUniqueExtension}, that there is a $C^{r+1}$ gauge
transformation $u$ of $E|_{B(x_0,r_0)}$ such that $u(A,\Phi)=(A,\Phi)$ on
$B(x_0,r_0)$. Let $(A,\Phi)$ again denote the induced pair defined by the
isomorphism $E|_{B(x_0,\varrho)}\simeq E_0\times{S^3}\times(0,\varrho)$
(given by a choice of unitary frame for $E_0=E_{x_0}$ and parallel, radial
translation via $A$) and let $u$ be the induced gauge transformation on
$E_0\times S^3\times(0,r_0)$. Then $u(A) = uAu^{-1}-(d_Au)u^{-1}$ and
$A=B+C\,dr$, where $C=A(\frac{\rd}{\rd r})=0$, and so
$du/dr = 0$. 
We now extend $u$ by parallel translation via $A$ along radial geodesics
emanating from $x_0$ to a gauge transformation $\hatu$ on 
$E_0\times S^3\times(0,\varrho)$.

Let $(\hatA,\hatPhi)=\hatu(A,\Phi)$ be the gauge-equivalent pair on
$S^3\times (0,\varrho)$, so $(\hatA,\hatPhi)$ is a $C^r$ solution to
\eqref{eq:PT} over $S^3\times (0,\varrho)$ with $\hatA=\hatu(A)$ in radial
gauge. In particular, $(\hatA,\hatPhi) = (A,\Phi)$ over $S^3\times (0,r_0)$: we
need to show that $(\hatA,\hatPhi) = (A,\Phi)$ over 
$S^3\times(0,r_1)$ in order to prove Proposition
\ref{prop:LocalUniqueExtension}. 

The one-parameter family of pairs $(\hatB(r),\hatPsi(r))$ on
$(\su(E_0),W^+_0\otimes E_0)$ over $S^3$ also satisfies the
ordinary differential equation in
\eqref{eq:PolarMonopole} so, subtracting these two pairs of ordinary
differential equations, we obtain 
for $r\geq r_0$,
\begin{align*}
{\frac{d(\hatB-B)}{dr}} + *_{g_r}(F_{\hatB}-F_{B})
&= \sigma\gamma^{-1}(\hatPsi\otimes\hatPsi^*
- \Psi\otimes\Psi^*)_{00}, \\
{\frac{d(\hatPsi-\Psi)}{dr}} 
&= D_{B,\theta,f}\Psi - D_{\hatB,\theta,f}\hatPsi.
\notag
\end{align*}
Since $F_{B}= d(B-\Gamma) + (B-\Gamma)\wedge(B-\Gamma)$ and $D_{B,\theta,f}
= D + \gamma(B-\Gamma) + \gamma(\theta) + f\id_{W^+_0\otimes E_0}$, we
obtain an ordinary differential equation for the difference pair
$$
(b,\psi) := (\hatB-B,\hatPsi-\Psi) \in \Om^1(S^3,\su(E_0))
\oplus \Om^0(S^3,W_0^+\otimes E_0),
$$
so that,
\begin{align}
{\frac{db}{dr}} &= -*_{g_r}db + *_{g_r}(\hatB\wedge b + b\wedge B) 
+ \sigma\gamma^{-1}(\hatPsi\otimes\psi^* 
+ \psi\otimes\Psi^*)_{00}, 
\label{eq:DiffPolarMonopole}\\
{\frac{d\psi}{dr}} &= D\psi + \gamma(\hatB-\Gamma)\psi 
 + \gamma(b)\Psi + \gamma(\theta)\psi + f\psi,
\notag
\end{align}
where $r\in(r_0,\varrho)$ and
$\Ga$ is the product connection on $\su(E_0)\times S^3$ over $S^3$
defined by our trivialization. The above system has the schematic form
\begin{equation}
{\frac{d(b,\psi)}{dr}} 
= \begin{pmatrix}-*_{g_r}d & 0 \\ 0 
& D\end{pmatrix}(b,\psi) 
+ \sZ_r(b,\psi),
\label{eq:SchemDiffPolarMonopole}
\end{equation}
where $\sZ_r$ is a one-parameter family of
zeroth-order operators with coefficients depending on
$g_r$, $B$, $\sigma$, and $\Psi$.

\subsubsection{Reduction to the Agmon-Nirenberg theorem}
We first observe that the operator $-*_{g_r}d$ on $\Omega^1(S^3,\su(E_0))$
is self-adjoint with respect to the $L^2$ inner product induced by the
metric $g_r$ on $S^3$. Indeed, as $*_{g_r}^2 = 1$ on $\Omega^1(S^3)$
and $d\tr(b\wedge b') = \tr(db\wedge b') - \tr(b\wedge db')$, we have
\begin{align*}
\int_{S^3}\langle b,-*_{g_r}d b'\rangle_r \,d\vol_r
&= 
- \int_{S^3}\tr(b\wedge *_{g_r}(-*_{g_r}db'))
\\
&= \int_{S^3}\tr(b\wedge db') = \int_{S^3}\tr(db\wedge b')
\\
&= 
\int_{S^3}\tr(b\wedge db') 
= \int_{S^3}\tr(b'\wedge *_{g_r}^2db)
\\
&= -\int_{S^3}\langle b',*_{g_r}d b\rangle_r \,d\vol_r,
\end{align*}
and therefore,
$$
\int_{S^3}\langle b,-*_{g_r}d b'\rangle_r \,d\vol_r
= 
\int_{S^3}\langle -*_{g_r}d b,b'\rangle_r \,d\vol_r.
$$
The Dirac operators $D:\Omega^1(S^3,W^+_0\otimes E_0) \to
\Omega^1(S^3,W^-_0\otimes E_0)$ are defined by the one-parameter family of 
$g_r$-compatible Clifford maps $\gamma_r:T^*S^3\to \End(W_0^+)$, the
one-parameter family of $\U(2)$ connections on $W^+_0\times S^3$,
one-parameter family of determinant connections on $\det E_0\times S^3$,
and the product $\SO(3)$ connection on $\su(E_0)\times S^3$. Then, by
\cite[Proposition II.5.3]{LM} we have
$$
\langle D\psi,\psi'\rangle
=
\langle \psi,D\psi'\rangle + \divg_{g_r}\xi,
$$
where $\xi(r)$ is the one-parameter family of vector fields on $S^3$ defined by
$\alpha(\xi) = -\langle\psi,\gamma(\alpha)\psi'\rangle$, for all
$\alpha\in\Omega^1(S^3)$.  Hence, the Divergence theorem for the metric
$g_r$ implies that $D$ is self-adjoint with respect the the family of $L^2$
inner products on $\Omega^1(S^3,W^+_0\otimes E_0)$ defined by $g_r$:
$$
\int_{S^3}\langle D\psi,\psi'\rangle \,d\vol_r
= 
\int_{S^3}\langle \psi,D\psi'\rangle \,d\vol_r
$$
The Agmon-Nirenberg theorem is not immediately applicable to the ordinary
differential equation in \eqref{eq:SchemDiffPolarMonopole} since the
differential operator $\sP_r := -*_{g_r}d\oplus D$ is only self-adjoint
on the Hilbert space
$\fH_r := L^2(S^3,\Lambda^1\otimes\su(E_0))\oplus L^2(S^3,W^+_0\otimes E_0)$
with $L^2$ inner products $(\cdot\ ,\ \cdot)_r$
defined by the family of metrics
$g_r$ on $S^3$, rather than a fixed inner product.

If $d\vol_r$ is the volume form on $S^3$ defined by the metric $g_r$ then
we may write
$$
d\vol_r  = h_r^2d\vol, \qquad r\in (0,\varrho),
$$
for some positive function $h_r$ on $S^3$, where $d\vol$ is the volume form
on $S^3$ defined by the standard metric. However,
\begin{equation}
\sQ_r 
:= 
h_r\sP_rh_r^{-1},
\qquad r\in (0,\varrho),
\label{eq:DefnQr}
\end{equation}
is a differentiable path of self-adjoint, first-order, elliptic
differential operators on the fixed, real Hilbert space underlying
$\fH := L^2(S^3,\Lambda^1\otimes\su(E_0))\oplus L^2(S^3,W^+_0\otimes E_0)$
with $L^2$ inner product $(\cdot\ ,\ \cdot)$
defined by the standard metric on $S^3$. Indeed, if we define a 
Hilbert-space isomorphism
$\fH_r\to\fH$ by $(b,\psi)\mapsto (\beta,\varphi):=h_r(b,\psi)$, then
\begin{align*}
\int_{S^3}\langle \sQ_r(\beta,\varphi),(\beta',\varphi')\rangle\,d\vol
&=
\int_{S^3}\langle h_r\sP_rh_r^{-1}h_r(b,\psi),h_r(b',\psi')\rangle\,d\vol
\\
&= 
\int_{S^3}\langle \sP_r(b,\psi),(b',\psi')\rangle \,d\vol_r
\\
&=
\int_{S^3}\langle (b,\psi),\sP_r(b',\psi')\rangle \,d\vol_r
\\
&=
\int_{S^3}\langle h_r(b,\psi),h_r\sP_rh_r^{-1}h_r(b',\psi')\rangle \,d\vol,
\end{align*}
and therefore,
$$
\int_{S^3}\langle \sQ_r(\beta,\varphi),(\beta',\varphi')\rangle\,d\vol
=
\int_{S^3}\langle (\beta,\varphi),\sQ_r(\beta',\varphi')\rangle \,d\vol,
\qquad r\in(0,\varrho).
$$
The operators $\sQ_r$ have dense domain $L^2_1$. 
Since $(b,\psi) = h_r^{-1}(\beta,\varphi)$, we see that
$$
{\frac{d(b,\psi)}{dr}}
= 
-h_r^{-2}{\frac{dh_r}{dr}}(\beta,\varphi) 
+ h_r^{-1}{\frac{d(\beta,\varphi)}{dr}},
$$
and so, subsituting into \eqref{eq:SchemDiffPolarMonopole} gives
\begin{equation}
\label{eq:COVSchemDiffPolarMonopole}
{\frac{d(\beta,\varphi)}{dr}} 
= \sQ_r(\beta,\varphi)+ 
h_r\sZ_rh_r^{-1}(\beta,\varphi) 
+ h_r^{-1}{\frac{dh_r}{dr}}(\beta,\varphi),
\qquad r\in (0,\varrho),
\end{equation}
with $(\beta,\varphi)$ a path in $\fH$.

\subsubsection{Verification of the Agmon-Nirenberg conditions}
We can now conclude the proof of our unique continuation result for
reducible $\PU(2)$ monopoles:

\begin{proof}[Completion of proof of
Proposition~\ref{prop:LocalUniqueExtension}] 
The estimates \eqref{eq:FlowBound} and \eqref{eq:RealOperatorBound} are, in
principle, straightforward to check since we only need them for $r$ varying in
the compact interval $[\half r_0,r_1]$; the second condition,
\eqref{eq:RealOperatorBound}, requires a little more explanation since we
need an estimate for $d\sQ/dr$.
We first check condition \eqref{eq:FlowBound}.  Comparing
\eqref{eq:DiffPolarMonopole},
\eqref{eq:SchemDiffPolarMonopole}, and \eqref{eq:COVSchemDiffPolarMonopole}
we find that
$$
\left\|{\frac{d(\beta,\varphi)}{dr}} - \sQ_r(\beta,\varphi)\right\|
\le c_1\|(\beta,\varphi)\|, \qquad \half r_0\le r\le \half\varrho,
$$
for some positive constant $c_1 = c_1(g,r_0,A,\Phi)$ and so the ordinary
differential equation in 
\eqref{eq:COVSchemDiffPolarMonopole} 
obeys the estimate \eqref{eq:FlowBound} on the interval $[\half r_0,r_1]$.
To check condition \eqref{eq:RealOperatorBound}, observe that the
definition of $\sQ_r$ in \eqref{eq:SchemDiffPolarMonopole} and
\eqref{eq:DefnQr} yields the pointwise bounds
$$
\left|{\frac{d\sQ_r}{dr}}(\beta,\varphi)\right|
\le 
 C\left(\left|\left(\cov\beta,\cov\varphi\right)\right|
+|(\beta,\varphi)|\right),
$$
where $\half r_0 \le r \le \half\varrho$, for some positive constant
$C=C(g,\tau,A_W,r_0)$, and $\cov$ denotes covariant derivatives on 
$\Lambda^1\otimes\su(E_0)$ and $W^+_0\otimes E_0\times S^3$ which are
independent of $r$. Thus, using the standard elliptic estimate for
$\sQ_r$ we obtain the $L^2$ bound
$$
\left\|{\frac{d\sQ_r}{dr}}(\beta,\varphi)\right\|
\le C\left(\|\sQ_r(\beta,\varphi)\| + \|(\beta,\varphi)\|\right), 
$$
and so, for $\half r_0 \le r \le \half\varrho$,
$$
\left|\left((\beta,\varphi),{\frac{d\sQ_r}{dr}}(\beta,\varphi)\right)\right|
\le C\left(\|\sQ_r(\beta,\varphi)\|(\beta,\varphi)\|
+ \|(\beta,\varphi)\|^2\right).
$$
Therefore, \eqref{eq:RealOperatorBound} is obeyed with $c_2=c_3=C$ on the
interval $[\half r_0,r_1]$. By Theorem
\ref{thm:UniqueFlow} and Remark \ref{rmk:ANComplexHilbert}, the solution
$(\beta(r),\varphi(r))$ vanishes for $r\in(0,r_1)$ since it vanishes for
$r\in(0,r_0)$. Hence, $(a(r),\phi(r))=0$ for $r\in(0,r_1)$ and
$\hatu(A,\Phi)=(A,\Phi)$ on the ball $B(x_0,r_1)$. This completes the proof
of Proposition \ref{prop:LocalUniqueExtension}.
\end{proof}


\appendix
\section{Holonomy perturbations and regularity}
\label{app:Appendix}
When defining our holonomy perturbations in 
\S \ref{subsubsec:HolonomyPerturbations} we deferred a detailed
discussion of several important regularity issues which arise
in their construction. The first concerns the regularity of sections of $E$
and $\su(E)$ which are constructed by parallel transport via Sobolev
connections and is described in \S \ref{app:Parallel}. The second concerns
the regularity of the $\sG_E$ equivariant maps
$\fh_{j,l,\al}:\sA_E(X)\to L^2_{k+1}(X,\su(E))$ and we show in 
\S \ref{app:HolonomyRegularity} that these maps are $C^\8$. (All of the
$\sG_E$ equivariant maps discussed here are $\ssG_E$ equivariant since
$S_Z^1$ acts trivially on connections in $\sA_E$; if the connection on
$\det E$ is not fixed, then all of the $\sG_E$ equivariant maps are $\Aut
E$ equivariant.)

The definition of the maps $\fm_{j,l,\al}$ also makes use of the
existence of a locally finite,
$C^\8$ `positive partition' on a paracompact manifold modelled
on a separable Hilbert space, in the sense of Proposition
\ref{prop:PositivePartition}. The differentials of the cutoff functions in
these partitions need not be bounded {\em a priori\/} for the usual reason
that in a Hilbert manifold, closed and bounded sets are not compact.
Nonetheless, as we shall see in 
\S \ref{app:Cutoff}, it is possible to modify the standard proof
\cite{Lang} of the existence
of a $C^\8$ partition of unity on a paracompact $C^\8$ manifold to produce
$C^\8$ cutoff functions all of whose differentials are bounded.  The sums
defining the perturbations $\vectau\cdot\vecfm$ and
$\vecvartheta\cdot\vecfm$ of \S \ref{subsubsec:HolonomyPerturbations} are
finite when restricted to any of the open neighborhoods $(\pi\circ
r_{Y_j})^{-1}(U_{j,\alpha}) \subset \sA_E^*(X)$, where
$r_{Y_j}:\sA_E(X)\to\sA_E(Y_j)$ and 
$\pi:\sA_E(X)\to\sB_E(Y_j)$ is the canonical projection.
However, the number of terms in these sums may be infinite in the
neighborhood of a reducible connection: It is for this reason that they
must converge with respect to a suitable choice of weights, as
described in
\S \ref{subsubsec:HolonomyPerturbations}. We specify the choice of weights 
in \S \ref{app:HolonomyConvergence} and explain 
why the $\sG_E$ equivariant maps $\vectau\cdot\vecfm$ and
$\vecvartheta\cdot\vecfm$ are $C^\8$ on $\sA_E(X)$.

\subsection{Parallel transport for Sobolev connections}
\label{app:Parallel}
The fact that parallel transport is well-defined for $L^2_2$ connections
has been pointed out in \cite{MorganGTNotes}. However, to construct
$L^2_{k+1}$ sections constructed by parallel transport via $L^2_k$
connections with $k\ge 2$, we need to regularize them, taking care to do this
gauge equivariantly. 
(Given an $L^2_k$ connection $A$ it does not follow that local sections
constructed by $A$-parallel transport are in $L^2_{k+1}$.) If we ignored
the issue of equivariance then it would suffice to use the 
standard smoothing kernel
described in \cite[\S 7.2 \& \S 7.3]{GT}, yielding a $C^\8$ section
irrespective of whether the connection is $L^2_k$ or $C^\8$. Instead, we
use the kernel $K_t(A)(x,y) \in \Hom(\su(E)|_y,\su(E)|_x)$, $t\in (0,\8)$,
of the heat operator
$$
K_t(A) = \exp(-td_A^*d_A): L^2(X,\su(E)) \to L^2_{k+1}(X,\su(E)),
$$
defined by the $C^\8$ metric $g$ on $TX$ and $L^2_k$ connection $A$ on $E$
\cite[\S 1.6]{Gilkey}; related constructions are described in
\cite[p. 339]{TauFrame}, \cite[p. 177]{TauStable}. The key well-known
properties \cite{Chavel}, \cite[\S 1.6]{Gilkey} 
of the heat kernel which we need are summarized below:

\begin{lem}
\label{lem:HeatKernel}
Continue the above notation and suppose $\zeta\in L^2(X,\su(E))$.
\begin{enumerate}
\item As $t\to 0$, the heat operator $\exp(-td_A^*d_A)$ converges to the $L^2$
orthogonal projection onto $\Ker(d_A^*d_A)^\perp$.
\item If $A$ on $\su(E)$ is irreducible then
$K_t(A)(x,\cdot)$ converges to the Dirac delta distribution $\de(x,\cdot)$.
\item If $A$ is $L^2_k$ then $K_t(A)\zeta$ is $L^2_{k+1}$; if $A$ is $C^\8$
then $K_t(A)\zeta$ is $C^\8$.
\item If $\zeta\in C^0(\tU,\su(E))$, for an open set $\tU\subset X$,
and $A$ is irreducible on $\su(E)$, then $K_t(A)\zeta\to \zeta$ 
in $C^0(U,\su(E))$ as $t\to 0$ for any open subset $U\Subset \tU$.
\end{enumerate}
\end{lem}

Lemma \ref{lem:HeatKernel} continues to hold for unitary connections $A$ on
Hermitian bundles $E$ over compact $C^\8$ manifolds with boundary $\barY =
Y\cup\rd Y$; we use the Neumann boundary conditions of 
\cite[p. 192]{DK}, \cite[Proposition 2.1]{TauL2}
to obtain an $L^2$ self-adjoint Laplacian
$d_A^*d_A$. Our main application will be to sections defined by parallel
translation:

\begin{lem}
\label{lem:HeatKernelRadialGauge}
Let $k\ge 2$ be an integer and let $A$ be an $L^2_k$ unitary connection on a
Hermitian two-plane bundle $E$ over $X$.  Let
$x_0$ be a point in $X$, let $\zeta_0\in \su(E)|_{x_0}$,
and let $B(x_0,r)$ be a geodesic ball with center $x_0$ and radius $r>0$.
\begin{enumerate}
\item Parallel transport with respect to $A$ is well-defined along $C^\8$
paths in $X$.
\item If $\zeta$ is the section of $\su(E)|_{B(x_0,r)}$ obtained by
$A$-parallel transport of $\zeta_0$ along radial geodesics originating at
$x_0$ and $A$ is $C^\8$ then $\zeta$ is $C^\8$ on $B(x_0,r)$. If $A$ is
$L^2_k$ then $\zeta$ is $C^0$ on $B(x_0,r)$ and its mollification
$K_t(A|_{B(x_0,r)})\zeta$ is $L^2_{k+1}$ on $B(x_0,r)$ for any $t>0$.
\end{enumerate}
\end{lem}

\begin{pf}
Let $\ga:[0,1]\to X$ be a $C^\8$ path such that $\ga(0)=x_0$ and let
$U\subset X$ be an open neighborhood of $\ga([0,1])$. We may suppose
without loss of generality that $A\in L^2_k(U,\Lambda^1\otimes\fu(2))$ is a
connection matrix and consider the parallel transport of $\eta_0\in\CC^2$.
Note that we have a continuous Sobolev embedding
$L^2_2(U)\to L^2([0,1])$ \cite[Theorem V.4]{Adams}. We wish to solve the
ordinary differential equation
\begin{equation}
\frac{d\eta}{dt} + a\eta = 0, \qquad t\in [0,1],
\label{eq:ODE}
\end{equation}
where $a = \ga^*A \in L^2([0,1],\fu(2))$, with initial condition $\eta(0) =
\eta_0$, for a solution $\eta \in (L^2_1\cap C^0)([0,1],\CC^2)$.

If $A$ is $C^\8$ then there is a unique solution $\eta\in C^\8([0,1],\CC^2)$ to
\eqref{eq:ODE} (see, for example, \cite{Hartman}) and it obeys the inequality
\cite[Lemma IV.4.1]{Hartman} 
\begin{align}
|\eta(t)| &\le |\eta_0|\exp\left(\int_0^t |a(s)|\,ds\right) 
\label{eq:ODEBound}\\
&\le |\eta_0|\exp\left(\|a\|_{L^1([0,1])}\right)
\le |\eta_0|\exp(c\|A\|_{L^2_2(U)}), \notag
\end{align}
for $t\in[0,1]$. If $A$ is only $L^2_2$, we may choose any sequence
$\{A_\al\}$ of $C^\8$ connections which converge to $A$ in $L^2_2$ and let
$\{\eta_\al\}$ be the corresponding sequence of $C^\8$ solutions to
\eqref{eq:ODE}. {}From \eqref{eq:ODEBound} 
it follows that $\{\eta_\al\}$ is $C^0$-Cauchy 
and then \eqref{eq:ODE} implies that $\{\eta_\al\}$ is
$L^2_1$-Cauchy with limit 
$\eta \in (L^2_1\cap C^0)([0,1],\CC^2)$ solving \eqref{eq:ODE}.
This proves Assertion (1).

{}From the proof of Assertion (1) we see that $\zeta\in C^0(U,\su(2))$. If the
connection $A$ is $C^\8$ on $B(x_0,r)$ then by differentiating
\eqref{eq:ODE} with respect to local coordinates (regarded as parameters
\cite[Chapter 5]{Hartman})
on $B(x_0,r)$ and solving the resulting first order ordinary differential
equations for the derivatives, it follows that the
section $\zeta$ is $C^\8$ on $B(x_0,r)$. The rest of Assertion (2) follows
from the properties of the smoothing kernel $K_t(A)(x,y)$ given in
Lemma \ref{lem:HeatKernel}.
\end{pf}

\subsection{Regularity of holonomy maps}
\label{app:HolonomyRegularity}
We consider the regularity of the holonomy maps from $A\in\sA_E$ to
$h_{\ga,x_0}(A)\in \fu(E)|_{x_0}$, to $\hat\fh_\ga(A)\in
C^0(B(x_0,r_0),\su(E))$, and to $\fh_\ga(A)\in
L^2_{k+1}(B(x_0,r_0),\su(E))$. See also \cite[\S 3.3]{BraamDon}, 
\cite[Lemma 1b.1]{Floer}, \cite[\S 8a]{TauCasson}, and \cite[\S
5A]{TauFloer} for related calculations.

Let $\ga:[0,1]\to X$ denote a parametrization for an oriented path
$\ga\subset X$ with $\gamma(0) = x_0$. The general Sobolev embedding
\cite[Theorem V.4]{Adams} implies that we have a continuous restriction map
$L^2_k(X,\La^1\otimes\su(E))\to L^2_{k-2}(\ga,\La^1\otimes\su(E))$, for
$k\ge 2$.  We fix an isomorphism $E|_\ga \simeq \ga\times\CC^2$ of $\U(2)$
bundles and denote the connection matrix corresponding to a unitary
connection $A$ on $E$ again by $A\in L^2_{k-2}(\ga,\La^1\otimes\fu(2))$. Let
$A+sa$, $s\in\RR$, be a one-parameter family of nearby connections on $E$,
with connection matrices $A(t)+sa(t)\in L^2_{k-2}(\ga,\La^1\otimes\fu(2))$. Let
$\xi_0\in\CC^2$ correspond to a point in the fiber $E|_{x_0}$ and, with
respect to the connection $A+sa$, let $\xi(t;s)$ be the parallel transport
of $\xi_0$ along $\ga(t)$, so $\xi(t;s)$ solves
\begin{equation}
\left(\frac{d}{dt} + A(t) + sa(t)\right)\xi(t;s) = 0, \qquad \xi(0;s) = \xi_0.
\label{eq:ParametrizedHolonomy}
\end{equation}
Thus, if $P_\gamma(t;A)\in\Isom(E|_{\gamma(0)},E_{\gamma(t)}) \simeq \U(2)$ 
denotes parallel transport along $\gamma$ from
$\gamma(0)$ to $\gamma(t)$ with respect to the connection $A$, we have
$$
\xi(t;s) = P_{\ga}(t;A+sa)\xi_0, \quad t\in[0,1].
$$
Therefore, our task is to compute
$$
(DP_{\gamma}(t;\cdot))_A(a)\xi_0 
= \left.\frac{d}{ds}P_{\ga}(t;A+sa)\xi_0\right|_{s=0}
= \left.\frac{d\xi}{ds}(t;s)\right|_{s=0}, \quad t\in[0,1].
$$
Differentiating \eqref{eq:ParametrizedHolonomy} with respect to $s$ we see
that $d\xi(t;s)/ds$ solves
\begin{equation}
\left(\frac{d}{dt} + A(t) + sa(t)\right)\frac{d\xi}{ds}
= -a(t)\xi(t;s), 
\qquad \frac{d\xi}{ds}(0;s) = 0.
\label{eq:DerivParametrizedHolonomy}
\end{equation}
Let $Y(t;s)$ be the fundamental matrix solution for the linear differential
operator on the left-hand side, so the solution $d\xi(t;s)/ds$ can be
written as \cite[Corollary IV.2.1]{Hartman}
$$
\frac{d\xi}{ds}(t;s) 
= -Y(t;s)\int_0^t Y^{-1}(\tau;s)a(\tau)\xi(\tau;s)\,d\tau.
$$
Since $\xi(t;s) = Y(t;s)\xi_0 = P_{\ga}(t;A+sa)\xi_0$
for any $\xi_0\in E|_{\gamma(0)}$, setting $s=0$ above
gives
\begin{equation}
(DP_{\gamma}(t;\cdot))_A(a)
= 
-P_{\ga}(t;A)\int_0^t P_{\ga}^{-1}(\tau;A)a(\tau)P_{\ga}(\tau;A)\,d\tau,
\label{eq:DifferentialParallelTransport}
\end{equation}
and so, by the Sobolev embedding theorem, a derivative bound
\begin{equation}
|(DP_{\gamma}(t;\cdot))_A(a)|
\le
\|a\|_{L^1(\ga)} \le c\|a\|_{L^2_{2,A}(X)}, \quad t\in [0,1].
\label{eq:DiffParTransportBound}
\end{equation}
The estimate \eqref{eq:DiffParTransportBound} implies that the
right-hand side of \eqref{eq:DifferentialParallelTransport} is
well-defined whenever $k\ge 2$.
In particular, when $t=1$ we have $P_{\ga}(1;A) = h_{\ga,x_0}(A)$, and so
\eqref{eq:DifferentialParallelTransport} gives
\begin{equation}
(Dh_{\gamma,x_0})_A(a) 
= - h_{\ga,x_0}(A)\int_\ga P_\ga^{-1}(A) a P_\ga(A).
\label{eq:DifferentialHolonomy}
\end{equation}
Thus, we have a well-defined differential 
$$
(Dh_{\gamma,x_0})_A: L^2_k(X,\La^1\otimes\fu(E)) 
\to T_{h_{\gamma,x_0}(A)}\U(E)|_{x_0} \simeq \fu(E)|_{x_0}
$$
of the holonomy map $h_{\gamma,x_0}:\sA_E(X)\to \U(E)|_{x_0}$ at
the point $A\in\sA_E(X)$.

A similar argument shows that all higher derivatives exist by repeatedly
applying the derivative formulas
\eqref{eq:DifferentialParallelTransport} and
\eqref{eq:DifferentialHolonomy} (see, for example, \cite[\S 8a]{TauCasson})
and so we may conclude:

\begin{lem}
\label{lem:HolonomyPointRegularity}
For $k\ge 2$, the holonomy map $h_{\gamma,x_0}:\sA_E(X)\to
\U(E)|_{x_0}$ is $C^\8$.
\end{lem}

In the same vein, the holonomy map $h_{\gamma,x_0}:\sA_E(X)\to
\SO(\su(E))|_{x_0}$ is $C^\8$ and we have:
we have:

\begin{lem}
\label{lem:HolonomySectionRegularity}
For $k\ge 2$, the following holonomy maps are $C^\8$:
\begin{align*}
&\sA_E(X)\ni A\mapsto\hat\fh_\ga(A)\in C^0(B(x_0,r_0),\su(E)), \\
&\sA_E(X)\ni A\mapsto\fh_\ga(A)\in L^2_{k+1}(B(x_0,r_0),\su(E)).
\end{align*}
\end{lem}

\begin{pf}
The fact that $\hat\fh_\ga$ is a $C^\8$ map follows from the proof of Lemma
\ref{lem:HolonomyPointRegularity}. The map $\fh_\ga$ is defined by
$\fh_\ga(A) = K_t(A)\hat\fh_\ga(A)$ and so the conclusion follows from the
fact that the map $\sA_E(X)\ni A\mapsto
K_t(A)\in\Hom(L^2(B(x_0,r_0),\su(E)),L^2_{k+1}(B(x_0,r_0),\su(E)))$ is $C^\8$.
\end{pf}

\subsection{Positive partitions and cutoff functions with bounded
differentials on Hilbert manifolds}
\label{app:Cutoff} 
We modify Lang's proof of the existence
of a $C^\8$ {\em partition of unity\/} on a paracompact $C^\8$ manifold
modelled on a separable Hilbert space
\cite[Corollary II.3.8]{Lang} to produce a $C^\8$ {\em positive
partition\/} whose cutoff functions have bounded differentials of all orders;
see Proposition \ref{prop:PositivePartition}. Recall:

\begin{prop}
\cite[Corollary II.3.8]{Lang}
\label{prop:PartitionUnity}
Let $\sX$ be a paracompact $C^\8$ manifold modelled on a separable Hilbert
space $\fH$. Then $\sX$ admits locally finite, $C^\8$ partitions of unity:
For any open cover of $\sX$ there is a countable, locally finite open
subcover $\{U_\alpha\}_{\alpha=1}^\8$ and a family of $C^\8$ functions
$\psi_\alpha:\sX\to [0,1]$ such that
\begin{itemize}
\item $\supp\psi_\alpha\subset U_\alpha$,
\item $\sum_{\alpha=1}^\8\psi_\alpha(x) = 1$ for all $x\in \sX$.
\end{itemize}
\end{prop}

We shall need the following generalizations to a Banach space setting of
the analogous familiar facts from analysis on finite-dimensional spaces.

\begin{lem}
\cite[Proposition 1.3.10]{AMR}
\label{lem:C0Complete}
Let $M$ be a topological space and let $(N,d)$ be a complete metric
space. Then the set $C(M,N)$ of all bounded continuous maps is a complete
metric space with respect to the metric $D(f,g) := \sup\{d(f(x),g(x)):
x\in M\}$.
\end{lem}

Let $\bE,\bF$ be Banach spaces and let $U\subset\bE$ be an open subset. By
analogy with the usual definitions in finite dimensions, we let
$C^s(U,\bF)$ be the set of $C^s$ maps $f:U\to\bF$ with norm 
$$
\|f\|_{C^s(U)} := \max_{0\le p \le s}\sup_{x\in U}\|(D^p f)_x\| < \8,
$$
where $\|(D^p f)_x\|$ is the norm of
$(D^pf)_x\in\Hom(\otimes^p\bE,\bF)$. Lemma \ref{lem:C0Complete}
implies that $C^0(U,\bF)$ is a Banach space. In general we have:

\begin{lem}
For any integer $s\ge 0$, the set $C^s(U,\bF)$ is a Banach space.
\end{lem}

\begin{pf}
For the case $s=1$ this follows from Lemma \ref{lem:C0Complete} and
\cite[Theorem 1.1.5]{Hormander}. The general case
\cite[pp. 113--114]{AMR} is proved in the same way using Taylor's Theorem.
\end{pf}

By analogy with the finite-dimensional case we set $C^\8(U,\bF) =
\cap_{s\ge 0}C^s(U,\bF)$. 

The $C^\8$ cutoff functions produced by Proposition
\ref{prop:PartitionUnity} may not necessarily have bounded differentials on
Hilbert manifolds, as their supports are not compact, so we now describe a
modified procedure which does produce bounded cutoff functions.

\begin{lem}
\cite[Lemma II.3.5]{Lang}
\label{lem:Scallop}
Let $M$ be a metric space and let $\{B_\alpha\}_{\alpha=1}^\8$ be a
covering of a subset $W\subset M$ by open balls. Then there exists a
locally finite open covering $\{V_\alpha\}_{\alpha=1}^\8$ of $W$ such that
$V_\alpha\subset B_\alpha$ for all $\alpha$ and 
$$
V_\alpha 
= 
B_\alpha\cap {}^c\barB(x_{\alpha,1},r_{\alpha,1})\cap\cdots
\cap {}^c\barB(x_{\alpha,\alpha-1},r_{\alpha,\alpha-1}),
$$
where ${}^c\barB = M - \barB$.
\end{lem}

\begin{lem}
\label{lem:Separation}
Let $B_0,B_1\dots,B_m$ be open balls in a Hilbert space $\fH$ and let $V$
be a {\em scalloped\/} open subset:
$$
V = B_0\cap {}^c\barB_1\cap\cdots\cap {}^c\barB_m.
$$
Then there exists a $C^\8$ function $\omega:\fH\to [0,1]$
such that $\om(x)>0$ if $x\in V$ and $\omega(x)=0$ otherwise, while the
differentials $D^p\omega$ are bounded for all $p\in\NN$ with bounds depending
only on $p,r_0,\dots,r_m$.
\end{lem}

\begin{pf}
Let $\varphi:\RR\to[0,1]$ be a $C^\8$ function such that $\varphi(t)=1$ for
$t\le 1$, $0 < \varphi(t) < 1$ for $1 < t < 2$, and $\varphi(t)=0$ for
$t\ge 2$. Suppose $B_\alpha = B(x_\alpha,r_\alpha)$ for $\alpha=0,\dots,m$. Set
$$
\varphi_\alpha(x) 
:= 
\varphi(r_\alpha^{-2}\|x-x_\alpha\|^2), 
\qquad x\in\fH, \quad \alpha = 1,\dots,m,
$$
so $0 \le \varphi_\alpha < 1$ on ${}^c\barB_\alpha$ and $\varphi_\alpha = 1$
on $\barB_\alpha$. Set $\varphi_0(x) = \varphi(2r_0^{-2}\|x-x_0\|^2)$, so $0 <
\varphi_0 \le 1$ on $B_0$ and $\varphi_0 = 0$ on ${}^cB_0$. Then
$$
\omega := \varphi_0\prod_{\alpha=1}^m(1-\varphi_\alpha)
$$
is positive on $V$ and zero on ${}^cV = \fH - V$, while 
$$
\|D^p\om\| \le C,
$$
for some positive constant $C = C(\varphi,p,r_0,\dots,r_m)$ and all $p\in\NN$.
\end{pf}

\begin{lem}
\label{lem:Cutoff}
Let $A_1,A_2$ be non-empty, closed disjoint subsets of a separable Hilbert
space $\fH$. Then there exists a $C^\8$ function $\chi:\fH\to [0,1]$ such that
$\chi(x) = 0$ if $x\in A_1$ and $\chi(x) > 0$ if $x\in A_2$, with bounded
differentials of all orders.
\end{lem}

\begin{pf}
By Lindel\"of's Theorem there is a countable collection of balls
$\{B_\alpha\}_{\alpha=1}^\8$ covering $A_2$ and such that $B_\alpha \subset
{}^cA_1 = \fH - A_1$. Let $W = \cup_\alpha B_\alpha$ and find a locally
finite refinement $\{V_\alpha\}$ of scalloped open subsets $V_\alpha\subset
B_\alpha$ using Lemma \ref{lem:Scallop} which covers $W$. Using Lemma
\ref{lem:Separation} we find a $C^\8$ function $\omega_\alpha:\fH\to\RR$
such that $\omega_\alpha$ is positive on $V_\alpha$, zero on
${}^cV_\alpha$, and has bounded differentials of all orders. The sum
$$
\chi := \sum_{\alpha=1}^\8 
\frac{2^{-\alpha}\omega_\alpha}{1 + \|\omega_\alpha\|_{C^\alpha(\fH)}}
$$
is finite on any neighborhood $V_\alpha$. The function $\chi:\fH\to\RR$ is
positive on $\cup_\alpha V_\alpha = W \supset A_2$ and zero on
$A_1$. Moreover, $\chi$ has bounded differentials of all orders, as desired.
\end{pf}

\begin{rmk}
It is at this point that the usual proof of existence of a partition of
unity encounters difficulties if we wish to ensure that the cutoff function
$\chi$ has bounded differentials of all orders. In the proof of Theorem
II.3.7 of
\cite{Lang}, Lang first constructs a function $\omega$ obeying the
conclusions of Lemma \ref{lem:Cutoff} and then a $C^\8$ function
$\sigma:\fH\to\RR$ such that $\sigma > 0$ on the closed set ${}^cU$
disjoint from $A_2$, where $U$ is the open set where $\omega>0$, and
$\sigma = 0$ on $A_2$. He then defines $\chi = \omega/(\omega+\sigma)$, so
$\chi = 1$ on $A_2$. In finite dimensions, the compactness of closed
bounded sets ensures that $\chi$ will have bounded differentials of all
orders if $A_2$ is bounded, but this obviously fails in infinite dimensions.
\end{rmk}

\begin{prop}
\label{prop:PositivePartition}
Let $\sX$ be a paracompact $C^\8$ manifold modelled on a separable Hilbert
space $\fH$. Then $\sX$ admits a locally finite, $C^\8$ {\em positive
partition:\/} For any open cover of $\sX$ there is a countable, locally
finite open subcover by parametrizations 
$\{U_\alpha,\pi_\alpha\}_{\alpha=1}^\8$, where
$\pi_\alpha:\fH\supset\bU_\alpha\simeq U_\alpha\subset\sX$, and a family of
$C^\8$ functions $\chi_\alpha:\sX\to [0,1]$ such that
\begin{itemize}
\item $\supp\chi_\alpha\subset U_\alpha$,
\item $\sum_{\alpha=1}^\8\chi_\alpha(x) > 0$ for all $x\in \sX$,
\item $\chi_\alpha\circ\pi_\alpha:\fH\to[0,1]$ has bounded differentials of
all orders. 
\end{itemize}
\end{prop}

\begin{pf}
Let $\{B_x\}$ be an open covering of $\sX$ by balls and, using
paracompactness and Lemma \ref{lem:Scallop}, let $\{U_\alpha\}$ be a
countable, locally finite refinement such that each open subset $U_\alpha$
is contained in some $B_{x(\alpha)}$, where $\pi_\alpha^{-1}:\sX\supset
B_{x(\alpha)}\to \fH$ is a coordinate chart and $\pi_\alpha^{-1}(U_\alpha)
= \bU_\alpha$. We now find open refinements $\{V_\alpha\}$ and then
$\{W_\alpha\}$ such that
$$
\barW_\alpha \subset V_\alpha \subset \barV_\alpha \subset U_\alpha,
$$
the bar denoting closure in $\sX$. For each $\alpha$, Lemma
\ref{lem:Cutoff} and the identification
$\pi_\alpha:\fH\supset\bU_\alpha\simeq U_\alpha\subset\sX$, provides a
$C^\8$ cutoff function $\chi_\alpha:\sX\to\RR$ such that $\chi_\alpha > 0$
on $\barW_\alpha$ and $\chi_\alpha = 0$ on $\sX-V_\alpha$. The $C^\8$ map
$\chi_\alpha\circ\pi_\alpha:\bU_\alpha\to[0,1]$ extends by zero to a $C^\8$
map $\chi_\alpha\circ\pi_\alpha:\fH\to[0,1]$ having bounded differentials of
all orders. Then $\{U_\alpha,\chi_\alpha\}$ is the desired positive
partition.
\end{pf}

Let $\barY = Y\cup\rd Y\subset X$ be a smooth submanifold-with-boundary.
For any open subset $U\subset \sA_E(Y)$ and $C^s$ function $f:U\to\CC$,
we define
\begin{equation}
\|f\|_{C^s(U)}
:=
\sup_{A\in U}
\sup_{\substack{1\le i\le s\\ \|a_i\|_{L^2_{k,A}(Y)}\le 1}}
|(D^s f)_A(a_1,\dots,a_s)|.
\label{eq:NormFunctionDefn}
\end{equation}
Let $\{B([A_\alpha],r_\alpha)\}$ be a countable covering of
$\sB_E^*(Y)$ by $L^2_k$ balls with subordinate locally finite
subcover $\{U_\alpha\}$ and $C^\8$ positive partition $\{\chi_\alpha\}$,
with $\supp\chi_\alpha\subset U_\alpha$.  We may suppose that
$B([A_\alpha],r_\alpha) = \pi(\bB(A_\alpha,r_\alpha))$, where
$\bB(A_\alpha,r_\alpha) \subset \sA_E^*(Y)$ is an open $L^2_k$ ball
in the Coulomb slice $A_\alpha+\Ker d_{A_\alpha}^*$ with center $A_\alpha$
and radius $r_\alpha$, and $\pi:\sA_E^*(Y) \to \sB_E^*(Y)$
is the canonical projection. Let $\bU_\alpha = \pi^{-1}(U_\alpha)\cap
\bB(A_\alpha,r_\alpha)$, so $\bU_\alpha\subset A_\alpha+\Ker
d_{A_\alpha}^*$ and
$$
\supp\chi_\alpha\circ\pi\subset \pi^{-1}(U_\alpha) 
= \sG_E(Y)\cdot\bU_\alpha \simeq
\sG_E(Y)\times\bU_\alpha. 
$$
Now $\chi_\alpha\circ\pi(u(A)) = \chi_\alpha\circ\pi(A)$
for all $u\in \sG_E(Y)$ and by the construction 
of Proposition \ref{prop:PositivePartition} the $C^\8$ map
\begin{equation}
\chi_\alpha\circ\pi: A_\alpha+\Ker d_{A_\alpha}^*
\subset L^2_{k+1}(Y,\su(E))\to [0,1]
\label{eq:CutoffSlice}
\end{equation}
has bounded differentials of all orders with respect to the {\em fixed\/}
$L^2_{k+1,A_\alpha}$ norm on $L^2_{k+1}(Y,\su(E))$, that is, 
it has bounded differentials of all orders in the
sense of Proposition \ref{prop:PositivePartition}.  Now for any $A\in
\bU_\alpha$ the $L^2_{k+1,A_\alpha}$ and $L^2_{k+1,A}$ norms on
$L^2_{k+1}(Y,\su(E))$ compare uniformly
with constants depending only on
$r_\alpha \ge \|A-A_\alpha\|_{L^2_{k;A_\alpha}}$
(since $k\ge 2$), and so the maps
\eqref{eq:CutoffSlice} have bounded differentials of all orders in the sense
of Equation \eqref{eq:NormFunctionDefn}. Moreover, the same holds for the
maps
$$
\chi_\alpha\circ\pi: \sA_E^*(Y)  \to [0,1].
$$
Finally, we have a restriction map
$r_Y:\sA_E^*(X)\to\sA_E^*(Y)$ given by $A\mapsto A|_Y$ and
the composition 
$$
\chi_\alpha\circ\pi\circ r_Y: \sA_E^*(X)\to [0,1]
$$ 
again has bounded differentials of all orders in the sense of Equation
\eqref{eq:NormFunctionDefn}. 

\subsection{Uniform convergence of holonomy perturbations on neighborhoods
of reducibles}
\label{app:HolonomyConvergence}
In \S \ref{subsubsec:HolonomyPerturbations} we constrained the sequences
perturbations $\vectau$ and $\vecvartheta$ to vary in certain weighted
$\ell^1_\delta$ spaces. We now show that a sequence of positive weights
$\delta\in\ell^\8(\RR^+)$ may be chosen in such a way that the sums
$\vectau\cdot\vecfm$ and $\vecvartheta\cdot\vecfm$ and all their differentials
converge uniformly on $\sA_E(X)$.

For any open subset $U\subset \sA_E(X)$ and $C^s$ map $\ff:U\to
L^2_{k+1}(X,\su(E))$, we define
\begin{equation}
\|\ff\|_{C^s(U)}
:=
\sup_{A\in U}
\sup_{\substack{1\le i\le s \\ \|a_i\|_{L^2_{k,A}(X)}\le 1}}
\|(D^s\ff)_A(a_1,\dots,a_s)\|_{L^2_{k+1,A}(X)}.
\label{eq:NormDefn}
\end{equation}
We then have:

\begin{prop}
\label{prop:SmoothAcrossReducibles}
Continue the notation of \S \ref{subsubsec:HolonomyPerturbations}
and let $k\ge 2$ be an integer. Then
there exists a sequence $\delta=(\delta_\alpha)_{\alpha=1}^\8 \in
\ell^\8((0,1])$ of positive weights
such that the $\sG_E$ equivariant maps
\begin{align*}
&\vectau\cdot\vecfm:\sA_E(X)\to 
L^2_{k+1}(X,\gl(\La^+)\otimes_\RR\so(\su(E))), \\
&\vecvartheta\cdot\vecfm:\sA_E(X)\to
L^2_{k+1}(X,\Hom(W^+,W^-)\otimes_\CC\fsl(E)),
\end{align*} 
are $C^\8$, with uniformly bounded differentials of all orders
in the sense of \eqref{eq:NormDefn}. In
particular, they satisfy the following $C^0$ estimates for $r\ge k+1$,
\begin{align*}
\sup_{A\in\sA_E}
\|\vecvartheta\cdot\vecfm(A)\|_{L^2_{k+1,A}(X)} 
&\le C\|\vecvartheta\|_{\ell^1_\delta(C^r(X))}, \\
\sup_{A\in\sA_E}
\|\vectau\cdot\vecfm(A)\|_{L^2_{k+1,A}(X)} 
&\le C\|\vectau\|_{\ell^1_\delta(C^r(X))},
\end{align*} 
for some positive constant $C=C(g,k)$, and more generally they satisfy the
$C^s$ estimates of \eqref{eq:CsBoundPert2} and \eqref{eq:CsBoundPert1}, for
every integer $s\ge 0$.
\end{prop}

\begin{rmk}
As should be clear from their definition, the maps
$\vectau\cdot\vecfm$ and $\vecvartheta\cdot\vecfm$ are not analytic,
although this will cause no difficulty in practice.
\end{rmk}

\begin{pf}
It suffices to consider $\vecvartheta\cdot\vecfm$, since the argument for
$\vectau\cdot\vecfm$ is obviously identical.  We first observe that the sum
$\vecvartheta\cdot\vecfm(A)$ is finite for each connection
$A\in\sA_E(X)$, and so defines an element of
$L^2_{k+1}(X,\Hom(W^+,W^-)\otimes\gl(E))$, while it is identically zero
if $A$ is reducible on $X$. However, on any open neighborhood of a
reducible connection $A\in\sA_E(X)$, the number of terms in the sum
$\vecvartheta\cdot\vecfm$ may be infinite and so our task is to choose a
sequence of weights $\delta$ such that this sum and all its differentials
converge uniformly on $\sA_E(X)$. The sum $\vecvartheta\cdot\vecfm$
will then define a $C^\8$ map. We begin with a couple of preparatory
lemmas.

\begin{lem}
\label{lem:BoundedHeatKernel}
Let $\barY = Y\cup\rd Y\subset X$ be a smooth submanifold with boundary,
let $B(A_0,r_0)\subset\sA_E^*(Y)$ be an open $L^2_k$ ball with
$k\ge 2$, and let $\eta\in L^2(Y,\su(E))$. Then the maps
$$
B(A_0,r_0)\ni A\mapsto \ff(A) :=
e^{-td_A^*d_A}\eta \in L^2_{k+1}(Y,\su(E))
$$
have bounded differentials of all orders, in the sense of \eqref{eq:NormDefn},
with constants $c=c(k,r_0,s,t)$:
$$
\|\ff\|_{C^s(B(A_0,r_0))} \le c\|\eta\|_{L^2(Y)}.
$$
\end{lem}

\begin{pf}
The Sobolev multiplication theorems and the fact that $k\ge 2$ imply that
the norms
$$
\|\zeta\|_{L^2_{k+1,A}(Y)}
\quad\text{and}\quad
\|\zeta\|'_{L^2_{k+1,A}(Y)} := \|(1+\De_A)^{(k+1)/2}\zeta\|_{L^2(Y)}
$$
on $L^2_{k+1}(Y,\su(E))$ are equivalent for 
$A\in B(A_0,r_0)\subset \sA_E^*(Y)$, with constants
depending at most on $k,r_0$. 

Since $\De_A = d_A^*d_A$, its derivative with respect to $A$ in the
direction $\de A=a$ is given by $\de\De_A = [a,\cdot]^*d_A +
d_A^*[a,\cdot]$; we use the abbreviation $\de\De_A =
(D\De_{(\cdot)})_A(a)$. 
The connection $A_0$ is irreducible by hypothesis and so for
a small enough open $L^2_k$ ball $B(A_0,r_0)$, there is a positive constant
$\la_0>0$ such that $\la[A]\ge \la_0 >0$ for all $A\in B(A_0,r_0)$, where
$\la[A]$ is the least eigenvalue of $\De_A$. Thus,
$$
\Spec(\De_A) \subset [\la_0,\8), \qquad A\in B(A_0,r_0),
$$
and, for any holomorphic function $f$ on an open neighborhood
$\Omega\subset\CC$ with $\Spec(\De_A)\subset\Omega$ and $\Gamma$ any
contour that surrounds $\Spec(\De_A)$ in $\Omega$, we have \cite{Rudin} 
\begin{align*}
f(\De_A) 
&= 
\frac{1}{2\pi i}\oint_\Gamma f(\la)(\la-\De_A)^{-1}\,d\la, \\
(Df)_{\De_A}(\de\De_A)
&=
\frac{1}{2\pi i}\oint_\Gamma f(\la)(\la-\De_A)^{-1}(\de\De_A)
(\la-\De_A)^{-1}\,d\la,
\end{align*}
and similarly for all higher-order derivatives. Note that
\begin{align*}
&(1+\De_A)^{(k+1)/2}f(\De_A) \\
&\qquad = 
\frac{1}{2\pi i}\oint_\Gamma (1+\la)^{(k+1)/2}
f(\la)(\la-\De_A)^{-1}\,d\la, \\
&(1+\De_A)^{(k+1)/2}(Df)_{\De_A}(\de\De_A) \\
&\qquad =
\frac{1}{2\pi i}\oint_\Gamma (1+\la)^{(k+1)/2}f(\la)
(\la-\De_A)^{-1}(\de\De_A)(\la-\De_A)^{-1}\,d\la.
\end{align*}
We can fix $\Omega$ and $\Gamma\subset\Omega$ such that
$\dist(\Gamma,\Spec(\De_A)) \ge d_0 > 0$, for some positive constant $d_0$
and all $A\in B(A_0,r_0)$ (see \cite[\S 1.6]{Gilkey}). Now choose $f(z) =
e^{-tz}$ and $f_{k+1}(z) = (1+z)^{(k+1)/2}f(z)$, and estimate as in
\cite[p. 53]{Gilkey}:
\begin{align*}
\|e^{-t\De_A}\eta\|_{L^2_{k+1,A}(Y)}
&\le
c\|(1+\De_A)^{(k+1)/2}e^{-t\De_A}\eta\|_{L^2(Y)} \\
&\le 
c\oint_\Gamma \left|(1+\la)^{(k+1)/2}e^{-t\la}\,d\la\right| 
\sup_{\la\in\Gamma}\|(\la-\De_A)^{-1}\eta\|_{L^2(Y)} \\
&\le 
c\|\eta\|_{L^2(Y)},
\end{align*}
which gives the desired $C^0$ bound
$$
\sup_{A\in B(A_0,r_0)}\|e^{-t\De_A}\eta\|_{L^2_{k+1,A}(Y)}
\le c\|\eta\|_{L^2(Y)},
$$
where $c=c(d_0,k,r_0,t)$. 

To obtain the $C^1$ bound, observe that
\begin{align*}
&\|(Df)_{\De_A}(\de\De_A)\eta\|_{L^2_{k+1,A}(Y)} \\
&\le
c\|(1+\De_A)^{(k+1)/2}(Df)_{\De_A}(\de\De_A)\eta\|_{L^2(Y)} \\
&\le 
c\oint_\Gamma \left|(1+\la)^{(k+1)/2}f(\la)\,d\la\right| 
\sup_{\la\in\Gamma}
\|(\la-\De_A)^{-1}(\de\De_A)(\la-\De_A)^{-1}\eta\|_{L^2(Y)} \\
&\le 
\sup_{\la\in\Gamma}
c\|(\la-\De_A)^{-1}(\de\De_A)(\la-\De_A)^{-1}\eta\|_{L^2(Y)}.
\end{align*}
Our expression for $\de \De_A$ gives
\begin{align*}
&\|(\la-\De_A)^{-1}(\de \De_A)(\la-\De_A)^{-1}\eta\|_{L^2(Y)} \\
&\le 
\|(\la-\De_A)^{-1}([a,\cdot]^*d_A + d_A^*[a,\cdot])
(\la-\De_A)^{-1}\eta\|_{L^2(Y)} \\
&\le 
c\left(\|[a,\cdot]^*d_A(\la-\De_A)^{-1}\eta\|_{L^2(Y)} 
+ \|d_A^*[a,\cdot](\la-\De_A)^{-1}\eta\|_{L^2(Y)}\right) \\
&\le 
c\|a\|_{L^2_{2;A}(Y)}\|(\la-\De_A)^{-1}\eta\|_{L^2_{2;A}(Y)} \\
&\le 
c\|a\|_{L^2_{k,A}(Y)}\|\eta\|_{L^2(Y)}, 
\end{align*}
where $c=c(r_0,d_0)$. Combining these estimates yields
$$
\|\de(e^{-t\De_A})\eta\|_{L^2_{k+1,A}(Y)} 
\le 
c\|a\|_{L^2_{k,A}(Y)}\|\eta\|_{L^2(Y)},
$$
and so we have the desired $C^1$ bound
$$
\sup_{A\in B(A_0,r_0)}
\sup_{\|a\|_{L^2_{k,A}(X)}\le 1}
\|(De^{-t\De_{(\cdot)}})_A(a)\eta\|_{L^2_{k+1,A}(Y)}
\le c\|\eta\|_{L^2(Y)}.
$$
for some $c=c(d_0,k,r_0,t)$. 
The analysis can be repeated, essentially unchanged, for all higher
differentials and is left to the reader. 
\end{pf}

\begin{lem}
\label{lem:BoundedSmoothHolonomyMap}
The $\sG_E$ equivariant holonomy maps
$$
\fm_{j,l,\alpha}:\sA_E(X)\to
L^2_{k+1}(X,\su(E))
$$ 
of \eqref{eq:HolonomySection} are $C^\8$ with bounded differentials of all
orders in the sense of \eqref{eq:NormDefn}.
\end{lem}

\begin{pf}
Recall from \eqref{eq:HeatKernelSmoothedHolonomy}
and \eqref{eq:HolonomySection} that
\begin{align*}
\fm_{j,l,\alpha}(A)
&= 
\beta_j[A]\chi_{j,\alpha}[A|_{B(x_j,2R_0)}]
\varphi_j\fh_{\gamma_{j,l,\alpha}}(A), \\
\varphi_j\fh_{\gamma_{j,l,\alpha}}(A) 
&= 
\varphi_jK_t(A|_{B(x_0,2R_0)})\hat\fh_{\gamma_{j,l,\alpha}}(A).
\end{align*}
The cutoff functions $\chi_{j,\alpha}\circ \pi\circ
r_{B(x_j,2R_0)}:\sA_E(X)\to[0,1]$
have bounded differentials of all orders, in the sense of
\eqref{eq:NormFunctionDefn},
by the remarks following Proposition \ref{prop:PositivePartition};
moreover, they are supported in $\sA_E^*(X)$. The
the functions $\beta_j:\sA_E(X)\to[0,1]$
have bounded differentials of all orders, as is clear from
their definition in \eqref{eq:ConnEnergyCutoff}. Finally, one can see from
the proofs of Lemmas \ref{lem:HolonomyPointRegularity}
and \ref{lem:HolonomySectionRegularity}, together with
Lemma \ref{lem:BoundedHeatKernel}, that the maps
$$
\varphi_jK_t(\cdot|_{B(x_0,2R_0)})\hat\fh_{\gamma_{j,l,\alpha}}:
\sA_E^*(X)\to L^2_{k+1}(X,\su(E))
$$
have bounded differentials of all orders in the sense of
\eqref{eq:NormDefn} on open subsets $\sG_E\cdot\bB(A_0,r_0)\subset
\sA_E^*(X)$, where $\bB(A_0,r_0)$ is an open $L^2_k$ ball in
$\bK_{A_0}\subset \sA_E^*(X)$.
\end{pf}

Given Lemma \ref{lem:BoundedSmoothHolonomyMap}, we have
$$
M_{\alpha,s} 
:=
\max_{j,l}\|\fm_{j,\alpha,l}\|_{C^s(\sA_E)}  < \8,
$$
where we recall that $1\le l\le 3$ and $1\le j\le N_b$.  Now choose a
sequence of positive weights $\delta=(\delta_\alpha)_{\alpha=1}^\8 \in
\ell^\8(\RR^+)$ by setting
$$
\delta_\alpha := 
\left(1 + \max_{0\le s\le\al}M_{\al,s}\right)^{-1},
\qquad\alpha \in \NN,
$$
and suppose $\vecvartheta\in\ell^1_\delta(\AAA,C^r(X))$. Then
$$
\|\vecvartheta\|_{\ell^1_\delta(C^r(X))}
= \sum_{j,l,\alpha}\delta_\alpha^{-1}\|\vartheta_{j,l,\alpha}\|_{C^r(X)} 
< \8,
$$
and therefore,
\begin{align*}
\|\vecvartheta\cdot\vecfm\|_{C^s(\sA_E)}
&\le \sum_{j,l,\alpha}\|\vartheta_{j,l,\alpha}
\fm_{j,l,\alpha}\|_{C^s(\sA_E)} \\
&\le c\sum_{j,l,\alpha}\|\vartheta_{j,l,\alpha}\|_{L^2_{k+1}(X)}
\|\fm_{j,l,\alpha}\|_{C^s(\sA_E)} 
\quad \text{by \eqref{eq:NormDefn}}\\
&\le c\sum_{j,l,\alpha}\|\vartheta_{j,l,\alpha}\|_{C^r(X)}M_{\al,s},
\end{align*}
since $r\ge k+1$. Here, $c=c(g,k)$ is a universal positive constant coming
from the continuous Sobolev multiplication $L^2_{k+1}\times L^2_{k+1,A} \to
L^2_{k+1,A}$.  Hence, using the facts that $M_{\al,s}\le 1 + \max_{0\le
t\le\al}M_{\al,t} =
\de_\al^{-1}$ for $\alpha\ge s$
and $1\le\de_\al^{-1}$ for $1\le\alpha\le s-1$ , we get
\begin{align}
\|\vecvartheta\cdot\vecfm\|_{C^s(\sA_E)}
&\le c\sum_{j,l}\left(\sum_{\alpha=1}^{s-1}M_{\al,s}
\|\vartheta_{j,l,\alpha}\|_{C^r(X)}
+ \sum_{\alpha\ge s}
\de_\al^{-1}\|\vartheta_{j,l,\alpha}\|_{C^r(X)}\right) \notag\\
&\le  c\left(1+\max_{1\le\alpha\le s-1}M_{\al,s}\right)\sum_{j,l}
\left(\sum_{\alpha=1}^{s-1}\|\vartheta_{j,l,\alpha}\|_{C^r(X)}
+ \sum_{\alpha\ge
s}\de_\al^{-1}\|\vartheta_{j,l,\alpha}\|_{C^r(X)}\right)
\notag\\
&\le  c\left(1+\max_{1\le\alpha\le s-1}M_{\al,s}\right)
\sum_{j,l,\alpha}\de_\al^{-1}\|\vartheta_{j,l,\alpha}\|_{C^r(X)} 
\notag\\
&= C\|\vecvartheta\|_{\ell^1_\delta(C^r(X))} <\8,
\label{eq:CsBoundPert2}
\end{align}
where $C=C(g,s,k)$ is defined by the last equality above. In
particular, we see that $\vecvartheta\cdot\vecfm$ is a $C^s$ map on
$\sA_E(X)$, for every integer $s\ge 0$. The same argument gives
\begin{equation}
\|\vectau\cdot\vecfm\|_{C^s(\sA_E)}
\le C\|\vectau\|_{\ell^1_\delta(C^r(X))} <\8,
\label{eq:CsBoundPert1}
\end{equation}
and $\vectau\cdot\vecfm$ is a $C^s$ map on
$\sA_E(X)$, for every integer $s\ge 0$.
\end{pf}




\end{document}